%% file: quasisplit.tex
\documentclass[11pt,reqno]{amsart}

\input{macros}

\title[Categorification of iquantum groups]{Categorification of quasi-split iquantum groups}

\author{Jonathan Brundan}
\address[J.B.]{Department of Mathematics, University of Oregon, Eugene, OR, USA
}
\urladdr{\href{https://pages.uoregon.edu/brundan}{https://pages.uoregon.edu/brundan}, \textrm{\textit{ORCiD}:} \href{https://orcid.org/0009-0009-2793-216X}{0009-0009-2793-216X}}
\email{brundan@uoregon.edu}

\author{Weiqiang Wang}
\address[W.W.]{Department of Mathematics, University of Virginia, Charlottesville, VA, USA}
\urladdr{\href{https://uva.theopenscholar.com/weiqiang-wang}{
https://uva.theopenscholar.com/weiqiang-wang}}
\email{ww9c@virginia.edu}

\author{Ben Webster}
\address[B.W.]{Department of Pure Mathematics, University of Waterloo \& Perimeter Institute for Theoretical Physics,
Waterloo, ON, Canada}
\urladdr{\href{https://uwaterloo.ca/scholar/b2webste}{https://uwaterloo.ca/scholar/b2webste},\textrm{\textit{ORCiD}:} \href{https://orcid.org/0000-0003-1896-5540}{0000-0003-1896-5540}}
\email{ben.webster@uwaterloo.ca}

\thanks{J.B. is supported in part by NSF grant DMS-2348840.
W.W. is supported in part by the NSF grant DMS-2401351. B.W. is supported by Discovery Grant RGPIN-2024-03760 from the Natural Sciences and Engineering Research Council of Canada. This research was also supported by Perimeter Institute for Theoretical Physics,
and it was completed whilst J.B. was visiting the Okinawa Institute of Science and Technology (OIST) through the Theoretical Sciences Visiting Program (TSVP)
}



\begin{document}

\begin{abstract}
We introduce a new family of graded 2-categories generalizing the 2-quantum groups introduced by Khovanov, Lauda and Rouquier. We use them to categorify quasi-split iquantum groups in all symmetric types.
\end{abstract}

\maketitle

\vspace{-11mm}
\tableofcontents

\input{s1-introduction}

\input{s2-quantum}
\input{s3-iquantum}

\input{s4-embedding}
\input{s5-nondegeneracy}
\input{s6-categorification}
\input{appendix-a}
\input{appendix-b}
\input{appendix-c}

\bibliographystyle{alphaurl}
\bibliography{quasisplit}

\end{document}

%% file: macros.tex
\usepackage{
    amsmath,
    amsfonts,
    amssymb,
    amsthm,
    amscd,
    comment,
    enumitem,
    etoolbox,
    textcomp,
    gensymb,
    mathtools,
    mathdots,
    accents
}
\usepackage[dvipsnames]{xcolor}
\usepackage[all]{xy}
\makeatletter
\@namedef{subjclassname@2020}{
  \textup{2020} Mathematics Subject Classification}
\makeatother

\usepackage[T1]{fontenc}
\usepackage{dsfont}
\usepackage[colorlinks=true, linkcolor=black, citecolor=black, urlcolor=black, breaklinks=true]{hyperref}

\usepackage{zref-clever}

\zcsetup{nameinlink=true,
countertype = {enumi=item, enumii=item, enumiii=item, enumiv=item},
counterresetby = {enumii=enumi, enumiii=enumii, enumiv=enumiii}}
    
\zcRefTypeSetup{item}{
    Name-sg = {} ,
    name-sg = {} ,
    Name-pl = {} ,
    name-pl = {} ,
}
\zcRefTypeSetup{cor}{
    Name-sg = Corollary ,
    name-sg = Corollary ,
    Name-pl = Corollaries ,
    name-pl = Corollaries ,
}
\zcRefTypeSetup{defin}{
    Name-sg = Definition ,
    name-sg = Definition ,
    Name-pl = Definitions ,
    name-pl = Definitions ,
}
\zcRefTypeSetup{eg}{
    Name-sg = Example ,
    name-sg = Example ,
    Name-pl = Examples ,
    name-pl = Examples ,
}
\zcRefTypeSetup{equation}{
    Name-sg = {} ,
    name-sg = {} ,
    Name-pl = {} ,
    name-pl = {} ,
}
\zcRefTypeSetup{lem}{
    Name-sg = Lemma ,
    name-sg = Lemma ,
    Name-pl = Lemmas ,
    name-pl = Lemmas ,
}
\zcRefTypeSetup{rem}{
    Name-sg = Remark ,
    name-sg = Remark ,
    Name-pl = Remarks ,
    name-pl = Remarks ,
}
\zcRefTypeSetup{conj}{
    Name-sg = Conjecture ,
    name-sg = Conjecture ,
    Name-pl = Conjectures ,
    name-pl = Conjectures ,
}
\zcRefTypeSetup{hyp}{
    Name-sg = Hypothesis ,
    name-sg = Hypothesis ,
    Name-pl = Hypotheses ,
    name-pl = Hypotheses ,
}
\zcsetup{countertype={subsection=subsection}}
\zcRefTypeSetup{subsection}{
    Name-sg = Subsection ,
    name-sg = subsection ,
    Name-pl = Subsections ,
    name-pl = subsections ,
}
\zcRefTypeSetup{table}{
    Name-sg = Table ,
    name-sg = Table ,
    Name-pl = Tables ,
    name-pl = Tables ,
}
\zcRefTypeSetup{theo}{
    Name-sg = Theorem ,
    name-sg = Theorem ,
    Name-pl = Theorems ,
    name-pl = Theorems ,
}

\newcommand{\cref}[1]{\zcref{#1}}
\newcommand{\Cref}[1]{\zcref[S]{#1}}

\AddToHook{env/conj/begin}{\zcsetup{countertype={theo=conj}}}
\AddToHook{env/hyp/begin}{\zcsetup{countertype={theo=hyp}}}
\AddToHook{env/cor/begin}{\zcsetup{countertype={theo=cor}}}
\AddToHook{env/defin/begin}{\zcsetup{countertype={theo=defin}}}
\AddToHook{env/eg/begin}{\zcsetup{countertype={theo=eg}}}
\AddToHook{env/lem/begin}{\zcsetup{countertype={theo=lem}}}
\AddToHook{env/rem/begin}{\zcsetup{countertype={theo=rem}}}
\AddToHook{env/table/begin}{\zcsetup{countertype={theo=table}}}

\newcommand\llbracket{[\![}
\newcommand\rrbracket{]\!]}

\newcommand\N{\mathbb{N}}

\newcommand\catQH{\mathbf{QH}}
\newcommand\Q{\mathbb{Q}}
\newcommand\Z{\mathbb{Z}}
\newcommand\kk{\Bbbk}
\newcommand{\KK}{\mathbb{K}}
\newcommand\one{\mathds{1}}

\newcommand{\ba}{{\eta}}
\newcommand{\bb}{{\varpi}}
\newcommand\bk{\text{\boldmath{$k$}}}
\newcommand\bi{\text{\boldmath{$i$}}}
\newcommand\bj{\text{\boldmath{$j$}}}
\newcommand\bl{\text{\boldmath{$\ell$}}}
\newcommand{\grrank}{\operatorname{rank}}
\newcommand\grdim{\operatorname{dim}}

\newcommand\contra{\operatorname{Contra}}
\newcommand\cov{\operatorname{Cov}}
\newcommand\B{\mathbf{B}}
\newcommand\catC{\mathbf{C}}
\newcommand{\CATC}{\mathfrak{C}}
\newcommand{\Frac}{\operatorname{Frac}}
\newcommand\catR{\mathbf{R}}
\newcommand\catS{\mathbf{S}}

\newcommand\catI{\mathbf{I}}
\newcommand\catJ{\mathbf{J}}
\newcommand\catH{\mathbf{H}}
\newcommand{\NB}{\mathbf{NB}}

\newcommand{\Kar}{\operatorname{Kar}}
\newcommand\eps{\varepsilon}
\newcommand\op{\mathrm{op}}
\newcommand{\f}{\mathbf{f}}
\newcommand{\End}{\operatorname{End}}
\newcommand{\Hom}{\operatorname{Hom}}

\newcommand\rev{\mathrm{rev}}
\newcommand\transpose{\textup{T}}
              \newcommand\lround{(\!(}
\newcommand\rround{)\!)}
\newcommand\llangle{\langle\!\langle}
\newcommand\rrangle{\rangle\!\rangle}
\newcommand{\qbinom}[2]{{\textstyle\genfrac{[}{]}{0pt}{}{#1}{#2}}}

\newcommand{\BSWW}{\mathtt{BSWW}}
\def\mod#1{{#1\text{-}\!\operatorname{mod}}}
\def\proj#1{{#1\text{-}\!\operatorname{proj}}}
\def\gmod#1{{#1\text{-}\!\operatorname{gmod}}}
\def\Nak{\operatorname{Nak}}
\def\costandard#1{{#1\text{-}\!\operatorname{gmod}_\nabla}}
\def\standard#1{{#1\text{-}\!\operatorname{gmod}_\Delta}}
\def\gproj#1{{#1\text{-}\!\operatorname{gproj}}}
\def\ginj#1{{#1\text{-}\!\operatorname{ginj}}}
\def\shape{\operatorname{sh}}
\def\oshape{\vec{\operatorname{sh}}}
\def\KL{\mathtt{K}}
\def\BWWI{\mathtt{I}}
\def\BWWJ{\mathtt{J}}

\def\del{\varsigma}
\newcommand{\atled}{\text{\rotatebox[origin=c]{180}{$\delta$}}}

\newcommand{\R}{{\mathrm R}}
\renewcommand{\H}{{\mathrm H}}
\newcommand{\HH}{{\H^\lambda}}
\newcommand{\U}{{\mathrm U}}
\newcommand{\QH}{{\mathrm{QH}}}
\newcommand{\NH}{{\mathrm{NH}}}

\newcommand\UU{\mathfrak{U}}
\newcommand\UUloc{\underline{\mathfrak{U}}}

\newcommand\xUUi{{\widetilde{\mathfrak{U}}^\imath}}
\newcommand\END{\mathbf{End}}
\newcommand\HOM{\mathbf{Hom}}

\DeclareMathOperator{\Add}{Add}

\DeclareMathOperator{\id}{id}
\DeclareMathOperator{\im}{im}

\DeclareMathOperator{\sgn}{sgn}

\DeclareMathOperator{\wt}{wt}

\usepackage{tikz}
\usetikzlibrary{arrows,patterns}
\usetikzlibrary{shapes.geometric}
\usetikzlibrary{decorations.markings,decorations.pathreplacing}
\usetikzlibrary{calc}
\usepackage{tikz-cd}

\tikzset{labl/.style={anchor=south, rotate=90, inner sep=.5mm}}

\newcommand*\getscale[1]{%
  \begingroup
    \pgfgettransformentries{\scaleA}{\scaleB}{\scaleC}{\scaleD}{\whatevs}{\whatevs}%
    \pgfmathsetmacro{#1}{sqrt(abs(\scaleA*\scaleD-\scaleB*\scaleC))}%
    \expandafter
  \endgroup
  \expandafter\def\expandafter#1\expandafter{#1}%
}

\tikzset{anchorbase/.style={>=To,baseline={([yshift=-0.5ex]current bounding box.center)}}}
\tikzset{ 
    centerzero/.style={>=To,baseline={([yshift=-0.5ex]#1)}},
    centerzero/.default={0,0}
}
\tikzset{wipe/.style={white,line width=4pt}}
\tikzset{rotateup/.style={anchor=south, rotate=90, inner sep=.8mm}}
\tikzset{rotatedown/.style={anchor=north, rotate=-90, inner sep=.8mm}}

\tikzset{gcolor/.style={green!60!black,semithick}}
\tikzset{pinhead/.style={thin,gray,fill=yellow!40!white}}
\tikzset{iQ/.style={semithick,black}}
\tikzset{bulb/.style={thin,black,fill=black}}
\tikzset{Fakebubble/.style={thin,black,fill=red!20!white}}
\tikzset{internalbubble/.style={thin,black,fill=blue!20!white}}
\tikzset{Internalbubble/.style={thin,black,fill=green}}
\tikzset{Q/.style={thin,black}}
\tikzset{Bulb/.style={thin,black,fill=white}}
\tikzset{fakebubble/.style={semithick,black,fill=red!20!white}}

\newcommand{\Circledbar}[2]{
\getscale{\scalefactor};
\node[circle,draw,Bulb,inner sep=.85pt] at (#1) 
{$\hspace{0.06pt}\color{black}\scriptscriptstyle #2$};\draw[line width=.7pt,black] ($(#1)+({-.06 / \scalefactor},{.075 / \scalefactor})$) to ($(#1)+({.06 / \scalefactor},{.075 / \scalefactor})$)}

\newcommand{\Circled}[2]{\node[circle,draw,Bulb,inner sep=.85pt] at (#1) {$\hspace{0.06pt}\color{black}\scriptscriptstyle #2$}}

\newcommand{\Pin}[3]{
    \path (#1) node[inner sep=1.6pt] (x) {} to (#2) node[rectangle,rounded corners,draw,pinhead,inner sep=2.5pt](y) {$\color{black}\scriptstyle#3$};
    \draw[Triangle Cap-,thick,gray] (x)--(y);
    \opendot{#1}
}

\newcommand{\Pinpin}[4]{
\path (#1) node[inner sep=1.6pt] (x) {} to (#3) node[rectangle,rounded corners,draw,pinhead,inner sep=2.5pt](y) {$\color{black}\scriptstyle#4$};
    \draw[Triangle Cap-,thick,gray] (x)--(#2)--(y);
\opendot{#2}; \opendot{#1}
}

\newcommand{\Pinpinpin}[5]{
\path (#1) node[inner sep=1.6pt] (x) {} to (#4) node[rectangle,rounded corners,draw,pinhead,inner sep=2.5pt](y) {$\color{black}\scriptstyle#5$};
    \draw[Triangle Cap-,thick,gray] (x)--(#2)--(#3)--(y);
\opendot{#3}; \opendot{#2}; \opendot{#1}
}

\newcommand{\Pinpinpinmid}[5]{
\path (#1) node[inner sep=1pt] (x) {} to (#3) node[inner sep=1pt] (z) {} to (#4) node[rectangle,rounded corners,draw,pinhead,inner sep=2.5pt](y) {$\color{black}\scriptstyle#5$};
\draw[-,thick,gray] (#2)--(y);
\draw[Triangle Cap-,thick,gray] (x)--(#2);
\draw[Triangle Cap-,thick,gray] (z)--(#2);
\opendot{#3}; \opendot{#2}; \opendot{#1}
}

\newcommand{\pin}[3]{
    \path (#1) node[inner sep=1.6pt] (x) {} to (#2) node[rectangle,rounded corners,draw,pinhead,inner sep=2.5pt](y) {$\color{black}\scriptstyle#3$};
    \draw[Triangle Cap-,thick,gray] (x)--(y);
    \closeddot{#1}
}

\newcommand{\pinpin}[4]{
\path (#1) node[inner sep=1.6pt] (x) {} to (#3) node[rectangle,rounded corners,draw,pinhead,inner sep=2.5pt](y) {$\color{black}\scriptstyle#4$};
    \draw[Triangle Cap-,thick,gray] (x)--(#2)--(y);
\closeddot{#2}; \closeddot{#1}
}

\newcommand{\pinpinpin}[5]{
\path (#1) node[inner sep=1.6pt] (x) {} to (#4) node[rectangle,rounded corners,draw,pinhead,inner sep=2.5pt](y) {$\color{black}\scriptstyle#5$};
    \draw[Triangle Cap-,thick,gray] (x)--(#2)--(#3)--(y);
\closeddot{#3}; \closeddot{#2}; \closeddot{#1}
}

\newcommand{\pinpinpinmid}[5]{
\path (#1) node[inner sep=1pt] (x) {} to (#3) node[inner sep=1pt] (z) {} to (#4) node[rectangle,rounded corners,draw,pinhead,inner sep=2.5pt](y) {$\color{black}\scriptstyle#5$};
\draw[-,thick,gray] (#2)--(y);
\draw[Triangle Cap-,thick,gray] (x)--(#2);
\draw[Triangle Cap-,thick,gray] (z)--(#2);
\closeddot{#3}; \closeddot{#2}; \closeddot{#1}
}

\newcommand{\circledbar}[2]{
\getscale{\scalefactor};
\node[circle,draw,bulb,inner sep=.85pt] at (#1) {$\hspace{0.06pt}\color{white}\scriptscriptstyle \pmb#2$};\draw[line width=1.1pt,white] ($(#1)+({-.06/ \scalefactor},{.08 / \scalefactor})$) to ($(#1)+({.06/ \scalefactor},{.08 / \scalefactor})$)}

\newcommand{\circled}[2]{\node[circle,draw,bulb,inner sep=.85pt] at (#1) {$\hspace{0.06pt}\color{white}\scriptscriptstyle\pmb #2$}}

\newcommand{\cross}[2][45]{
\getscale{\scalefactor};
\draw[line width=1.3pt] ($(#2)+({-.12*cos(#1)/\scalefactor},{-.12*sin(#1)/\scalefactor})$) -- ($(#2)+({.12*cos(#1)/\scalefactor},{.12*sin(#1)/\scalefactor})$);
\draw[line width=1.3pt] ($(#2)+({-.12*sin(#1)/\scalefactor},{.12*cos(#1)/\scalefactor})$) -- ($(#2)+({.12*sin(#1)/\scalefactor},{-.12*cos(#1)/\scalefactor})$);
}

\newcommand\braidup{to[out=up,in=down]}
\newcommand\braiddown{to[out=down,in=up]}

\colorlet{darkg}{green!50!black}

\newcommand\dotlabel[1]{$\color{teal}\scriptstyle{#1}$}
\newcommand\regionlabel[1]{$\color{gray}\scriptstyle{#1}$}
\newcommand\strandlabel[1]{$\color{olive}\scriptstyle{#1}$}
\newcommand\stringlabel[2]{\node at (#1) {\strandlabel{#2}}} 
 
\newcommand\catlabel[2]{\node at (#1) {$\color{darkg}\scriptstyle{#2}$}} 
\newcommand\botlabel[1]{node[inner sep=0.5pt,anchor=north] {\strandlabel{#1}}}
\newcommand\toplabel[1]{node[anchor=south,inner sep=0.5pt] {\strandlabel{#1}}}
\newcommand\dottybubblelabel[2]{\node at (#1) {$\color{teal}\scriptstyle{#2}$}} 
\newcommand\bubblelabel[2]{\node at (#1) {$\scriptstyle{#2}$}} 
\newcommand\region[2]{\node at (#1) {\regionlabel{#2}}}

\newcommand\closeddot[2][black]{
    \node at (#2) {\color{#1}$\scriptstyle{\bullet}$}
}

\newcommand\opendot[2][black]{
    \node at (#2) {\color{white} $\scriptstyle{\bullet}$};\node at (#2) {\color{#1}$\scriptstyle{\circ}$}
}

\newcommand\Opendot[2][black]{
    \node at (#2) {\color{white} $\bullet$};\node at (#2) {\color{#1}$\circ$}
}

\newcommand\multcloseddot[3]{
    \closeddot{#1};
    \draw (#1) node[anchor=#2] {\dotlabel{#3}}
}

\newcommand\multopendot[3]{
    \opendot{#1};
    \draw (#1) node[anchor=#2] {\dotlabel{#3}}
}

\newcommand{\teleporter}[3][black]{
\draw[Triangle Cap-Triangle Cap,ultra thick,blue!40!white] (#2)--(#3);
\opendot{#2}; \opendot{#3}}

\newcommand{\Teleporter}[3][black]{
\draw[Triangle Cap-Triangle Cap,ultra thick,green] (#2)--(#3);
\opendot{#2}; \opendot{#3}}

\newcommand{\bentteleporter}[4][black]{
\draw[Triangle Cap-Triangle Cap,ultra thick,blue!40!white] (#2)..controls#4..(#3);
\opendot{#2}; \opendot{#3}
}

\newcommand\clockwisebubble[2][black]{%
\draw[-to,thin,#1] (#2)++(0,.2) arc(90:-270:0.2);
}
\newcommand\anticlockwisebubble[2][black]{%
  \draw[-to,thin,#1] (#2)++(0,.2) arc(90:450:0.2);
}
\newcommand\filledclockwisebubble[2][black]{%
\draw[-to,#1,Fakebubble] (#2)++(0,.2) arc(90:-270:0.2);
}
\newcommand\filledanticlockwisebubble[2][black]{%
  \draw[-to,#1,Fakebubble] (#2)++(0,.2) arc(90:450:0.2);
}

\newcommand\anticlockwiseinternalbubble[2][black]{%
 \getscale{\scalefactor};
 \draw[-to,#1,internalbubble] ($(#2)+({-.12 / \scalefactor},0)$)  arc(-180:180:{.12/\scalefactor});
}
\newcommand\clockwiseinternalbubble[2][black]{%
 \getscale{\scalefactor};
 \draw[-to,#1,internalbubble] ($(#2)+({.12 / \scalefactor},0)$)  arc(0:-360:{.12/\scalefactor});
}

\newcommand\babyanticlockwiseinternalbubble[2][black]{%
 \getscale{\scalefactor};
 \draw[-to,#1,internalbubble] ($(#2)+({-.1 / \scalefactor},0)$)  arc(-180:180:{.1/\scalefactor});
}
\newcommand\babyclockwiseinternalbubble[2][black]{%
 \getscale{\scalefactor};
 \draw[-to,#1,internalbubble] ($(#2)+({.1 / \scalefactor},0)$)  arc(0:-360:{.1/\scalefactor});
}

\newcommand\Anticlockwiseinternalbubble[2][black]{%
 \getscale{\scalefactor};
 \draw[-to,#1,Internalbubble] ($(#2)+({-.12 / \scalefactor},0)$)  arc(-180:180:{.12/\scalefactor});
}
\newcommand\Clockwiseinternalbubble[2][black]{%
 \getscale{\scalefactor};
 \draw[-to,#1,Internalbubble] ($(#2)+({.12 / \scalefactor},0)$)  arc(0:-360:{.12/\scalefactor});
}

\newcommand\Babyanticlockwiseinternalbubble[2][black]{%
 \getscale{\scalefactor};
 \draw[-to,#1,Internalbubble] ($(#2)+({-.1 / \scalefactor},0)$)  arc(-180:180:{.1/\scalefactor});
}
\newcommand\Babyclockwiseinternalbubble[2][black]{%
 \getscale{\scalefactor};
 \draw[-to,#1,Internalbubble] ($(#2)+({.1 / \scalefactor},0)$)  arc(0:-360:{.1/\scalefactor});
}

\newtheorem{theo}{Theorem}[section]

\newtheorem{lem}[theo]{Lemma}
\newtheorem{cor}[theo]{Corollary}
\newtheorem{conj}[theo]{Conjecture}

\theoremstyle{definition}
\newtheorem{defin}[theo]{Definition}
\newtheorem{hyp}[theo]{Hypothesis}
\newtheorem{rem}[theo]{Remark}
\newtheorem{eg}[theo]{Example}

\numberwithin{equation}{section}
\allowdisplaybreaks

%% file: s1-introduction.tex
\setcounter{section}{0}

\section{Introduction}\label{s1-introduction}

Drinfeld-Jimbo quantum groups were categorified by Khovanov and Lauda \cite{KL3} and Rouquier \cite{Rou}. They introduced a family of graded 2-categories whose Grothendieck rings are isomorphic to Lusztig's modified integral form for the 
underlying quantum groups.
This was proved 
for $\mathfrak{sl}_n$ and conjectured for other types in \cite{KL3}, and proved in general in \cite{canonical,unfurling}. 
The ``2-quantum groups'' arising from 
this work play an important role in higher 
representation theory and the related program to categorify Reshetikhin-Turaev invariants in low-dimensional topology.

In this paper, we introduce a more general family of graded 2-categories which are expected to categorify {\em iquantum groups} arising from quasi-split quantum symmetric pairs. We prove this whenever the underlying Cartan matrix is symmetric. For simplicity in the rest of the introduction, we restrict our attention just 
to this symmetric case.

\vspace{2mm}

Let $Q$ be a loop-free quiver with vertex set $I$
and $\#(i\rightarrow j)$ arrows from vertex $i$ to vertex $j$. 
The corresponding Cartan matrix $A = (a_{i,j})_{i, j \in I}$ is defined by
$$
a_{i,j} :=
\begin{cases}
2&\text{if $i=j$}\\
-\#(i \rightarrow j) - \#(j \rightarrow i)&\text{if $i \neq j$.}
\end{cases}
$$
Also fix a realization of this Cartan matrix, that is, free Abelian groups $X$ and $Y$
with a perfect pairing
$Y \times X \rightarrow \Z, (h,\lambda) \mapsto h(\lambda)$, 
and choices of {\em simple coroots} 
$h_i \in Y$ and {\em simple roots} $\alpha_j \in X$
such that $h_i(\alpha_j) = a_{i,j}$ for $i,j \in I$.
We assume the simple coroots 
are linearly independent.
Let $\tau:X \rightarrow X$ be an involution 
with dual $\tau^*:Y\rightarrow Y$, such that
$\tau(\alpha_i) = \alpha_{\tau i}$
and $\tau^*(h_i) = h_{\tau i}$
for an induced involution $\tau:I \rightarrow I$.
It follows that
$a_{\tau i, \tau j} = a_{i,j}$
for all $i,j \in I$.
Let
$Y^\imath := \ker(\id + \tau^*)$.
Also define $\del \in \Z^I$
by setting
$$
\del_i := \begin{cases}
-1&\text{if $i = \tau i$}\\
\#(i \rightarrow \tau i)&\text{if $i \neq \tau i$.}
\end{cases}
$$
To this data, there is associated a {\em quasi-split iquantum group} $\U^\imath$. It is one of a general family of $\Q(q)$-algebras introduced 
by Kolb \cite{Kolb}, following Letzter's ground-breaking work \cite{letzter1} in finite type.
Letzter and Kolb
constructed $\U^\imath$ as a certain coideal subalgebra of the quantized enveloping algebra $\U
= \langle q^h, e_i, f_i\:|\:h \in Y,i \in I\rangle$ of this Cartan type generated by $q^{h}\:(h \in Y^\imath)$ and
$$
b_i :=  f_i + q^{\del_i} e_{\tau i} k_i^{-1}
$$ 
for $i \in I$, where $k_i := q^{h_i}$.
The pair $(\U, \U^\imath)$ is called a 
{\em quantum symmetric pair}. 
Ordinary quantum groups can be viewed as a special case: they are the quasi-split iquantum groups arising from the quantum symmetric pairs of {\em diagonal type}; see \cite[Rem.~4.10]{BW18QSP} or \cref{iexisting} below. 

We need here the modified form $\dot\U^\imath$
of the iquantum group $\U^\imath$
introduced in \cite{BW18QSP,BW21iCB}. 
This is a locally unital algebra with mutually orthogonal distinguished idempotents $1_\lambda$ indexed by {\em iweights} $\lambda$, that is, elements of the
quotient Abelian group $X^\imath := X / \im(\id+\tau)$.  
Denoting a pre-image of $\lambda \in X^\imath$ in $X$ by $\hat\lambda$, we use 
$\lambda \pm \alpha_i$ to denote the image of $\hat\lambda\pm \alpha_i$ in $X^\imath$, and set
$$
\lambda_i := 
(h_i-h_{\tau i})(\hat\lambda) 
\in \Z.
$$
Also let $b_i^{(n)} 1_\lambda$ be 
the usual divided power $b_i^n 1_\lambda / [n]_q^! $ if $i \neq \tau i$, or the {\em idivided power} from \cite{BW18KL,BeW18} 
if $i =\tau i$ (see \cref{idividedpowerrelation} below). 
Then the algebra $\dot\U^\imath$ can be defined by generators and relations 
as the
locally unital $\Q(q)$-algebra with distinguished idempotents
$1_\lambda\:(\lambda \in X^\imath)$ and
generators
$b_i 1_\lambda = 1_{\lambda - \alpha_i} b_i$
for all $i \in I$ and $\lambda \in X^\imath$,
subject to the relation
\begin{align*}
\sum_{n=0}^{1-a_{i,j}}(-1)^n
b_i^{(n)} b_j b_i^{(1-a_{i,j}-n)} 1_\lambda
=
\delta_{i,\tau j}\!
\prod_{r=1}^{-a_{i,j}}(q^r-q^{-r})\times
\frac{(-1)^{a_{i,j}}q^{\lambda_i-\del_i-\binom{a_{i,j}}{2}}-q^{
\binom{a_{i,j}}{2}+\del_i-\lambda_i}}{q-q^{-1}}
b_i^{(-a_{i,j})}
1_\lambda
\end{align*}
for all $i \neq j$ in $I$ and $\lambda \in X^\imath$.
The mysterious expression on the 
right hand side of this relation in the case $i = \tau j$ was worked out originally in works of Letzter \cite{letzzonal} and
Balagović and Kolb \cite{BK}.
The case $i = \tau i$ is also difficult as idivided powers are more complicated; the relation in this situation first appeared in \cite{serre}.

\vspace{2mm}

Now we come to the new definition of the
2-iquantum group $\UU^\imath$. 
The full formulation of this given in \cref{def2iqg}
applies to symmetrizable (not merely symmetric) Cartan matrices but
requires a mild additional hypothesis: {\em if $i,j \in I$ are $\tau$-fixed points then
$a_{i,j} \equiv a_{j,i}\pmod{2}$} (see \cref{unfortunately}).
The definition is a bit tidier if we make the following stronger assumption:
{\em if $i,j \in I$ are $\tau$-fixed points then $a_{i,j}$ is even}.
This makes it possible to choose the orientation of edges in $Q$ in such a way that 
$\#(i \rightarrow j) = \#(\tau j \rightarrow \tau i)$; 
in other words, $\tau$ defines an isomorphism from $Q$ to the opposite quiver.
We assume this is the case for the remainder of the introduction.
We also fix an algebraically closed field $\kk$,
which should not be of characteristic 2
if $a_{i,\tau i} \neq 0$ for some $i$.
For each $i$, let
$c_i:X \rightarrow \kk^\times$ be a group homomorphism
such that
\begin{itemize}
\item
$c_i(\alpha_j) = (-1)^{\#(j \rightarrow i)}$ for all $i,j \in I$;
\item
$c_{\tau i}(\tau(\lambda)) = (-1)^{h_i(\lambda)} c_i(\lambda)$
for all $\lambda \in X$ and $i \in I$ with $i \neq \tau i$.
\end{itemize}
Assuming that the simple roots are linearly independent, such a family of what we call {\em normalization homomorphisms} always exists.
We also fix the choice of a decomposition $I = I_1 \sqcup I_0 \sqcup I_{-1}$
so that $I_0 = \{i \in I\:|\:\tau i = i\}$
and $I_{-1} = \tau(I_1)$.
Then, for $i,j \in I$ and $\lambda \in X^\imath$
with pre-image $\hat\lambda \in X$, we let
\begin{align*}
Q_{i,j}(x,y)&:=\begin{cases}
 (x-y)^{\#(i \rightarrow j)}(y-x)^{\#(j \rightarrow i)}
&\text{if $i \neq j$}\\
0&\text{if $i=j$,}
\end{cases}\\
Q^\imath_{i,j}(x,y)&:=
(-1)^{\delta_{i,\tau j}} Q_{i,j}(x,y),\\
\zeta_i(\lambda)&:=
\begin{cases}
2^{\del_i} (-1)^{h_{\tau i}(\hat\lambda)}c_i(\hat\lambda)/c_{\tau i}(\hat\lambda)&\text{if $i \in I_1$}\\
-\frac{1}{2} (-1)^{h_{i}(\hat\lambda)}
&\text{if $i \in I_0$}\\
-2^{\del_i} (-1)^{h_{\tau i}(\hat\lambda)}c_i(\hat\lambda)/c_{\tau i}(\hat\lambda)&\text{if $i \in I_{-1}$.}
\end{cases}.
\end{align*}

The {\em 2-iquantum group} $\UU^\imath$ with these parameters is the graded $\kk$-linear 2-category
with object set $X^\imath$, 
generating 1-morphisms 
$B_i \one_\lambda = \one_{\lambda -\alpha_i} B_i
: \lambda \rightarrow \lambda-\alpha_i$
for $\lambda \in X^\imath$ and $i \in I$,
with identity 2-endomorphisms
denoted by unoriented strings
$\begin{tikzpicture}[iQ,anchorbase]
\draw[-] (0,-0.2)\botlabel{i} -- (0,0.2);
\region{0.2,0}{\lambda};
\region{-0.4,0}{\lambda-\alpha_i};
\end{tikzpicture}$,
and the generating 2-morphisms listed in \cref{table0}.
\begin{table}
\begin{align*}
\begin{array}{|l|c|}
\hline
\hspace{13mm}\text{Generator}&\text{Degree}\\
\hline
\,\ \ \begin{tikzpicture}[iQ,centerzero]
\draw[-] (0,-0.3) \botlabel{i} -- (0,0.3) \toplabel{i};
\closeddot{0,0};
\region{0.2,0}{\lambda};
\end{tikzpicture}\ \, 
\colon B_i \one_\lambda \Rightarrow B_i \one_\lambda&2\\
\:\begin{tikzpicture}[iQ,centerzero]
\draw[-] (-0.25,-0.15) \botlabel{\tau i} to [out=90,in=90,looseness=3](0.25,-0.15) \botlabel{i};
\region{0.45,0.1}{\lambda};
\node at (0,.3) {$\phantom.$};
\node at (0,-.4) {$\phantom.$};
\end{tikzpicture}\hspace{-.03pt}
\colon B_{\tau i} B_{i} \one_\lambda \Rightarrow \one_\lambda
&1+\del_i-\lambda_i\\
\:\begin{tikzpicture}[iQ,centerzero]
\draw[-] (-0.25,0.15) \toplabel{\tau i} to[out=-90,in=-90,looseness=3] (0.25,0.15) \toplabel{i};
\region{0.45,-.1}{\lambda};
\node at (0,.2) {$\phantom.$};\node at (0,-.3) {$\phantom.$};
\end{tikzpicture}
\colon \one_\lambda \Rightarrow B_{\tau i} B_i \one_\lambda
&1+\del_i-\lambda_i\\
\!\!\!\begin{tikzpicture}[baseline=-1mm,iQ]
\draw[-] (-0.25,0) arc(180:-180:0.25);
\node at (-0.42,0) {\strandlabel{\tau i}};
\region{0.45,0}{\lambda};
\dottybubblelabel{0,0}{n};
\node at (0,-.3) {$\phantom{.}$};
\end{tikzpicture}
\colon \one_\lambda \Rightarrow\one_\lambda
\text{ for }0 \leq n \leq \del_i-\lambda_i
&2 n\\
\;\begin{tikzpicture}[iQ,centerzero,scale=.9]
\draw[-] (-0.3,-0.3) \botlabel{i} -- (0.3,0.3) \toplabel{i};
\draw[-] (0.3,-0.3) \botlabel{j} -- (-0.3,0.3) \toplabel{j};
\region{0.35,0}{\lambda};
\end{tikzpicture}
\;\,\colon B_{i} B_{j} \one_\lambda \Rightarrow B_j B_i \one_\lambda
&-a_{i,j}\\
\hline
\end{array}
\end{align*}
\caption{Generating 2-morphisms of $\UU^\imath$}\label{table0}
\end{table}
The generating 2-morphisms are subject to certain relations.
We write these down using a generating function formalism which has proven very useful in previous work, for example, \cite{HKM, K0}; it will be explained fully later in the text. 
Briefly, our generating functions are formal Laurent series in a variable $u^{-1}$, and the variables $x,y,z$ denote dots on strings in order from left to right. 
We omit 1- and 2-cell labels which can be inferred from context, and we omit them altogether if we are writing something which is true for all possible labels.
There are also a couple of shorthands for important generating functions:
\begin{align*}
\begin{tikzpicture}[iQ,centerzero,scale=1]
\draw[-] (0,-0.3) -- (0,0.3);
\circled{0,0}{u};
\end{tikzpicture}
&:=
\begin{tikzpicture}[iQ,centerzero,scale=1]
\draw[-] (0,-0.3) -- (0,0.3);
\pin{0,0}{.8,0}{\frac{1}{u-x}};
\end{tikzpicture}\ =
\sum_{n \geq 0} 
\begin{tikzpicture}[iQ,centerzero,scale=1]
\draw[-] (0,-0.3) -- (0,0.3);
\multcloseddot{0,0}{east}{n};
\end{tikzpicture}
u^{-n-1} 
,\\
\begin{tikzpicture}[iQ,centerzero,scale=1]
\draw[-] (0,-0.3) -- (0,0.3);
\circledbar{0,0}{u};
\end{tikzpicture}
&:=
\begin{tikzpicture}[iQ,centerzero,scale=1]
\draw[-] (0,-0.3) -- (0,0.3);
\pin{0,0}{.8,0}{\frac{1}{u+x}};
\end{tikzpicture}\ =
\sum_{n \geq 0} 
\begin{tikzpicture}[iQ,centerzero,scale=1]
\draw[-] (0,-0.3) -- (0,0.3);
\multcloseddot{0,0}{east}{n};
\end{tikzpicture}
(-1)^n u^{-n-1} 
,\\
\begin{tikzpicture}[baseline=-1mm,iQ]
\draw[fakebubble] (-0.25,0) arc(180:-180:0.25);
\node at (-0.42,0) {\strandlabel{\tau i}};
\bubblelabel{0,0}{u};
\region{0.4,0}{\lambda};
\end{tikzpicture}
&:=
\begin{dcases}
-\tfrac{1}{2u} \id_{\one_\lambda}&\text{if $i  =\tau i$}\\
\sum_{n = 0}^{\del_i-\lambda_i}
\begin{tikzpicture}[baseline=-1mm,iQ]
\draw[-] (-0.25,0) arc(180:-180:0.25);
\node at (-0.42,0) {\strandlabel{\tau i}};
\region{0.4,0}{\lambda};
\dottybubblelabel{0,0}{n};
\end{tikzpicture}\ 
u^{\del_i-\lambda_i-n}&\text{if $i \neq \tau i$},
\end{dcases}\\
\begin{tikzpicture}[baseline=-1mm,iQ]
\draw[-] (-0.25,0) arc(180:-180:0.25);
\node at (-0.42,0) {\strandlabel{\tau i}};
\region{0.95,0}{\lambda};
\node at (.55,0) {$(u)$};
\end{tikzpicture}
&:= 
\zeta_i(\lambda)^{-1}\left(
\begin{tikzpicture}[baseline=-1mm,iQ]
\draw[-] (-0.25,0) arc(180:-180:0.25);
\node at (-0.42,0) {\strandlabel{\tau i}};
\region{0.6,0}{\lambda};
\circled{0.25,0}{u};
\end{tikzpicture}+
\begin{tikzpicture}[baseline=-1mm,iQ]
\draw[fakebubble] (-0.25,0) arc(180:-180:0.25);
\node at (-0.42,0) {\strandlabel{\tau i}};
\region{0.4,0}{\lambda};
\bubblelabel{0,0}{u};
\node at (0,.3) {$\phantom{.}$};
\node at (0,-.3) {$\phantom{.}$};
\end{tikzpicture}\right).
\end{align*}
Here are the defining relations:
\begin{align*}
\begin{tikzpicture}[baseline=-1mm,iQ]
\draw[-] (-0.25,0) arc(180:-180:0.25);
\node at (-0.42,0) {\strandlabel{\tau i}};
\region{0.95,0}{\lambda};
\node at (.55,0) {$(u)$};
\end{tikzpicture}&\in
u^{\del_i-\lambda_i} \id_{\one_\lambda}
+ u^{\del_i-\lambda_i-1} \End_{\UU^\imath}(\one_\lambda)\llbracket u^{-1}\rrbracket,\!\!\!\!\!\!\!\!\!\!\!\!\!\!\!\\
\begin{tikzpicture}[baseline=-1mm,iQ]
\draw[-] (-0.25,0) arc(180:-180:0.25);
\node at (-0.35,0) {\strandlabel{i}};
\node at (.55,0) {$(u)$};
\end{tikzpicture}\!\!
\begin{tikzpicture}[baseline=-1mm,iQ]
\draw[-] (-0.25,0) arc(180:-180:0.25);
\node at (-0.44,0) {\strandlabel{\tau i}};
\node at (.7,0) {$(-u)$};
\region{-.8,0}{\lambda};
\end{tikzpicture}
&= (-1)^{\del_i-\lambda_i} u^{-a_{i,\tau i}}\id_{\one_\lambda},
\\
\begin{tikzpicture}[centerzero,iQ]
\draw[-] (-0.25,0) arc(180:-180:0.25);
\node at (-.42,0) {\strandlabel{\tau i}};
\node at (.55,0) {$(u)$};
\draw[-] (1,-0.5)\botlabel{j} to (1,0.5);
\pin{1,0}{2.1,0}{(u-x)^{a_{i,j}}};
\end{tikzpicture}
&=\begin{tikzpicture}[centerzero,iQ]
\draw[-] (-0.7,-0.5)\botlabel{j} to (-0.7,0.5);
\pin{-.7,0}{-1.75,0}{(u+x)^{a_{\tau i, j}}};
\draw[-] (-0.25,0) arc(180:-180:0.25);
\node at (-0.42,0) {\strandlabel{\tau i}};
\node at (.55,0) {$(u)$};
\end{tikzpicture},&\!\!\!\!\!\!\!\!\!\!\!\!\!\!\!\!\!
\begin{tikzpicture}[iQ,baseline=-3pt,scale=1]
\draw[-] (-.3,.25)\toplabel{i} to [out=-90,in=-90,looseness=3] (.3,.25);
\closeddot{-.27,.03};
\end{tikzpicture}\ =
-\ 
\begin{tikzpicture}[iQ,baseline=-3pt,scale=1]
\draw[-] (-.3,.25)\toplabel{i} to [out=-90,in=-90,looseness=3] (.3,.25);
\closeddot{.27,.03};
\end{tikzpicture}
\ ,\quad
\begin{tikzpicture}[iQ,centerzero,scale=1]
\draw[-] (-.3,-.25) \botlabel{i} to [out=90,in=90,looseness=3] (.3,-.25);
\closeddot{-.27,-.03};
\end{tikzpicture}\ &=
-\begin{tikzpicture}[iQ,centerzero,scale=1]
\draw[-] (-.3,-.25)\botlabel{i} to [out=90,in=90,looseness=3] (.3,-.25);
\closeddot{.27,-.03};
\end{tikzpicture}\ ,
\end{align*}

\vspace{-4mm}

\begin{align*}
&
\begin{tikzpicture}[iQ,centerzero={0,-2pt},scale=1.25]
\draw[-] (-0.3,0.4) -- (-0.3,0) arc(180:360:0.15) arc(180:0:0.15) -- (0.3,-0.4) \botlabel{i};
\end{tikzpicture}
=
\begin{tikzpicture}[iQ,centerzero={0,-2pt},scale=1.25]
\draw[-] (0,-0.4)\botlabel{i} -- (0,0.4);
\end{tikzpicture}\ =
\begin{tikzpicture}[iQ,centerzero={0,-2pt},scale=1.25]
\draw[-] (-0.3,-0.4) \botlabel{i}-- (-0.3,0) arc(180:0:0.15) arc(180:360:0.15) -- (0.3,0.4);
\end{tikzpicture}\ ,&
\begin{tikzpicture}[anchorbase,scale=1,iQ]
\draw[-] (0,-0.5)\botlabel{i} to[out=up,in=180] (0.3,0.2) to[out=0,in=up] (0.45,0) to[out=down,in=0] (0.3,-0.2) to[out=180,in=down] (0,0.5);
\end{tikzpicture}
=
\left[\ 
\begin{tikzpicture}[anchorbase,scale=1,iQ]
\draw[-] (-0.8,-0.5)\botlabel{i} -- (-0.8,0.5);
\circledbar{-0.8,0}{u};
\draw[fakebubble] (-0.4,0) arc(180:-180:0.25);
\node at (-0.48,0) {\strandlabel{i}};
\bubblelabel{-.15,0}{u};
\end{tikzpicture}\ 
\right]_{u:-1}\!\!\!\!\!\!\!,\qquad\quad\!
\begin{tikzpicture}[iQ,baseline=-5pt,scale=1.1]
\draw[-] (-0.25,0.15)\toplabel{i}  to[out=-90,in=-90,looseness=3] (0.25,0.15);
\draw[-] (0.3,-0.4)to[out=up,in=down] (0,0.15)\toplabel{j};
\end{tikzpicture}=
\begin{tikzpicture}[iQ,baseline=-5pt,scale=1.1]
\draw[-] (-0.25,0.15) \toplabel{i} to[out=-90,in=-90,looseness=3] (0.25,0.15);
\draw[-] (-0.3,-0.4) to[out=up,in=down] (0,0.15)\toplabel{j};
\end{tikzpicture}\ ,\qquad\quad\:
\begin{tikzpicture}[iQ,baseline=.5pt,scale=1.1]
\draw[-] (-0.25,-0.15)  to[out=90,in=90,looseness=3] (0.25,-0.15)\botlabel{i};
\draw[-] (-0.3,0.4) \braiddown (0,-0.15)\botlabel{j};
\end{tikzpicture}
=\begin{tikzpicture}[iQ,baseline=.5pt,scale=1.1]
\draw[-] (-0.25,-0.15)  to[out=90,in=90,looseness=3] (0.25,-0.15)\botlabel{i};
\draw[-] (0.3,0.4) \braiddown (0,-0.15)\botlabel{j};
\end{tikzpicture}\ ,
\end{align*}

\vspace{-4.5mm}

\begin{align*}
\begin{tikzpicture}[iQ,centerzero,scale=1]
\draw[-] (-0.3,-0.3) \botlabel{i} -- (0.3,0.3);
\draw[-] (0.3,-0.3) \botlabel{j} -- (-0.3,0.3);
\closeddot{-0.15,-0.15};
\end{tikzpicture}
-
\begin{tikzpicture}[iQ,centerzero,scale=1]
\draw[-] (-0.3,-0.3) \botlabel{i} -- (0.3,0.3);
\draw[-] (0.3,-0.3) \botlabel{j} -- (-0.3,0.3);
\closeddot{0.15,0.15};
\end{tikzpicture}
&= \delta_{i,j}\ 
\begin{tikzpicture}[iQ,centerzero,scale=1]\draw[-] (-0.2,-0.3) \botlabel{i} -- (-0.2,0.3);
\draw[-] (0.2,-0.3) \botlabel{j} -- (0.2,0.3);
\end{tikzpicture}-
\delta_{i,\tau j}\ 
\begin{tikzpicture}[iQ,centerzero,scale=1]
\draw[-] (-0.2,-0.3) \botlabel{i} to [looseness=2,out=90,in=90] (0.2,-0.3);
\draw[-] (0.2,0.3) to [looseness=2,out=-90,in=-90] (-0.2,0.3)\toplabel{j};
\end{tikzpicture}=
\begin{tikzpicture}[iQ,centerzero,scale=1]
\draw[-] (-0.3,-0.3) \botlabel{i} -- (0.3,0.3);
\draw[-] (0.3,-0.3) \botlabel{j} -- (-0.3,0.3);
\closeddot{-0.15,0.15};
\end{tikzpicture}
- 
\begin{tikzpicture}[iQ,centerzero,scale=1]
\draw[-] (-0.3,-0.3) \botlabel{i} -- (0.3,0.3);
\draw[-] (0.3,-0.3) \botlabel{j} -- (-0.3,0.3);
\closeddot{0.15,-0.15};
\end{tikzpicture}\ ,\\
\begin{tikzpicture}[iQ,centerzero,scale=1.08]
\draw[-] (-0.2,-0.4) \botlabel{i} to[out=45,in=down] (0.15,0) to[out=up,in=-45] (-0.2,0.4);
\draw[-] (0.2,-0.4) \botlabel{j} to[out=135,in=down] (-0.15,0) to[out=up,in=225] (0.2,0.4);
\end{tikzpicture}
&=
\begin{tikzpicture}[iQ,centerzero,scale=1.26]
\draw[-] (-0.2,-0.3) \botlabel{i} -- (-0.2,0.3);
\draw[-] (0.1,-0.3) \botlabel{j} -- (0.1,0.3);
\pinpin{0.1,0}{-0.2,0}{-.9,0}{Q^\imath_{i,j}(x,y)};
\end{tikzpicture}
+
\delta_{i,\tau j}
\left[\ \begin{tikzpicture}[iQ,centerzero,scale=1.26]
\draw[-] (-0.2,-0.3) \botlabel{i} to [looseness=2.2,out=90,in=90] (0.2,-0.3);
\draw[-] (0.2,0.3) to [looseness=2.2,out=-90,in=-90] (-0.2,0.3)\toplabel{i};
\circled{-.17,.14}{u};
\circled{-.17,-.14}{u};
\draw[fakebubble] (0.4,0) arc(180:-180:0.2);
\node at (0.32,0) {\strandlabel{j}};
\bubblelabel{.6,0}{u};
\end{tikzpicture}\ \right]_{u:-1}\!\!,
\\
\begin{tikzpicture}[iQ,centerzero,scale=1.1]
\draw[-] (-0.4,-0.4) \botlabel{i} -- (0.4,0.4);
\draw[-] (0,-0.4) \botlabel{j} to[out=135,in=down] (-0.32,0) to[out=up,in=225] (0,0.4);
\draw[-] (0.4,-0.4) \botlabel{k} -- (-0.4,0.4);
\end{tikzpicture}
-
\begin{tikzpicture}[iQ,centerzero,scale=1.1]
\draw[-] (-0.4,-0.4) \botlabel{i} -- (0.4,0.4);
\draw[-] (0,-0.4) \botlabel{j} to[out=45,in=down] (0.32,0) to[out=up,in=-45] (0,0.4);
\draw[-] (0.4,-0.4) \botlabel{k} -- (-0.4,0.4);
\end{tikzpicture}
&= \delta_{i,k} 
\begin{tikzpicture}[iQ,centerzero,scale=1.1]
\draw[-] (-0.28,-0.4) \botlabel{i} -- (-0.28,0.4);
\draw[-] (0,-0.4) \botlabel{j} -- (0,0.4);
\draw[-] (0.28,-0.4) \botlabel{k} -- (0.28,0.4);
\pinpinpin{-.28,0}{0,0}{.28,0}{1.7,0}{
\frac{Q^\imath_{i,j}(x,y)-Q^\imath_{i,j}(z,y)}{x-z}};
\end{tikzpicture}
+
\delta_{i,\tau j}\delta_{j,\tau k}\!
\left[\begin{tikzpicture}[iQ,centerzero,scale=1.4]
\draw[-] (-0.2,-0.3) \botlabel{i} to [looseness=2.2,out=90,in=90] (0.2,-0.3);
\draw[-] (0.2,0.3) to [looseness=2.2,out=-90,in=-90] (-0.2,0.3)\toplabel{i};
\draw[-] (.85,-.3) \botlabel{k} to (.85,.3);
\circled{-.17,.14}{u};
\circled{-.17,-.14}{u};
\circled{.85,0}{u};
\draw[fakebubble] (0.33,0) arc(180:-180:0.16);
\node at (0.26,0) {\strandlabel{j}};
\bubblelabel{.49,0}{u};
\end{tikzpicture}
-
\begin{tikzpicture}[iQ,centerzero,scale=1.4]
\draw[-] (0.3,-0.3) \botlabel{j} to [looseness=2.2,out=90,in=90] (0.7,-0.3);
\draw[-] (0.7,0.3) to [looseness=2.2,out=-90,in=-90] (0.3,0.3)\toplabel{j};
\draw[-] (-.4,-.3) \botlabel{i} to (-.4,.3);
\circledbar{-.4,0}{u};
\circledbar{.65,.14}{u};
\circledbar{.65,-.14}{u};
\draw[fakebubble] (-0.14,0) arc(180:-180:0.16);
\node at (-0.2,0) {\strandlabel{i}};
\bubblelabel{.02,0}{u};
\end{tikzpicture}
\right]_{\!u:-1}\notag\end{align*}

\vspace{-7mm}
$$
\hspace{65mm}
-\delta_{i,\tau j} 
\begin{tikzpicture}[iQ,centerzero,scale=1.4]
\draw[-] (-0.3,0.3)  -- (0.3,-0.3)\botlabel{k};
\draw[-] (0,-0.3) [out=90,in=90,looseness=2.5] to  (-0.3,-0.3)\botlabel{i};
\draw[-] (0.3,0.3) [out=-90,in=-90,looseness=2.5] to (0,0.3) \toplabel{j};
\pinpinpinmid{-.08,-.1}{0,0}{.08,.1}{1.35,0}{
\frac{Q^\imath_{j,k}(x,y)-Q^\imath_{j, k}(z,y)}{x-z}};
\end{tikzpicture}
-\delta_{j, \tau k}
\begin{tikzpicture}[iQ,centerzero,scale=1.4]
\draw[-] (0.3,0.3)  -- (-0.3,-0.3)\botlabel{i};
\draw[-] (0,-0.3) \botlabel{j} [out=90,in=90,looseness=2.5] to  (0.3,-0.3);
\draw[-] (-0.3,0.3) \toplabel{k}[out=-90,in=-90,looseness=2.5] to (0,0.3) ;
\pinpinpinmid{.08,-.1}{0,0}{-.08,.1}{1.35,0}{
\frac{Q^\imath_{j, i}(x,y)-Q^\imath_{j, i}(z,y)}{x-z}};
\end{tikzpicture}\ .
$$

The 2-quantum groups of Khovanov, Lauda and Rouquier are isomorphic to 2-iquantum groups of diagonal type. 
There are two other special cases already in the literature: the categorification of iquantum groups of quasi-split type AIII with an even number of nodes from \cite{BSWW}, and the nil-Brauer category which categorifies the split iquantum group of rank one from \cite{BWWa1,BWWbasis}.
The dictionaries to all of these special cases are explained in \cref{rankone,iexisting,bsww}.

\vspace{2mm}

The main results proved about this new family of 2-categories are as follows. 
\begin{itemize}
\item {\bf\cref{ind}}: We prove that $\UU^\imath$ is non-degenerate.
By this, we mean that the spanning sets for its 2-morphism spaces which are obtained by an obvious straightening rule are actually bases. This result generalizes the non-degeneracy of 2-quantum groups conjectured and proved for $\mathfrak{sl}_n$ by Khovanov and Lauda in \cite{KL3}, and proved in general in \cite{unfurling}.
\item {\bf\cref{bcor}}:
Then we explain how the graded ranks of 2-morphism spaces in $\UU^\imath$ can be computed using the non-degenerate symmetric bilinear form $(\cdot,\cdot)^\imath$ on the iquantum group $\dot\U^\imath$ defined in \cite{BW18QSP,BW21iCB}.
This extends \cite[Th.~2.7]{KL3}.
\item {\bf \cref{TBA}}: We classify indecomposables in
the appropriate idempotent completions
of the morphism categories in $\UU^\imath$.
This is done by equipping the path algebras of these categories
with {\em graded triangular bases} in the sense of \cite{GTB}, with
underlying Cartan algebras that are certain 
quiver Hecke algebras (see \cref{iQHA}).
These bases also lead to the definition of
{\em standard modules} for the morphism categories of $\UU^\imath$, which are standardizations of the indecomposable projective modules for these Cartan algebras.
\item 
{\bf\cref{cthm}}:
We prove that the split Grothendieck ring
of $\UU^\imath$ is isomorphic to 
the $\Z[q,q^{-1}]$-form
$\dot\U^\imath_\Z$ 
of $\dot\U^\imath$ generated by the divided/idivided powers $b_i^{(n)} 1_\lambda$.
This generalizes the known result for 2-quantum groups which was
proved assuming non-degeneracy by Khovanov and Lauda in \cite{KL3}.
A consequence is that isomorphism classes of 
indecomposable objects in the categorification give rise to a new basis for $\dot\U^\imath_\Z$
which we call the {\em iorthodox basis}.
Standard modules correspond to another basis, but over $\Z\lround
q\rround$ rather than over $\Z[q,q^{-1}]$, which we call the 
{\em standardized orthodox basis}.
\end{itemize}
In the main body of the text, we work in a more general setup than in this introduction, allowing Cartan matrices that are merely symmetrizable,
and more general choices of $Q_{i,j}(x,y)$.
We conjecture that $\UU^\imath$ is non-degenerate in the general setup.
However, at present, we are only able to prove it for ``geometric'' parameters as in \cref{tojapan}; this is more general than the situation described in the introduction but still only applies to 
symmetric Cartan matrices. The other results referenced above are proved in complete generality but depend on non-degeneracy.

\vspace{2mm}

The main technical tool at the heart of the paper is a 2-functor $\Xi^\imath$ from $\UU^\imath$ to a localization of the 2-quantum group $\UU$ associated to $\U$; see \cref{boxingday}.
We view this as a categorification of the standard embedding $\U^\imath\hookrightarrow\U$.
The insights gained from the calculations establishing its existence were essential
when we were working out some of the finer details in the defining relations of $\UU^\imath$. 
The 2-functor $\Xi^\imath$ also plays a key role in the proof of non-degeneracy. We point out a special case of independent interest, which is a variant of this functor for ordinary 2-quantum groups related to the comultiplication $\Delta:\U \rightarrow \U\otimes \U$; see \cref{valentinesday}.
We use this in \cref{nondeg} to give a self-contained proof of the non-degeneracy of 2-quantum groups for any symmetrizable Cartan matrix and any choice of parameters, making this article independent of \cite{unfurling}.

Another important ingredient in our approach is an
isometry $\jmath$ between $\dot\U^\imath 1_\lambda$ and the algebra $\f$ 
that is the negative (or positive) half of $\U$; see \cref{jbversion} which extends \cite[Th.~2.8]{PBW} and \cite[Th.~2.1]{BWWa1}.
The iSerre relation is equivalent to the relation obtained by applying $\jmath^{-1}$ to the ordinary Serre relation in $\f$.
This point of view leads to a new proof of Balagović-Kolb-Letzter's relation; see \cref{another}. In the categorification,
the role of the isometry $\jmath$ is played by
the {\em standardization functor} arising from the theory of graded triangular bases. We use this to prove a categorical analog of the BKL relation; see \cref{wheniserre2}. Our proof is quite different and much shorter than the one in \cite{BSWW} which treated the special case $a_{i,\tau i}=-1$.
The other cases of the iSerre relation are categorified by a more explicit split exact complex, which is the same as in \cite[Cor.~7]{KL2} and \cite[Prop.~4.2]{Rou} when $i \neq \tau i$, but there are additional subtleties involving the nil-Brauer category in the difficult case $i = \tau i$; see \cref{wheniserre}.

\vspace{2mm}

Moving forward, a basic open problem is to find situations in which the iorthodox basis coincides with the icanonical basis for $\dot\U^\imath_\Z$ from \cite{BW18QSP, BW21iCB}. This is closely related to the positivity of the icanonical basis, which was proved already for finite and affine Satake types AIII in \cite{LiW18, AffineIIIMemoirs},
and conjectured for all locally finite symmetric quasi-split types in
\cite[\S9]{WangICM}. In view of this, it seems reasonable to conjecture that the iorthodox basis in characteristic 0 coincides with the icanonical basis for these types. This is already known to be the case
in finite ADE diagonal types \cite[Th.~8.7]{canonical}, for the iquantum group of split rank one \cite{BWWa1}, and for the iquantum group of quasi-split type AIII with an even number of nodes \cite[Th.~A]{BSWW}. 

The icanonical bases arising from iquantum groups of quasi-split AIII types are intimately connected to Kazhdan-Lusztig theory in types B/C/D. In view of this, the 2-iquantum groups of these types introduced in \cite{BSWW} and the present paper 
are widely expected to be related to the affine Brauer category introduced in \cite{RS}, and the affine Brauer algebras appearing in many previous works such as \cite{Naz,AMR,ES}. It would be interesting to find a conceptual explanation of this
which is similar to the way that 2-quantum groups of types $A_{p-1}^{(1)}$ and $A_\infty$ are related to Heisenberg categories in \cite{HKM}. 

\vspace{2mm}
\noindent
{\em Conventions.} We work over an $\N$-graded 
commutative ring
$\kk = \bigoplus_{n \geq 0} \kk_n$; for \cref{TBA,cthm} we need $\kk_0$ to be a field. We use $\kk^\times$ to denote the homogeneous units in $\kk_0$. 
Let $[n]_q$ denote the
quantum integer $\frac{q^n-q^{-n}}{q-q^{-1}}$,  $[n]_q^!$ be the quantum factorial, and $\qbinom{m}{n}_q$ be the quantum binomial coefficient.
The bar involution on $\Q(q)$ is the field automorphism 
$f(q) \mapsto \overline{f(q)} := f(q^{-1})$. 
We also use $q$ to denote the {\em upward} grading shift functor on the monoidal category of graded $\kk$-modules, i.e., $(qV)_n = V_{n-1}$. A free graded $\kk$-module $V$
is {\em locally finite-dimensional and bounded below}
if it is isomorphic to 
$\kk^{\oplus f(q)} := \bigoplus_{i \in \Z} (q^i \kk)^{\oplus n_i}$
for a formal Laurent series 
$f(q) = \sum_{i \in \Z} n_i q^i \in \N\lround q\rround$. In this situation, we denote $f(q)$ by $\grrank_q V$ (or $\grdim_q V$ if $\kk$ is a field).

By a {\em graded category}, we mean a category enriched in graded $\kk$-modules. 
A {\em graded 2-category} is a category enriched in graded categories.
If $\catC$ is a graded category, the dual graded category $\catC^\op$ has the same objects as $\catC$. For objects $X$ and $Y$ in $\catC$, the morphism space $\Hom_{\catC^\op}(X,Y)$
is a copy $\{f^\op\:|\:f \in \Hom_{\catC}(Y,X)\}$ of $\Hom_{\catC}(Y,X)$ with
$\deg(f^\op) = \deg(f)$ and $f^\op \circ g^\op = (g \circ f)^\op$.
The {\em opposite} $\CATC^\op$ of a graded 2-category $\CATC$ is the graded 2-category defined by reversing vertical composition (its morphism categories are the duals of the morphism categories in $\CATC$). The
{\em reverse} $\CATC^\rev$ is defined by reversing horizontal composition (its morphism categories are the same as in $\CATC$).

\vspace{2mm}
\noindent
{\em Acknowledgements.} The authors would like to thank Alistair Savage and Liron Speyer for helpful discussions.

%% file: s2-quantum.tex
\setcounter{section}{1}

\section{Reminders about the categorification of quantum groups}\label{s2-quantum}

In this section, we fix a symmetrizable Cartan datum, briefly recall the definition of the corresponding quantized enveloping algebra and Lusztig's modified integral form.
Then we write down the definition of the 2-quantum group introduced by Khovanov, Lauda and Rouquier \cite{KL3,Rou}, and record some further relations working in terms of generating functions.

\subsection{The quantum group \texorpdfstring{$\U$}{}}\label{data}

Let $A = (a_{i,j})_{i,j \in I}$ be a symmetrizable generalized Cartan matrix. This means that $a_{i,i} = 2$ for all $i \in I$,
$a_{i,j} \in -\N$ for $i\neq j$ in $I$,
$a_{i,j} = 0\Leftrightarrow a_{j,i}=0$, and
there are given positive integers $d_i\:(i \in I)$
such that $d_i a_{i,j} = d_j a_{j,i}$ 
for all $i,j \in I$.
We do not insist that the set $I$ is finite, but the number of non-zero entries in each row and column of the Cartan matrix should be finite.
Let $X$, the {\em weight lattice},
and $Y$, the {\em coweight lattice},
be free Abelian groups
with a given perfect pairing
$Y\times X \rightarrow \Z, (h,\lambda) \mapsto h(\lambda)$.
We assume that $X$ contains
elements $\alpha_i\:(i \in I)$, called {\em simple roots}, and $Y$ contains elements
$h_i\:(i \in I)$, called
{\em simple coroots},
such that $h_i(\alpha_j) = a_{i,j}$ for all $i,j \in I$.
We will assume that the simple coroots are linearly independent (``$Y$-regular'') but do 
not require linear independence of the simple roots.
Let 
\begin{equation}\label{dominantweights}
X^+ := \{\lambda \in X\:|\:h_i(\lambda) \geq 0\text{ for all }i \in I\text{, and }h_i(\lambda)=0\text{ for all but finitely many $i$}\}
\end{equation}
be the set of {\em dominant weights}.

The {\em quantized enveloping algebra}
$\U$ attached to this root datum 
is a $\Q(q)$-algebra generated by elements
$(q^h)_{h \in Y}$ and
$(e_i, f_i)_{i \in I}$
satisfying the usual relations. 
In particular, $q^h q^{h'} = q^{h+h'}$, and
$$
[e_i, f_j] = \delta_{i, j}\frac{k_i-k_i^{-1}}{q_i-q_i^{-1}}
$$
where 
$k_i := q^{d_i h_i}$ and $q_i := q^{d_i}$.
Our notational choices for $\U$ follow \cite{Lubook} closely, the only departure being the use of $e_i, f_i, k_i$ for its generators, which Lusztig denotes by
$E_i, F_i, \widetilde{K}_i$. We do this so
that we can use the upper case letters 
to denote corresponding {\em functors} when we pass to the categorification of $\U$.
We view $\U$ as a Hopf algebra with the comultiplication $\Delta$ from \cite{Lubook}:
\begin{align}
\Delta(e_i) &= e_i \otimes 1 + k_i \otimes e_i,
&\Delta(f_i) &= 1 \otimes f_i + f_i \otimes k_i^{-1},
&\Delta(q^h) &= q^h \otimes q^h.
\end{align}

Actually, more natural for categorification
is Lusztig's {\em modified form}
$\dot\U$. 
This is the locally unital $\Q(q)$-algebra 
$\bigoplus_{\lambda, \mu \in X} 1_\lambda \U 1_{\mu}$
with the (mutually orthogonal) distinguished idempotents
$\{1_\lambda\:|\:\lambda \in X\}$ and
generators
\begin{align*}
e_i 1_\lambda &= 1_{\lambda+\alpha_i} e_i,&
f_i 1_\lambda &= 1_{\lambda-\alpha_i} f_i
\end{align*}
for $\lambda \in X$ and $i \in I$.
The relations in $\dot\U$ are derived from the ones in $\U$
using the rule that 
$q^h 1_\lambda = 1_\lambda q^h = q^{h(\lambda)} 1_\lambda$
for $\lambda \in X, h \in Y$. In particular,
$k_i 1_\lambda = q_i^{h_i(\lambda)} 1_\lambda$.
For example, the relation displayed in the previous paragraph becomes
$[e_i, f_j] 1_\lambda
= \delta_{i,j} [h_i(\lambda)]_{q_i} 1_\lambda$,
interpreting the commutator $[e_i, f_j] 1_\lambda$ as $e_i f_j 1_\lambda - f_j e_i 1_\lambda$.
There is also a natural completion
$\widehat{\U}$ of $\dot\U$
consisting of matrices $(u_{\lambda,\mu})_{\lambda,\mu \in X}$ with $u_{\lambda,\mu} \in 1_\lambda \U 1_{\mu}$
such that there are only finitely many non-zero entries in each row and column.
We identify $\dot\U$ with a subalgebra of $\widehat{\U}$
so that $1_\lambda u 1_\mu$ is the matrix
with this as its $(\lambda,\mu)$-entry
and 0 in all other positions.
A general element
$(u_{\lambda,\mu})_{\lambda,\mu \in X} \in \widehat{\U}$
may also be denoted as the infinite sum
$\sum_{\lambda,\mu \in X} u_{\lambda,\mu}$.

Finally, we pass to $\Z[q,q^{-1}]$-forms.
Let $\dot\U_\Z$ be the integral form for $\dot\U$
generated by the divided powers $e_i^{(n)} 1_\lambda := e_i^n 1_\lambda\, \big/ \, [n]_{q_i}^!$ and
$f_i^{(n)}1_\lambda :=
f_i^n 1_\lambda\,\big/\, [n]_{q_i}^!$ for all $\lambda \in X, i \in I$ and $n \geq 0$.
There are some useful symmetries:
\begin{align}
\label{omegainv}
\text{(linear involution)}\quad
\omega:\dot\U_\Z&\stackrel{\sim}{\rightarrow}\dot\U_\Z,
&
e_i^{(n)} 1_\lambda&\mapsto
f_i^{(n)} 1_{-\lambda},
&
f_i^{(n)} 1_\lambda&\mapsto
e_i^{(n)} 1_{-\lambda},
\\
\text{(anti-linear involution)}\quad
\psi:\dot\U_\Z &\stackrel{\sim}{\rightarrow} \dot\U_\Z,
&e_i^{(n)} 1_\lambda &\mapsto e_i^{(n)} 1_{\lambda},
&f_i^{(n)} 1_\lambda &\mapsto f_i^{(n)} 1_{\lambda},
\label{psiinv}
\\
\text{(linear anti-involution)}\quad
\sigma:\dot\U_\Z&\stackrel{\sim}{\rightarrow}\dot\U_\Z,
&
e_i^{(n)} 1_\lambda&\mapsto
1_{-\lambda} e_{i}^{(n)}, &
f_i^{(n)} 1_\lambda&\mapsto
1_{-\lambda} f_{i}^{(n)}.\label{sigmainv}
\end{align}
Let $\widehat{\U}_\Z$ be the completion of $\dot\U_\Z$ defined like in the previous paragraph.

\subsection{The 2-quantum group \texorpdfstring{$\UU$}{}}\label{existing}

In addition to the root datum just chosen, we need some additional parameters.
Assume that $\kk$ is an $\N$-graded commutative ring
as in {\em Conventions}.
For $i, j$ in $I$,
let $Q_{i,j}(x,y) \in \kk[x,y]$ be 
a polynomial which is $0$ if $i = j$, 
and which is
homogeneous of degree $-2 d_i a_{i,j}$ if $x$ and $y$ are assigned the
degrees $2d_i$ and $2 d_j$, respectively. We assume that
\begin{equation}\label{hermitian}
Q_{i,j}(x,y) = Q_{j,i}(y,x).
\end{equation}
When $i \neq j$, we have that
\begin{equation}\label{generic}
Q_{i,j}(x,y)
= t_{i,j} x^{-a_{i,j}}
+ \sum_{\substack{0 \leq r < -a_{i,j}\\0 \leq s < -a_{j,i}}} t_{i,j;r,s} x^r y^s + t_{j,i}
y^{-a_{j,i}}
\end{equation}
for $t_{i,j}, t_{j,i} \in \kk_0$ and
$t_{i,j;r,s}=t_{j,i;s,r} \in \kk_{-2 d_i a_{i,j} - 2d_i r - 2d_j s}$.
We assume moreover that
$t_{i,j}$ is invertible, 
so it is an element of $\kk^\times$, and set
$t_{i,i} := 1$.

\begin{defin}[2-quantum groups]\label{def2qg}
The {\em 2-quantum group}
$\UU$ with these parameters is the
graded 2-category 
with object set $X$, generating $1$-morphisms
\begin{align}
E_i \one_\lambda =\one_{\lambda+\alpha_i} E_i &\colon \lambda \to \lambda + \alpha_i,&
\one_{\lambda-\alpha_i} F_i=  F_i \one_{\lambda} &\colon \lambda \to \lambda - \alpha_i
\end{align}
for $\lambda \in X$ and $i \in I$,
with identity 2-endomorphisms represented by oriented strings
$\begin{tikzpicture}[Q,anchorbase]
\draw[-to] (0,-0.2)\botlabel{i} -- (0,0.2);
\region{0.2,0}{\lambda};
\end{tikzpicture}$ and
$\begin{tikzpicture}[Q,anchorbase]
\draw[to-] (0,-0.2) \botlabel{i}-- (0,0.2);
\region{0.2,0.1}{\lambda};
\end{tikzpicture}$,
and the four families of generating 2-morphisms
displayed in the top half of \cref{table1}.
\begin{table}
\begin{align*}
\begin{array}{|l|c|}
\hline
\hspace{13mm}\text{Generator}&\text{Degree}\\
\hline
\ \ \begin{tikzpicture}[Q,centerzero]
\draw[-to] (0,-0.3) \botlabel{i} -- (0,0.3) \toplabel{i};
\opendot{0,0};
\region{0.2,0}{\lambda};
\end{tikzpicture}\ 
\colon E_i \one_\lambda \Rightarrow E_i \one_\lambda&2 d_i\\
\,\begin{tikzpicture}[Q,centerzero,scale=.9]
\draw[-to] (-0.3,-0.3) \botlabel{i} -- (0.3,0.3) \toplabel{i};
\draw[-to] (0.3,-0.3) \botlabel{j} -- (-0.3,0.3) \toplabel{j};
\region{0.35,0}{\lambda};
\end{tikzpicture}
\;\colon E_i E_j \one_\lambda \Rightarrow E_j E_i \one_\lambda
&-d_i a_{i,j}\\
\:\begin{tikzpicture}[Q,centerzero]
\draw[-to] (-0.25,-0.15) \botlabel{i} to [out=90,in=90,looseness=3](0.25,-0.15) \botlabel{i};
\region{0.45,0.1}{\lambda};
\node at (0,.3) {$\phantom.$};
\node at (0,-.4) {$\phantom.$};
\end{tikzpicture}
\colon E_i F_i \one_\lambda \Rightarrow \one_\lambda
&d_i(1-h_i(\lambda))\\
\:\begin{tikzpicture}[Q,centerzero]
\draw[-to] (-0.25,0.15) \toplabel{i} to[out=-90,in=-90,looseness=3] (0.25,0.15) \toplabel{i};
\region{0.45,-.1}{\lambda};
\node at (0,.2) {$\phantom.$};\node at (0,-.3) {$\phantom.$};
\end{tikzpicture}
\colon \one_\lambda \Rightarrow F_i E_i \one_\lambda
&d_i(1+h_i(\lambda))
\\
\hline
\ \ \begin{tikzpicture}[Q,centerzero]
\draw[to-] (0,-0.3) \botlabel{i} -- (0,0.3) \toplabel{i};
\opendot{0,0};
\region{0.2,0}{\lambda};
\end{tikzpicture}\ 
\colon F_i \one_\lambda \Rightarrow F_i \one_\lambda&2 d_i\\
\,\begin{tikzpicture}[Q,centerzero,scale=.9]
\draw[to-] (-0.3,-0.3) \botlabel{i} -- (0.3,0.3) \toplabel{i};
\draw[to-] (0.3,-0.3) \botlabel{j} -- (-0.3,0.3) \toplabel{j};
\region{0.35,0}{\lambda};
\end{tikzpicture}
\;\colon F_i F_j \one_\lambda \Rightarrow F_j F_i \one_\lambda
&-d_i a_{i,j}\\
\,\begin{tikzpicture}[Q,centerzero,scale=.9]
\draw[-to] (-0.3,-0.3) \botlabel{i} -- (0.3,0.3) \toplabel{i};
\draw[to-] (0.3,-0.3) \botlabel{j} -- (-0.3,0.3) \toplabel{j};
\region{0.35,0}{\lambda};
\end{tikzpicture}
\;\colon E_i F_j \one_\lambda \Rightarrow F_j E_i \one_\lambda
&0\\
\,\begin{tikzpicture}[Q,centerzero,scale=.9]
\draw[to-] (-0.3,-0.3) \botlabel{i} -- (0.3,0.3) \toplabel{i};
\draw[-to] (0.3,-0.3) \botlabel{j} -- (-0.3,0.3) \toplabel{j};
\region{0.35,0}{\lambda};
\end{tikzpicture}
\;\colon F_i E_j \one_\lambda \Rightarrow E_j F_i \one_\lambda
&0\\
\:\begin{tikzpicture}[Q,centerzero]
\draw[to-] (-0.25,-0.15) \botlabel{i} to [out=90,in=90,looseness=3](0.25,-0.15) \botlabel{i};
\region{0.45,0.1}{\lambda};
\node at (0,.3) {$\phantom.$};
\node at (0,-.4) {$\phantom.$};
\end{tikzpicture}
\colon F_i E_i \one_\lambda \Rightarrow \one_\lambda
&d_i(1+h_i(\lambda))\\
\:\begin{tikzpicture}[Q,centerzero]
\draw[to-] (-0.25,0.15) \toplabel{i} to[out=-90,in=-90,looseness=3] (0.25,0.15) \toplabel{i};
\region{0.45,-.1}{\lambda};
\node at (0,.2) {$\phantom.$};\node at (0,-.3) {$\phantom.$};
\end{tikzpicture}
\colon \one_\lambda \Rightarrow E_i F_i \one_\lambda
&d_i(1-h_i(\lambda))\\
\hline
\end{array}
\end{align*}
\caption{Generating 2-morphisms of $\UU$}\label{table1}
\end{table}
The degrees of the generating 2-morphisms are also recorded in this table.
The 2-morphisms in the bottom part of \cref{table1}
are not needed for the initial definition; they
will be introduced in the next subsection.
The generating 2-morphisms are subject to relations still to be explained. 
Before writing these down, we explain some conventions.
\begin{itemize}
\item
We will usually only label one of the 2-cells in a string diagram with a weight---the others can then be worked out implicitly. When we omit
{\em all} labels in 2-cells, it should be understood that we are discussing something that holds
for all possible labels.
\item
We will label strings just at one end. If we omit a label or orientation, it means that we are discussing something that holds for all possibilities.
\item
When a dot is labelled by a multiplicity,
we mean to take its power under vertical composition.
For a polynomial
$f(x) = \sum_{r=0}^n c_r x^r$, we use the
shorthand
\begin{equation*}\label{singlepin}
    \begin{tikzpicture}[Q,centerzero]
        \draw[-] (0,-0.3) -- (0,0.3);
             \Pin{0,0}{-.7,0}{f(x)};
    \end{tikzpicture}
    = \begin{tikzpicture}[Q,centerzero]
        \draw[-] (0,-0.3) -- (0,0.3);
             \Pin{0,0}{.7,0}{f(x)};
    \end{tikzpicture}\
    :=
    \sum_{r=0}^n c_r\
    \begin{tikzpicture}[Q,centerzero]
        \draw[-] (0,-0.3) -- (0,0.3);
        \multopendot{0,0}{west}{r};
    \end{tikzpicture}
\end{equation*}
to ``pin'' $f(x)$ to a dot on a string.
Similarly, for
$f(x,y) = \sum_{r=0}^n\sum_{s=0}^m c_{r,s} x^r y^s$, we use
\begin{align*}
\begin{tikzpicture}[Q,centerzero]
\draw[-] (0,-0.3) -- (0,0.3);
\draw[-] (0.4,-0.3) -- (0.4,0.3);
\Pinpin{.4,0}{0,0}{-.8,0}{f(x,y)};
\end{tikzpicture}
=
\begin{tikzpicture}[Q,centerzero]
\draw[-] (0,-0.3) -- (0,0.3);
\draw[-] (0.4,-0.3) -- (0.4,0.3);
\Pinpin{0,0}{0.4,0}{1.2,0}{f(x,y)};
\end{tikzpicture}\ 
& :=
\sum_{r=0}^n\sum_{s=0}^m c_{r,s}
\begin{tikzpicture}[Q,centerzero]
\draw[-] (0,-0.3) -- (0,0.3);
\draw[-] (0.4,-0.3) -- (0.4,0.3);
\multopendot{.4,0}{west}{s};
\multopendot{0,0}{east}{r};
\end{tikzpicture}
\end{align*}
This notation extends to polynomials $f(x,y,z)$ in three variables pinned to three dots, with
the convention that the variables in alphabetic order correspond to the dots ordered by the lexicographic order on their Cartesian coordinates.
Thus, $x$ corresponds to the leftmost dot, and the lowest one if there are
several such dots in the same vertical line.
\item
We use the following shorthand
to denote the composite 2-morphism obtained by ``rotating'' the generating 2-morphism:
\begin{align}
\begin{tikzpicture}[Q,centerzero]
\draw[to-] (0.3,-0.3) \botlabel{j} -- (-0.3,0.3);
\draw[-to] (-0.3,-0.3) \botlabel{i} -- (0.3,0.3);
\region{0.4,0}{\lambda};        
\end{tikzpicture}\ 
&:=
\begin{tikzpicture}[Q,centerzero,scale=1.2]
\draw[-to] (0.1,-0.3)\botlabel{i} \braidup (-0.1,0.3);
\draw[-to] (-0.4,0.3) -- (-0.4,0.1) to[out=down,in=left] (-0.2,-0.2) to[out=right,in=left] (0.2,0.2) to[out=right,in=up] (0.4,-0.1)  -- (0.4,-0.3)\botlabel{j};
\region{0.6,0}{\lambda};        
\end{tikzpicture}
\ .\label{rightpivot}
\end{align}
\end{itemize}
Now for the relations.
There are three families.
First, we have
the \emph{quiver Hecke algebra} relations:
\begin{align}
\begin{tikzpicture}[Q,centerzero]
\draw[-to] (-0.3,-0.3) \botlabel{i} -- (0.3,0.3);
\draw[-to] (0.3,-0.3) \botlabel{j} -- (-0.3,0.3);
\opendot{-0.15,-0.15};
\end{tikzpicture}
-
\begin{tikzpicture}[Q,centerzero]
\draw[-to] (-0.3,-0.3) \botlabel{i} -- (0.3,0.3);
\draw[-to] (0.3,-0.3) \botlabel{j} -- (-0.3,0.3);
\opendot{0.15,0.15};
\end{tikzpicture}
&= \delta_{i,j}  \ 
\begin{tikzpicture}[Q,centerzero]
\draw[-to] (-0.2,-0.3) \botlabel{i} -- (-0.2,0.3);
\draw[-to] (0.2,-0.3) \botlabel{i} -- (0.2,0.3);
\end{tikzpicture} =    \begin{tikzpicture}[Q,centerzero]
\draw[-to] (-0.3,-0.3) \botlabel{i} -- (0.3,0.3);
\draw[-to] (0.3,-0.3) \botlabel{j} -- (-0.3,0.3);
\opendot{-0.15,0.15};
\end{tikzpicture}
-
\begin{tikzpicture}[Q,centerzero]
\draw[-to] (-0.3,-0.3) \botlabel{i} -- (0.3,0.3);
\draw[-to] (0.3,-0.3) \botlabel{j} -- (-0.3,0.3);
\opendot{0.15,-0.15};
\end{tikzpicture} \ ,\label{dotslide}\\\label{quadratic}
\begin{tikzpicture}[Q,centerzero,scale=1.1]
\draw[-to] (-0.2,-0.4) \botlabel{i} to[out=45,in=down] (0.15,0) to[out=up,in=-45] (-0.2,0.4);
\draw[-to] (0.2,-0.4) \botlabel{j} to[out=135,in=down] (-0.15,0) to[out=up,in=225] (0.2,0.4);
\end{tikzpicture}
&=
\begin{tikzpicture}[Q,centerzero,scale=1.1]
\draw[-to] (-0.2,-0.4) \botlabel{i} -- (-0.2,0.4);
\draw[-to] (0.2,-0.4) \botlabel{j} -- (0.2,0.4);
\Pinpin{0.2,0}{-0.2,0}{-1.1,0}{Q_{i,j}(x,y)};
\end{tikzpicture}\ ,\\\label{braid}
\begin{tikzpicture}[Q,centerzero,scale=1.1]
\draw[-to] (-0.4,-0.4) \botlabel{i} -- (0.4,0.4);
\draw[-to] (0,-0.4) \botlabel{j} to[out=135,in=down] (-0.32,0) to[out=up,in=225] (0,0.4);
\draw[-to] (0.4,-0.4) \botlabel{k} -- (-0.4,0.4);
\end{tikzpicture}
\ -\
\begin{tikzpicture}[Q,centerzero,scale=1.1]
\draw[-to] (-0.4,-0.4) \botlabel{i} -- (0.4,0.4);
\draw[-to] (0,-0.4) \botlabel{j} to[out=45,in=down] (0.32,0) to[out=up,in=-45] (0,0.4);
\draw[-to] (0.4,-0.4) \botlabel{k} -- (-0.4,0.4);
\end{tikzpicture}
&=\delta_{i,k}\ 
\begin{tikzpicture}[Q,centerzero,scale=1.1]
\draw[-to] (-0.3,-0.4) \botlabel{i} -- (-0.3,0.4);
\draw[-to] (0,-0.4) \botlabel{j} -- (0,0.4);
\draw[-to] (0.3,-0.4) \botlabel{i} -- (0.3,0.4);
\Pinpinpin{.3,0}{0,0}{-.3,0}{-1.7,0}{
\frac{Q_{i,j}(x,y)-Q_{i,j}(z,y)}{x-z}};
\end{tikzpicture}\ .
\end{align}
Next, the \emph{right adjunction relations}:
\begin{align}\label{rightadj}
\begin{tikzpicture}[Q,centerzero,scale=1.2]
\draw[-to] (-0.3,0.4) -- (-0.3,0) arc(180:360:0.15) arc(180:0:0.15) -- (0.3,-0.4) \botlabel{i};
\end{tikzpicture}
&=
\begin{tikzpicture}[Q,centerzero,scale=1.2]
\draw[to-] (0,-0.4) \botlabel{i} -- (0,0.4);
\end{tikzpicture}
\ ,&
\begin{tikzpicture}[Q,centerzero,scale=1.2]
\draw[-to] (-0.3,-0.4) \botlabel{i}-- (-0.3,0) arc(180:0:0.15) arc(180:360:0.15) -- (0.3,0.4);
\end{tikzpicture}
&=
\begin{tikzpicture}[Q,centerzero,scale=1.2]
\draw[-to] (0,-0.4) \botlabel{i} -- (0,0.4);
\end{tikzpicture}\ ,
\end{align}
Finally, we have the 
\emph{inversion relations} which assert that
\begin{equation}\label{thismatrix}
\begin{tikzpicture}[Q,centerzero]
\draw[-to] (-0.3,-0.3) \botlabel{i}-- (0.3,0.3);
\draw[to-] (0.3,-0.3) \botlabel{j} -- (-0.3,0.3);
\region{0.38,0.02}{\lambda};
\end{tikzpicture}:
E_i F_j \one_\lambda \Rightarrow F_j E_i \one_\lambda\end{equation}
is an isomorphism for all $\lambda \in X$ and $i \neq j$, as are the following matrices for all $\lambda$ and $i$:
\begin{equation}
M_{\lambda;i} :=
\begin{dcases}
\begin{pmatrix}
\begin{tikzpicture}[Q,centerzero]
\draw[-to] (-0.25,-0.25) \botlabel{i}-- (0.25,0.25);
\draw[to-] (0.25,-0.25) \botlabel{i} -- (-0.25,0.25);
\region{0.35,0.02}{\lambda};
\end{tikzpicture} &
\begin{tikzpicture}[Q,centerzero]
\draw[-to] (-0.25,0.15) \toplabel{i} to[out=-90,in=-90,looseness=3] (0.25,0.15);
\region{-0.45,-.1}{\lambda};
\node at (0,.2) {$\phantom.$};\node at (0,-.3) {$\phantom.$};
\end{tikzpicture}
&
\begin{tikzpicture}[Q,centerzero]
\draw[-to] (-0.25,0.15) \toplabel{i} to[out=-90,in=-90,looseness=3] (0.25,0.15);
\region{-0.45,-.1}{\lambda};
\node at (0,.2) {$\phantom.$};\node at (0,-.3) {$\phantom.$};
\opendot{0.23,-0.03};
\end{tikzpicture}
&\!\!\cdots\!\!
&
\begin{tikzpicture}[Q,centerzero]
\draw[-to] (-0.25,0.15) \toplabel{i} to[out=-90,in=-90,looseness=3] (0.25,0.15);
\region{-0.45,-.1}{\lambda};
\node at (0,.2) {$\phantom.$};\node at (0,-.3) {$\phantom.$};
\multopendot{0.23,-0.03}{west}{-h_i(\lambda)-1};
\end{tikzpicture}
\end{pmatrix}\phantom{_T}
&\text{if } h_i(\lambda) \leq 0\\
\begin{pmatrix}   
\begin{tikzpicture}[Q,centerzero]
\draw[-to] (-0.25,-0.25) \botlabel{i}-- (0.25,0.25);
\draw[to-] (0.25,-0.25)\botlabel{i}  -- (-0.25,0.25);
\region{0.38,0.02}{\lambda};
\end{tikzpicture} &
\begin{tikzpicture}[Q,centerzero]
\draw[-to] (-0.25,-0.15) \botlabel{i} to [out=90,in=90,looseness=3](0.25,-0.15);
\region{0.45,0.1}{\lambda};
\node at (0,.3) {$\phantom.$};
\node at (0,-.4) {$\phantom.$};
\end{tikzpicture}
&
\begin{tikzpicture}[Q,centerzero]
\draw[-to] (-0.25,-0.15) \botlabel{i} to [out=90,in=90,looseness=3](0.25,-0.15);
\region{0.45,0.1}{\lambda};
\node at (0,.3) {$\phantom.$};
\node at (0,-.4) {$\phantom.$};
\opendot{-0.23,.03};
\end{tikzpicture}
&
\!\!\!\cdots\!\!\!
&
\begin{tikzpicture}[Q,centerzero]
\draw[-to] (-0.25,-0.15) \botlabel{i} to [out=90,in=90,looseness=3](0.25,-0.15);
\region{0.45,0.1}{\lambda};
\node at (0,.3) {$\phantom.$};
\node at (0,-.4) {$\phantom.$};
\multopendot{-0.23,.03}{east}{h_i(\lambda)-1};
\end{tikzpicture}
\end{pmatrix}^\transpose\phantom{_T}
&\text{if } h_i(\lambda) \geq 0.
\end{dcases}
\end{equation}
\end{defin}

\begin{rem}
The relations \cref{dotslide} are equivalent to the formula
\begin{equation}\label{demazure}
\begin{tikzpicture}[Q,centerzero,scale=1.2]
\draw[-to] (-0.3,-0.3) \botlabel{i} -- (0.3,0.3);
\draw[-to] (0.3,-0.3) \botlabel{j} -- (-0.3,0.3);
\Pinpin{.15,.15}{-.15,.15}{-.8,.15}{f(x,y)};
\end{tikzpicture}
-
\begin{tikzpicture}[Q,centerzero,scale=1.3]
\draw[-to] (-0.3,-0.3) \botlabel{i} -- (0.3,0.3);
\draw[-to] (0.3,-0.3) \botlabel{j} -- (-0.3,0.3);
\Pinpin{.15,-.15}{-.15,-.15}{-.8,-.15}{f(y,x)};
\end{tikzpicture}
= \delta_{i,j}  \ 
\begin{tikzpicture}[Q,centerzero,scale=1.3]
\draw[-to] (-0.2,-0.3) \botlabel{i} -- (-0.2,0.3);
\draw[-to] (0.2,-0.3) \botlabel{i} -- (0.2,0.3);
\Pinpin{.2,0}{-.2,0}{-1.15,0}{\frac{f(x,y)-f(y,x)}{x-y}};
\end{tikzpicture}
\end{equation}
for any $f(x,y) \in \kk[x,y]$.
The fraction on the right hand side
is the image of $f(x,y)$ under the 
Demazure operator.
\end{rem}

\subsection{More generators and generating functions}\label{stoptheworld}

The presentation for $\UU$ just explained is due to Rouquier \cite{Rou}. In \cite{KL3},
Khovanov and Lauda found a different presentation with more generators and relations.
The two approaches were reconciled in \cite{Brundan}.
We will explain this here in a form which also incorporates a renormalization of the extra generators.
This depends on the additional choice of homomorphisms of Abelian groups
\begin{equation}\label{great}
c_i:X \rightarrow \kk^\times
\end{equation}
with $c_i(\alpha_i) = 1$ for all $i \in I$,
which we call {\em normalization homomorphisms}.
Let
\begin{align}\label{gollygosh}
r_{i,j} &:= \frac{c_i(\alpha_j)}{t_{i,j}}
\end{align}
for $i,j \in I$.
To match the conventions of \cite{Brundan}, take $c_i = 1$ for all $i$.
Another choice proposed in \cite{BHLW} is to choose
$c_i$ so that $c_i(\alpha_j) = t_{i,j}$.
It is very convenient when possible to find such functions because it implies that $r_{i,j} = 1$ for all $i,j \in I$, hence, by further relations recorded shortly, string diagrams are invariant 
under planar isotopy. Note we are requiring that $c_i$ is a group homomorphism, whereas \cite{BHLW} 
only assumes that $c_i(\lambda+\alpha) = c_i(\lambda)c_i(\alpha)$
for any $\lambda \in X$ and $\alpha$ in the root lattice.
It seems to us to be reasonable to impose this extra requirement, 
and it will be helpful later in the article.

\begin{eg}[Geometric parameters]\label{tojapan}
Suppose that $d_i = 1$ for all
$i \in I$. The Cartan matrix is symmetric. 
Let $Q$ be a quiver with vertex set $I$, no
loops, and $\#(i
\rightarrow j)$ directed edges from $i$ to $j$
such that $\#(i\rightarrow j) + \#(j \rightarrow i) = -a_{i,j}$ for
all $i \neq j$.
Then, for $i \neq j$, one can take
$$
Q_{i,j}(x,y) := (x-y)^{\#(i \rightarrow j)} (y-x)^{\#(j \rightarrow i)}.
$$
We call this a {\em geometric choice} of parameters for a symmetric
Cartan matrix. 
In the introduction, we worked with geometric parameters
for a quiver admitting an anti-involution $\tau$, i.e., an involution
$\tau$ of the vertex set such that
$\#(i \rightarrow j) = \#(\tau i \leftarrow \tau j)$, 
so that $Q_{\tau i,\tau j}(x,y) = Q_{i,j}(-x,-y)$.
\end{eg}

Now we define the additional 2-morphisms:
\begin{itemize}
\item
Let
$\begin{tikzpicture}[Q,centerzero]
\draw[-to] (0.25,-0.25) \botlabel{j} -- (-0.25,0.25);
\draw[to-] (-0.25,-0.25) \botlabel{i} -- (0.25,0.25);
\region{0.33,0}{\lambda};        
\end{tikzpicture}$
be
$\Big(\begin{tikzpicture}[Q,centerzero]
\draw[to-] (0.25,-0.25)\botlabel{i}  -- (-0.25,0.25);
\draw[-to] (-0.25,-0.25) \botlabel{j} -- (0.25,0.25);
\region{0.33,0}{\lambda};        
\end{tikzpicture}\Big)^{-1}$ if $i \neq j$, 
or the first entry of the matrix $-M_{\lambda;i}^{-1}$ if $i=j$.
\item
Let
$\begin{tikzpicture}[Q,centerzero]
\draw[to-] (-0.25,-0.15) \botlabel{i} to [out=90,in=90,looseness=3](0.25,-0.15);
\region{0.45,0.1}{\lambda};
\node at (0,.3) {$\phantom.$};
\node at (0,-.4) {$\phantom.$};
\end{tikzpicture}$
be the last entry of $c_i(\lambda)^{-1}
M_{\lambda;i}^{-1}$ if $h_i(\lambda) < 0$ or $-c_i(\lambda)^{-1}
 \begin{tikzpicture}[Q,centerzero,scale=1.2]
 \draw[to-] (-.2,-.3)\botlabel{i} to [out=90,in=-90,looseness=1] (.15,.15) to [out=90,in=90,looseness=1.7] (-.15,.15) to [out=-90,in=90,looseness=1] (.2,-.3);
 \multopendot{-0.14,0.15}{east}{h_i(\lambda)};
\region{0.4,0}{\lambda};
\end{tikzpicture}$ if $h_i(\lambda) \geq 0$.
\item Let $\begin{tikzpicture}[Q,centerzero]
\draw[to-] (-0.25,0.15) \toplabel{i} to[out=-90,in=-90,looseness=3] (0.25,0.15);
\region{0.45,-.1}{\lambda};
\node at (0,.2) {$\phantom.$};\node at (0,-.3) {$\phantom.$};
\end{tikzpicture}$
be the last entry of $c_i(\lambda)
M_{\lambda;i}^{-1}$ if $h_i(\lambda) > 0$ or
$c_i(\lambda) \begin{tikzpicture}[Q,centerzero,scale=1.2]
\draw[to-] (-.2,.3)\toplabel{i} to [out=-90,in=90,looseness=1] (.15,-.15) to [out=-90,in=-90,looseness=1.7] (-.15,-.15) to [out=90,in=-90,looseness=1] (.2,.3);
\multopendot{0.14,-0.15}{west}{-h_i(\lambda)};
\region{-0.4,0}{\lambda};
\end{tikzpicture}$ if $h_i(\lambda) \leq 0$.
\item
Let $\begin{tikzpicture}[Q,centerzero,scale=1.1]
\draw[to-] (0,-0.4)\botlabel{i} -- (0,0.4);
\opendot{0,0};
\region{0.2,0}{\lambda};
\end{tikzpicture}
:=
\begin{tikzpicture}[Q,centerzero,scale=1.1]
\draw[to-] (0.3,-0.4) \botlabel{i}-- (0.3,0) arc(0:180:0.15) arc(360:180:0.15) -- (-0.3,0.4);
\opendot{0,0};
\region{0.5,0}{\lambda};
\end{tikzpicture}=
 \begin{tikzpicture}[Q,centerzero,scale=1.1]
\draw[-to] (0.3,0.4) -- (0.3,0) arc(360:180:0.15) arc(0:180:0.15) -- (-0.3,-0.4)\botlabel{i};
\opendot{0,0};
\region{0.45,0}{\lambda};
\end{tikzpicture}$.
\item Let
$\begin{tikzpicture}[Q,centerzero]
\draw[to-] (0.3,-0.3) \botlabel{j}-- (-0.3,0.3) ;
\draw[to-] (-0.3,-0.3) \botlabel{i}-- (0.3,0.3) ;
\region{0.4,0}{\lambda};        
\end{tikzpicture} :=r_{i,j}\ 
\begin{tikzpicture}[Q,anchorbase,scale=.7]
\draw[-to] (1.3,.4) to (1.3,-1.2)\botlabel{j};
\draw[-] (-1.3,-.4) to (-1.3,1.2);
\draw[-] (.5,1.1) to [out=0,in=90,looseness=1] (1.3,.4);
\draw[-] (-.35,.4) to [out=90,in=180,looseness=1] (.5,1.1);
\draw[-] (-.5,-1.1) to [out=180,in=-90,looseness=1] (-1.3,-.4);
\draw[-] (.35,-.4) to [out=-90,in=0,looseness=1] (-.5,-1.1);
\draw[-] (.35,-.4) to [out=90,in=-90,looseness=1] (-.35,.4);
\draw[-] (-0.35,-.5) to[out=0,in=180,looseness=1] (0.35,.5);
\draw[-to] (0.35,.5) to[out=0,in=90,looseness=1] (0.8,-1.2)\botlabel{i};
\draw[-] (-0.35,-.5) to[out=180,in=-90,looseness=1] (-0.8,1.2);
\region{1.55,0}{\lambda};
\end{tikzpicture}=
r_{j,i}\ 
\begin{tikzpicture}[Q,anchorbase,scale=.7]
\draw[-to] (-1.3,.4) to (-1.3,-1.2) \botlabel{i};
\draw[-] (1.3,-.4) to (1.3,1.2);
\draw[-] (-.5,1.1) to [out=180,in=90,looseness=1] (-1.3,.4);
\draw[-] (.35,.4) to [out=90,in=0,looseness=1] (-.5,1.1);
\draw[-] (.5,-1.1) to [out=0,in=-90,looseness=1] (1.3,-.4);
\draw[-] (-.35,-.4) to [out=-90,in=180,looseness=1] (.5,-1.1);
\draw[-] (-.35,-.4) to [out=90,in=-90,looseness=1] (.35,.4);
\draw[-] (0.35,-.5) to[out=180,in=0,looseness=1] (-0.35,.5);
\draw[-to] (-0.35,.5) to[out=180,in=90,looseness=1] (-0.8,-1.2) \botlabel{j};
\draw[-] (0.35,-.5) to[out=0,in=-90,looseness=1] (0.8,1.2);
\region{1.55,0}{\lambda};
\end{tikzpicture}.$
\end{itemize}
The equalities in the definitions of the downward dot and crossing just given are by no means obvious; they follow from the next relations \cref{leftpivots,lotsmore,ruby,wax}. 
By \cite[Th.~4.3]{Brundan},
the leftward cups and caps satisfy the zig-zag relations:
\begin{align}\label{leftpivots}
\begin{tikzpicture}[centerzero,Q,scale=1.2]
\draw[to-] (-0.3,0.4) -- (-0.3,0) arc(180:360:0.15) arc(180:0:0.15) -- (0.3,-0.4)\botlabel{i};
\end{tikzpicture}
&=
\begin{tikzpicture}[centerzero,Q,scale=1.2]
\draw[-to] (0,-0.4)\botlabel{i} -- (0,0.4);
\end{tikzpicture}
\ ,&
\begin{tikzpicture}[centerzero,Q,scale=1.2]
\draw[to-] (-0.3,-0.4)\botlabel{i} -- (-0.3,0) arc(180:0:0.15) arc(180:360:0.15) -- (0.3,0.4);
\end{tikzpicture}
&=
\begin{tikzpicture}[centerzero,Q,scale=1.2]
\draw[to-] (0,-0.4) \botlabel{i} -- (0,0.4);
\end{tikzpicture}\ .
\end{align}
We also have that
\begin{align}\label{lotsmore}
\begin{tikzpicture}[Q,centerzero,scale=1.2]
\draw[-to] (-0.25,0.15) \toplabel{i} to[out=-90,in=-90,looseness=3] (0.25,0.15);
\opendot{-0.22,-0.05};
\end{tikzpicture}
&=
\begin{tikzpicture}[Q,centerzero,scale=1.2]
\draw[-to] (-0.25,0.15) \toplabel{i} to[out=-90,in=-90,looseness=3] (0.25,0.15);
\opendot{0.23,-0.05};
\end{tikzpicture}\ ,&
\begin{tikzpicture}[Q,centerzero,scale=1.2]
\draw[to-] (-0.25,0.15) \toplabel{i} to[out=-90,in=-90,looseness=3] (0.25,0.15);
\opendot{-0.22,-0.05};
\end{tikzpicture}
&=
\begin{tikzpicture}[Q,centerzero,scale=1.2]
\draw[to-] (-0.25,0.15) \toplabel{i} to[out=-90,in=-90,looseness=3] (0.25,0.15);
\opendot{0.23,-0.05};
\end{tikzpicture}\ ,&
  \begin{tikzpicture}[Q,centerzero,scale=1.2]
\draw[-to] (-0.25,-0.15) \botlabel{i} to[out=90,in=90,looseness=3] (0.25,-0.15);
\opendot{-0.22,0.05};
\end{tikzpicture}
&=
\begin{tikzpicture}[Q,centerzero,scale=1.2]
\draw[-to] (-0.25,-0.15) \botlabel{i} to[out=90,in=90,looseness=3] (0.25,-0.15);
\opendot{0.23,0.05};
\end{tikzpicture}\ ,&
\begin{tikzpicture}[Q,centerzero,scale=1.2]
\draw[to-] (-0.25,-0.15) \botlabel{i} to[out=90,in=90,looseness=3] (0.25,-0.15);
\opendot{-0.22,0.05};
\end{tikzpicture}
&=
\begin{tikzpicture}[Q,centerzero,scale=1.2]
\draw[to-] (-0.25,-0.15) \botlabel{i} to[out=90,in=90,looseness=3] (0.25,-0.15);
\opendot{0.23,0.05};
\end{tikzpicture}\ ,\\\label{ruby}
\begin{tikzpicture}[Q,anchorbase,scale=1.2]
\draw[-to] (-0.25,0.15) \toplabel{i} to[out=-90,in=-90,looseness=3] (0.25,0.15);
\draw[-to] (-0.3,-0.4) to[out=up,in=down] (0,0.15)\toplabel{j};
\end{tikzpicture}
&=
\begin{tikzpicture}[Q,anchorbase,scale=1.2]
\draw[-to] (-0.25,0.15) \toplabel{i} to[out=-90,in=-90,looseness=3] (0.25,0.15);
\draw[-to] (0.3,-0.4)to[out=up,in=down] (0,0.15)\toplabel{j};
\end{tikzpicture}\ ,&
\begin{tikzpicture}[Q,anchorbase,scale=1.2]
\draw[to-] (-0.25,0.15) \toplabel{i} to[out=-90,in=-90,looseness=3] (0.25,0.15);
\draw[to-] (-0.3,-0.4) to[out=up,in=down] (0,0.15)\toplabel{j};
\end{tikzpicture}
&=
\begin{tikzpicture}[Q,anchorbase,scale=1.2]
\draw[to-] (-0.25,0.15) \toplabel{i} to[out=-90,in=-90,looseness=3] (0.25,0.15);
\draw[to-] (0.3,-0.4)to[out=up,in=down] (0,0.15)\toplabel{j};
\end{tikzpicture}\ ,
&
\begin{tikzpicture}[Q,anchorbase,scale=1.2]
\draw[-to] (-0.25,0.15) \toplabel{i} to[out=-90,in=-90,looseness=3] (0.25,0.15);
\draw[to-] (-0.3,-0.4) to[out=up,in=down] (0,0.15)\toplabel{j};
\end{tikzpicture}
&=r_{j,i}\ 
\begin{tikzpicture}[Q,anchorbase,scale=1.2]
\draw[-to] (-0.25,0.15) \toplabel{i} to[out=-90,in=-90,looseness=3] (0.25,0.15);
\draw[to-] (0.3,-0.4)to[out=up,in=down] (0,0.15)\toplabel{j};
\end{tikzpicture}\ ,&
\begin{tikzpicture}[Q,anchorbase,scale=1.2]
\draw[to-] (-0.25,0.15) \toplabel{i} to[out=-90,in=-90,looseness=3] (0.25,0.15);
\draw[-to] (-0.3,-0.4) to[out=up,in=down] (0,0.15)\toplabel{j};
\end{tikzpicture}
&=
r_{i,j}^{-1}\ 
\begin{tikzpicture}[Q,anchorbase,scale=1.2]
\draw[to-] (-0.25,0.15) \toplabel{i} to[out=-90,in=-90,looseness=3] (0.25,0.15);
\draw[-to] (0.3,-0.4)to[out=up,in=down] (0,0.15)\toplabel{j};
\end{tikzpicture}\ ,\\\label{wax}
\begin{tikzpicture}[Q,anchorbase,scale=1.2]
\draw[-to] (-0.25,-0.15) \botlabel{i} to[out=90,in=90,looseness=3] (0.25,-0.15);
\draw[to-] (-0.3,0.4) \braiddown (0,-0.15)\botlabel{j};
\end{tikzpicture}
&=
\begin{tikzpicture}[Q,anchorbase,scale=1.2]
\draw[-to] (-0.25,-0.15) \botlabel{i} to[out=90,in=90,looseness=3] (0.25,-0.15);
\draw[to-] (0.3,0.4) \braiddown (0,-0.15)\botlabel{j};
\end{tikzpicture}\ ,&
\begin{tikzpicture}[Q,baseline=-1mm,scale=1.2]
\draw[to-] (-0.25,-0.15) \botlabel{i} to[out=90,in=90,looseness=3] (0.25,-0.15);
\draw[-to] (-0.3,0.4) \braiddown (0,-0.15)\botlabel{j};
\end{tikzpicture}
&=
\begin{tikzpicture}[Q,baseline=-1mm,scale=1.2]
\draw[to-] (-0.25,-0.15) \botlabel{i} to[out=90,in=90,looseness=3] (0.25,-0.15);
\draw[-to] (0.3,0.4) \braiddown (0,-0.15)\botlabel{j};
\end{tikzpicture}\ ,
&\begin{tikzpicture}[Q,baseline=-1mm,scale=1.2]
\draw[-to] (-0.25,-0.15) \botlabel{i} to[out=90,in=90,looseness=3] (0.25,-0.15);
\draw[-to] (-0.3,0.4) \braiddown (0,-0.15)\botlabel{j};
\end{tikzpicture}
&=r_{j,i}^{-1}\ 
\begin{tikzpicture}[Q,baseline=-1mm,scale=1.2]
\draw[-to] (-0.25,-0.15) \botlabel{i} to[out=90,in=90,looseness=3] (0.25,-0.15);
\draw[-to] (0.3,0.4) \braiddown (0,-0.15)\botlabel{j};
\end{tikzpicture}\ ,&
\begin{tikzpicture}[Q,anchorbase,scale=1.2]
\draw[to-] (-0.25,-0.15) \botlabel{i} to[out=90,in=90,looseness=3] (0.25,-0.15);
\draw[to-] (-0.3,0.4) \braiddown (0,-0.15)\botlabel{j};
\end{tikzpicture}
&=r_{i,j}\ 
\begin{tikzpicture}[Q,anchorbase,scale=1.2]
\draw[to-] (-0.25,-0.15) \botlabel{i} to[out=90,in=90,looseness=3] (0.25,-0.15);
\draw[to-] (0.3,0.4)\braiddown (0,-0.15)\botlabel{j};
\end{tikzpicture}\ .
\end{align}
These follow from
\cite[Lem.~2.1, Cor.~2.4, Lem.~5.1, Th.~5.3]{Brundan}.

\begin{rem}\label{snow}
We emphasize that the 2-category $\UU$ depends only on the Cartan datum and the parameters $Q_{i,j}(x,y)$ chosen earlier.
It is merely the normalization of the leftward cups and caps
and the downward crossings that depends on the choice of the normalization homomorphisms $c_i$.
\end{rem}

In most subsequent calculations, we will work systematically with generating functions, which in general will be formal Laurent series in auxiliary variables $u^{-1}, v^{-1},\dots$.
For such a generating function $f(u)$,
we use the notation
$[f(u)]_{u:n}$ to denote its $u^{n}$-coefficient,
$[f(u)]_{u;\geq 0}$ for its polynomial part, 
and so on.
For any polynomial $f(x)$, we have that
\begin{align}\label{trick}
\left[\frac{f(u)}{u-x}\right]_{u:-1}&=f(x),&
\left[\frac{f(u)}{u-x}\right]_{u:<0}&=\frac{f(x)}{u-x}.
\end{align}
We view the series
\begin{align}\label{often}
\frac{1}{u-x} &= \sum_{r \geq 0} x^r u^{-r-1}.
\end{align}
as a generating function for multiple dots on a string, introducing the shorthand
\begin{align}\label{dgf}
\begin{tikzpicture}[Q,centerzero,scale=1.1]
\draw[-] (0,-0.3) -- (0,0.3);
\Circled{0,0}{u};
\end{tikzpicture}
&:=
\begin{tikzpicture}[Q,centerzero,scale=1.1]
\draw[-] (0,-0.3) -- (0,0.3);
\Pin{0,0}{.7,0}{\frac{1}{u-x}};
\end{tikzpicture}\
.
\end{align}
The following is a consequence of \cref{dotslide}:
\begin{align}\label{gendotslide}
\begin{tikzpicture}[Q,centerzero,scale=1.2]
\draw[-to] (-0.3,-0.3) \botlabel{i} -- (0.3,0.3);
\draw[-to] (0.3,-0.3) \botlabel{j} -- (-0.3,0.3);
\Circled{-0.15,-0.15}{u};
\end{tikzpicture}
-
\begin{tikzpicture}[Q,centerzero,scale=1.3]
\draw[-to] (-0.3,-0.3) \botlabel{i} -- (0.3,0.3);
\draw[-to] (0.3,-0.3) \botlabel{j} -- (-0.3,0.3);
\Circled{0.15,0.15}{u};
\end{tikzpicture}
&= \delta_{i,j}  \ 
\begin{tikzpicture}[Q,centerzero,scale=1.3]
\draw[-to] (-0.2,-0.3) \botlabel{i} -- (-0.2,0.3);
\draw[-to] (0.2,-0.3) \botlabel{i} -- (0.2,0.3);
\Circled{-0.2,0}{u};
\Circled{0.2,0}{u};
\end{tikzpicture} =    
\begin{tikzpicture}[Q,centerzero,scale=1.3]
\draw[-to] (-0.3,-0.3) \botlabel{i} -- (0.3,0.3);
\draw[-to] (0.3,-0.3) \botlabel{j} -- (-0.3,0.3);
\Circled{-0.15,0.15}{u};
\end{tikzpicture}
-
\begin{tikzpicture}[Q,centerzero,scale=1.3]
\draw[-to] (-0.3,-0.3) \botlabel{i} -- (0.3,0.3);
\draw[-to] (0.3,-0.3) \botlabel{j} -- (-0.3,0.3);
\Circled{0.15,-0.15}{u};
\end{tikzpicture} \ .
\end{align}

The most interesting relations in $\UU$ involve bubbles.
To formulate them, we first introduce the
{\em fake bubbles}\footnote{In the literature, our fake bubbles
are usually denoted by negatively dotted bubbles, a convention introduced by Lauda \cite{Lauda}.
This would be too confusing in this context---later on we will invert dots on strings labelled by $i=\tau i$, so that negatively dotted bubbles will have a natural meaning.}
\begin{align}\label{fake1}
\begin{tikzpicture}[baseline=-1mm,Q]
\draw[to-] (-0.25,0) arc(180:-180:0.25);
\node at (0,-.4) {\strandlabel{i}};
\region{0.5,0}{\lambda};
\dottybubblelabel{0,0}{n};
\end{tikzpicture}&:=
c_i(\lambda)^{-n-1}
\det \Big[-\ \begin{tikzpicture}[baseline=-1mm,Q]
\draw[-to] (0.25,0) arc(360:0:0.25);
\node at (0,-.4) {\strandlabel{i}};
\region{0.45,0}{\lambda};
\multopendot{-0.25,0}{east}{r-s+h_i(\lambda)};
\end{tikzpicture}
\Big]_{r,s=1,\dots,n}
\\\intertext{for $0 \leq n \leq h_i(\lambda)$, and}
\begin{tikzpicture}[baseline=-1mm,Q]
\draw[-to] (0.25,0) arc(360:0:0.25);
\node at (0,-.4) {\strandlabel{i}};
\region{0.5,0}{\lambda};
\dottybubblelabel{0,0}{n};
\end{tikzpicture}&:=
c_i(\lambda)^{n+1}
\det \Big[-\ \begin{tikzpicture}[baseline=-1mm,Q]
\draw[to-] (-0.25,0) arc(180:-180:0.25);
\node at (0,-.4) {\strandlabel{i}};
\region{-0.5,0}{\lambda};
\multopendot{0.25,0}{west}{r-s-h_i(\lambda)};
\end{tikzpicture} \Big]_{r,s=1,\dots,n}
\label{fake2}
\end{align}
for $0 \leq n \leq -h_i(\lambda)$.
When $n=0$, these are $0 \times 0$ matrices, whose determinants should be interpreted as $\id_{\one_\lambda}$.
We put the fake bubbles together to form the
{\em fake bubble polynomials}
\begin{align}\label{fakebubblepolynomial}
\begin{tikzpicture}[baseline=-1mm,Q]
\draw[Fakebubble,to-] (-0.25,0) arc(180:-180:0.25);
\node at (0,-0.4) {\strandlabel{i}};
\region{0.5,0}{\lambda};
\bubblelabel{0,0}{u};
\end{tikzpicture}
&:=
\sum_{n = 0}^{h_i(\lambda)}
\begin{tikzpicture}[baseline=-1mm,Q]
\draw[to-] (-0.25,0) arc(180:-180:0.25);
\node at (0,-.4) {\strandlabel{i}};
\region{0.5,0}{\lambda};
\dottybubblelabel{0,0}{n};
\end{tikzpicture}\ 
u^{h_i(\lambda)-n},&
\begin{tikzpicture}[baseline=-1mm,Q]
\draw[Fakebubble,-to] (0.25,0) arc(360:0:0.25);
\node at (0,-0.4) {\strandlabel{i}};
\region{0.5,0}{\lambda};
\bubblelabel{0,0}{u};
\end{tikzpicture}
&:=
\sum_{n = 0}^{-h_i(\lambda)}
\begin{tikzpicture}[baseline=-1mm,Q]
\draw[-to] (0.25,0) arc(360:0:0.25);
\node at (0,-.4) {\strandlabel{i}};
\region{0.5,0}{\lambda};
\dottybubblelabel{0,0}{n};
\end{tikzpicture}\ 
u^{-h_i(\lambda)-n},
\end{align}
then define the {\em bubble generating functions}
\begin{align}
\begin{tikzpicture}[baseline=-1mm,Q]
\draw[to-] (-0.25,0) arc(180:-180:0.25);
\node at (0,-.4) {\strandlabel{i}};
\region{0.95,0}{\lambda};
\node at (.55,0) {$(u)$};
\end{tikzpicture}
&:= 
c_i(\lambda)\left(
\begin{tikzpicture}[baseline=-1mm,Q]
\draw[to-] (-0.25,0) arc(180:-180:0.25);
\node at (0,-0.4) {\strandlabel{i}};
\region{0.6,0}{\lambda};
\Circled{0.25,0}{u};
\end{tikzpicture}+
\begin{tikzpicture}[baseline=-1mm,Q]
\draw[Fakebubble,to-] (-0.25,0) arc(180:-180:0.25);
\node at (0,-0.4) {\strandlabel{i}};
\region{0.4,0}{\lambda};
\bubblelabel{0,0}{u};
\node at (0,.3) {$\phantom{.}$};
\node at (0,-.3) {$\phantom{.}$};
\end{tikzpicture}\right),&
\label{bubblegeneratingfunctions}
\begin{tikzpicture}[baseline=-1mm,Q]
\draw[-to] (-0.25,0) arc(180:-180:0.25);
\node at (0,-.4) {\strandlabel{i}};
\region{0.95,0}{\lambda};
\node at (.55,0) {$(u)$};
\end{tikzpicture}
&:= 
c_i(\lambda)^{-1}\left(
\begin{tikzpicture}[baseline=-1mm,Q]
\draw[-to] (0.25,0) arc(360:0:0.25);
\node at (0,-0.4) {\strandlabel{i}};
\region{0.5,0}{\lambda};
\Circled{-0.25,0}{u};
\end{tikzpicture}+
\begin{tikzpicture}[baseline=-1mm,Q]
\draw[Fakebubble,-to] (0.25,0) arc(360:0:0.25);
\node at (0,-0.4) {\strandlabel{i}};
\region{0.5,0}{\lambda};
\bubblelabel{0,0}{u};
\node at (0,.3) {$\phantom{.}$};
\node at (0,-.3) {$\phantom{.}$};
\end{tikzpicture}\right).
\end{align}
The determinants in \cref{fake1,fake2} arise from the solution of
equations implied by the relations
\begin{align}\label{infgrass0}
\begin{tikzpicture}[baseline=-1mm,Q]
\draw[to-] (-0.25,0) arc(180:-180:0.25);
\node at (0,-.4) {\strandlabel{i}};
\region{0.95,0}{\lambda};
\node at (.55,0) {$(u)$};
\end{tikzpicture}& 
\in u^{h_i(\lambda)} \id_{\one_\lambda} + u^{h_i(\lambda)-1} \End_{\UU}(\one_\lambda)\llbracket u^{-1}\rrbracket,\\\label{infgrass1}
\begin{tikzpicture}[baseline=-1mm,Q]
\draw[-to] (-0.25,0) arc(180:-180:0.25);
\node at (0,-.4) {\strandlabel{i}};
\region{0.95,0}{\lambda};
\node at (.55,0) {$(u)$};
\end{tikzpicture} &\in 
u^{-h_i(\lambda)} + u^{-h_i(\lambda)-1} \End_{\UU}(\one_\lambda)\llbracket u^{-1}\rrbracket,\\
\label{infgrass}\qquad&
\begin{tikzpicture}[Q,centerzero,scale=1]
\draw[to-] (-0.68,0) arc(180:-180:0.25);
\node[black] at (0.1,0) {$(u)$};
\node at (-.4,-.43) {\strandlabel{i}};
\region{.6,0}{\lambda};
\end{tikzpicture} \ 
\begin{tikzpicture}[Q,centerzero,scale=1]
\draw[-to] (-.25,0) arc(180:-180:0.25);
\node[black] at (0.54,0) {$(u)$};
\node at (0,-0.4) {\strandlabel{i}};
\end{tikzpicture} = \id_{\one_\lambda}.
\end{align}
Again, these are derived from the defining relations in \cite{Brundan}. 
More generally, for any $n \geq 0$, we can define
the {\em bubbles of degree $2 d_i n$} by the same determinants
\cref{fake1,fake2}. When $n \geq \max(0,h_i(\lambda)+1)$ and $n \geq \max(0, -h_i(\lambda)+1)$, respectively,
the resulting 2-endomorphisms of $\one_\lambda$ are equal to the actual dotted bubbles
\begin{align}
\begin{tikzpicture}[baseline=-1mm,Q]
\draw[to-] (-0.25,0) arc(180:-180:0.25);
\node at (0,-.4) {\strandlabel{i}};
\region{0.5,0}{\lambda};
\dottybubblelabel{0,0}{n};
\end{tikzpicture} &= \begin{tikzpicture}[baseline=-1mm,Q]
\draw[to-] (-0.25,0) arc(180:-180:0.25);
\node at (0,-.4) {\strandlabel{i}};
\region{-0.5,0}{\lambda};
\multopendot{0.25,0}{west}{n-h_i(\lambda)-1};
\end{tikzpicture}=c_i(\lambda)^{-1} \bigg[\ 
\begin{tikzpicture}[baseline=-1mm,Q]
\draw[to-] (-0.25,0) arc(180:-180:0.25);
\node at (0,-.4) {\strandlabel{i}};
\region{0.95,0}{\lambda};
\node at (.55,0) {$(u)$};
\end{tikzpicture}\bigg]_{u:h_i(\lambda)-n},\\
\begin{tikzpicture}[baseline=-1mm,Q]
\draw[-to] (-0.25,0) arc(180:-180:0.25);
\node at (0,-.4) {\strandlabel{i}};
\region{0.5,0}{\lambda};
\dottybubblelabel{0,0}{n};
\end{tikzpicture} &= 
\begin{tikzpicture}[baseline=-1mm,Q]
\draw[-to] (0.25,0) arc(360:0:0.25);
\node at (0,-.4) {\strandlabel{i}};
\region{0.45,0}{\lambda};
\multopendot{-0.25,0}{east}{n+h_i(\lambda)-1};
\end{tikzpicture}
=c_i(\lambda) \bigg[\ 
\begin{tikzpicture}[baseline=-1mm,Q]
\draw[-to] (-0.25,0) arc(180:-180:0.25);
\node at (0,-.4) {\strandlabel{i}};
\region{0.95,0}{\lambda};
\node at (.55,0) {$(u)$};
\end{tikzpicture}\bigg]_{u:-h_i(\lambda)-n}.
\end{align}

To further illustrate the usefulness of the generating functions, we next record the remaining relations 
from \cite{KL3} written using 
them. They all follow from relations
derived in \cite{KL3,Brundan} by equating coefficients. Most of these generating function forms were originally worked out in \cite{HKM,unfurling}.
Let
\begin{align}\label{boshgosh}
R_{i,j}(x,y) &:= \begin{dcases}
t_{i,j}^{-1} Q_{i,j}(x,y)
&\text{if $i \neq j$}\\
\frac{1}{(x-y)^2}&\text{if $i=j$.}
\end{dcases}
\end{align}
Then:
\begin{align}
\label{bubslide}
\begin{tikzpicture}[anchorbase,Q,scale=.9]
\draw[-to] (-0.6,-0.5)\botlabel{j} to (-0.6,0.5);
\draw[to-] (-0.25,0) arc(180:-180:0.25);
\node at (0,-.4) {\strandlabel{i}};
\node at (.55,0) {$(u)$};
\end{tikzpicture}
&=\ 
\begin{tikzpicture}[anchorbase,Q,scale=.9]
\draw[to-] (-0.25,0) arc(180:-180:0.25);
\node at (0,-.4) {\strandlabel{i}};
\node at (.55,0) {$(u)$};
\draw[-to] (1,-0.5)\botlabel{j} to (1,0.5);
\Pin{1,0}{2.1,0}{R_{i, j}(u,x)};
\end{tikzpicture}\ ,
&\begin{tikzpicture}[anchorbase,Q,scale=.9]
\draw[-to] (-0.25,0) arc(180:-180:0.25);
\node at (0,-0.4) {\strandlabel{i}};
\node at (.55,0) {$(u)$};
\draw[-to] (1,-0.5)\botlabel{j} to (1,0.5);
\end{tikzpicture}\ 
&=\ 
\begin{tikzpicture}[anchorbase,Q,scale=.9]
\draw[-to] (-0.6,-0.5)\botlabel{j} to (-0.6,0.5);
\draw[-to] (-0.25,0) arc(180:-180:0.25);
\node at (0,-.4) {\strandlabel{i}};
\node at (.55,0) {$(u)$};
\Pin{-.6,0}{-1.7,0}{R_{i, j}(u,x)};\end{tikzpicture}
\ ,\\
\begin{tikzpicture}[anchorbase,scale=1.1,Q]
\draw[-to] (0,-0.5)\botlabel{i} to[out=up,in=180] (0.3,0.2) to[out=0,in=up] (0.45,0) to[out=down,in=0] (0.3,-0.2) to[out=180,in=down] (0,0.5);
\Circled{.42,0}{u};
\end{tikzpicture}
&=
-\left[\ 
\begin{tikzpicture}[anchorbase,scale=1.1,Q]
\draw[-to] (-0.8,-0.5)\botlabel{i} -- (-0.8,0.5);
\Circled{-0.8,0}{u};
\draw[Fakebubble,to-] (-0.4,0) arc(-180:180:0.2);
\bubblelabel{-.2,0}{u};
\node at (-0.2,-.33) {\strandlabel{i}};
\end{tikzpicture}
\ \right]_{u:< 0},
&
\begin{tikzpicture}[anchorbase,scale=1.1,Q]
\draw[-to] (0,-0.5)\botlabel{i} to[out=up,in=0] (-0.3,0.2) to[out=180,in=up] (-0.45,0) to[out=down,in=180] (-0.3,-0.2) to[out=0,in=down] (0,0.5);
\Circled{-.42,0}{u};
\end{tikzpicture}
&=
\left[
\ \begin{tikzpicture}[anchorbase,scale=1.1,Q]
\draw[-to] (1.2,-0.5)\botlabel{i} -- (1.2,0.5);
\Circled{1.2,0}{u};
\draw[Fakebubble,to-] (0.4,0) arc(180:-180:0.2);
\bubblelabel{.6,0}{u};
\node at (.6,-.33) {\strandlabel{i}};
\end{tikzpicture}\ 
\right]_{u:< 0},
\label{curlrels}\\
\label{altquadratic}
\begin{tikzpicture}[Q,centerzero,scale=1.2]
\draw[-to] (-0.2,-0.4) \botlabel{i} to[out=45,in=down] (0.15,0) to[out=up,in=-45] (-0.2,0.4);
\draw[to-] (0.2,-0.4) \botlabel{j} to[out=135,in=down] (-0.15,0) to[out=up,in=225] (0.2,0.4);
\Circled{-.15,0}{u};
\end{tikzpicture}
&=
(-1)^{\delta_{i,j}}\begin{tikzpicture}[Q,centerzero,scale=1.4]
\draw[-to] (-0.14,-0.3) \botlabel{i} -- (-0.14,0.3);
\draw[to-] (0.14,-0.3) \botlabel{j} -- (0.14,0.3);
\Circled{.14,0}{u};
\end{tikzpicture}
+
\delta_{i,j}
\left[\,\begin{tikzpicture}[Q,centerzero,scale=1.6]
\draw[-to] (-0.2,-0.3) \botlabel{i} to [looseness=2.2,out=90,in=90] (0.2,-0.3);
\draw[-to] (0.2,0.3) to [looseness=2.2,out=-90,in=-90] (-0.2,0.3)\toplabel{i};
\Circled{-.15,.12}{u};
\Circled{-.15,-.12}{u};
\draw[Fakebubble,to-] (0.3,0) arc(180:-180:0.142);
\bubblelabel{.442,0}{u};
\node at (0.442,-.23) {\strandlabel{i}};
\end{tikzpicture}\,\right]_{\!u:<0}\!\!\!\!\!,
&
\begin{tikzpicture}[Q,centerzero,scale=1.2]
\draw[to-] (-0.2,-0.4) \botlabel{i} to[out=45,in=down] (0.15,0) to[out=up,in=-45] (-0.2,0.4);
\draw[-to] (0.2,-0.4) \botlabel{j} to[out=135,in=down] (-0.15,0) to[out=up,in=225] (0.2,0.4);
\Circled{.15,0}{u};
\end{tikzpicture}
&=
(-1)^{\delta_{i,j}}\begin{tikzpicture}[Q,centerzero,scale=1.4]
\draw[to-] (-0.14,-0.3) \botlabel{i} -- (-0.14,0.3);
\draw[-to] (0.14,-0.3) \botlabel{j} -- (0.14,0.3);
\Circled{-.14,0}{u};
\end{tikzpicture}
+
\delta_{i,j}
\left[\,\begin{tikzpicture}[Q,centerzero,scale=1.6]
\draw[-to] (0.2,-0.3)  to [looseness=2.2,out=90,in=90] (-0.2,-0.3)\botlabel{i};
\draw[-to] (-0.2,0.3)\toplabel{i} to [looseness=2.2,out=-90,in=-90] (0.2,0.3);
\Circled{.15,.12}{u};
\Circled{.15,-.12}{u};
\draw[Fakebubble,to-] (-0.6,0) arc(-180:180:0.142);
\bubblelabel{-.458,0}{u};
\node at (-0.458,-.23) {\strandlabel{i}};
\end{tikzpicture}\,\right]_{\!u:<0}\!\!\!\!\!\!\!,
\end{align}

\vspace{-2mm}

\begin{align}
\label{altbraid}
\begin{tikzpicture}[Q,centerzero,scale=1.5]
\draw[-to] (-0.3,-0.3) \botlabel{i} -- (0.3,0.3);
\draw[to-] (0,-0.3) \botlabel{j} to[out=135,in=down] (-0.22,0) to[out=up,in=225] (0,0.3);
\draw[-to] (0.3,-0.3) \botlabel{k} -- (-0.3,0.3);
\end{tikzpicture}
\!-\!
\begin{tikzpicture}[Q,centerzero,scale=1.5]
\draw[-to] (-0.3,-0.3) \botlabel{i} -- (0.3,0.3);
\draw[to-] (0,-0.3) \botlabel{j} to[out=45,in=down] (0.22,0) to[out=up,in=-45] (0,0.3);
\draw[-to] (0.3,-0.3) \botlabel{k} -- (-0.3,0.3);
\end{tikzpicture}
&=
\delta_{i,j}\delta_{j,k}
\left[\ \begin{tikzpicture}[Q,centerzero,scale=1.5]
\draw[-to] (-0.2,-0.3) \botlabel{i} to [looseness=2.2,out=90,in=90] (0.2,-0.3);
\draw[-to] (0.2,0.3) to [looseness=2.2,out=-90,in=-90] (-0.2,0.3)\toplabel{i};
\draw[-to] (.85,-.3) \botlabel{i} to (.85,.3);
\Circled{-.15,.12}{u};
\Circled{-.15,-.12}{u};
\Circled{.85,0}{u};
\draw[Fakebubble,to-] (0.3,0) arc(180:-180:0.142);
\bubblelabel{.442,0}{u};
\node at (0.442,-.23) {\strandlabel{i}};
\end{tikzpicture}
+
\begin{tikzpicture}[Q,centerzero,scale=1.5]
\draw[to-] (0,-0.3)  to [looseness=2.2,out=90,in=90] (0.4,-0.3)\botlabel{i};
\draw[to-] (0.4,0.3) \toplabel{i} to [looseness=2.2,out=-90,in=-90] (0,0.3);
\draw[-to] (-.65,-.3) \botlabel{i} to (-.65,.3);
\Circled{-.65,0}{u};
\Circled{.34,.12}{u};
\Circled{.34,-.12}{u};
\draw[Fakebubble,to-] (-0.4,0) arc(-180:180:0.142);
\bubblelabel{-.258,0}{u};
\node at (-0.258,-.23) {\strandlabel{i}};
\end{tikzpicture}
\ \right]_{\!u:<0}\!\!\!\!\!.
\end{align}
There are also some more
general versions of some of \cref{curlrels,altquadratic}:
\begin{align}
\begin{tikzpicture}[anchorbase,scale=1,Q]
\draw[-to] (0,-0.5)\botlabel{i} to[out=up,in=180] (0.3,0.2) to[out=0,in=up] (0.45,0) to[out=down,in=0] (0.3,-0.2) to[out=180,in=down] (0,0.5);
\Circled{.42,0}{u};
\region{0.7,0}{\lambda};
\end{tikzpicture}\!\!\!
&=
-c_i(\lambda)\left[\ 
\begin{tikzpicture}[anchorbase,scale=1,Q]
\draw[-to] (-0.8,-0.5)\botlabel{i} -- (-0.8,0.5);
\Circled{-0.8,0}{u};
\draw[-to] (-.4,0) arc(180:-180:0.2);
\node at (0.26,0) {$(u)$};
\node at (-.17,-.33) {\strandlabel{i}};
\region{0.7,0}{\lambda};
\end{tikzpicture}
\right]_{u:< 0},
&\!\!\!\!
\begin{tikzpicture}[anchorbase,scale=1,Q]
\draw[-to] (0,-0.5)\botlabel{i} to[out=up,in=0] (-0.3,0.2) to[out=180,in=up] (-0.45,0) to[out=down,in=180] (-0.3,-0.2) to[out=0,in=down] (0,0.5);
\Circled{-.42,0}{u};
\region{-0.7,0}{\lambda};
\end{tikzpicture}\!
&=
c_i(\lambda)^{-1}
\left[
\ \begin{tikzpicture}[anchorbase,scale=1,Q]
\draw[-to] (1.2,-0.5)\botlabel{i} -- (1.2,0.5);
\Circled{1.2,0}{u};
\draw[to-] (0,0) arc(180:-180:0.2);
\node at (0.65,0) {$(u)$};
\node at (0.2,-.33) {\strandlabel{i}};
\region{-0.25,0}{\lambda};
\end{tikzpicture}\ 
\right]_{u:< 0},
\label{Curlrels}\\
\label{Altquadratic}
\begin{tikzpicture}[Q,centerzero,scale=1.4]
\draw[-to] (-0.15,-0.3) \botlabel{i} to[out=45,in=down] (0.11,0) to[out=up,in=-45] (-0.15,0.3);
\draw[to-] (0.15,-0.3) \botlabel{j} to[out=135,in=down] (-0.11,0) to[out=up,in=225] (0.15,0.3);
\Circled{-.11,0}{u};
\region{0.25,0}{\lambda};
\end{tikzpicture}\!
&=
(-1)^{\delta_{i,j}}\begin{tikzpicture}[Q,centerzero,scale=1.4]
\draw[-to] (-0.14,-0.3) \botlabel{i} -- (-0.14,0.3);
\draw[to-] (0.14,-0.3) \botlabel{j} -- (0.14,0.3);
\Circled{.14,0}{u};
\region{0.35,0}{\lambda};
\end{tikzpicture}\!\!
+\!
\delta_{i,j}c_i(\lambda)^{-1}\!\!
\left[\,\begin{tikzpicture}[Q,centerzero,scale=1.4]
\draw[-to] (-0.2,-0.3) \botlabel{i} to [looseness=2.2,out=90,in=90] (0.2,-0.3);
\draw[-to] (0.2,0.3) to [looseness=2.2,out=-90,in=-90] (-0.2,0.3)\toplabel{i};
\Circled{-.15,.12}{u};
\Circled{-.15,-.12}{u};
\draw[to-] (0.27,0) arc(180:-180:0.142);
\node at (0.73,0) {$(u)$};
\node at (0.42,-.23) {\strandlabel{i}};
\region{0.45,.3}{\lambda};
\end{tikzpicture}\!\!\right]_{\!u:<0}\!\!\!\!,\,
&\!\!\!\!
\begin{tikzpicture}[Q,centerzero,scale=1.4]
\draw[to-] (-0.15,-0.3) \botlabel{i} to[out=45,in=down] (0.11,0) to[out=up,in=-45] (-0.15,0.3);
\draw[-to] (0.15,-0.3) \botlabel{j} to[out=135,in=down] (-0.11,0) to[out=up,in=225] (0.15,0.3);
\Circled{.11,0}{u};
\region{-.25,0}{\lambda};
\end{tikzpicture}
&=
(-1)^{\delta_{i,j}}\!\!\begin{tikzpicture}[Q,centerzero,scale=1.4]
\draw[to-] (-0.14,-0.3) \botlabel{i} -- (-0.14,0.3);
\draw[-to] (0.14,-0.3) \botlabel{j} -- (0.14,0.3);
\Circled{-.14,0}{u};
\region{-.35,0}{\lambda};
\end{tikzpicture}\!
+\!
\delta_{i,j}c_i(\lambda)\!\!
\left[\,\begin{tikzpicture}[Q,centerzero,scale=1.4]
\draw[-to] (0.2,-0.3)  to [looseness=2.2,out=90,in=90] (-0.2,-0.3)\botlabel{i};
\draw[-to] (-0.2,0.3)\toplabel{i} to [looseness=2.2,out=-90,in=-90] (0.2,0.3);
\Circled{.15,.12}{u};
\Circled{.15,-.12}{u};
\draw[to-] (-0.82,0) arc(-180:180:0.142);
\node at (-0.34,0) {$(u)$};
\node at (-0.658,-.23) {\strandlabel{i}};
\region{-0.65,.3}{\lambda};
\end{tikzpicture}\,\right]_{\!u:<0}\!\!\!\!.
\end{align}

The following theorem is a version of the main result of \cite{Brundan}. It gives another presentation for $\UU$
which is more symmetric but
less efficient than Rouquier's presentation formulated above. It is essentially the presentation 
of Khovanov and Lauda from \cite{KL3}.

\begin{theo}\label{klpres}
The 2-category $\UU$ can be defined equivalently
as the graded 2-category with objects $X$, generating 1-morphisms $E_i \one_\lambda = \one_{\lambda+\alpha_i} E_i, F_i \one_\lambda = \one_{\lambda-\alpha_i} F_i$\:($\lambda \in X, i \in I$), and all of the generating 2-morphism listed in \cref{table1} (both halves), subject to the relations
\cref{dotslide,quadratic,braid,rightadj,leftpivots,lotsmore,ruby,wax,infgrass0,infgrass,curlrels,altquadratic} (interpreted using the shorthands \cref{fake1,fake2,fakebubblepolynomial,bubblegeneratingfunctions}).
\end{theo}

\begin{rem}\label{notevendone}
Yet more relations can be obtained by partially or fully rotating 
\cref{dotslide,quadratic,braid,bubslide,curlrels,altquadratic,altbraid,Curlrels,Altquadratic}, i.e., attaching some cups to the bottom and some caps to the top of these relations then simplifying using \cref{leftpivots,lotsmore,ruby,wax}.
For example, the full rotations of \cref{dotslide,quadratic,braid}
produce
\begin{align}
\begin{tikzpicture}[Q,centerzero]
\draw[to-] (-0.3,-0.3) \botlabel{i} -- (0.3,0.3);
\draw[to-] (0.3,-0.3) \botlabel{j} -- (-0.3,0.3);
\opendot{0.15,0.15};
\end{tikzpicture}
-\begin{tikzpicture}[Q,centerzero]
\draw[to-] (-0.3,-0.3) \botlabel{i} -- (0.3,0.3);
\draw[to-] (0.3,-0.3) \botlabel{j} -- (-0.3,0.3);
\opendot{-0.15,-0.15};
\end{tikzpicture}
&= \delta_{i,j}  \ 
\begin{tikzpicture}[Q,centerzero]
\draw[to-] (-0.2,-0.3) \botlabel{i} -- (-0.2,0.3);
\draw[to-] (0.2,-0.3) \botlabel{i} -- (0.2,0.3);
\end{tikzpicture} =
\begin{tikzpicture}[Q,centerzero]
\draw[to-] (-0.3,-0.3) \botlabel{i} -- (0.3,0.3);
\draw[to-] (0.3,-0.3) \botlabel{j} -- (-0.3,0.3);
\opendot{0.15,-0.15};
\end{tikzpicture}- \begin{tikzpicture}[Q,centerzero]
\draw[to-] (-0.3,-0.3) \botlabel{i} -- (0.3,0.3);
\draw[to-] (0.3,-0.3) \botlabel{j} -- (-0.3,0.3);
\opendot{-0.15,0.15};
\end{tikzpicture}
\ ,\label{downwarddotslide}\\
\begin{tikzpicture}[Q,centerzero,scale=1.1]
\draw[to-] (-0.2,-0.4) \botlabel{i} to[out=45,in=down] (0.15,0) to[out=up,in=-45] (-0.2,0.4);
\draw[to-] (0.2,-0.4) \botlabel{j} to[out=135,in=down] (-0.15,0) to[out=up,in=225] (0.2,0.4);
\end{tikzpicture}
&=
\begin{tikzpicture}[Q,centerzero,scale=1.1]
\draw[to-] (-0.2,-0.3) \botlabel{i} -- (-0.2,0.3);
\draw[to-] (0.2,-0.3) \botlabel{j} -- (0.2,0.3);
\Pinpin{0.2,0}{-0.2,0}{-1.1,0}{{'}Q_{i,j}(x,y)};
\end{tikzpicture}\ ,\label{downwardscrossing}\\
\begin{tikzpicture}[Q,centerzero,scale=1.1]
\draw[to-] (-0.4,-0.4) \botlabel{i} -- (0.4,0.4);
\draw[to-] (0,-0.4) \botlabel{j} to[out=45,in=down] (0.32,0) to[out=up,in=-45] (0,0.4);
\draw[to-] (0.4,-0.4) \botlabel{k} -- (-0.4,0.4);
\region{0.5,0}{\lambda};
\end{tikzpicture}
-
\begin{tikzpicture}[Q,centerzero,scale=1.1]
\draw[to-] (-0.4,-0.4) \botlabel{i} -- (0.4,0.4);
\draw[to-] (0,-0.4) \botlabel{j} to[out=135,in=down] (-0.32,0) to[out=up,in=225] (0,0.4);
\draw[to-] (0.4,-0.4) \botlabel{k} -- (-0.4,0.4);
\region{0.5,0}{\lambda};
\end{tikzpicture}
&=\delta_{i,k}\ 
\begin{tikzpicture}[Q,centerzero,scale=1.2]
\draw[to-] (-0.3,-0.3) \botlabel{i} -- (-0.3,0.3);
\draw[to-] (0,-0.3) \botlabel{j} -- (0,0.3);
\draw[to-] (0.3,-0.3) \botlabel{i} -- (0.3,0.3);
\Pinpinpin{.3,0}{0,0}{-.3,0}{-1.7,0}{
\frac{{'}Q_{i,j}(x,y)-{'}Q_{i,j}(z,y)}{x-z}};
\end{tikzpicture}\label{downwardsbraid}
\end{align}
where 
\begin{equation}\label{altparams}
{'}Q_{i,j}(x,y) := r_{i,j} r_{j,i} Q_{i,j}(x,y).
\end{equation}
Thus, the downward dots and crossings also satisfy quiver Hecke algebra relations, but for 
a different matrix of parameters and some different sign conventions.
\end{rem}

\subsection{Symmetries}\label{symmetries}

Next we explain how to lift the symmetries $\omega,\psi$ and $\sigma$ from \cref{omegainv,sigmainv,psiinv}
to the 2-quantum group.

There is an isomorphism of graded 2-categories
\begin{equation}\label{Omegainv}
\bar\Omega:\UU^\op \stackrel{\sim}{\rightarrow} \UU
\end{equation}
defined on objects by $\lambda \mapsto -\lambda$,
on generating 1-morphisms by
$E_i \one_\lambda \mapsto F_i \one_{-\lambda}$
and
$F_i \one_\lambda \mapsto E_i \one_{-\lambda}$,
and on generating 2-morphisms by
\begin{align*}
\begin{tikzpicture}[Q,centerzero]
\draw[-to] (0,-0.3) \botlabel{i} -- (0,0.3);
\opendot{0,0};
\region{0.3,0}{\lambda};
\end{tikzpicture}
&\mapsto \begin{tikzpicture}[Q,centerzero]
\draw[to-] (0,-0.3) \botlabel{i} -- (0,0.3);
\opendot{0,0};
\region{0.4,0}{-\lambda};
\end{tikzpicture},&
\begin{tikzpicture}[Q,centerzero]
\draw[to-] (0,-0.3) \botlabel{i} -- (0,0.3);
\opendot{0,0};
\region{0.3,0}{\lambda};
\end{tikzpicture}
&\mapsto\begin{tikzpicture}[Q,centerzero]
\draw[-to] (0,-0.3) \botlabel{i} -- (0,0.3);
\opendot{0,0};
\region{0.4,0}{-\lambda};
\end{tikzpicture},\\
\begin{tikzpicture}[Q,centerzero,scale=.9]
\draw[-to] (-0.3,-0.3) \botlabel{i} -- (0.3,0.3);
\draw[-to] (0.3,-0.3) \botlabel{j} -- (-0.3,0.3);
\region{0.4,0}{\lambda};
\end{tikzpicture}
&\mapsto-r_{j,i}^{-1}\ \begin{tikzpicture}[Q,centerzero,scale=.9]
\draw[to-] (-0.3,-0.3) \botlabel{j} -- (0.3,0.3);
\draw[to-] (0.3,-0.3) \botlabel{i} -- (-0.3,0.3);
\region{0.5,0}{-\lambda};
\end{tikzpicture},&
\begin{tikzpicture}[Q,centerzero,scale=.9]
\draw[to-] (-0.3,-0.3) \botlabel{i} -- (0.3,0.3);
\draw[to-] (0.3,-0.3) \botlabel{j} -- (-0.3,0.3);
\region{0.4,0}{\lambda};
\end{tikzpicture}
&\mapsto-r_{i,j}\ \begin{tikzpicture}[Q,centerzero,scale=.9]
\draw[-to] (-0.3,-0.3) \botlabel{j} -- (0.3,0.3);
\draw[-to] (0.3,-0.3) \botlabel{i} -- (-0.3,0.3);
\region{0.5,0}{-\lambda};
\end{tikzpicture},\\
\begin{tikzpicture}[Q,centerzero,scale=.9]
\draw[-to] (-0.3,-0.3) \botlabel{i} -- (0.3,0.3);
\draw[to-] (0.3,-0.3) \botlabel{j} -- (-0.3,0.3);
\region{0.4,0}{\lambda};
\end{tikzpicture}
&\mapsto -\ \begin{tikzpicture}[Q,centerzero,scale=.9]
\draw[-to] (-0.3,-0.3) \botlabel{j} -- (0.3,0.3);
\draw[to-] (0.3,-0.3) \botlabel{i} -- (-0.3,0.3);
\region{0.5,0}{-\lambda};
\end{tikzpicture},&
\begin{tikzpicture}[Q,centerzero,scale=.9]
\draw[to-] (-0.3,-0.3) \botlabel{i} -- (0.3,0.3);
\draw[-to] (0.3,-0.3) \botlabel{j} -- (-0.3,0.3);
\region{0.4,0}{\lambda};
\end{tikzpicture}
&\mapsto-\ \begin{tikzpicture}[Q,centerzero,scale=.9]
\draw[to-] (-0.3,-0.3) \botlabel{j} -- (0.3,0.3);
\draw[-to] (0.3,-0.3) \botlabel{i} -- (-0.3,0.3);
\region{0.5,0}{-\lambda};
\end{tikzpicture},\\
\begin{tikzpicture}[Q,centerzero]
\draw[-to] (-0.25,0.15) \toplabel{i} to[out=-90,in=-90,looseness=3] (0.25,0.15);
\region{0.45,-.1}{\lambda};
\end{tikzpicture}
&\mapsto \begin{tikzpicture}[Q,centerzero]
\draw[-to] (-0.25,-0.15) \botlabel{i} to [out=90,in=90,looseness=3](0.25,-0.15);
\region{0.55,0.1}{-\lambda};
\end{tikzpicture},&
\begin{tikzpicture}[Q,centerzero]
\draw[-to] (-0.25,-0.15) \botlabel{i} to [out=90,in=90,looseness=3](0.25,-0.15);
\region{0.45,0.1}{\lambda};
\end{tikzpicture}
&\mapsto \begin{tikzpicture}[Q,centerzero]
\draw[-to] (-0.25,0.15) \toplabel{i} to[out=-90,in=-90,looseness=3] (0.25,0.15);
\region{0.55,-.1}{-\lambda};
\end{tikzpicture},\\
\begin{tikzpicture}[Q,centerzero]
\draw[to-] (-0.25,-0.15) \botlabel{i} to [out=90,in=90,looseness=3](0.25,-0.15);
\region{0.45,0.1}{\lambda};
\end{tikzpicture}
&\mapsto \begin{tikzpicture}[Q,centerzero]
\draw[to-] (-0.25,0.15) \toplabel{i} to[out=-90,in=-90,looseness=3] (0.25,0.15);
\region{0.55,-.1}{-\lambda};
\end{tikzpicture},&
\begin{tikzpicture}[Q,centerzero]
\draw[to-] (-0.25,0.15) \toplabel{i} to[out=-90,in=-90,looseness=3] (0.25,0.15);
\region{0.45,-.1}{\lambda};
\end{tikzpicture}
&\mapsto \begin{tikzpicture}[Q,centerzero]
\draw[to-] (-0.25,-0.15) \botlabel{i} to [out=90,in=90,looseness=3](0.25,-0.15);
\region{0.55,0.1}{-\lambda};
\end{tikzpicture},\\
\begin{tikzpicture}[baseline=-1mm,Q]
\draw[to-] (-0.25,0) arc(180:-180:0.25);
\node at (0,-.4) {\strandlabel{i}};
\region{0.95,0}{\lambda};
\node at (.55,0) {$(u)$};
\end{tikzpicture}
&\mapsto 
\begin{tikzpicture}[baseline=-1mm,Q]
\draw[-to] (-0.25,0) arc(180:-180:0.25);
\node at (0,-.4) {\strandlabel{i}};
\region{1.05,0}{-\lambda};
\node at (.55,0) {$(u)$};
\end{tikzpicture},&
\begin{tikzpicture}[baseline=-1mm,Q]
\draw[-to] (-0.25,0) arc(180:-180:0.25);
\node at (0,-.4) {\strandlabel{i}};
\region{0.95,0}{\lambda};
\node at (.55,0) {$(u)$};
\end{tikzpicture}
&\mapsto 
\begin{tikzpicture}[baseline=-1mm,Q]
\draw[to-] (-0.25,0) arc(180:-180:0.25);
\node at (0,-.4) {\strandlabel{i}};
\region{1.05,0}{-\lambda};
\node at (.55,0) {$(u)$};
\end{tikzpicture}.
\end{align*}
Roughly, this reflects string diagrams in a horizontal axis and negates all 2-cell labels (there are some additional scalars arising from crossings).
The proof of existence of $\bar\Omega$ is an easy relations check, using the original ``minimal'' presentation for $\UU$; see \cite[Th.~2.3]{Brundan}.
It is the 2-categorical analog not of $\omega$ but of the anti-linear involution
$\bar\omega := \psi\circ\omega$. 

Corresponding to the bar involution $\psi$, there is an isomorphism of graded 2-categories
\begin{equation}\label{Psiinv}
\Psi:\UU^\op \stackrel{\sim}{\rightarrow} \UU
\end{equation}
defined on objects by $\lambda \mapsto \lambda$,
on generating 1-morphisms by
$E_i \one_\lambda \mapsto E_i \one_{\lambda}$
and
$F_i \one_\lambda \mapsto F_i \one_{\lambda}$,
and on generating 2-morphisms by
\begin{align*}
\begin{tikzpicture}[Q,centerzero]
\draw[-to] (0,-0.3) \botlabel{i} -- (0,0.3);
\opendot{0,0};
\region{0.3,0}{\lambda};
\end{tikzpicture}
&\mapsto \begin{tikzpicture}[Q,centerzero]
\draw[-to] (0,-0.3) \botlabel{i} -- (0,0.3);
\opendot{0,0};
\region{0.3,0}{\lambda};
\end{tikzpicture},&
\begin{tikzpicture}[Q,centerzero]
\draw[to-] (0,-0.3) \botlabel{i} -- (0,0.3);
\opendot{0,0};
\region{0.3,0}{\lambda};
\end{tikzpicture}
&\mapsto\begin{tikzpicture}[Q,centerzero]
\draw[to-] (0,-0.3) \botlabel{i} -- (0,0.3);
\opendot{0,0};
\region{0.3,0}{\lambda};
\end{tikzpicture},\\
\begin{tikzpicture}[Q,centerzero,scale=.9]
\draw[-to] (-0.3,-0.3) \botlabel{i} -- (0.3,0.3);
\draw[-to] (0.3,-0.3) \botlabel{j} -- (-0.3,0.3);
\region{0.4,0}{\lambda};
\end{tikzpicture}
&\mapsto \begin{tikzpicture}[Q,centerzero,scale=.9]
\draw[-to] (-0.3,-0.3) \botlabel{j} -- (0.3,0.3);
\draw[-to] (0.3,-0.3) \botlabel{i} -- (-0.3,0.3);
\region{0.4,0}{\lambda};
\end{tikzpicture},&
\begin{tikzpicture}[Q,centerzero,scale=.9]
\draw[to-] (-0.3,-0.3) \botlabel{i} -- (0.3,0.3);
\draw[to-] (0.3,-0.3) \botlabel{j} -- (-0.3,0.3);
\region{0.4,0}{\lambda};
\end{tikzpicture}
&\mapsto \begin{tikzpicture}[Q,centerzero,scale=.9]
\draw[to-] (-0.3,-0.3) \botlabel{j} -- (0.3,0.3);
\draw[to-] (0.3,-0.3) \botlabel{i} -- (-0.3,0.3);
\region{0.4,0}{\lambda};
\end{tikzpicture},\\
\begin{tikzpicture}[Q,centerzero,scale=.9]
\draw[-to] (-0.3,-0.3) \botlabel{i} -- (0.3,0.3);
\draw[to-] (0.3,-0.3) \botlabel{j} -- (-0.3,0.3);
\region{0.4,0}{\lambda};
\end{tikzpicture}
&\mapsto r_{j,i}^{-1}\ \begin{tikzpicture}[Q,centerzero,scale=.9]
\draw[to-] (-0.3,-0.3) \botlabel{j} -- (0.3,0.3);
\draw[-to] (0.3,-0.3) \botlabel{i} -- (-0.3,0.3);
\region{0.4,0}{\lambda};
\end{tikzpicture},&
\begin{tikzpicture}[Q,centerzero,scale=.9]
\draw[to-] (-0.3,-0.3) \botlabel{i} -- (0.3,0.3);
\draw[-to] (0.3,-0.3) \botlabel{j} -- (-0.3,0.3);
\region{0.4,0}{\lambda};
\end{tikzpicture}
&\mapsto r_{i,j}\ \begin{tikzpicture}[Q,centerzero,scale=.9]
\draw[-to] (-0.3,-0.3) \botlabel{j} -- (0.3,0.3);
\draw[to-] (0.3,-0.3) \botlabel{i} -- (-0.3,0.3);
\region{0.4,0}{\lambda};
\end{tikzpicture},\\
\begin{tikzpicture}[Q,centerzero]
\draw[-to] (-0.25,0.15) \toplabel{i} to[out=-90,in=-90,looseness=3] (0.25,0.15);
\region{0.45,-.1}{\lambda};
\end{tikzpicture}
&\mapsto \begin{tikzpicture}[Q,centerzero]
\draw[to-] (-0.25,-0.15) \botlabel{i} to [out=90,in=90,looseness=3](0.25,-0.15);
\region{0.45,0.1}{\lambda};
\end{tikzpicture},&
\begin{tikzpicture}[Q,centerzero]
\draw[-to] (-0.25,-0.15) \botlabel{i} to [out=90,in=90,looseness=3](0.25,-0.15);
\region{0.45,0.1}{\lambda};
\end{tikzpicture}
&\mapsto \begin{tikzpicture}[Q,centerzero]
\draw[to-] (-0.25,0.15) \toplabel{i} to[out=-90,in=-90,looseness=3] (0.25,0.15);
\region{0.45,-.1}{\lambda};
\end{tikzpicture},\\
\begin{tikzpicture}[Q,centerzero]
\draw[to-] (-0.25,-0.15) \botlabel{i} to [out=90,in=90,looseness=3](0.25,-0.15);
\region{0.45,0.1}{\lambda};
\end{tikzpicture}
&\mapsto \begin{tikzpicture}[Q,centerzero]
\draw[-to] (-0.25,0.15) \toplabel{i} to[out=-90,in=-90,looseness=3] (0.25,0.15);
\region{0.45,-.1}{\lambda};
\end{tikzpicture},&
\begin{tikzpicture}[Q,centerzero]
\draw[to-] (-0.25,0.15) \toplabel{i} to[out=-90,in=-90,looseness=3] (0.25,0.15);
\region{0.45,-.1}{\lambda};
\end{tikzpicture}
&\mapsto \begin{tikzpicture}[Q,centerzero]
\draw[-to] (-0.25,-0.15) \botlabel{i} to [out=90,in=90,looseness=3](0.25,-0.15);
\region{0.45,0.1}{\lambda};
\end{tikzpicture},\\
\begin{tikzpicture}[baseline=-1mm,Q]
\draw[to-] (-0.25,0) arc(180:-180:0.25);
\node at (0,-.4) {\strandlabel{i}};
\region{0.95,0}{\lambda};
\node at (.55,0) {$(u)$};
\end{tikzpicture}
&\mapsto \begin{tikzpicture}[baseline=-1mm,Q]
\draw[to-] (-0.25,0) arc(180:-180:0.25);
\node at (0,-.4) {\strandlabel{i}};
\region{.95,0}{\lambda};
\node at (.55,0) {$(u)$};
\end{tikzpicture},&
\begin{tikzpicture}[baseline=-1mm,Q]
\draw[-to] (-0.25,0) arc(180:-180:0.25);
\node at (0,-.4) {\strandlabel{i}};
\region{0.95,0}{\lambda};
\node at (.55,0) {$(u)$};
\end{tikzpicture}
&\mapsto\begin{tikzpicture}[baseline=-1mm,Q]
\draw[-to] (-0.25,0) arc(180:-180:0.25);
\node at (0,-.4) {\strandlabel{i}};
\region{.95,0}{\lambda};
\node at (.55,0) {$(u)$};
\end{tikzpicture}.
\end{align*}
The proof of existence of $\Psi$ is again a relations check, but one needs to use the more symmetric presentation from \cref{klpres}.

The composition
\begin{equation}\label{Omega}
\Omega := \Psi \circ \bar\Omega^\op :\UU\stackrel{\sim}{\rightarrow}\UU
\end{equation}
is a categorical analog of the Chevalley involution $\omega$. It maps $\lambda \mapsto -\lambda,
E_i \one_\lambda \mapsto F_i \one_{-\lambda},
F_i \one_\lambda \mapsto E_i \one_{-\lambda}$, and
\begin{align*}
\begin{tikzpicture}[Q,centerzero]
\draw[-to] (0,-0.3) \botlabel{i} -- (0,0.3);
\opendot{0,0};
\region{0.3,0}{\lambda};
\end{tikzpicture}
&\mapsto \begin{tikzpicture}[Q,centerzero]
\draw[to-] (0,-0.3) \botlabel{i} -- (0,0.3);
\opendot{0,0};
\region{0.4,0}{-\lambda};
\end{tikzpicture},&
\begin{tikzpicture}[Q,centerzero]
\draw[to-] (0,-0.3) \botlabel{i} -- (0,0.3);
\opendot{0,0};
\region{0.3,0}{\lambda};
\end{tikzpicture}
&\mapsto\begin{tikzpicture}[Q,centerzero]
\draw[-to] (0,-0.3) \botlabel{i} -- (0,0.3);
\opendot{0,0};
\region{0.4,0}{-\lambda};
\end{tikzpicture},\\
\begin{tikzpicture}[Q,centerzero,scale=.9]
\draw[-to] (-0.3,-0.3) \botlabel{i} -- (0.3,0.3);
\draw[-to] (0.3,-0.3) \botlabel{j} -- (-0.3,0.3);
\region{0.4,0}{\lambda};
\end{tikzpicture}
&\mapsto -r_{j,i}^{-1}\ \begin{tikzpicture}[Q,centerzero,scale=.9]
\draw[to-] (-0.3,-0.3) \botlabel{i} -- (0.3,0.3);
\draw[to-] (0.3,-0.3) \botlabel{j} -- (-0.3,0.3);
\region{0.5,0}{-\lambda};
\end{tikzpicture},&
\begin{tikzpicture}[Q,centerzero,scale=.9]
\draw[to-] (-0.3,-0.3) \botlabel{i} -- (0.3,0.3);
\draw[to-] (0.3,-0.3) \botlabel{j} -- (-0.3,0.3);
\region{0.4,0}{\lambda};
\end{tikzpicture}
&\mapsto -r_{i,j}\ \begin{tikzpicture}[Q,centerzero,scale=.9]
\draw[-to] (-0.3,-0.3) \botlabel{i} -- (0.3,0.3);
\draw[-to] (0.3,-0.3) \botlabel{j} -- (-0.3,0.3);
\region{0.5,0}{-\lambda};
\end{tikzpicture},\\
\begin{tikzpicture}[Q,centerzero,scale=.9]
\draw[-to] (-0.3,-0.3) \botlabel{i} -- (0.3,0.3);
\draw[to-] (0.3,-0.3) \botlabel{j} -- (-0.3,0.3);
\region{0.4,0}{\lambda};
\end{tikzpicture}
&\mapsto -r_{i,j}^{-1}\ \begin{tikzpicture}[Q,centerzero,scale=.9]
\draw[to-] (-0.3,-0.3) \botlabel{i} -- (0.3,0.3);
\draw[-to] (0.3,-0.3) \botlabel{j} -- (-0.3,0.3);
\region{0.5,0}{-\lambda};
\end{tikzpicture},&
\begin{tikzpicture}[Q,centerzero,scale=.9]
\draw[to-] (-0.3,-0.3) \botlabel{i} -- (0.3,0.3);
\draw[-to] (0.3,-0.3) \botlabel{j} -- (-0.3,0.3);
\region{0.4,0}{\lambda};
\end{tikzpicture}
&\mapsto -r_{j,i}\ \begin{tikzpicture}[Q,centerzero,scale=.9]
\draw[-to] (-0.3,-0.3) \botlabel{i} -- (0.3,0.3);
\draw[to-] (0.3,-0.3) \botlabel{j} -- (-0.3,0.3);
\region{0.5,0}{-\lambda};
\end{tikzpicture},\\
\begin{tikzpicture}[Q,centerzero]
\draw[-to] (-0.25,0.15) \toplabel{i} to[out=-90,in=-90,looseness=3] (0.25,0.15);
\region{0.45,-.1}{\lambda};
\end{tikzpicture}
&\mapsto \begin{tikzpicture}[Q,centerzero]
\draw[to-] (-0.25,0.15) \toplabel{i} to[out=-90,in=-90,looseness=3] (0.25,0.15);
\region{0.55,-.1}{-\lambda};
\end{tikzpicture},&
\begin{tikzpicture}[Q,centerzero]
\draw[-to] (-0.25,-0.15) \botlabel{i} to [out=90,in=90,looseness=3](0.25,-0.15);
\region{0.45,0.1}{\lambda};
\end{tikzpicture}
&\mapsto\begin{tikzpicture}[Q,centerzero]
\draw[to-] (-0.25,-0.15) \botlabel{i} to [out=90,in=90,looseness=3](0.25,-0.15);
\region{0.55,0.1}{-\lambda};
\end{tikzpicture},\\
\begin{tikzpicture}[Q,centerzero]
\draw[to-] (-0.25,-0.15) \botlabel{i} to [out=90,in=90,looseness=3](0.25,-0.15);
\region{0.45,0.1}{\lambda};
\end{tikzpicture}
&\mapsto \begin{tikzpicture}[Q,centerzero]
\draw[-to] (-0.25,-0.15) \botlabel{i} to [out=90,in=90,looseness=3](0.25,-0.15);
\region{0.55,0.1}{-\lambda};
\end{tikzpicture},&
\begin{tikzpicture}[Q,centerzero]
\draw[to-] (-0.25,0.15) \toplabel{i} to[out=-90,in=-90,looseness=3] (0.25,0.15);
\region{0.45,-.1}{\lambda};
\end{tikzpicture}
&\mapsto \begin{tikzpicture}[Q,centerzero]
\draw[to-] (-0.25,-0.15) \botlabel{i} to [out=90,in=90,looseness=3](0.25,-0.15);
\region{0.55,0.1}{-\lambda};
\end{tikzpicture},\\
\begin{tikzpicture}[baseline=-1mm,Q]
\draw[to-] (-0.25,0) arc(180:-180:0.25);
\node at (0,-.4) {\strandlabel{i}};
\region{0.95,0}{\lambda};
\node at (.55,0) {$(u)$};
\end{tikzpicture}
&\mapsto \begin{tikzpicture}[baseline=-1mm,Q]
\draw[to-] (-0.25,0) arc(-180:180:0.25);
\node at (0,-.4) {\strandlabel{i}};
\region{1.05,0}{-\lambda};
\node at (.55,0) {$(u)$};
\end{tikzpicture},&
\begin{tikzpicture}[baseline=-1mm,Q]
\draw[-to] (-0.25,0) arc(180:-180:0.25);
\node at (0,-.4) {\strandlabel{i}};
\region{0.95,0}{\lambda};
\node at (.55,0) {$(u)$};
\end{tikzpicture}
&\mapsto\begin{tikzpicture}[baseline=-1mm,Q]
\draw[-to] (-0.25,0) arc(-180:180:0.25);
\node at (0,-.4) {\strandlabel{i}};
\region{1.05,0}{-\lambda};
\node at (.55,0) {$(u)$};
\end{tikzpicture}.
\end{align*}
Roughly, this reverses orientations of strings and negates all 2-cell labels (with some additional scalars arising from crossings).
Beware that $\Omega \neq \Omega^{-1}$ unless
$r_{i,j} = r_{j,i}$ for all $i,j \in I$.

There is also an isomorphism
\begin{equation}\label{Sigmainv}
\Sigma:\UU^\rev\stackrel{\sim}{\rightarrow} \UU
\end{equation}
which is the 2-categorical analog of  the symmetry $\sigma$ from \cref{sigmainv}.
This is defined on objects by $\lambda \mapsto -\lambda$, on generating 1-morphisms by
$E_i \one_\lambda \mapsto \one_{-\lambda} E_i$
and $F_i \one_\lambda \mapsto \one_{-\lambda} F_i$, and on generating 2-morphisms by
\begin{align*}
\begin{tikzpicture}[Q,centerzero]
\draw[-to] (0,-0.3) \botlabel{i} -- (0,0.3);
\opendot{0,0};
\region{0.3,0}{\lambda};
\end{tikzpicture}
&\mapsto \begin{tikzpicture}[Q,centerzero]
\draw[-to] (0,-0.3) \botlabel{i} -- (0,0.3);
\opendot{0,0};
\region{-0.4,0}{-\lambda};
\end{tikzpicture},&
\begin{tikzpicture}[Q,centerzero]
\draw[to-] (0,-0.3) \botlabel{i} -- (0,0.3);
\opendot{0,0};
\region{0.3,0}{\lambda};
\end{tikzpicture}
&\mapsto\begin{tikzpicture}[Q,centerzero]
\draw[to-] (0,-0.3) \botlabel{i} -- (0,0.3);
\opendot{0,0};
\region{-0.4,0}{-\lambda};
\end{tikzpicture},\\
\begin{tikzpicture}[Q,centerzero,scale=.9]
\draw[-to] (-0.3,-0.3) \botlabel{i} -- (0.3,0.3);
\draw[-to] (0.3,-0.3) \botlabel{j} -- (-0.3,0.3);
\region{0.4,0}{\lambda};
\end{tikzpicture}
&\mapsto-\ \begin{tikzpicture}[Q,centerzero,scale=.9]
\draw[-to] (-0.3,-0.3) \botlabel{j} -- (0.3,0.3);
\draw[-to] (0.3,-0.3) \botlabel{i} -- (-0.3,0.3);
\region{-0.55,0}{-\lambda};
\end{tikzpicture},&
\begin{tikzpicture}[Q,centerzero,scale=.9]
\draw[to-] (-0.3,-0.3) \botlabel{i} -- (0.3,0.3);
\draw[to-] (0.3,-0.3) \botlabel{j} -- (-0.3,0.3);
\region{0.4,0}{\lambda};
\end{tikzpicture}
&\mapsto-\ \begin{tikzpicture}[Q,centerzero,scale=.9]
\draw[to-] (-0.3,-0.3) \botlabel{j} -- (0.3,0.3);
\draw[to-] (0.3,-0.3) \botlabel{i} -- (-0.3,0.3);
\region{-0.55,0}{-\lambda};
\end{tikzpicture},\\
\begin{tikzpicture}[Q,centerzero,scale=.9]
\draw[-to] (-0.3,-0.3) \botlabel{i} -- (0.3,0.3);
\draw[to-] (0.3,-0.3) \botlabel{j} -- (-0.3,0.3);
\region{0.4,0}{\lambda};
\end{tikzpicture}
&\mapsto -r_{j,i}^{-1}\ \begin{tikzpicture}[Q,centerzero,scale=.9]
\draw[to-] (-0.3,-0.3) \botlabel{j} -- (0.3,0.3);
\draw[-to] (0.3,-0.3) \botlabel{i} -- (-0.3,0.3);
\region{-0.55,0}{-\lambda};
\end{tikzpicture},&
\begin{tikzpicture}[Q,centerzero,scale=.9]
\draw[to-] (-0.3,-0.3) \botlabel{i} -- (0.3,0.3);
\draw[-to] (0.3,-0.3) \botlabel{j} -- (-0.3,0.3);
\region{0.4,0}{\lambda};
\end{tikzpicture}
&\mapsto-r_{i,j}\ \begin{tikzpicture}[Q,centerzero,scale=.9]
\draw[-to] (-0.3,-0.3) \botlabel{j} -- (0.3,0.3);
\draw[to-] (0.3,-0.3) \botlabel{i} -- (-0.3,0.3);
\region{-0.55,0}{-\lambda};
\end{tikzpicture}),\\
\begin{tikzpicture}[Q,centerzero]
\draw[-to] (-0.25,0.15) \toplabel{i} to[out=-90,in=-90,looseness=3] (0.25,0.15);
\region{0.45,-.1}{\lambda};
\end{tikzpicture}
&\mapsto \begin{tikzpicture}[Q,centerzero]
\draw[to-] (-0.25,0.15)  to[out=-90,in=-90,looseness=3] (0.25,0.15)\toplabel{i};
\region{-0.55,-.1}{-\lambda};
\end{tikzpicture},&
\begin{tikzpicture}[Q,centerzero]
\draw[-to] (-0.25,-0.15) \botlabel{i} to [out=90,in=90,looseness=3](0.25,-0.15);
\region{0.45,0.1}{\lambda};
\end{tikzpicture}
&\mapsto 
\begin{tikzpicture}[Q,centerzero]
\draw[to-] (-0.25,-0.15)  to [out=90,in=90,looseness=3](0.25,-0.15)\botlabel{i};
\region{-0.55,0.1}{-\lambda};
\end{tikzpicture},\\
\begin{tikzpicture}[Q,centerzero]
\draw[to-] (-0.25,-0.15) \botlabel{i} to [out=90,in=90,looseness=3](0.25,-0.15);
\region{0.45,0.1}{\lambda};
\end{tikzpicture}
&\mapsto \begin{tikzpicture}[Q,centerzero]
\draw[-to] (-0.25,-0.15)  to [out=90,in=90,looseness=3](0.25,-0.15)\botlabel{i};
\region{-0.55,0.1}{-\lambda};
\end{tikzpicture},&
\begin{tikzpicture}[Q,centerzero]
\draw[to-] (-0.25,0.15) \toplabel{i} to[out=-90,in=-90,looseness=3] (0.25,0.15);
\region{0.45,-.1}{\lambda};
\end{tikzpicture}&\mapsto
\begin{tikzpicture}[Q,centerzero]
\draw[-to] (-0.25,0.15) to[out=-90,in=-90,looseness=3] (0.25,0.15)\toplabel{i} ;
\region{-0.55,-.1}{-\lambda};
\end{tikzpicture},\\
\begin{tikzpicture}[baseline=-1mm,Q]
\draw[to-] (-0.25,0) arc(180:-180:0.25);
\node at (0,-.4) {\strandlabel{i}};
\region{0.95,0}{\lambda};
\node at (.55,0) {$(u)$};
\end{tikzpicture}
&\mapsto
\begin{tikzpicture}[baseline=-1mm,Q]
\draw[-to] (-0.25,0) arc(180:-180:0.25);
\node at (0,-.4) {\strandlabel{i}};
\region{-.6,0}{-\lambda};
\node at (.55,0) {$(u)$};
\end{tikzpicture},&
\begin{tikzpicture}[baseline=-1mm,Q]
\draw[-to] (-0.25,0) arc(180:-180:0.25);
\node at (0,-.4) {\strandlabel{i}};
\region{0.95,0}{\lambda};
\node at (.55,0) {$(u)$};
\end{tikzpicture}
&\mapsto
\begin{tikzpicture}[baseline=-1mm,Q]
\draw[to-] (-0.25,0) arc(180:-180:0.25);
\node at (0,-.4) {\strandlabel{i}};
\region{-.6,0}{-\lambda};
\node at (.55,0) {$(u)$};
\end{tikzpicture}.
\end{align*}
Roughly, this reflects string diagrams in a vertical axis  (with some additional scalars arising from crossings).

Finally, let $\bar\UU$ be the 2-quantum group defined from the same Cartan datum as $\UU$ but replacing
$Q_{i,j}(x,y)$ with $\bar Q_{i,j}(x,y) := r_{i,j} r_{j,i} Q_{i,j}(-x,-y)$ and 
$c_i$ with $\bar c_i:X \rightarrow \kk^\times,\lambda \mapsto (-1)^{h_i(\lambda)} c_i(\lambda)$. The scalars $r_{i,j}$ 
become $\bar r_{i,j} := r_{j,i}^{-1}$.
Then there is an isomorphism of graded 2-categories
\begin{equation}\label{newyearseve}
J:\bar\UU \stackrel{\sim}{\rightarrow}\UU
\end{equation}
which is the identity on objects and 1-morphisms, and is defined on 2-morphisms by
\begin{align*}
\begin{tikzpicture}[Q,centerzero]
\draw[-to] (0,-0.3) \botlabel{i} -- (0,0.3);
\opendot{0,0};
\region{0.3,0}{\lambda};
\end{tikzpicture}
&\mapsto -\ \begin{tikzpicture}[Q,centerzero]
\draw[-to] (0,-0.3) \botlabel{i} -- (0,0.3);
\opendot{0,0};
\region{0.3,0}{\lambda};
\end{tikzpicture},&
\begin{tikzpicture}[Q,centerzero]
\draw[to-] (0,-0.3) \botlabel{i} -- (0,0.3);
\opendot{0,0};
\region{0.3,0}{\lambda};
\end{tikzpicture}
&\mapsto
-\ \begin{tikzpicture}[Q,centerzero]
\draw[to-] (0,-0.3) \botlabel{i} -- (0,0.3);
\opendot{0,0};
\region{0.3,0}{\lambda};
\end{tikzpicture},\\
\begin{tikzpicture}[Q,centerzero,scale=.9]
\draw[-to] (-0.3,-0.3) \botlabel{i} -- (0.3,0.3);
\draw[-to] (0.3,-0.3) \botlabel{j} -- (-0.3,0.3);
\region{0.4,0}{\lambda};
\end{tikzpicture}
&\mapsto -r_{i,j}\ \begin{tikzpicture}[Q,centerzero,scale=.9]
\draw[-to] (-0.3,-0.3) \botlabel{i} -- (0.3,0.3);
\draw[-to] (0.3,-0.3) \botlabel{j} -- (-0.3,0.3);
\region{0.4,0}{\lambda};
\end{tikzpicture},&
\begin{tikzpicture}[Q,centerzero,scale=.9]
\draw[to-] (-0.3,-0.3) \botlabel{i} -- (0.3,0.3);
\draw[to-] (0.3,-0.3) \botlabel{j} -- (-0.3,0.3);
\region{0.4,0}{\lambda};
\end{tikzpicture}
&\mapsto -r_{j,i}^{-1}\ \begin{tikzpicture}[Q,centerzero,scale=.9]
\draw[to-] (-0.3,-0.3) \botlabel{i} -- (0.3,0.3);
\draw[to-] (0.3,-0.3) \botlabel{j} -- (-0.3,0.3);
\region{0.4,0}{\lambda};
\end{tikzpicture},\\
\begin{tikzpicture}[Q,centerzero,scale=.9]
\draw[-to] (-0.3,-0.3) \botlabel{i} -- (0.3,0.3);
\draw[to-] (0.3,-0.3) \botlabel{j} -- (-0.3,0.3);
\region{0.4,0}{\lambda};
\end{tikzpicture}
&\mapsto -r_{j,i}\ \begin{tikzpicture}[Q,centerzero,scale=.9]
\draw[-to] (-0.3,-0.3) \botlabel{i} -- (0.3,0.3);
\draw[to-] (0.3,-0.3) \botlabel{j} -- (-0.3,0.3);
\region{0.4,0}{\lambda};
\end{tikzpicture},&
\begin{tikzpicture}[Q,centerzero,scale=.9]
\draw[to-] (-0.3,-0.3) \botlabel{i} -- (0.3,0.3);
\draw[-to] (0.3,-0.3) \botlabel{j} -- (-0.3,0.3);
\region{0.4,0}{\lambda};
\end{tikzpicture}
&\mapsto -r_{i,j}^{-1}\ \begin{tikzpicture}[Q,centerzero,scale=.9]
\draw[to-] (-0.3,-0.3) \botlabel{i} -- (0.3,0.3);
\draw[-to] (0.3,-0.3) \botlabel{j} -- (-0.3,0.3);
\region{0.4,0}{\lambda};
\end{tikzpicture},\\
\begin{tikzpicture}[Q,centerzero]
\draw[-to] (-0.25,0.15) \toplabel{i} to[out=-90,in=-90,looseness=3] (0.25,0.15);
\region{0.45,-.1}{\lambda};
\end{tikzpicture}
&\mapsto \begin{tikzpicture}[Q,centerzero]
\draw[-to] (-0.25,0.15) \toplabel{i} to[out=-90,in=-90,looseness=3] (0.25,0.15);
\region{0.45,-.1}{\lambda};
\end{tikzpicture},&
\begin{tikzpicture}[Q,centerzero]
\draw[-to] (-0.25,-0.15) \botlabel{i} to [out=90,in=90,looseness=3](0.25,-0.15);
\region{0.45,0.1}{\lambda};
\end{tikzpicture}
&\mapsto \begin{tikzpicture}[Q,centerzero]
\draw[-to] (-0.25,-0.15) \botlabel{i} to [out=90,in=90,looseness=3](0.25,-0.15);
\region{0.45,0.1}{\lambda};
\end{tikzpicture},\\
\begin{tikzpicture}[Q,centerzero]
\draw[to-] (-0.25,-0.15) \botlabel{i} to [out=90,in=90,looseness=3](0.25,-0.15);
\region{0.45,0.1}{\lambda};
\end{tikzpicture}
&\mapsto -\ \begin{tikzpicture}[Q,centerzero]
\draw[to-] (-0.25,-0.15) \botlabel{i} to [out=90,in=90,looseness=3](0.25,-0.15);
\region{0.45,0.1}{\lambda};
\end{tikzpicture},&
\begin{tikzpicture}[Q,centerzero]
\draw[to-] (-0.25,0.15) \toplabel{i} to[out=-90,in=-90,looseness=3] (0.25,0.15);
\region{0.45,-.1}{\lambda};
\end{tikzpicture}&\mapsto-\ \begin{tikzpicture}[Q,centerzero]
\draw[to-] (-0.25,0.15) \toplabel{i} to[out=-90,in=-90,looseness=3] (0.25,0.15);
\region{0.45,-.1}{\lambda};
\end{tikzpicture}
,\\
\begin{tikzpicture}[baseline=-1mm,Q]
\draw[to-] (-0.25,0) arc(180:-180:0.25);
\node at (0,-.4) {\strandlabel{i}};
\region{0.95,0}{\lambda};
\node at (.54,0) {$(u)$};
\end{tikzpicture}
&\mapsto(-1)^{h_i(\lambda)}\ \begin{tikzpicture}[baseline=-1mm,Q]
\draw[to-] (-0.25,0) arc(180:-180:0.25);
\node at (0,-.4) {\strandlabel{i}};
\region{1.3,0}{\lambda};
\node at (.68,0) {$(-u)$};
\end{tikzpicture}
,&
\begin{tikzpicture}[baseline=-1mm,Q]
\draw[-to] (-0.25,0) arc(180:-180:0.25);
\node at (0,-.4) {\strandlabel{i}};
\region{0.95,0}{\lambda};
\node at (.54,0) {$(u)$};
\end{tikzpicture}
&\mapsto (-1)^{h_i(\lambda)}\ 
\begin{tikzpicture}[baseline=-1mm,Q]
\draw[-to] (-0.25,0) arc(180:-180:0.25);
\node at (0,-.4) {\strandlabel{i}};
\region{1.3,0}{\lambda};
\node at (.68,0) {$(-u)$};
\end{tikzpicture}.
\end{align*}

%% file: s3-iquantum.tex
\setcounter{section}{2}

\section{New 2-categories for quasi-split
  \texorpdfstring{$\mathrm{i}$}{}quantum groups}\label{s3-iquantum}

We continue with the Cartan datum and parameters fixed in the previous section.
We are now going to propose an analogous framework for categorification of quasi-split iquantum groups arising from quantum symmetric pairs. 
 
\subsection{The iquantum group \texorpdfstring{$\U^\imath(\del)$}{}}\label{iqg}

From a quasi-split Satake diagram of the same Cartan type as $\U$, we obtain an
involution $\tau$ of the set $I$
with $a_{i,j} = a_{\tau i,\tau j}$ and $d_i = d_{\tau i}$. 
We allow $\tau =\id$. 
Occasionally, it will be useful to have chosen a decomposition
\begin{equation}\label{usefuldecomposition}
I = I_1 \sqcup I_0 \sqcup I_{-1}
\end{equation}
such that $I_1$ is a set of representatives for the $\tau$-orbits of size 2, $I_0$ consists of the $\tau$-fixed points, and
$I_{-1} = \tau(I_1)$.
We pick
$\del = (\del_i)_{i \in I} \in \Z^I$
so that
\begin{align}
\label{recap}
\del_i + \del_{\tau i} &= -a_{i,\tau i}
\qquad\text{and}\qquad 
\text{$\del_i \geq 0$ if $i \neq \tau i$.}
\end{align}
If $i = \tau i$ then $\del_i=-1$, and if $a_{i,\tau i} = 0$ then $\del_i=0$.

We also assume there are 
involutions $\tau:X \rightarrow X$ and $\tau^*:Y\rightarrow Y$ such that
$\tau^*(h)(\lambda) = h (\tau(\lambda))$
for all $\lambda \in X$ and $h \in Y$,
and also
$\tau (\alpha_i) = \alpha_{\tau i}$ 
and $\tau^*(h_i) = h_{\tau i}$ for each $i \in I$. 
The {\em iweight lattice} and {\em icoweight lattice} are the Abelian groups
\begin{align*}
X^\imath &:= X \big/ 
\im (\id+\tau),
&
Y^\imath &:= \ker (\id+\tau^*).
\end{align*} 
When $\lambda,\mu,\dots$ are iweights in the quotient $X^\imath$, we will use the notation
$\hat \lambda,\hat \mu,\dots$ to denote a pre-image in $X$.
For $\lambda \in X^\imath$ and $\pi \in X$, 
we write $\lambda +\pi$
for the sum of $\lambda$ and the canonical image of $\pi$ in $X^\imath$.
For $\lambda \in X^\imath$ with pre-image $\hat\lambda \in X$ and $i \in I$, we define
\begin{align}
\lambda_i &:= (h_i - h_{\tau i})(\hat\lambda) \in \Z.
\end{align}
This is well defined independent of the choice of the pre-image $\hat\lambda$.
It is also useful to note that
\begin{align}\label{somewhatannoying}
(\lambda+\alpha_j)_i &= \lambda_i + a_{i,j} - a_{i,\tau j},\\
\label{dontmentionit}
\del_{\tau i}-\lambda_{\tau i} &= 
\lambda_i-\del_i
-a_{i,\tau i}.
\end{align}
If $\tau i = i$ then $\lambda_i = 0$. In this case, the {\em parity} of $h_i(\hat\lambda)$ is well defined independent of the choice of the pre-image $\hat\lambda$, so it makes sense to set
\begin{equation}\label{stillhaveparity}
t_i(\lambda) := \begin{cases}
0&\text{if $h_i(\hat\lambda)$ is even}\\
1&\text{if $h_i(\hat\lambda)$ is odd,}
\end{cases}
\end{equation}
for $i \in I$ with $i=\tau i$.

The {\em iquantum group}
associated to the above data is the subalgebra
$\U^\imath = \U^\imath(\del)$ of $\U$
generated by the elements
$q^h\:(h \in Y^\imath)$
and 
\begin{equation}\label{unpredictable}
b_i :=  f_i + q_i^{\del_i} e_{\tau i} k_i^{-1}
\end{equation}
for each $i \in I$.
It is a (right) coideal subalgebra of $\U$, and 
the pair $(\U, \U^\imath)$ is a 
{\em quantum symmetric pair} of quasi-split type. 

\begin{rem}\label{notthebest}
The formula \cref{unpredictable} is essentially the ``standard embedding'' from \cite{Kolb}, but 
it is not the most general embedding
$\U^\imath(\del)\hookrightarrow \U$: one can take
$$
b_i :=  f_i + q_i^{\del_i} e_{\tau i} k_i^{-1}
+
\textstyle{\frac{q_i^{\mu_i} - q_i^{-\mu_i}}{q_i-q_i^{-1}}} k_i^{-1}
$$
for some $\mu = (\mu_i)_{i \in I} \in \Z^I$
such that $\mu_i = 0$ either if $i \neq \tau i$,
or if there exists $j \in I$ with $j = \tau j$ and $a_{i,j}\equiv 1 \pmod{2}$.
\end{rem}

We will need the {\em modified form}
$\dot\U^\imath=\dot\U^\imath(\del)$
for $\U^\imath$, which was introduced in \cite{BW18QSP} (see also \cite[Sec.~3.5]{BW21iCB}).
This can be viewed as a
subalgebra of the completion 
$\hat \U$ of the modified form of $\U$.
Given $\lambda \in X^\imath$, we let
$1_\lambda := \sum_{\hat\lambda}1_{\hat\lambda} \in \hat \U$
summing over all pre-images
$\hat\lambda \in X$ of $\lambda$.
Then 
$$
\dot\U^\imath := \bigoplus_{\lambda,\mu \in X^\imath} 1_\lambda \U^\imath 1_\mu \subset \hat \U.
$$ 
Although rarely a unital subalgebra of $\hat \U$,
it is a locally unital algebra in its own right with distinguished idempotents $\{1_\lambda\:|\:\lambda \in X^\imath\}$. As such,
it is generated by the elements
$b_i 1_\lambda = 1_{\lambda - \alpha_i} b_i$
for all $i \in I$ and $\lambda \in X^\imath$; this depends on the observation that
$-\alpha_i$ and $\alpha_{\tau i}$ have the same image in $X^\imath$.

Following \cite{BW18KL,BeW18},
we define elements $b_i^{(n)} 1_\lambda =1_{\lambda-n\alpha_i} b_i^{(n)}\in \dot\U^\imath$, 
which we call {\em divided powers} if $i \neq \tau i$
or {\em idivided powers} if $i = \tau i$:
\begin{equation}\label{idividedpowerrelation}
b_i^{(n)} 1_\lambda :=
\begin{dcases}
\frac{1}{[n]_{q_i}^!}
\:b_i^n 1_\lambda&\text{if $i \neq \tau i$}\\
\frac{1}{[n]_{q_i}^!}
\prod_{\substack{m=0\\m\equiv h_i(\hat\lambda)\!\!\!\!\!\pmod{2}}}^{n-1}
\!\!\!\left(b_i^2 - [m]_{q_i}^2\right) 1_\lambda
&\text{if $i = \tau i$ and $n$ is even}\\
\frac{1}{[n]_{q_i}^!}
\:b_i\!\!\!\!\!\!\!
\prod_{\substack{m=1\\m\equiv h_i(\hat\lambda)\!\!\!\!\!\pmod{2}}}^{n-1}
\!\!\!\left(b_i^2 - [m]_{q_i}^2\right) 1_\lambda
&\text{if $i = \tau i$ and $n$ is odd.}\\
\end{dcases}
\end{equation}
Let $\dot\U_\Z^\imath = \dot\U_\Z^\imath(\del) \subset \widehat{\U}_\Z$ be the $\Z[q,q^{-1}]$-form
for $\dot\U^\imath$ generated by 
$b_i^{(n)} 1_\lambda\:(\lambda \in X^\imath, i \in I, n \geq 0)$. 
Again, this was introduced originally in \cite{BW18QSP}.

There is an explicit presentation for
$\dot\U^\imath$, which is implied by the presentation for $\U^\imath$ derived in 
\cite[Th.~3.1]{serre} (see also \cite[Th.~3.6]{BK}):
it is the locally unital algebra with
(mutually orthogonal) distinguished
idempotents $\{1_\lambda\:|\:\lambda \in X^\imath\}$ and generators $b_i 1_\lambda = 1_{\lambda- \alpha_i} b_i\:(\lambda \in X^\imath, i \in I)$
subject just to the {\em iSerre relations}
\begin{align}
\sum_{n=0}^{1-a_{i,j}}(-1)^n
b_i^{(n)} b_j b_i^{(1-a_{i,j}-n)} 1_\lambda
=
\delta_{i,\tau j}\!
\prod_{r=1}^{-a_{i,j}}
(q_i^r-q_i^{-r})\times
\frac{(-1)^{a_{i,j}}q_i^{\lambda_i-\del_i-\binom{a_{i,j}}{2}}-q_i^{
\binom{a_{i,j}}{2}+\del_i-\lambda_i}}{q_i-q_i^{-1}}
b_i^{(-a_{i,j})}
1_\lambda
\label{secretary}
\end{align}
for all $i \neq j$ in $I$ and $\lambda \in X^\imath$,
where $b_i^{(n)} 1_\lambda$ is as defined by \cref{idividedpowerrelation}.
If $a_{i,\tau i}=0$ (when $\del_i = 0$) then this relation
implies that
$[b_{\tau i}, b_i] 1_\lambda = [\lambda_i]_{q_i} 1_\lambda$.
If $a_{i,\tau i} = -1$ then it gives that
$(b_{\tau i} b_i^{(2)}
- b_i b_{\tau i} b_i
+ 
b_i^{(2)} b_{\tau i})1_\lambda
= - (q_i^{\lambda_i-\del_i-1}+q_i^{1+\del_i-\lambda_i}) b_i 1_\lambda.$

Next, we introduce some symmetries.
Define ${'}\del = ({'}\del_i)_{i \in I}$ 
by ${'}\del_i := \del_{\tau i}$.
Then there are maps
\begin{align}
\label{newyearsday2}
\text{(linear isomorphism)}\quad
\omega^\imath:\dot\U_\Z^\imath(\del)&\stackrel{\sim}{\rightarrow}\dot\U_\Z^\imath({'}\del),
&
b_i^{(n)} 1_\lambda&\mapsto
b_{\tau i}^{(n)} 1_{-\lambda},\\\label{snoring}
\text{(anti-linear involution)}\quad
\psi^\imath:\dot\U_\Z^\imath(\del) &\stackrel{\sim}{\rightarrow} \dot\U_\Z^\imath(\del),
&b_i^{(n)} 1_\lambda &\mapsto b_i^{(n)} 1_{\lambda},\\\label{newyearsday1}
\text{(linear anti-isomorphism)}\quad
\sigma^\imath:\dot\U_\Z^\imath(\del)&\stackrel{\sim}{\rightarrow}\dot\U_\Z^\imath({'}\del),
&
b_i^{(n)} 1_\lambda&\mapsto
1_{-\lambda} b_{i}^{(n)}.
\end{align}
The existence of these is a relations check using \cref{secretary}.

\begin{lem}\label{lazydog}
Suppose we are given two families of parameters
$(\del_i)_{i \in I}$
and $(\del^\dagger_i)_{i \in I}$
satisfying \cref{recap} and an element $\pi \in X^\imath$
such that $\pi_i = \del^\dagger_i-\del_i$ for each $i \in I$.  
Then
there is a $\Q(q)$-algebra isomorphism
\begin{align*}
s:\dot\U^\imath(\del) 
&\stackrel{\sim}{\rightarrow} \dot\U^\imath(\del^\dagger),&
b_i 1_\lambda & \mapsto b_i 1_{\lambda+\pi}.
\end{align*}
It restricts to an isomorphism
$\dot\U_\Z^\imath(\del)\stackrel{\sim}{\rightarrow}
\dot\U_\Z^\imath(\del^\dagger)$ between the integral forms.
\end{lem}

\begin{proof}
Note that such an iweight $\pi$ exists providing $\del_i = \del^\dagger_i$ for all but finitely many $i \in I$;
for example, we could take $\pi$ to be the image in $X^\imath$ of $\sum_{i \in I_1} (\del^\dagger_i-\del_i) \varpi_i$.
To construct the homomorphism $s$, we use the 
description by generators and relations.
There is only the one relation \cref{secretary}, 
and the argument reduces easily to checking that
\begin{equation}\label{whatpiisfor}
\lambda_i-\del_i =  (\lambda+\pi)_i-\del^\dagger_i,
\end{equation}
which is clear.
The homomorphism is an isomorphism because it
has a two-sided inverse mapping
$b_i 1_\lambda \mapsto b_i 1_{\lambda-\pi}$.
It restricts to an isomorphism between integral forms because it takes
$b_i^{(n)} 1_\lambda \mapsto b_i^{(n)} 1_{\lambda+\pi}$ for all $n \geq 1$.
\end{proof}

\cref{lazydog}
justifies omitting $\del$ from the notation $\dot\U_\Z^\imath(\del)$, although its embedding into $\widehat{\U}_\Z$ does depend on this choice.

\subsection{Definition of the 2-iquantum group \texorpdfstring{$\UU^\imath(\del,\zeta)$}{}}\label{new2cat}

As usual, $\kk$ is an $\N$-graded commutative ring.
If $a_{i,\tau i} \neq 0$ for some $i \in I$ then we also require that 2 is invertible in $\kk$.
In the previous subsection, we fixed a choice of Satake datum 
and parameters $\del = (\del_i)_{i \in I} \in \Z^I$ satisfying \cref{recap}, as needed to define the
quasi-split iquantum group $\U^\imath$.
We also choose parameters
$Q_{i,j}(x,y) \in \kk[x,y]$ as in \cref{generic}
such that $t_{i,j} := Q_{i,j}(1,0) \in \kk^\times$ for all $i \neq j$ in $I$, and normalization homomorphisms 
$c_i:X\rightarrow \kk^\times$ 
as in \cref{great} such that $c_i(\alpha_i) = 1$ for all $i \in I$.
Recall also the definitions of $r_{i,j} \in \kk^\times$
and $R_{i,j}(x,y) \in \kk(x,y)$ from \cref{gollygosh,boshgosh}.

\begin{hyp}\label{hypotheses}
We make four additional assumptions about the parameter choices, which
ensure some compatibility with the involution $\tau$:
we require that
\begin{align}\label{lizka}
c_{\tau i}(\tau \lambda) &= (-1)^{h_i(\lambda)} c_i(\lambda),\\
c_{\tau i}(\alpha_i) &= t_{\tau i, i}\quad\text{(equivalently, $r_{\tau i, i}=1$)},\label{luckya}\\
\label{qid}
x^{a_{i, \tau i}} 
Q_{\tau i, i}(x,-x) &\in \kk^\times\\\intertext{for all $\lambda \in X$ and $i \in I$ with $i \neq \tau i$,
and that}
\label{tauantisymmetric}
Q_{\tau i,\tau j}(x,y) &= 
r_{i,j} r_{j,i} r_{i,\tau j} r_{j,\tau i} Q_{i, j}(-x,-y)
\end{align}
for all $i,j \in I$.
\end{hyp}

\begin{rem}\label{unfortunately}
The assumption \cref{tauantisymmetric}
imposes some 
restrictions on the root datum.
Suppose that there exist $i \neq j$ in $I$ with $i = \tau i$ and $j = \tau j$. Taking $x^{-a_{i,j}}$-coefficients in \cref{tauantisymmetric} gives 
\begin{equation}\label{rocketman}
(r_{i,j} r_{j,i})^2 = (-1)^{a_{i,j}}.
\end{equation}
Since the left hand side is symmetric in $i$ and $j$, this implies that
$(-1)^{a_{i,j}} = (-1)^{a_{j,i}}$, hence, we must have that $a_{i,j} \equiv a_{j,i}\pmod{2}$.
\end{rem}

We also fix a choice of invertible scalars
$\zeta = (\zeta_i)_{i \in I}$ such that
\begin{align}\label{being}
\zeta_i&=-\tfrac{1}{2}\text{ if $i = \tau i$,}
&
Q_{\tau i,i}(x,-x) &= (-1)^{\del_i+1} \zeta_i \zeta_{\tau i}
x^{-a_{i,\tau i}}
\text{ if $i \neq \tau i$.}
\end{align}
Such scalars exist due to the assumption \cref{qid}.
All choices have now been made, but there is a little bit of additional notation.
Let $\gamma_i:X^\imath \rightarrow \kk^\times$ be the group homomorphism defined by
\begin{equation}\label{purer}
\gamma_i(\lambda) := 
(-1)^{h_{\tau i}(\hat\lambda)}
\frac{c_i(\hat\lambda)}{c_{\tau i}(\hat\lambda)}
\end{equation}
for $\lambda \in X^\imath$ and any pre-image $\hat \lambda \in X$.
The assumption that $c_i$ is a homomorphism implies that
this is well defined
independent of the choice of $\hat\lambda$.
Note also that
\begin{equation}\label{howthingswork}
\gamma_i(\lambda)\gamma_{\tau i}(\lambda) = (-1)^{\lambda_i}
\end{equation}
for $i \in I$ and $\lambda \in X^\imath$,
and
\begin{equation}\label{morethingsworking}
\gamma_i(\lambda) = (-1)^{t_i(\lambda)}
\end{equation}
if $i = \tau i$.
We use the shorthand
\begin{equation}\label{childhood}
\zeta_i(\lambda) := \zeta_i \gamma_i(\lambda).
\end{equation}
The function $\zeta_i:X^\imath \rightarrow \kk^\times$ is
{\em not} a homomorphism of Abelian groups.
Finally, we use the additional shorthand
\begin{equation}\label{modifiedsigns}
Q^\imath_{i,j}(x,y) := (-1)^{\delta_{i,\tau j}} Q_{i,j}(x,y)
\end{equation}
in order to draw attention to an important but annoying sign\footnote{There was an egregious typo in the first $\tt{arXiv}$ version of the article in this definition!}.

Now comes the main definition.
We emphasize that the 2-category $\UU^\imath(\del,\zeta)$ that it introduces
depends not only on $\del$ and $\zeta$, but also on the Satake datum, the choices of $Q_{i,j}(x,y)$ and the normalization homomorphisms $c_i$ satisfying \cref{hypotheses}.

\begin{defin}[2-iquantum groups]\label{def2iqg}
The {\em 2-iquantum group}
$\UU^\imath=\UU^\imath(\del,\zeta)$ with the parameters chosen above
is the graded 2-category
with object set $X^\imath$, 
generating 1-morphisms 
$$
B_i \one_\lambda = \one_{\lambda -\alpha_i} B_i
: \lambda \rightarrow \lambda-\alpha_i
$$
for $\lambda \in X^\imath$ and $i \in I$,
with identity 2-endomorphisms
denoted by unoriented strings
$\begin{tikzpicture}[iQ,anchorbase]
\draw[-] (0,-0.2)\botlabel{i} -- (0,0.2);
\region{0.2,0}{\lambda};
\region{-0.4,0}{\lambda-\alpha_i};
\end{tikzpicture}$,
and generating 2-morphisms for all $i,j \in I$ and $\lambda \in X^\imath$ as \cref{table2}.
\begin{table}
\begin{align*}
\begin{array}{|l|c|}
\hline
\hspace{13mm}\text{Generator}&\text{Degree}\\
\hline
\,\ \ \begin{tikzpicture}[iQ,centerzero]
\draw[-] (0,-0.3) \botlabel{i} -- (0,0.3) \toplabel{i};
\closeddot{0,0};
\region{0.2,0}{\lambda};
\end{tikzpicture}\ \, 
\colon B_i \one_\lambda \Rightarrow B_i \one_\lambda&2 d_i\\
\:\begin{tikzpicture}[iQ,centerzero]
\draw[-] (-0.25,-0.15) \botlabel{\tau i} to [out=90,in=90,looseness=3](0.25,-0.15) \botlabel{i};
\region{0.45,0.1}{\lambda};
\node at (0,.3) {$\phantom.$};
\node at (0,-.4) {$\phantom.$};
\end{tikzpicture}\hspace{-.03pt}
\colon B_{\tau i} B_{i} \one_\lambda \Rightarrow \one_\lambda
&d_i(1+\del_i-\lambda_i)\\
\:\begin{tikzpicture}[iQ,centerzero]
\draw[-] (-0.25,0.15) \toplabel{\tau i} to[out=-90,in=-90,looseness=3] (0.25,0.15) \toplabel{i};
\region{0.45,-.1}{\lambda};
\node at (0,.2) {$\phantom.$};\node at (0,-.3) {$\phantom.$};
\end{tikzpicture}
\colon \one_\lambda \Rightarrow B_{\tau i} B_i \one_\lambda
&d_i(1+\del_i-\lambda_i)\\
\!\!\!\begin{tikzpicture}[baseline=-1mm,iQ]
\draw[-] (-0.25,0) arc(180:-180:0.25);
\node at (-0.42,0) {\strandlabel{\tau i}};
\region{0.45,0}{\lambda};
\dottybubblelabel{0,0}{n};
\node at (0,-.3) {$\phantom{.}$};
\end{tikzpicture}
\colon \one_\lambda \Rightarrow\one_\lambda
\text{ for }0 \leq n \leq \del_i-\lambda_i
&2d_i n\\
\;\begin{tikzpicture}[iQ,centerzero,scale=.9]
\draw[-] (-0.3,-0.3) \botlabel{i} -- (0.3,0.3) \toplabel{i};
\draw[-] (0.3,-0.3) \botlabel{j} -- (-0.3,0.3) \toplabel{j};
\region{0.35,0}{\lambda};
\end{tikzpicture}
\;\,\colon B_{i} B_{j} \one_\lambda \Rightarrow B_j B_i \one_\lambda
&-d_i a_{i,j} \\
\hline
\end{array}
\end{align*}
\caption{Generating 2-morphisms of $\UU^\imath$}\label{table2}
\end{table}
The degrees of the generating 2-morphisms are also listed in this table.
The generating 2-morphisms are subject to certain relations.
Before writing these down, we explain some further conventions.
\begin{itemize}
\item 
As in the previous section, we will usually 
label strings just at one end, but note for unoriented cups and caps that the string label switches from $i$ to $\tau i$ at the critical point.
\item
As before, we will only label one of the 2-cells by the iweight $\lambda$, and if we omit all labels we are writing something that is true for all possible labels.
\item
We will now be using black (closed) dots rather than the white (open) dots used in the previous section. Pins are defined in just the same way as before.
In place of \cref{dgf}, we now need {\em two}
shorthands:
\begin{align}\label{idgf}
\ .
\end{align}
Recall also our convention when multi-variable polynomials are pinned to several dots---the alphabetical ordering of the variables $x,y,z$ corresponds to the lexicographical ordering of the Cartesian coordinates of the dots.
\item
We use the following generating functions for dotted bubbles:
\begin{align}\label{ifakebubblepolynomial}
%
\right).
\end{align}
The relation \cref{ibubblerel} below
implies that
\begin{align}\label{kitty}
 %
 \bigg]_{u:\del_i-\lambda_i-n}.
\end{align}
Thus, we have bubbles $%
$
for $1 \leq n \leq \del_i-\lambda_i$. We do not count 
the bubble $%
$ of lowest degree
as a fake bubble because it is just a 
scalar---it equals $\zeta_i(\lambda) \id_{\one_\lambda}$.
\end{itemize}
The defining relations are as follows: 
\begin{align}
\label{ibubblerel}
\left[\!\!
%
\ .
\end{align}
\end{defin}

The additional generating function in \cref{idgf} is needed because we now have that
\begin{align}\label{whatfun}
%
\ .
\end{align}
The following are consequences of \cref{idotslide}:
\begin{align}\label{genidotslide1}
%
\ .
\end{align}
Using one of these, the following general form of the curl relation \cref{icurl} can be deduced:
\begin{equation}\label{fancyicurl}
%
.
\end{equation}
In the case $i = \tau i$, the calculation proving this is explained in
the proof of \cite[Th.~2.5]{BWWbasis}; this depends on the connection
to nil-Brauer explained in \cref{rankone}.
Since it is instructive to see at this point, here is the proof in the easier case that $i \neq \tau i$:
\begin{align*}
%
\right]_{\!u:<0}.
\end{align*}

\begin{rem}\label{teapot}
Using \cref{izigzag,ipitchfork},
one can show that
\begin{align}
%
\ .\label{baldy}
\end{align}
Now the relation \cref{iquadratic} can be rotated through 180$^\circ$ by attaching a nested pair of cups and caps,
using \cref{baldy} to simplify the result, 
to obtain the identity \cref{tauantisymmetric}.
This is the reason for this assumption---it ensures that the quadratic relation is invariant under rotation.
Rotating the curl relation \cref{fancyicurl} in a similar way gives the new relation
\begin{equation}\label{icurlrotate}
%
.
\end{equation}
We leave it to the reader to check that the other defining relations of $\UU^\imath$ are invariant under rotation.
When checking the rotation invariance of \cref{ibraid}, the following are helpful:
\begin{align}
%
.
\end{align}
\end{rem}

\begin{lem}\label{lazydog2}
Let $\del$ and $\zeta$ be as above, and suppose that
$\del^\dagger = (\del^\dagger_i)_{i \in I}$ and
$\zeta^\dagger=(\zeta^\dagger_i)_{i \in I}$
is another choice for these parameters satisfying
\cref{recap,being}.
Let $\UU^\imath(\del,\zeta)$ and
$\UU^\imath(\del^\dagger,\zeta^\dagger)$
be the corresponding 2-iquantum groups
defined using the same choices of $Q_{i,j}(x,y)$ and $c_i$.
Let $\pi$ be an
iweight such that
$\pi_i = \del^\dagger_i - \del_i$ for each $i$,
and choose scalars
$\kappa_i \in \kk^\times$
such that $\zeta_i = \kappa_i^2 \zeta^\dagger_i \gamma_i(\pi)$ 
and $\kappa_i \kappa_{\tau i} = 1$
for each $i \in I$.
Then there is an isomorphism of graded 2-categories
\begin{align*}
S:\UU^\imath(\del,\zeta) 
&\stackrel{\sim}{\rightarrow} \UU^\imath(\del^\dagger,\zeta^\dagger)
\end{align*}
taking the object $\lambda$ to $\lambda +\pi$,
the 1-morphism $B_i \one_\lambda$ to $B_i \one_{\lambda+\pi}$, and defined on 
a string diagram representing a 2-morphism by 
adding $\pi$ to the labels of every 2-cell
then multiplying
by $\kappa_i$ for each
occurrence of the generator
$\begin{tikzpicture}[iQ,centerzero,scale=.8]
\draw[-] (-0.25,-0.15) \botlabel{\tau i} to [out=90,in=90,looseness=3](0.25,-0.15);
\region{0.45,0.1}{\lambda};
\end{tikzpicture}$
or the generator
$\begin{tikzpicture}[iQ,centerzero,scale=.8]
\draw[-] (-0.25,0.15) \toplabel{\tau i} to [out=-90,in=-90,looseness=3](0.25,0.15);
\region{0.45,-0.1}{\lambda};
\end{tikzpicture}$ in the diagram,
and by $\kappa_i^2$ for each fake bubble
$\begin{tikzpicture}[baseline=-.5mm,iQ]
\draw[-] (-0.25,0) arc(180:-180:0.25);
\node at (-0.42,0) {\strandlabel{\tau i}};
\region{0.45,0}{\lambda};
\dottybubblelabel{0,0}{n};
\end{tikzpicture}$.
\end{lem}

\begin{proof}
It is well defined by a straightforward verification of the defining relations.
One needs to use \cref{whatpiisfor} to see that everything makes sense.
It is an isomorphism because it has a two-sided inverse defined in a similar way, replacing each $\kappa_i$ by $\kappa_i^{-1}$.
\end{proof}

\cref{lazydog2}
justifies omitting $\del$ and $\zeta$ from the notation $\UU^\imath(\del,\zeta)$. 
Note though that the degrees of cups, caps and bubbles depend on the choice of $\del$.

\begin{defin}[Weak 2-iquantum groups]\label{extended}
We refer to the relation \cref{ibraid} in the case that
the string labels $(i,j,k)$ are $(i,\tau i,i)$ for $i \neq \tau i$
as the {\em ibraid relation}. 
The {\em weak 2-iquantum group}
$\xUUi = 
\xUUi(\del,\zeta)$ is the graded 2-category defined in almost the same way as in \cref{def2iqg} but omitting all instances of the ibraid relation.
\end{defin}

The following theorem shows that the ibraid relation holds
automatically if $a_{i,\tau i} = 0$, and it holds
up to multiplication by a polynomial when $a_{i,\tau i} < 0$.
Consequently, although it cannot 
be omitted in the definition of $\UU^\imath$, the precise form of 
the ibraid relation is determined by the other relations.

\begin{theo}\label{necessary}
In the weak 2-iquantum group $\xUUi(\del,\zeta)$, we have 
for $i \neq \tau i$ that
\begin{align}
&\begin{tikzpicture}[iQ,centerzero,scale=1.2]
\draw[-] (-0.4,-0.4) \botlabel{i} -- (0.4,0.4);
\draw[-] (0,-0.4) \botlabel{\tau i} to[out=135,in=down] (-0.32,0) to[out=up,in=225] (0,0.4);
\draw[-] (0.4,-0.4) \botlabel{i} -- (-0.4,0.4);
\pinpin{-0.12,.3}{.3,.3}{1.15,.3}{Q_{i,\tau i}(y,x)};
\end{tikzpicture}
\!-\!\!\!
\begin{tikzpicture}[iQ,centerzero,scale=1.2]
\draw[-] (-0.4,-0.4) \botlabel{i} -- (0.4,0.4);
\draw[-] (0,-0.4) \botlabel{\tau i} to[out=45,in=down] (0.32,0) to[out=up,in=-45] (0,0.4);
\draw[-] (0.4,-0.4) \botlabel{i} -- (-0.4,0.4);
\pinpin{0.13,.3}{.3,.3}{1.15,.3}{Q_{i,\tau i}(y,x)};
\end{tikzpicture}
= \notag\\
&\qquad
\begin{tikzpicture}[iQ,centerzero,scale=1.2]
\draw[-] (-0.28,-0.4) \botlabel{i} -- (-0.28,0.4);
\draw[-] (0,-0.4) \botlabel{\tau i} -- (0,0.4);
\draw[-] (0.28,-0.4) \botlabel{i} -- (0.28,0.4);
\pinpinpin{.28,0}{0,0}{-.28,0}{-1.6,0}{
\frac{Q^\imath_{i,\tau i}(x,y)-Q^\imath_{i,\tau i}(z,y)}{x-z}};
\pinpin{0,.2}{.28,.2}{1.05,.2}{Q_{i,\tau i}(y,x)};
\end{tikzpicture}+
\left[\begin{tikzpicture}[iQ,centerzero,scale=1.2]
\draw[-] (-0.2,-0.4) \botlabel{i} to [looseness=2.5,out=90,in=90] (0.2,-0.4);
\draw[-] (0.2,0.4) to [looseness=2.5,out=-90,in=-90] (-0.2,0.4)\toplabel{i};
\draw[-] (.84,-.4) \botlabel{i} to (.84,.4);
\circled{-.17,.2}{u};
\circled{-.17,-.2}{u};
\circled{.84,0}{u};
\draw[fakebubble] (0.33,0) arc(180:-180:0.15);
\node at (0.19,0) {\strandlabel{\tau i}};
\bubblelabel{0.48,0}{u};
\pinpin{.18,.25}{.85,.25}{1.6,.25}{Q_{i,\tau i}(y,x)};
\end{tikzpicture}-
\begin{tikzpicture}[iQ,centerzero,scale=1.2]
\draw[-] (-0.65,-0.4) \botlabel{i} -- (-0.65,0.4);
\draw[-] (0,0.4) \toplabel{\tau i} to[out=down,in=down,looseness=3] (0.4,0.4);
\draw[-] (0.4,-.4) to[out=up,in=up,looseness=3] (0,-.4) \botlabel{\tau i};
\circledbar{-.65,0}{u};
\circledbar{.35,.15}{u};
\circledbar{.35,-.15}{u};
\draw[fakebubble] (-0.05,0) arc(360:0:0.15);
\node at (-0.42,0) {\strandlabel{i}};
\bubblelabel{-0.2,0}{u};
\pinpin{0.005,.33}{.395,.33}{1.15,.33}{Q_{i,\tau i}(y,x)};
\end{tikzpicture}\right]_{u:-1}\notag
\\&\qquad
-\,
\begin{tikzpicture}[iQ,centerzero,scale=1.6]
\draw[-] (-0.3,0.3) -- (0.3,-0.3)\botlabel{i};
\draw[-] (0,0.3) \toplabel{\tau i}[out=-90,in=-90,looseness=2.5] to  (0.3,0.3);
\draw[-] (-0.3,-0.3) \botlabel{i}[out=90,in=90,looseness=2.5] to (0,-0.3);
\pinpinpinmid{.08,.1}{0,0}{-.08,-.1}{-1.3,0}{
\frac{Q^\imath_{\tau i, i}(x,y)-Q^\imath_{\tau i, i}(z,y)}{x-z}};
\pinpin{0.01,.2}{.29,.2}{1,.2}{Q_{i,\tau i}(y,x)};
\end{tikzpicture}
-
\begin{tikzpicture}[iQ,centerzero,scale=1.6]
\draw[-] (0.3,0.3)  -- (-0.3,-0.3)\botlabel{i};
\draw[-] (0,-0.3) \botlabel{\tau i} [out=90,in=90,looseness=2.5] to  (0.3,-0.3);
\draw[-] (-0.3,0.3) \toplabel{i}[out=-90,in=-90,looseness=2.5] to (0,0.3);
\pinpinpinmid{-.08,.1}{0,0}{.08,-.1}{-1.3,0}{
\frac{Q^\imath_{\tau i, i}(x,y)-Q^\imath_{\tau i, i}(z,y)}{x-z}};
\pinpin{-0.01,.23}{.23,.23}{1,.23}{Q_{i,\tau i}(y,x)};
\end{tikzpicture}.\label{dandy}
\end{align}
This relation is the ibraid relation 
vertically composed on top with 
$\begin{tikzpicture}[iQ,anchorbase]
\draw(-.3,-.2) \botlabel{i} to (-.3,.2);
\draw(0,-.2) \botlabel{\tau i} to (0,.2);
\draw(.3,-.2) \botlabel{i} to (.3,.2);
\pinpin{0,0}{.3,0}{1.3,0}{Q_{i,\tau i}(y,x)};
\end{tikzpicture}\ $.
\end{theo}

\begin{proof}
See \cref{Appendix-A}.
\end{proof}

\subsection{Symmetries}\label{isymmetries}

Suppose we are given parameters $Q_{i,j}(x,y)$ and normalization functions $c_i:X \rightarrow \kk^\times$
satisfying \cref{hypotheses}.
Another admissible choice of parameters
is obtained by
replacing $Q_{i,j}(x,y)$ with ${'}Q_{i,j}(x,y) := r_{i,j} r_{j,i} Q_{i,j}(x,y)$, 
hence $r_{i,j}$ with ${'}r_{i,j} := r_{j,i}^{-1}$,
keeping the normalization functions $c_i:X \rightarrow \kk^\times$, hence $\gamma_i:X^\imath \rightarrow \kk^\times$,
unchanged. 
Let ${'}\UU^\imath(\del,\zeta)$ be the 2-iquantum group defined using these primed parameters.
Also define ${'}\del = ({'}\del_i)_{i \in I}$ and ${'}\zeta = ({'}\zeta_i)_{i \in I}$
so that
\begin{align}\label{sinead}
{'}\del_i &= \del_{\tau i},&
\del_i &= {'}\del_{\tau i},&
{'}\zeta_{i} &= 
(-1)^{\del_{\tau i}+1} \zeta_{\tau i},&
\zeta_i &= (-1)^{{'}\del_{\tau i}+1}\ {'}\zeta_{\tau i}.
\end{align}
These also satisfy \cref{recap,being},
so there are two more 2-iquantum groups
$\UU^\imath({'}\del,{'}\zeta)$ and
${'}\UU^\imath({'}\del,{'}\zeta)$.

There is an isomorphism of graded 
2-categories
\begin{equation}\label{elbowsup}
\Omega^\imath:\UU^\imath(\del,\zeta) \stackrel{\sim}{\rightarrow} \UU^\imath({'}\del,{'}\zeta)
\end{equation}
defined on objects by $\lambda \mapsto -\lambda$, on generating 1-morphisms by
$B_i \one_\lambda \mapsto B_{\tau i}\one_{-\lambda}$, and on generating 2-morphisms
by
\begin{align*}
\begin{tikzpicture}[iQ,centerzero]
\draw[-] (0,-0.3) \botlabel{i} -- (0,0.3);
\closeddot{0,0};
\region{0.25,0}{\lambda};
\end{tikzpicture}&\mapsto
-\ \begin{tikzpicture}[iQ,centerzero]
\draw[-] (0,-0.3) \botlabel{\tau i} -- (0,0.3);
\closeddot{0,0};
\region{0.35,0}{-\lambda};
\end{tikzpicture},&
\begin{tikzpicture}[iQ,centerzero,scale=.9]
\draw[-] (-0.3,-0.3) \botlabel{i} -- (0.3,0.3);
\draw[-] (0.3,-0.3) \botlabel{j} -- (-0.3,0.3);
\region{0.4,0}{\lambda};
\end{tikzpicture}&\mapsto 
- r_{\tau i, j}r_{\tau j, \tau i}\ \begin{tikzpicture}[iQ,centerzero,scale=.9]
\draw[-] (-0.3,-0.3) \botlabel{\tau i} -- (0.3,0.3);
\draw[-] (0.3,-0.3) \botlabel{\tau j} -- (-0.3,0.3);
\region{0.5,0}{-\lambda};
\end{tikzpicture}\\
\begin{tikzpicture}[iQ,centerzero]
\draw[-] (-0.25,-0.15) \botlabel{\tau i} to [out=90,in=90,looseness=3](0.25,-0.15);
\region{0.45,0.1}{\lambda};
\end{tikzpicture}&\mapsto
\begin{tikzpicture}[iQ,centerzero]
\draw[-] (-0.25,-0.15) \botlabel{i} to [out=90,in=90,looseness=3](0.25,-0.15);
\region{0.55,0.1}{-\lambda};
\end{tikzpicture}
,&
\begin{tikzpicture}[iQ,centerzero]
\draw[-] (-0.25,0.15) \toplabel{\tau i} to[out=-90,in=-90,looseness=3] (0.25,0.15);
\region{0.45,-.1}{\lambda};
\end{tikzpicture}&\mapsto
\begin{tikzpicture}[iQ,centerzero]
\draw[-] (-0.25,0.15) \toplabel{i} to[out=-90,in=-90,looseness=3] (0.25,0.15);
\region{0.55,-.1}{-\lambda};
\end{tikzpicture},\\
\begin{tikzpicture}[baseline=-1mm,iQ]
\draw[-] (-0.25,0) arc(180:-180:0.25);
\node at (-0.42,0) {\strandlabel{\tau i}};
\region{0.95,0}{\lambda};
\node at (.55,0) {$(u)$};
\end{tikzpicture}&\mapsto
(-1)^{\del_{i}-\lambda_i}\ \begin{tikzpicture}[baseline=-1mm,iQ]
\draw[-] (-0.25,0) arc(180:-180:0.25);
\node at (-0.37,0) {\strandlabel{i}};.
\region{1.36,0}{-\lambda};
\node at (.7,0) {$(-u)$};
\end{tikzpicture}.
\end{align*}
Also there is an isomorphism of graded 
2-categories
\begin{equation}\label{eastchinasea}
\Psi^\imath:\UU^\imath(\del,\zeta)^\op \stackrel{\sim}{\rightarrow} {'}\UU^\imath(\del,\zeta)
\end{equation}
defined on objects by $\lambda \mapsto \lambda$, on generating 1-morphisms by
$B_i \one_\lambda \mapsto B_i\one_\lambda$, and on generating 2-morphisms
by
\begin{align*}
\begin{tikzpicture}[iQ,centerzero]
\draw[-] (0,-0.3) \botlabel{i} -- (0,0.3);
\closeddot{0,0};
\region{0.25,0}{\lambda};
\end{tikzpicture}&\mapsto
\begin{tikzpicture}[iQ,centerzero]
\draw[-] (0,-0.3) \botlabel{i} -- (0,0.3);
\closeddot{0,0};
\region{0.25,0}{\lambda};
\end{tikzpicture},&
\begin{tikzpicture}[iQ,centerzero,scale=.9]
\draw[-] (-0.3,-0.3) \botlabel{i} -- (0.3,0.3);
\draw[-] (0.3,-0.3) \botlabel{j} -- (-0.3,0.3);
\region{0.4,0}{\lambda};
\end{tikzpicture}&\mapsto 
r_{j,i}^{-1}\ \begin{tikzpicture}[iQ,centerzero,scale=.9]
\draw[-] (-0.3,-0.3) \botlabel{j} -- (0.3,0.3);
\draw[-] (0.3,-0.3) \botlabel{i} -- (-0.3,0.3);
\region{0.4,0}{\lambda};
\end{tikzpicture}\\
\begin{tikzpicture}[iQ,centerzero]
\draw[-] (-0.25,-0.15) \botlabel{\tau i} to [out=90,in=90,looseness=3](0.25,-0.15);
\region{0.45,0.1}{\lambda};
\end{tikzpicture}&\mapsto
\begin{tikzpicture}[iQ,centerzero]
\draw[-] (-0.25,0.15) \toplabel{\tau i} to[out=-90,in=-90,looseness=3] (0.25,0.15);
\region{0.45,-.1}{\lambda};
\end{tikzpicture},&
\begin{tikzpicture}[iQ,centerzero]
\draw[-] (-0.25,0.15) \toplabel{\tau i} to[out=-90,in=-90,looseness=3] (0.25,0.15);
\region{0.45,-.1}{\lambda};
\end{tikzpicture}&\mapsto
\begin{tikzpicture}[iQ,centerzero]
\draw[-] (-0.25,-0.15) \botlabel{\tau i} to [out=90,in=90,looseness=3](0.25,-0.15);
\region{0.55,0.1}{\lambda};
\end{tikzpicture},\\
\begin{tikzpicture}[baseline=-1mm,iQ]
\draw[-] (-0.25,0) arc(180:-180:0.25);
\node at (-0.42,0) {\strandlabel{\tau i}};
\region{0.95,0}{\lambda};
\node at (.55,0) {$(u)$};
\end{tikzpicture}&\mapsto
\begin{tikzpicture}[baseline=-1mm,iQ]
\draw[-] (-0.25,0) arc(180:-180:0.25);
\node at (-0.42,0) {\strandlabel{\tau i}};.
\region{.95,0}{\lambda};
\node at (.55,0) {$(u)$};
\end{tikzpicture}.
\end{align*}
Finally,
there is an isomorphism
\begin{equation}
\Sigma^\imath:\UU^\imath(\del,\zeta)^\rev
\stackrel{\sim}{\rightarrow}
{'}\UU^\imath({'}\del,{'}\zeta)
\end{equation}
taking $\lambda \mapsto -\lambda$,
$B_i \one_\lambda \mapsto \one_{-\lambda} B_i$, and
\begin{align*}
\begin{tikzpicture}[iQ,centerzero]
\draw[-] (0,-0.3) \botlabel{i} -- (0,0.3);
\closeddot{0,0};
\region{0.25,0}{\lambda};
\end{tikzpicture}&\mapsto
\begin{tikzpicture}[iQ,centerzero]
\draw[-] (0,-0.3) \botlabel{i} -- (0,0.3);
\closeddot{0,0};
\region{-0.35,0}{-\lambda};
\end{tikzpicture},&
\begin{tikzpicture}[iQ,centerzero,scale=.9]
\draw[-] (-0.3,-0.3) \botlabel{i} -- (0.3,0.3);
\draw[-] (0.3,-0.3) \botlabel{j} -- (-0.3,0.3);
\region{0.4,0}{\lambda};
\end{tikzpicture}&\mapsto -r_{j,i}^{-1} \ 
\begin{tikzpicture}[iQ,centerzero,scale=.9]
\draw[-] (-0.3,-0.3) \botlabel{j} -- (0.3,0.3);
\draw[-] (0.3,-0.3) \botlabel{i} -- (-0.3,0.3);
\region{-0.5,0}{-\lambda};
\end{tikzpicture},\\
\begin{tikzpicture}[iQ,centerzero]
\draw[-] (-0.25,-0.15) \botlabel{\tau i} to [out=90,in=90,looseness=3](0.25,-0.15);
\region{0.45,0.1}{\lambda};
\end{tikzpicture}&\mapsto
\begin{tikzpicture}[iQ,centerzero]
\draw[-] (-0.25,-0.15) \botlabel{i} to [out=90,in=90,looseness=3](0.25,-0.15);
\region{-0.55,0.1}{-\lambda};
\end{tikzpicture}\ ,&
\begin{tikzpicture}[iQ,centerzero]
\draw[-] (-0.25,0.15) \toplabel{\tau i} to[out=-90,in=-90,looseness=3] (0.25,0.15);
\region{0.45,-.1}{\lambda};
\end{tikzpicture}&\mapsto
\begin{tikzpicture}[iQ,centerzero]
\draw[-] (-0.25,0.15) \toplabel{ i}to[out=-90,in=-90,looseness=3] (0.25,0.15);
\region{-0.55,-.1}{-\lambda};
\end{tikzpicture}\ ,\\
\begin{tikzpicture}[baseline=-1mm,iQ]
\draw[-] (-0.25,0) arc(180:-180:0.25);
\node at (-0.42,0) {\strandlabel{\tau i}};
\region{0.95,0}{\lambda};
\node at (.55,0) {$(u)$};
\end{tikzpicture}&\mapsto
(-1)^{\del_i-\lambda_i}\ \begin{tikzpicture}[baseline=-1mm,iQ]
\draw[-] (-0.25,0) arc(180:-180:0.25);
\node at (-0.37,0) {\strandlabel{i}};
\region{-.75,0}{-\lambda};
\node at (.68,0) {$(-u)$};
\end{tikzpicture}.
\end{align*}
All of these things follow by elementary calculations with the defining relations.
We regard $\Omega^\imath$, $\Psi^\imath$ and $\Sigma^\imath$
as categorical analogs of the maps $\omega^\imath$, $\psi^\imath$ and $\sigma^\imath$ from \cref{newyearsday2,snoring,newyearsday1}.

\subsection{Example: split rank one}\label{rankone}

Suppose that $i \in I$ is fixed by $\tau$.
Fix also $\lambda \in X^\imath$
and consider the strict graded 
monoidal category 
$\END_{\UU^\imath}(\lambda)$.
In view of \cref{morethingsworking,ibubblegeneratingfunction}, 
for $i$ with $\tau i = i$, 
the relation \cref{ibubblerel}
is equivalent to
\begin{equation}\label{equivform}
\begin{tikzpicture}[baseline=-1mm,iQ]
\draw[-] (-0.25,0) arc(180:-180:0.25);
\node at (-0.34,0) {\strandlabel{i}};
\region{0.5,0}{\lambda};
\end{tikzpicture}=t_i(\lambda)\id_{\one_\lambda}.
\end{equation}
For $t \in \{0,1\}$, let $\NB_{t}$ be
the {\em nil-Brauer category} from
\cite[Def.~2.1]{BWWbasis}.

\begin{theo}\label{aging}
For any element $i \in I$ that is fixed by $\tau$,
there is a strict 
graded monoidal functor
$\Upsilon:\NB_{t_i(\lambda)} \rightarrow \END_{\UU^\imath}(\lambda)$
such that
\begin{align*}
B &\mapsto B_i \one_\lambda,\\\notag
\begin{tikzpicture}[iQ,centerzero]
\draw[-] (0,-0.3) -- (0,0.3);
\closeddot{0,0};
\end{tikzpicture}&\mapsto \begin{tikzpicture}[iQ,centerzero]
\draw[-] (0,-0.3) \botlabel{i} -- (0,0.3);
\closeddot{0,0};
\region{0.2,0}{\lambda};
\end{tikzpicture},\quad
\begin{tikzpicture}[iQ,centerzero,scale=.9]
\draw[-] (-0.3,-0.3) -- (0.3,0.3);
\draw[-] (0.3,-0.3)  -- (-0.3,0.3);
\end{tikzpicture}\mapsto
\begin{tikzpicture}[iQ,centerzero,scale=.9]
\draw[-] (-0.3,-0.3) \botlabel{i} -- (0.3,0.3);
\draw[-] (0.3,-0.3) \botlabel{i} -- (-0.3,0.3);
\region{0.35,0}{\lambda};
\end{tikzpicture},\quad
\begin{tikzpicture}[iQ,centerzero]
\draw[-] (-0.25,-0.15) to [out=90,in=90,looseness=3](0.25,-0.15);
\end{tikzpicture}
\mapsto \begin{tikzpicture}[iQ,centerzero]
\draw[-] (-0.25,-0.15) \botlabel{i} to [out=90,in=90,looseness=3](0.25,-0.15);
\region{0.45,0.1}{\lambda};
\end{tikzpicture},\quad
\begin{tikzpicture}[iQ,centerzero]
\draw[-] (-0.25,0.15) to[out=-90,in=-90,looseness=3] (0.25,0.15);
\end{tikzpicture}\mapsto \begin{tikzpicture}[iQ,centerzero]
\draw[-] (-0.25,0.15) \toplabel{i} to[out=-90,in=-90,looseness=3] (0.25,0.15);
\region{0.45,-.1}{\lambda};
\end{tikzpicture},
\quad u^{-1} \mathbb{O}(u)\mapsto \begin{tikzpicture}[baseline=-1mm,iQ]
\draw[-] (-0.25,0) arc(180:-180:0.25);
\node at (-0.39,0) {\strandlabel{i}};
\node at (.55,0) {$(u)$};
\region{1,0}{\lambda};
\end{tikzpicture}.
\end{align*}
\end{theo}

\begin{proof}
The relations 
\cref{equivform,izigzag,ipivotal,ipitchfork,idotslide,icurl,iquadratic,ibraid} with all strings labelled by $i$
are equivalent to the defining relations of $\NB_t$.
\end{proof}

\subsection{Example: diagonal embedding}\label{iexisting}

In this subsection, we go back to the setup of \cref{s2-quantum},
dropping the additional hypotheses on parameters
from \cref{hypotheses}.
Let $\U$ be the quantum group from \cref{data}.
Let 
${\pmb I} := I^+\sqcup I^-$
be the disjoint union of two copies
$I^+ := \{i^+\:|\:i \in I\}$ and $I^- := \{i^-\:|\:i \in I\}$
of $I$.
This is the index set for the following ``doubled''
Cartan datum:
\begin{itemize}
\item
The Cartan matrix is
defined by $a_{i^+,j^+} = a_{i^-,j^-} := a_{i,j}$
and $a_{i^+,j^-} = a_{i^-,j^+} := 0$, and it is symmetrized by
$d_{i^+} = d_{i^-} := d_i$ for each $i,j \in I$.
\item
The weight lattice ${\pmb X}$ is $X \oplus X$,
and the coweight lattice ${\pmb Y}$
is $Y \oplus Y$ paired with ${\pmb X}$ in the obvious way, i.e., $(h^+,h^-) \in {\pmb Y}$ maps $(\lambda^+,\lambda^-) \in {\pmb X}$ to $h^+(\lambda^+)+h^-(\lambda^-)$.
\item
The simple roots are defined by $\alpha_{i^+} := (\alpha_i,0)$
and $\alpha_{i^-} := (0,\alpha_i)$
for $i \in I$.
\item
The simple coroots
are defined by $h_{i^+} := (h_i,0)$
and $h_{i^-} := (0,h_i)$ for $i \in I$.
\end{itemize}
We denote the resulting quantum group by ${\pmb\U}$.
It is identified with
the tensor product $\U \otimes_{\Q(q)} \U$
so that $q^{(h^+,h^-)}$, 
$e_{i^+}, e_{i^-}, f_{i^+}$ and $f_{i^-}$
correspond to $q^{h^+} \otimes q^{h^-}$, $e_i \otimes 1, 1 \otimes e_i, f_i \otimes 1$ and $1 \otimes f_i$, respectively.
The {\em diagonal embedding} of $\U$ into ${\pmb\U}$ is the $\Q(q)$-algebra homomorphism
\begin{align}\label{diagonalembedding}
\begin{split}
& (\omega \otimes \id) \circ \Delta:
\U \hookrightarrow {\pmb\U},
\\
q^h \mapsto q^{-h} \otimes q^h, \qquad
e_i &\mapsto b_{i^+} :=f_i \otimes 1 +  k_i^{-1} \otimes e_i, \qquad
f_i \mapsto b_{i^-} :=1 \otimes f_i + e_i \otimes k_{i}^{-1}
\end{split}
\end{align}
for $h \in Y, i \in I$.

Let ${\pmb\tau} :{\pmb I}\rightarrow{\pmb I}$ be the involution which switches $i^+$ and $i^-$ for each $i \in I$.
There are also involutions
${\pmb\tau}:{\pmb X}\rightarrow{\pmb X}$
and ${\pmb\tau}^*:{\pmb Y}\rightarrow {\pmb Y}$
switching the summands in the obvious way.
Let ${\pmb\U}^\imath$ be 
the corresponding iquantum group with iweight lattice
${\pmb X}^\imath$ and icoweight lattice ${\pmb Y}^\imath$. We omit $\del$ from the notation here because, in the present situation there is only one choice: for all $i \in I$ we must have that
\begin{align}
\del_{i^+}&:=0,&\del_{i^-}&:=0.
\end{align}
Thus, ${\pmb\U}^\imath$ is the subalgebra of ${\pmb\U}$
generated by $q^{-h}\otimes q^h\:(h \in Y)$ and $b_{i^+}$, $b_{i^-}$ for $i \in I$.
Comparing these formulae 
with \cref{diagonalembedding}, we see that
the diagonal embedding restricts to an algebra isomorphism
\begin{align}\label{nico}
d:\U &\stackrel{\sim}{\rightarrow} {\pmb\U}^\imath,&
q^h&\mapsto q^{-h} \otimes q^h,&
e_i &\mapsto b_{i^+},&
f_i &\mapsto b_{i^-}
\end{align}
for $h \in Y, i \in I$.
Thus, $d$ identifies $\U$ with ${\pmb\U}^\imath$, which is a
quasi-split iquantum group 
of diagonal type.
Note also that there is an isomorphism
\begin{align}
X &\stackrel{\sim}{\rightarrow} {\pmb X}^\imath,&
\lambda &\mapsto \bar\lambda
\end{align}
with $\bar\lambda$ being 
the image of $(0,\lambda) \in {\pmb X}$
in ${\pmb X}^\imath$.
This satisfies $\overline{\lambda+\alpha_i} = \overline{\lambda}-\alpha_{i^+} = \overline{\lambda}+\alpha_{i^-}$ for $i \in I$.

Our next theorem is a categorification of the isomorphism \cref{nico}.
Let $\UU$ be the 2-quantum group from \cref{existing} corresponding to $\U$.
So $\UU$
has object set $X$ and generating 1-morphisms $E_i \one_\lambda, F_i \one_\lambda\:(\lambda \in X, i \in I)$.
Its definition involves parameters
$Q_{i,j}(x,y)$. We also assume $c_i$ has been chosen as in \cref{great}, and let $r_{i,j}$ be as in \cref{gollygosh}. 
Let ${\pmb\UU}^\imath(\zeta)$ be the 2-iquantum group 
corresponding to ${\pmb\U}^\imath$
from \cref{new2cat} defined from the parameters
\begin{align}\label{this}
Q_{i^+,j^+}(x,y) &:= Q_{i,j}(x,y),&
Q_{i^-,j^-}(x,y)&:= r_{i,j} r_{j,i} Q_{i,j}(-x,-y),\\\label{that}
Q_{i^+,j^-}(x,y)&:= 1,&
Q_{i^-,j^+}(x,y)&:= 1,\\
\label{theother}
c_{i^+}\big((\lambda^+,\lambda^-)\big) &:= c_i(\lambda^+),&
c_{i^-}\big((\lambda^+,\lambda^-)\big)
&:= (-1)^{h_i(\lambda^-)} c_i(\lambda^-),\\
\zeta_{i^+}&:=1,& \zeta_{i^-} &:= -1.\label{bloppy}\\\intertext{for $i,j \in I$ and $\lambda^+,\lambda^- \in X$. According to the definitions \cref{gollygosh,purer}, we have that}
r_{i^+,j^+} &:= r_{i,j},&\label{aargh}
r_{i^-,j^-} &:= r_{j,i}^{-1},\\\label{aaargh}
r_{i^+,j^-}&:=1,
&r_{i^-,j^+}&:= 1,\\\label{blippy}
\gamma_{i^+}\big(\bar\lambda\big)
&=
c_i(\lambda)^{-1},&
\gamma_{i^-}\big(\bar\lambda\big)&=
(-1)^{h_i(\lambda)} c_i(\lambda)
\end{align}
for $i \in I$ and $\lambda \in X$.
Now it is straightforward to check that \cref{hypotheses}
is satisfied.

\begin{theo}\label{earphones}
\begin{samepage}
There is an
isomorphism of graded 2-categories $D:\UU \stackrel{\sim}{\rightarrow} {\pmb\UU}^\imath(\zeta)$ such that
\begin{align*}
\lambda &\mapsto \bar \lambda,&&&
E_i \one_\lambda &\mapsto B_{i^+} \one_{\bar \lambda},&
F_i \one_\lambda &\mapsto B_{i^-} \one_{\bar \lambda},\\
\begin{tikzpicture}[Q,centerzero]
\draw[-to] (0,-0.3) \botlabel{i} -- (0,0.3);
\opendot{0,0};
\region{0.25,0}{\lambda};
\end{tikzpicture}&\mapsto \begin{tikzpicture}[iQ,centerzero]
\draw[-] (0,-0.3) \botlabel{i^+} -- (0,0.3);
\closeddot{0,0};
\region{0.25,0}{\bar\lambda};
\end{tikzpicture}
\!,\!&
\begin{tikzpicture}[Q,centerzero,scale=.9]
\draw[-to] (-0.3,-0.3) \botlabel{i} -- (0.3,0.3);
\draw[-to] (0.3,-0.3) \botlabel{j} -- (-0.3,0.3);
\region{0.4,0}{\lambda};
\end{tikzpicture}&\mapsto 
\begin{tikzpicture}[iQ,centerzero,scale=.9]
\draw[-] (-0.3,-0.3) \botlabel{i^+} -- (0.3,0.3);
\draw[-] (0.3,-0.3) \botlabel{j^+} -- (-0.3,0.3);
\region{0.4,0}{\bar\lambda};
\end{tikzpicture},\!&
\begin{tikzpicture}[Q,centerzero]
\draw[-to] (-0.25,-0.15) \botlabel{i} to [out=90,in=90,looseness=3](0.25,-0.15);
\region{0.45,0.1}{\lambda};
\end{tikzpicture}&\mapsto \begin{tikzpicture}[iQ,centerzero]
\draw[-] (-0.25,-0.15) \botlabel{i^+} to [out=90,in=90,looseness=3](0.25,-0.15);
\region{0.45,0.1}{\bar\lambda};
\end{tikzpicture}\!,\!\!&
\begin{tikzpicture}[Q,centerzero]
\draw[-to] (-0.25,0.15) \toplabel{i} to[out=-90,in=-90,looseness=3] (0.25,0.15);
\region{0.45,-.1}{\lambda};
\end{tikzpicture}&\mapsto 
\begin{tikzpicture}[iQ,centerzero]
\draw[-] (-0.25,0.15) \toplabel{i^-} to[out=-90,in=-90,looseness=3] (0.25,0.15);
\region{0.45,-.1}{\bar\lambda};
\end{tikzpicture}\\\intertext{for $\lambda \in X, i,j \in I$.
It maps}
\begin{tikzpicture}[Q,centerzero]
\draw[to-] (0,-0.3) \botlabel{i} -- (0,0.3);
\opendot{0,0};
\region{0.25,0}{\lambda};
\end{tikzpicture}
&\mapsto
-\begin{tikzpicture}[iQ,centerzero]
\draw[-] (0,-0.3) \botlabel{i^-} -- (0,0.3);
\closeddot{0,0};
\region{0.25,0}{\bar \lambda};
\end{tikzpicture},&
\begin{tikzpicture}[Q,centerzero,scale=.9]
\draw[to-] (-0.3,-0.3) \botlabel{i} -- (0.3,0.3);
\draw[to-] (0.3,-0.3) \botlabel{j} -- (-0.3,0.3);
\region{0.4,0}{\lambda};
\end{tikzpicture}&\mapsto\!\!
\ \begin{tikzpicture}[iQ,centerzero,scale=.9]
\draw[-] (-0.3,-0.3) \botlabel{i^-} -- (0.3,0.3);
\draw[-] (0.3,-0.3) \botlabel{j^-} -- (-0.3,0.3);
\region{0.4,0}{\bar\lambda};
\end{tikzpicture},&
\begin{tikzpicture}[Q,centerzero,scale=.9]
\draw[-to] (-0.3,-0.3) \botlabel{i} -- (0.3,0.3);
\draw[to-] (0.3,-0.3) \botlabel{j} -- (-0.3,0.3);
\region{0.4,0}{\lambda};
\end{tikzpicture}
&\mapsto \begin{tikzpicture}[iQ,centerzero,scale=.9]
\draw[-] (-0.3,-0.3) \botlabel{i^+} -- (0.3,0.3);
\draw[-] (0.3,-0.3) \botlabel{j^-} -- (-0.3,0.3);
\region{0.4,0}{\bar\lambda};
\end{tikzpicture}\ ,&
\begin{tikzpicture}[Q,centerzero,scale=.9]
\draw[to-] (-0.3,-0.3) \botlabel{i} -- (0.3,0.3);
\draw[-to] (0.3,-0.3) \botlabel{j} -- (-0.3,0.3);
\region{0.4,0}{\lambda};
\end{tikzpicture}
&\mapsto 
\begin{tikzpicture}[iQ,centerzero,scale=.9]
\draw[-] (-0.3,-0.3) \botlabel{i^-} -- (0.3,0.3);
\draw[-] (0.3,-0.3) \botlabel{j^+} -- (-0.3,0.3);
\region{0.4,0}{\bar\lambda};
\end{tikzpicture}\ ,\\
\begin{tikzpicture}[Q,centerzero]
\draw[to-] (-0.25,0.15) \toplabel{i} to[out=-90,in=-90,looseness=3] (0.25,0.15);
\region{0.45,-.1}{\lambda};
\end{tikzpicture}\!\!
&\mapsto \begin{tikzpicture}[iQ,centerzero]
\draw[-] (-0.25,0.15) \toplabel{i^+} to[out=-90,in=-90,looseness=3] (0.25,0.15);
\region{0.45,-.1}{\bar\lambda};
\end{tikzpicture},\!\!
&
\begin{tikzpicture}[Q,centerzero]
\draw[to-] (-0.25,-0.15) \botlabel{i} to [out=90,in=90,looseness=3](0.25,-0.15);
\region{0.45,0.1}{\lambda};
\end{tikzpicture}\!\!&
\mapsto \begin{tikzpicture}[iQ,centerzero]
\draw[-] (-0.25,-0.15) \botlabel{i^-} to [out=90,in=90,looseness=3](0.25,-0.15);
\region{0.45,0.1}{\bar\lambda};
\end{tikzpicture},\!\!&
\begin{tikzpicture}[baseline=-1mm,Q]
\draw[to-] (-0.25,0) arc(180:-180:0.25);
\node at (0,-.4) {\strandlabel{i}};
\region{0.95,0}{\lambda};
\node at (.55,0) {$(u)$};
\end{tikzpicture}\!&\mapsto\!
\begin{tikzpicture}[baseline=-1mm,iQ]
\draw[-] (-0.25,0) arc(180:-180:0.25);
\node at (-.48,0) {\strandlabel{i^-}};
\region{0.95,0}{\bar\lambda};
\node at (.53,0) {$(u)$};
\end{tikzpicture}\!,&
\begin{tikzpicture}[baseline=-1mm,Q]
\draw[-to] (-0.25,0) arc(180:-180:0.25);
\node at (0,-.4) {\strandlabel{i}};
\region{0.95,0}{\lambda};
\node at (.55,0) {$(u)$};
\end{tikzpicture}\!\!&\mapsto
(-1)^{h_i(\lambda)}\!
\begin{tikzpicture}[baseline=-1mm,iQ]
\draw[-] (-0.25,0) arc(180:-180:0.25);
\node at (-.44,0) {\strandlabel{i^+}};
\region{1.32,0}{\bar\lambda};
\node at (.69,0) {$(-u)$};
\end{tikzpicture}\!.
\end{align*}
\end{samepage}
\end{theo}

\begin{proof}
One first checks from \cref{table1,table2} 
that the degrees of 2-morphisms
on either side are the same.
Then the theorem follows by comparing the presentation
of ${\pmb\UU}^\imath(\zeta)$ 
with the one for $\UU$ from \cref{klpres}, using also the additional relations mentioned in \cref{notevendone}.
\end{proof}

\subsection{Example: quasi-split type AIII with an even number of nodes}\label{bsww}

There is one more special case which already appears in the existing literature: the categorification of
the quasi-split iquantum group of type
AIII that was introduced in \cite{BSWW}.
Fix a choice of $r \geq 1$.
The Satake diagram is
\begin{equation}\label{AIIIeven}
\begin{tikzpicture}[anchorbase,scale=3]
\draw[to-to,red,dotted] (.7,.15) to (.7,-.15);
\draw[to-to,red,dotted] (.3,.15) to (.3,-.15);
\draw[to-to,red,dotted] (-.4,.15) to (-.4,-.15);
\draw[to-to,red,dotted] (-.8,.15) to (-.8,-.15);
\node at (-.8,.32) {$\frac{1}{2}$};
\node at (-.4,.32) {$\frac{3}{2}$};
\node at (.3,.32) {$r\!-\!\frac{3}{2}$};
\node at (.7,.32) {$r\!-\!\frac{1}{2}$};
\node at (-.8,-.32) {$-\frac{1}{2}$};
\node at (-.4,-.32) {$-\frac{3}{2}$};
\node at (.3,-.32) {$\frac{3}{2}\!-\!r$};
\node at (.7,-.32) {$\frac{1}{2}\!-\!r$};
\draw[thick] (-.8,.2) to (-0.4,.2);
\draw[thick,-] (-.8,.2) to[out=180,in=180,looseness=2](-.8,-.2);
\draw[thick] (.7,.2) to (0.3,.2);
\draw[thick,densely dashed] (-.4,.2) to (0.3,.2);
\draw[thick] (-.8,-.2) to (-0.4,-.2);
\draw[thick] (.7,-.2) to (0.3,-.2);
\draw[thick,densely dashed] (-.4,-.2) to (0.3,-.2);
\Opendot{.7,.2};
\Opendot{.3,.2};
\Opendot{-.4,.2};
\Opendot{-.8,.2};
\Opendot{.7,-.2};
\Opendot{.3,-.2};
\Opendot{-.4,-.2};
\Opendot{-.8,-.2};
\end{tikzpicture}
\end{equation}
So
$I = \textstyle\left\{\frac{1}{2}-r,\dots,-\frac{1}{2},\frac{1}{2},\dots, r-\frac{1}{2}\right\}$
and $\tau i = -i$.
As in \cite{BSWW}, we realize the weight lattice $X$ of the 
root system of type $A_{2r}$
as the Abelian group 
with generators $\eps_{-r},\dots,\eps_{-1},\eps_0,
\eps_1,\dots,\eps_r$
subject to the relation asserting that their sum is 0.
The simple roots are
$\alpha_i:=\eps_{i-\frac{1}{2}} - \eps_{i+\frac{1}{2}}$ for each $i \in I$. Let $Y$ be the dual lattice spanned by the simple coroots
$h_i\:(i \in I)$ such that $h_i(\alpha_j) = 2\delta_{i,j}-\delta_{i+1,j}-\delta_{i-1,j}$.
The involution $\tau:X\rightarrow X$ maps
$\eps_i \mapsto -\eps_{-i}$ and 
$\alpha_i \mapsto \alpha_{-i}$. Similarly, $\tau^*:Y \rightarrow Y$ takes $h_i \mapsto h_{-i}$.
Define $\del=(\del_i)_{i \in I}$ by
\begin{equation}
\del_i:=\begin{cases}
1&\text{if $i=\frac{1}{2}$}\\
0&\text{otherwise,}
\end{cases}
\end{equation}
and let $\U^\imath(\del)$ be the corresponding quasi-split iquantum group as in \cref{iqg}.

Assuming that $2$ is invertible in $\kk$, the categorification
of $\U^\imath(\del)$ discovered in \cite{BSWW} is a certain graded 2-category 
$\UU^{\BSWW}$ with object set $X^\imath$, and generating 1-morphisms $E_i \one_\lambda = \one_{\lambda+\alpha_i} E_i$ and $F_i \one_\lambda = \one_{\lambda-\alpha_i} F_i$ for $i \in I_1:= \{i \in I\:|\:i > 0\}$
and $\lambda \in X^\imath$. The identity 2-endomorphisms of $E_i \one_\lambda$ and $F_i \one_\lambda$ are denoted by an upward and a downward string labelled by $i$, respectively, just like we did for the ordinary 2-quantum group $\UU$. Then there are generating 2-morphisms denoted by dots, crossings, cups and caps subject to the relations in \cite[Def.~3.1, Def.~3.3]{BSWW}.

Also let $\UU^\imath(\del,\zeta)$ be the 2-iquantum group from \cref{new2cat} associated to the Satake diagram \cref{AIIIeven} and root datum just described, with parameters chosen as follows. First,
\begin{equation}\label{badchoice}
Q_{i,j}(x,y) := \begin{cases}
0&\text{if $i=j$}\\
x-y&\text{if $i > 0$ and $j=i+1$, or $i < 1$ and $j = i-1$}\\
y-x&\text{if $i < 0$ and $j =i+1$, or $i > 1$ and $j =i-1$}\\
1&\text{otherwise}
\end{cases}
\end{equation}
for $i,j \in I$, 
which are the geometric parameters
arising from the quiver
$\ \begin{tikzpicture}[anchorbase,scale=1.5]
\draw[to-to,red,densely dotted] (.7,.15) to (.7,-.15);
\draw[to-to,red,densely dotted] (.3,.15) to (.3,-.15);
\draw[to-to,red,densely dotted] (-.4,.15) to (-.4,-.15);
\draw[to-to,red,densely dotted] (-.8,.15) to (-.8,-.15);
\draw[thick,-to] (-.8,.2) to[out=200,in=90,looseness=1] (-1,-.03);\draw[thick,-] (-1,-.02) to [out=-90,in=160,looseness=1](-0.8,-.2);
\draw[thick,-] (.7,-.2) to (.52,-.2);  \draw[thick,to-] (.53,-.2) to (0.3,-.2);
\draw[thick,densely dashed] (-.4,.2) to (0.3,.2);
\draw[thick,-] (-0.58,-.2) to (-.3,-.2);
\draw[thick,-to] (-.8,-.2) to (-.57,-.2); 
\draw[thick,densely dashed] (-.4,-.2) to (0.3,-.2);
\draw[thick,-to] (-.8,.2) to (-.57,.2);\draw[thick,-](-.58,.2) to (-0.4,.2);
\draw[thick,to-] (.53,.2) to (0.3,.2);\draw[thick,-] (.7,.2) to (.52,.2);
\opendot{.7,.2};
\opendot{.3,.2};
\opendot{-.4,.2};
\opendot{-.8,.2};
\opendot{.7,-.2};
\opendot{.3,-.2};
\opendot{-.4,-.2};
\opendot{-.8,-.2};
\end{tikzpicture}
$.
For $i \in I$ and $\lambda \in X$, we let
\begin{align}
c_i(\lambda) &:= \begin{cases}
1&\text{if $i > 0$}\\
(-1)^{h_i(\lambda)}&\text{if $i < 0$,}
\end{cases}\\
\intertext{so that}
\gamma_i(\lambda)&=
\begin{cases}
1&\text{if $i > 0$}\\
(-1)^{\lambda_i}&\text{if $i < 0$}
\end{cases}\\\intertext{for $\lambda \in X^\imath$.
According to \cref{gollygosh}, we have that}
r_{i,j} &= \begin{cases}
-1&\text{if $j=i-1$ and $i \neq \frac{1}{2}$}\\
1&\text{otherwise}.
\end{cases}\\\intertext{Finally, let}
\zeta_i &:= \begin{cases}
1&\text{if $i > 0$}\\ 
-2&\text{if $i=-\frac{1}{2}$}\\
-1&\text{if $i < -1$.}
\end{cases}
\end{align}
\cref{hypotheses} holds.

\begin{theo}\label{pacific}
There is an
isomorphism of graded 2-categories
$\Phi:\UU^\BSWW \stackrel{\sim}{\rightarrow} \UU^\imath(\del,\zeta)$ such that
\begin{align*}
\lambda &\mapsto -\lambda,&&&
E_i \one_\lambda &\mapsto B_i \one_{-\lambda},&
F_i \one_\lambda &\mapsto B_{-i} \one_{-\lambda},\\
\begin{tikzpicture}[Q,centerzero]
\draw[-to] (0,-0.3) \botlabel{i} -- (0,0.3);
\opendot{0,0};
\region{0.2,0}{\lambda};
\end{tikzpicture}&\mapsto \begin{tikzpicture}[iQ,centerzero]
\draw[-] (0,-0.3) \botlabel{i} -- (0,0.3);
\closeddot{0,0};
\region{0.3,0}{-\lambda};
\end{tikzpicture}
,&
\begin{tikzpicture}[Q,centerzero]
\draw[to-] (0,-0.3) \botlabel{i} -- (0,0.3);
\opendot{0,0};
\region{0.2,0}{\lambda};
\end{tikzpicture}&\mapsto -\ \begin{tikzpicture}[iQ,centerzero]
\draw[-] (0,-0.3) \botlabel{-i} -- (0,0.3);
\closeddot{0,0};
\region{0.3,0}{-\lambda};
\end{tikzpicture},&
\begin{tikzpicture}[baseline=-1mm,Q]
\draw[to-] (-0.25,0) arc(180:-180:0.25);
\node at (0,-.4) {\strandlabel{i}};
\region{0.95,0}{\lambda};
\node at (.55,0) {$(u)$};
\end{tikzpicture}\!&\mapsto\!\!\!
\begin{tikzpicture}[baseline=-1mm,iQ]
\draw[-] (-0.25,0) arc(180:-180:0.25);
\node at (-.45,0) {\strandlabel{-i}};
\region{1,0}{-\lambda};
\node at (.53,0) {$(u)$};
\end{tikzpicture},&
\!\begin{tikzpicture}[baseline=-1mm,Q]
\draw[-to] (-0.25,0) arc(180:-180:0.25);
\node at (0,-.4) {\strandlabel{i}};
\region{0.95,0}{\lambda};
\node at (.55,0) {$(u)$};
\end{tikzpicture}\!&\mapsto
(-1)^{\del_{-i}-\lambda_i}\!\!
\begin{tikzpicture}[baseline=-1mm,iQ]
\draw[-] (-0.25,0) arc(180:-180:0.25);
\node at (-.36,0) {\strandlabel{i}};
\region{1.37,0}{-\lambda};
\node at (.69,0) {$(-u)$};
\end{tikzpicture},\\
\begin{tikzpicture}[Q,centerzero,scale=.9]
\draw[-to] (-0.3,-0.3) \botlabel{i} -- (0.3,0.3);
\draw[-to] (0.3,-0.3) \botlabel{j} -- (-0.3,0.3);
\region{0.35,0}{\lambda};
\end{tikzpicture}&\mapsto 
\begin{tikzpicture}[iQ,centerzero,scale=.9]
\draw[-] (-0.3,-0.3) \botlabel{i} -- (0.3,0.3);
\draw[-] (0.3,-0.3) \botlabel{j} -- (-0.3,0.3);
\region{0.45,0}{-\lambda};
\end{tikzpicture},\!\!\!&
\begin{tikzpicture}[Q,centerzero,scale=.9]
\draw[to-] (-0.3,-0.3) \botlabel{i} -- (0.3,0.3);
\draw[to-] (0.3,-0.3) \botlabel{j} -- (-0.3,0.3);
\region{0.35,0}{\lambda};
\end{tikzpicture}&\mapsto
\begin{tikzpicture}[iQ,centerzero,scale=.9]
\draw[-] (-0.3,-0.3) \botlabel{-i} -- (0.3,0.3);
\draw[-] (0.3,-0.3) \botlabel{-j} -- (-0.3,0.3);
\region{0.45,0}{-\lambda};
\end{tikzpicture},&
\begin{tikzpicture}[Q,centerzero,scale=.9]
\draw[-to] (-0.3,-0.3) \botlabel{i} -- (0.3,0.3);
\draw[to-] (0.3,-0.3) \botlabel{j} -- (-0.3,0.3);
\region{0.35,0}{\lambda};
\end{tikzpicture}
&\mapsto 
\begin{tikzpicture}[iQ,centerzero,scale=.9]
\draw[-] (-0.3,-0.3) \botlabel{i} -- (0.3,0.3);
\draw[-] (0.3,-0.3) \botlabel{-j} -- (-0.3,0.3);
\region{0.45,0}{-\lambda};
\end{tikzpicture},&
\begin{tikzpicture}[Q,centerzero,scale=.9]
\draw[to-] (-0.3,-0.3) \botlabel{i} -- (0.3,0.3);
\draw[-to] (0.3,-0.3) \botlabel{j} -- (-0.3,0.3);
\region{0.35,0}{\lambda};
\end{tikzpicture}
&\mapsto 
\begin{tikzpicture}[iQ,centerzero,scale=.9]
\draw[-] (-0.3,-0.3) \botlabel{-i} -- (0.3,0.3);
\draw[-] (0.3,-0.3) \botlabel{j} -- (-0.3,0.3);
\region{0.45,0}{-\lambda};
\end{tikzpicture},\\
\begin{tikzpicture}[Q,centerzero]
\draw[-to] (-0.25,-0.15) \botlabel{i} to [out=90,in=90,looseness=3](0.25,-0.15);
\region{0.45,0.1}{\lambda};
\end{tikzpicture}&\mapsto \begin{tikzpicture}[iQ,centerzero]
\draw[-] (-0.25,-0.15) \botlabel{i} to [out=90,in=90,looseness=3](0.25,-0.15);
\region{0.55,0.1}{-\lambda};
\end{tikzpicture},\!\!\!\!&
\begin{tikzpicture}[Q,centerzero]
\draw[-to] (-0.25,0.15) \toplabel{i} to[out=-90,in=-90,looseness=3] (0.25,0.15);
\region{0.45,-.1}{\lambda};
\end{tikzpicture}&\mapsto 
\begin{tikzpicture}[iQ,centerzero]
\draw[-] (-0.25,0.15) \toplabel{-i} to[out=-90,in=-90,looseness=3] (0.25,0.15);
\region{0.55,-.1}{-\lambda};
\end{tikzpicture},
&
\begin{tikzpicture}[Q,centerzero]
\draw[to-] (-0.25,0.15) \toplabel{i} to[out=-90,in=-90,looseness=3] (0.25,0.15);
\region{0.45,-.1}{\lambda};
\end{tikzpicture}
&\mapsto \begin{tikzpicture}[iQ,centerzero]
\draw[-] (-0.25,0.15) \toplabel{i} to[out=-90,in=-90,looseness=3] (0.25,0.15);
\region{0.55,-.1}{-\lambda};
\end{tikzpicture},&
\begin{tikzpicture}[Q,centerzero]
\draw[to-] (-0.25,-0.15) \botlabel{i} to [out=90,in=90,looseness=3](0.25,-0.15);
\region{0.45,0.1}{\lambda};
\end{tikzpicture}
&\mapsto \begin{tikzpicture}[iQ,centerzero]
\draw[-] (-0.25,-0.15) \botlabel{-i} to [out=90,in=90,looseness=3](0.25,-0.15);
\region{0.55,0.1}{-\lambda};
\end{tikzpicture}
\end{align*}
for $\lambda \in X^\imath, i,j \in I_1$.
\end{theo}

\begin{proof}
This is another elementary comparison of relations---all of the defining relations for $\UU^{\BSWW}$ hold in $\UU^\imath$ so that the 2-functor makes sense, and there is also a
2-functor $\UU^\imath\rightarrow \UU^{\BSWW}$
that is the 2-sided inverse of $\Phi$ defined by
checking that all of the defining relations for $\UU^\imath$ hold in $\UU^{\BSWW}$.
\end{proof}

\begin{rem}
There is another choice of parameters 
for this 2-iquantum group
which has the advantage that
$r_{i,j}=1$ for all $i, j \in I$, although the normalization functions $c_i$ are more complicated. It is a special case of the situation described in the introduction: we take the geometric parameters defined by
the quiver
$\ \begin{tikzpicture}[anchorbase,scale=1.5]
\draw[to-to,red,densely dotted] (.7,.15) to (.7,-.15);
\draw[to-to,red,densely dotted] (.3,.15) to (.3,-.15);
\draw[to-to,red,densely dotted] (-.4,.15) to (-.4,-.15);
\draw[to-to,red,densely dotted] (-.8,.15) to (-.8,-.15);
\draw[thick,] (-.8,-.2) to (-.62,-.2);\draw[thick,to-](-.63,-.2) to (-0.4,-.2);
\draw[thick,-] (-.8,.2) to[out=200,in=90,looseness=1] (-1,.03); \draw[thick,to-] (-1,.02) to [out=-90,in=160,looseness=1](-0.8,-.2);
 \draw[thick,-] (.47,-.2) to (0.3,-.2);\draw[thick,-to] (.7,-.2) to (.46,-.2);
\draw[thick,densely dashed] (-.4,-.2) to (0.3,-.2);
\draw[thick,-] (-.57,.2) to (-0.43,.2);\draw[thick,-to] (-.8,.2) to (-.56,.2);
\draw[thick,-] (.7,.2) to (.53,.2); \draw[thick,to-] (.54,.2) to (0.3,.2);
\draw[thick,densely dashed] (-.4,.2) to (0.3,.2);
\opendot{.7,.2};
\opendot{.3,.2};
\opendot{-.4,.2};
\opendot{-.8,.2};
\opendot{.7,-.2};
\opendot{.3,-.2};
\opendot{-.4,-.2};
\opendot{-.8,-.2};
\end{tikzpicture}
$, and define the normalization functions in such a way that
$c_i(\alpha_j) = (-1)^{\#(j\rightarrow i)}$ for this new quiver.
To construct $c_i$ explicitly,
for $i \in \Z+\frac{1}{2}$,
let $$
N_{>i} := \big|\{k \in \Z\:|\:-r \leq k \leq r\text{ and }k > i\}\big|.
$$
There is a homomorphism of Abelian 
groups $h_i^\vee:X \rightarrow \{\pm 1\}$
such that $h_i^\vee(\eps_k) = -1$ if $k > i$ and $N_{>i}$ is
even or if $k < i$ and $N_{>i}$ is odd,
and $h_i^\vee(\eps_k) = 1$ otherwise.
For $i \in I$, this function has the property that
$h_i^\vee(\alpha_j) = (-1)^{\delta_{i,j}}$ for all $j \in I$.
Then we can define $c_i:X \rightarrow\kk^\times$ by $c_i := h^\vee_{i-1}$.
We leave it to the reader to check that $\UU^\imath(\del,\zeta)$ defined with these new
parameters
is isomorphic to the one in \cref{pacific} via an isomorphism which
rescales caps, cups and crossings.
\end{rem}

%% file: s4-embedding.tex
\setcounter{section}{3}

\section{Categorification of the standard embedding}\label{s4-embedding}

In this section, we construct a 2-functor from the 2-iquantum group $\UU^\imath$ to a 2-category obtained by localizing the 2-quantum group $\UU$.
One can view this as a categorification of the standard embedding of $\dot\U^\imath_\Z$ into $\widehat{\U}_\Z$, although this statement should be taken with a pinch of salt since we are only able to do this after some localization---our formulae require certain polynomials in the dots to be invertible.
This 2-functor will be used in a critical way to prove the main results of the article in the final section.

\begin{table}
\begin{align*}
\begin{array}{|l|c|}
\hline
\hspace{13mm}\text{Generator}&\text{Degree}\\
\hline
\ \ \begin{tikzpicture}[Q,centerzero]
\draw[-to] (0,-0.3) \botlabel{i} -- (0,0.3)\toplabel{\phantom{i}};
\opendot{0,0};
\region{0.2,0}{\lambda};
\end{tikzpicture}\ 
\colon q_i^{\del_{\tau i}-h_{\tau i}(\lambda)}E_i \one_\lambda \Rightarrow q_i^{\del_{\tau i}-h_{\tau i}(\lambda)}E_i \one_\lambda&2 d_i\\
\ \ \begin{tikzpicture}[Q,centerzero]
\draw[to-] (0,-0.3) \botlabel{i} -- (0,0.3) \toplabel{\phantom{i}};
\opendot{0,0};
\region{0.2,0}{\lambda};
\end{tikzpicture}\ 
\colon F_i \one_\lambda \Rightarrow F_i \one_\lambda&2 d_i\\
\,\begin{tikzpicture}[Q,centerzero,scale=.9]
\draw[-to] (-0.3,-0.3) \botlabel{i} -- (0.3,0.3)\toplabel{\phantom{i}};
\draw[-to] (0.3,-0.3) \botlabel{j} -- (-0.3,0.3);
\region{0.35,0}{\lambda};
\end{tikzpicture}
\;\colon q_i^{\del_{\tau i}-h_{\tau i}(\lambda+\alpha_j)}
q_j^{\del_{\tau j}-h_{\tau j}(\lambda)}E_i E_j\one_\lambda \Rightarrow 
&-d_i a_{i,j}\\
\hspace{30mm}
q_i^{\del_{\tau i}-h_{\tau i}(\lambda)} q_j^{\del_{\tau j}-h_{\tau j}(\lambda+\alpha_i)}E_j E_i \one_\lambda&\\
\,\begin{tikzpicture}[Q,centerzero,scale=.9]
\draw[to-] (-0.3,-0.3) \botlabel{i} -- (0.3,0.3)\toplabel{\phantom{i}};
\draw[to-] (0.3,-0.3) \botlabel{j} -- (-0.3,0.3);
\region{0.35,0}{\lambda};
\end{tikzpicture}
\;\colon F_i F_j \one_\lambda \Rightarrow F_j F_i \one_\lambda
&-d_i a_{i,j}\\
\,\begin{tikzpicture}[Q,centerzero,scale=.9]
\draw[-to] (-0.3,-0.3) \botlabel{j} -- (0.3,0.3)\toplabel{\phantom{i}};
\draw[to-] (0.3,-0.3) \botlabel{i} -- (-0.3,0.3);
\region{0.35,0}{\lambda};
\end{tikzpicture}
\;\colon q_j^{\del_{\tau j}-h_{\tau j}(\lambda-\alpha_i)}E_j F_i \one_\lambda \Rightarrow q_j^{\del_{\tau j}-h_{\tau j}(\lambda)}F_i E_j \one_\lambda
&-d_i a_{i,\tau j}\\
\,\begin{tikzpicture}[Q,centerzero,scale=.9]
\draw[to-] (-0.3,-0.3) \botlabel{j} -- (0.3,0.3) \toplabel{\phantom{i}};
\draw[-to] (0.3,-0.3) \botlabel{i} -- (-0.3,0.3);
\region{0.35,0}{\lambda};
\end{tikzpicture}
\;\colon q_i^{\del_{\tau i}-h_{\tau i}(\lambda)}F_j E_i \one_\lambda \Rightarrow q_i^{\del_{\tau i}-h_{\tau i}(\lambda-\alpha_j)}E_i F_j \one_\lambda
&d_{i} a_{i, \tau j}\\
\:\begin{tikzpicture}[Q,centerzero]
\draw[-to] (-0.25,-0.15) \botlabel{i} to [out=90,in=90,looseness=3](0.25,-0.15);
\region{0.45,0.1}{\lambda};
\node at (0,.3) {$\phantom.$};
\node at (0,-.4) {$\phantom.$};
\end{tikzpicture}
\colon q_i^{\del_{\tau i}-h_{\tau i}(\lambda-\alpha_i)}E_i F_i \one_\lambda \Rightarrow \one_\lambda
&d_i(1\!+\!\del_i\!-\!h_i(\lambda)\!+\!h_{\tau i}(\lambda))\\
\:\begin{tikzpicture}[Q,centerzero]
\draw[-to] (-0.25,0.15) \toplabel{i} to[out=-90,in=-90,looseness=3] (0.25,0.15);
\region{0.45,-.1}{\lambda};
\node at (0,.2) {$\phantom.$};\node at (0,-.3) {$\phantom.$};
\end{tikzpicture}
\colon \one_\lambda \Rightarrow q_i^{\del_{\tau i}-h_{\tau i}(\lambda)}F_i E_i \one_\lambda
&d_i(1\!+\!\del_{\tau i}\!+\!h_i(\lambda)\!-\!h_{\tau i}(\lambda))
\\
\:\begin{tikzpicture}[Q,centerzero]
\draw[to-] (-0.25,-0.15) \botlabel{i} to [out=90,in=90,looseness=3](0.25,-0.15);
\region{0.45,0.1}{\lambda};
\node at (0,.3) {$\phantom.$};
\node at (0,-.4) {$\phantom.$};
\end{tikzpicture}
\colon q_i^{\del_{\tau i}-h_{\tau i}(\lambda)}F_i E_i \one_\lambda \Rightarrow \one_\lambda
&d_i(1\!-\!\del_{\tau i}\!+\!h_{i}(\lambda)\!+\!h_{\tau i}(\lambda))\\
\:\begin{tikzpicture}[Q,centerzero]
\draw[to-] (-0.25,0.15) \toplabel{i} to[out=-90,in=-90,looseness=3] (0.25,0.15);
\region{0.45,-.1}{\lambda};
\node at (0,.2) {$\phantom.$};\node at (0,-.3) {$\phantom.$};
\end{tikzpicture}
\colon \one_\lambda \Rightarrow q_i^{\del_{\tau i}-h_{\tau i}(\lambda-\alpha_i)}E_i F_i \one_\lambda
&d_i(1\!-\!\del_i\!-\!h_{i}(\lambda)\!-\!h_{\tau i}(\lambda))\\
\hline
\:\:\begin{tikzpicture}[Q,anchorbase,scale=.7]
\draw[-to,thin] (0,-.4)\botlabel{i} to (0,.5);
\clockwiseinternalbubble{0,.05};
\region{.45,0.05}{\lambda};
\node at (0,.45) {$\phantom{h}$};
\end{tikzpicture}\::q_i^{\del_{\tau i}-h_{\tau i}(\lambda)}E_i \one_\lambda\Rightarrow q_i^{\del_{\tau i}-h_{\tau i}(\lambda)}E_i \one_\lambda&
2d_i (\del_{\tau i}-h_{\tau i}(\lambda))\\
\:\:\begin{tikzpicture}[Q,anchorbase,scale=.7]
\node at (0,.4) {$\phantom{h}$};
\node at (0,-.42) {$\phantom{h}$};
\draw[-to,thin] (0,-.4)\botlabel{i} to (0,.5);
\anticlockwiseinternalbubble{0,.05};
\region{.4,0.05}{\lambda};
\end{tikzpicture}\:: q_i^{\del_{\tau i}-h_{\tau i}(\lambda)}E_i \one_\lambda\Rightarrow q_i^{\del_{\tau i}-h_{\tau i}(\lambda)}E_i \one_\lambda&-2 d_i(\del_{\tau i} - h_{\tau i}(\lambda))\\
\:\:\begin{tikzpicture}[Q,anchorbase,scale=.7]
\node at (0,.4) {$\phantom{h}$};
\node at (0,-.42) {$\phantom{h}$};
\draw[to-,thin] (0,-.4)\botlabel{i} to (0,.5);
\anticlockwiseinternalbubble{0,.05};
\region{.4,0.05}{\lambda};
\end{tikzpicture}\::F_i \one_\lambda\Rightarrow F_i \one_\lambda&2 d_i(\del_i+h_{\tau i}(\lambda))\\
\:\:\begin{tikzpicture}[Q,anchorbase,scale=.7]
\draw[to-,thin] (0,-.4)\botlabel{i} to (0,.5);
\clockwiseinternalbubble{0,.05};
\region{.45,0.05}{\lambda};
\node at (0,.45) {$\phantom{h}$};
\end{tikzpicture}\::F_i \one_\lambda\Rightarrow F_i \one_\lambda&
-2d_i (\del_{i}+h_{\tau i}(\lambda))\\\hline
\end{array}
\end{align*}
\caption{Adapted degrees of 2-morphisms in $\UU(\del,\zeta)$}\label{table3}
\end{table}

\subsection{Teleporters and internal bubbles}\label{newgrading}

We begin with a few reminders about $q$-envelopes; our conventions are the same as in \cite[Def.~6.10]{BE}. For a graded category $\catC$, 
its {\em $q$-envelope} (also known as {\em graded envelope}) 
is the graded category $\catC_q$ 
with objects that are formal symbols $q^n X$ for $n \in \Z$ and all objects
$X$ of $\catC$, with
$\Hom_{\catC_q}(q^n X, q^m Y)$ being the graded vector space $q^{m-n} \Hom_{\catC}(X,Y)$. Thus, a morphism $f:X \Rightarrow Y$ of degree $d$ in $\catC$ can also be viewed as a 2-morphism $f:q^n X \Rightarrow q^m Y$
in $\catC_q$ of degree $d+m-n$. Composition in 
$\catC_q$ is induced by composition in $\catC$ in the obvious way.

The category $\catC_q$ admits an invertible {\em grading shift functor}
$q:\catC_q\rightarrow \catC_q$ sending $q^n X$ to $q^{n+1} X$
and defined in the obvious way on morphisms.
The original category
$\catC$ is naturally identified with a full subcategory of $\catC_q$
so that object $X$ of $\catC$ is $q^0 X$ in $\catC_q$,
and the inclusion $\catC \rightarrow \catC_q$ is universal amongst all graded functors from $\catC$ to graded categories admitting a grading shift functor.
We also note that
\begin{equation}\label{awkwardinverse}
(\catC^\op)_q  = \left(\catC_{q^{-1}}\right)^\op,
\end{equation}
that is, the grading shift functor on $(\catC^\op)_q$ is the inverse of the grading shift functor of $\catC_q$ viewed as an endofunctor of 
$(\catC_q)^\op$.

If $\CATC$ is a graded 2-category, its $q$-envelope
$\CATC_q$
is the 2-category with the same objects as $\CATC$ and morphism categories that are the $q$-envelopes of the morphism categories of $\CATC$.
Horizontal composition of 1-morphisms in $\CATC_q$ 
is defined by $(q^n X) (q^m Y) := q^{n+m} (XY)$, and the definition of horizontal composition of 2-morphisms is also obvious. 
The original graded 2-category $\CATC$ is naturally identified with 
a wide ($=$ all objects) and full ($=$ all 2-morphisms) sub-2-category of $\CATC_q$.

Now we return to the setup of \cref{s3-iquantum}. In particular, $\tau:I
\rightarrow I$ is an involution, and
$\del = (\del_i)_{i \in I} \in \Z^I$ and $\zeta = (\zeta_i)_{i \in I}$
satisfy \cref{recap,being}.
We are going to work with a slightly different 
version $\UU(\del,\zeta)$ of the 2-quantum group $\UU$ which is graded in a way that is adapted to the standard embedding of $\dot\U^\imath$ into $\widehat{\U}$.
By definition, $\UU(\del,\zeta)$ is the wide and full sub-2-category of 
$\UU_q$ generated by the 1-morphisms $F_i \one_\lambda$ and
$q_i^{\del_{\tau i}-h_{\tau i}(\lambda)}E_i \one_\lambda$
for $i \in I$ and $\lambda \in X$.
The grading shift $q_i^{\del_{\tau i}-h_{\tau i}(\lambda)}$ 
matches the second term in the standard embedding \cref{unpredictable}.
The $q$-envelopes of $\UU$
and $\UU(\del,\zeta)$ are equal.
Also if we forget the grading, i.e., we view $\UU$ as a $\kk$-linear 2-category rather than a graded 2-category, then $\UU(\del,\zeta)$ is simply equal to $\UU$.

All string diagrams in the remainder of the 
section will represent 
2-morphisms in $\UU(\del,\zeta)$ rather than $\UU$, that is,
we will be using the upward string labelled $i$ to denote the
identity 2-endomorphism of $
q_i^{\del_{\tau i}-h_{\tau i}(\lambda)} E_i \one_\lambda$
(whereas before it was the identity 2-endomorphism of $E_i \one_\lambda$). The effect of this new convention is that the degrees of the 2-morphisms represented by cups and caps in this section are different to the degrees of the 2-morphisms represented by the same cups and caps in \cref{s2-quantum}.
We have listed the new degrees of all of the generating 2-morphisms of $\UU(\del,\zeta)$ in \cref{table3}.

Next we pass to a localization $\UUloc(\del,\zeta)$ of $\UU(\del,\zeta)$ which is obtained by adjoining two-sided inverses to the
2-morphisms\footnote{The invertibility of the morphisms in \cref{bensidea} is used in one place only---see the penultimate paragraph of the proof of \cref{boxingday}.}
\begin{equation}\label{bensidea}
\begin{tikzpicture}[Q,anchorbase,scale=1]
\draw[to-,thin] (0,-.3)\botlabel{i} to (0,.3);
\draw[-to,thin] (0.4,-.3)\botlabel{i} to (0.4,.3);
\Pinpin{0,0}{.4,0}{1.5,0}{Q_{i,\tau i}(x,-y)};
\region{2.65,0}{\lambda};
\end{tikzpicture}
\qquad\text{ and }\qquad
\begin{tikzpicture}[Q,anchorbase,scale=1]
\draw[-to,thin] (0,-.3)\botlabel{i} to (0,.3);
\draw[to-,thin] (0.4,-.3)\botlabel{i} to (0.4,.3);
\Pinpin{0,0}{.4,0}{1.5,0}{Q_{i,\tau i}(x,-y)};
\region{2.65,0}{\lambda};
\end{tikzpicture}
\end{equation}
for all $\lambda \in X$, $i \in I$ with $i \neq \tau i$,
and two-sided inverses to the 2-morphisms
\begin{equation}\label{teletubby}
\begin{tikzpicture}[Q,anchorbase,scale=1]
\draw[-to,thin] (-.5,-.3)\botlabel{i} to (-.5,.3);
\draw[-,ultra thick] (-.05,-.3) to (-.05,.3);
\draw[-to,thin] (0.4,-.3)\botlabel{\tau i} to (0.4,.3);
\Pinpin{.4,0}{-.5,0}{-1.1,0}{x+y};
\region{.6,0}{\lambda};
\end{tikzpicture}
\end{equation}
for all $\lambda \in X$ and $i \in I$,
where the thick string in the middle denotes some number (possibly zero) of other strings separating the outer strings.
We denote the inverse of \cref{teletubby} by a solid line which we refer to
as a {\em teleporter}:
\begin{equation}
\begin{tikzpicture}[Q,anchorbase,scale=1]
\draw[-to,thin] (-.5,-.3)\botlabel{i} to (-.5,.3);
\draw[-,ultra thick] (-.05,-.3) to (-.05,.3);
\draw[-to,thin] (0.4,-.3)\botlabel{\tau i} to (0.4,.3);
\teleporter{-.5,0}{0.4,0};
\region{.6,0}{\lambda};
\end{tikzpicture}
:=
\left(\:\:
\begin{tikzpicture}[Q,anchorbase,scale=1]
\draw[-,ultra thick] (-.05,-.3) to (-.05,.3);
\draw[-to,thin] (-.5,-.3)\botlabel{i} to (-.5,.3);
\draw[-to,thin] (0.4,-.3)\botlabel{\tau i} to (0.4,.3);
\Pinpin{.4,0}{-.5,0}{-1.1,0}{x+y};
\region{.6,0}{\lambda};
\end{tikzpicture}\right)^{-1}.\label{teleporter}
\end{equation}
We extend the teleporter notation further so that the endpoints do not need to be at the same horizontal
level, and the strings these points lie on do not need to be oriented in
the same direction. There is no ambiguity because of the relations \cref{lotsmore}.
We refer to \cite[$\S$3]{BWWbasis} for further
explanations in the special case that $i = \tau i$,
in which case
the dot is invertible in $\UUloc(\del,\zeta)$ with inverse
\begin{equation}
\begin{tikzpicture}[Q,centerzero]
\draw[-to] (0,-0.3) \botlabel{i} -- (0,0.3)\toplabel{\phantom{i}};
\multopendot{0,0}{east}{-1};
\region{0.2,0}{\lambda};
\end{tikzpicture}
:= 
\left(\ \begin{tikzpicture}[Q,centerzero]
\draw[-to] (0,-0.3) \botlabel{i} -- (0,0.3)\toplabel{\phantom{i}};
\opendot{0,0};
\region{0.2,0}{\lambda};
\end{tikzpicture}\!\right)^{-1}
=
{\textstyle\frac{1}{2}}
\begin{tikzpicture}[Q,centerzero]
\draw[-to] (0,-0.3) \botlabel{i} -- (0,0.3)\toplabel{\phantom{i}};
\region{0.6,0}{\lambda};
\bentteleporter{0,.15}{0,-.15}{(.35,.07) and (.35,-.07)};
\end{tikzpicture}\ .
\end{equation}

To work efficiently with teleporters,
we need the following variation on \cref{dgf}:
\begin{align}\label{dgfrpt}
\begin{tikzpicture}[Q,centerzero,scale=1.1]
\draw[-] (0,-0.3) -- (0,0.3);
\Circledbar{0,0}{u};
\end{tikzpicture}
&:=
\begin{tikzpicture}[Q,centerzero,scale=1.1]
\draw[-] (0,-0.3) -- (0,0.3);
\Pin{0,0}{.7,0}{\frac{1}{u+x}};
\end{tikzpicture}\
.
\end{align}
In fact, this can be obtained from \cref{dgf} by replacing $u$ by $-u$ then multiplying by $-1$. With this trick, the following is a consequence of \cref{gendotslide}:
\begin{align}
\label{gendotslidebar}
\begin{tikzpicture}[Q,centerzero,scale=1.3]
\draw[-to] (-0.3,-0.3) \botlabel{i} -- (0.3,0.3);
\draw[-to] (0.3,-0.3) \botlabel{j} -- (-0.3,0.3);
\Circledbar{0.15,0.15}{u};
\end{tikzpicture}-
\begin{tikzpicture}[Q,centerzero,scale=1.2]
\draw[-to] (-0.3,-0.3) \botlabel{i} -- (0.3,0.3);
\draw[-to] (0.3,-0.3) \botlabel{j} -- (-0.3,0.3);
\Circledbar{-0.15,-0.15}{u};
\end{tikzpicture}
&= \delta_{i,j}  \ 
\begin{tikzpicture}[Q,centerzero,scale=1.3]
\draw[-to] (-0.2,-0.3) \botlabel{i} -- (-0.2,0.3);
\draw[-to] (0.2,-0.3) \botlabel{i} -- (0.2,0.3);
\Circledbar{-0.2,0}{u};
\Circledbar{0.2,0}{u};
\end{tikzpicture} =    
\begin{tikzpicture}[Q,centerzero,scale=1.3]
\draw[-to] (-0.3,-0.3) \botlabel{i} -- (0.3,0.3);
\draw[-to] (0.3,-0.3) \botlabel{j} -- (-0.3,0.3);
\Circledbar{0.15,-0.15}{u};
\end{tikzpicture}-\begin{tikzpicture}[Q,centerzero,scale=1.3]
\draw[-to] (-0.3,-0.3) \botlabel{i} -- (0.3,0.3);
\draw[-to] (0.3,-0.3) \botlabel{j} -- (-0.3,0.3);
\Circledbar{-0.15,0.15}{u};
\end{tikzpicture}\ .
\end{align}
Here are a couple of useful identities involving teleporters:
\begin{align}
\begin{tikzpicture}[Q,centerzero,scale=1.4]
\draw[-to] (-0.3,-0.3) \botlabel{i} -- (0.3,0.3);
\draw[-to] (0.3,-0.3) \botlabel{j} -- (-0.3,0.3);
\draw[-to] (.6,-.3) \botlabel{\tau i}--(.6,.3);
\teleporter{.15,.15}{.6,.15};
\end{tikzpicture}
-
\begin{tikzpicture}[Q,centerzero,scale=1.4]
\draw[-to] (-0.3,-0.3) \botlabel{i} -- (0.3,0.3);
\draw[-to] (0.3,-0.3) \botlabel{j} -- (-0.3,0.3);
\draw[-to] (.6,-.3) \botlabel{\tau i}--(.6,.3);
\teleporter{-.15,-.15}{.6,-.15};
\end{tikzpicture}
&= \delta_{i,j}  \ 
\begin{tikzpicture}[Q,centerzero,scale=1.4]
\draw[-to] (.6,-.3) \botlabel{\tau i}--(.6,.3);
\draw[-to] (-0.2,-0.3) \botlabel{i} -- (-0.2,0.3);
\draw[-to] (0.2,-0.3) \botlabel{j} -- (0.2,0.3);
\teleporter{-.2,-.15}{.6,-.15};
\teleporter{.2,.15}{.6,.15};
\end{tikzpicture}\ ,\label{teleslide}&
\begin{tikzpicture}[Q,centerzero,scale=1.4]
\draw[-to] (-0.3,-0.3) \botlabel{i} -- (0.3,0.3);
\draw[-to] (0.3,-0.3) \botlabel{j} -- (-0.3,0.3);
\draw[-to] (.6,-.3) \botlabel{\tau j}--(.6,.3);
\teleporter{.15,-.15}{.6,-.15};
\end{tikzpicture}
-
\begin{tikzpicture}[Q,centerzero,scale=1.4]
\draw[-to] (-0.3,-0.3) \botlabel{i} -- (0.3,0.3);
\draw[-to] (0.3,-0.3) \botlabel{j} -- (-0.3,0.3);
\draw[-to] (.6,-.3) \botlabel{\tau j}--(.6,.3);
\teleporter{-.15,.15}{.6,.15};
\end{tikzpicture}
&= \delta_{i,j}  \ 
\begin{tikzpicture}[Q,centerzero,scale=1.4]
\draw[-to] (.6,-.3) \botlabel{\tau j}--(.6,.3);
\draw[-to] (-0.2,-0.3) \botlabel{i} -- (-0.2,0.3);
\draw[-to] (0.2,-0.3) \botlabel{j} -- (0.2,0.3);
\teleporter{-.2,.15}{.6,.15};
\teleporter{.2,-.15}{.6,-.15};
\end{tikzpicture}\ ,\\
\begin{tikzpicture}[Q,centerzero,scale=1.4]
\draw[-to] (-.2,-.3) \botlabel{i}--(-.2,.3);
\draw[-to] (.2,-.3) \botlabel{\tau i}--(.2,.3);
\Circled{-.2,.1}{u};
\teleporter{-.2,-.15}{.2,-.15};
\end{tikzpicture}
-
\begin{tikzpicture}[Q,centerzero,scale=1.4]
\draw[-to] (-.2,-.3) \botlabel{i}--(-.2,.3);
\draw[-to] (.2,-.3) \botlabel{\tau i}--(.2,.3);
\Circledbar{.2,.1}{u};
\teleporter{-.2,-.15}{.2,-.15};
\end{tikzpicture}
&=
\begin{tikzpicture}[Q,centerzero,scale=1.4]
\draw[-to] (-.2,-.3) \botlabel{i}--(-.2,.3);
\draw[-to] (.2,-.3) \botlabel{\tau i}--(.2,.3);
\Circled{-.2,0}{u};
\Circledbar{.2,0}{u};
\end{tikzpicture}\ ,
&
\begin{tikzpicture}[Q,centerzero,scale=1.4]
\draw[-to] (-.2,-.3) \botlabel{i}--(-.2,.3);
\draw[-to] (.2,-.3) \botlabel{\tau i}--(.2,.3);
\Circled{.2,.1}{u};
\teleporter{-.2,-.15}{.2,-.15};
\end{tikzpicture}
-
\begin{tikzpicture}[Q,centerzero,scale=1.4]
\draw[-to] (-.2,-.3) \botlabel{i}--(-.2,.3);
\draw[-to] (.2,-.3) \botlabel{\tau i}--(.2,.3);
\Circledbar{-.2,.1}{u};
\teleporter{-.2,-.15}{.2,-.15};
\end{tikzpicture}
&=
\begin{tikzpicture}[Q,centerzero,scale=1.4]
\draw[-to] (-.2,-.3) \botlabel{i}--(-.2,.3);
\draw[-to] (.2,-.3) \botlabel{\tau i}--(.2,.3);
\Circled{.2,0}{u};
\Circledbar{-.2,0}{u};
\end{tikzpicture}\ .\label{partialfractions}
\end{align}
The first equality in \cref{partialfractions} follows because
$\frac{1}{u-x}-\frac{1}{u+y}=\frac{x+y}{(u-x)(u+y)}$, and the second
equality is similar.
There are plenty of variations
of \cref{teleslide,partialfractions}
obtained by rotating or adding additional strings passed over by the teleporter.
We will apply these as needed later on.

Finally, we introduce one more 
shorthand for the following 2-morphisms in
$\UUloc(\del,\zeta)$ which we refer to as
{\em internal bubbles}. The definition of these involves 
$\zeta=(\zeta_i)_{i \in I}$ for the first time. 
In the special case that $i = \tau i$, the internal bubbles were defined already in \cite[(3.38)--(3.39)]{BWWbasis}
by exactly the same formula as below.
In general, for $i \in I$ and $\lambda \in X$, we let
\begin{align}\label{internal1}
\begin{tikzpicture}[Q,anchorbase,scale=1]
\draw[-to,thin] (0,-.4)\botlabel{i} to (0,.5);
\clockwiseinternalbubble{0,.05};
\end{tikzpicture}
&:=(-1)^{\del_{\tau i} + 1}\zeta_{\tau i}
\left[\begin{tikzpicture}[Q,anchorbase,scale=1]
	\draw[-to,thin] (0.3,-0.4)\botlabel{i} to (.3,.5);
	\multopendot{.3,.2}{west}{\del_{\tau i}};
\filledclockwisebubble{1.1,0};	
\stringlabel{1.1,-.3}{\tau i};
\Circledbar{.3,-.1}{u};
\bubblelabel{1.1,0}{u};
\end{tikzpicture}\right]_{u:-1}\!\!
-
(-1)^{\del_{\tau i}+1}\zeta_{\tau i}
\begin{tikzpicture}[Q,anchorbase,scale=1]
	\draw[-to,thin] (0.28,-0.4)\botlabel{i} to (.28,.5);
	\multopendot{.28,.2}{west}{\del_{\tau i}};
\stringlabel{1.2,-.4}{\tau i};
\clockwisebubble{1.2,-.1};
\teleporter{.27,-.1}{1,-.1};
\end{tikzpicture}
,\\\label{internal2}
\begin{tikzpicture}[Q,anchorbase,scale=1]
\draw[-to,thin] (0,-.4)\botlabel{i} to (0,.5);
\anticlockwiseinternalbubble{0,.05};
\end{tikzpicture}
&:=\zeta_i
\left[\begin{tikzpicture}[Q,anchorbase,scale=1]
	\draw[-to,thin] (-0.3,-0.4)\botlabel{i} to (-.3,.5);
	\multopendot{-.3,.2}{east}{\del_i};
\filledanticlockwisebubble{-1.1,0};	
\bubblelabel{-1.1,0}{u};
\stringlabel{-1.1,-.3}{\tau i};
\Circledbar{-.3,-.1}{u};
\end{tikzpicture}\right]_{u:-1}\!\!
-
\zeta_i\ 
\begin{tikzpicture}[Q,anchorbase,scale=1]
	\draw[-to,thin] (-0.28,-0.4)\botlabel{i} to (-.28,.5);
	\multopendot{-.28,.2}{east}{\del_i};
\anticlockwisebubble{-1.2,-.1};
\stringlabel{-1.2,-.4}{\tau i};
\teleporter{-.27,-.1}{-1,-.1};
\end{tikzpicture},\\\label{internal3}
\begin{tikzpicture}[Q,anchorbase,scale=1]
\draw[to-,thin] (0,-.4)\botlabel{i} to (0,.5);
\anticlockwiseinternalbubble{0,.05};
\end{tikzpicture}
&:= \zeta_{i}
\left[\begin{tikzpicture}[Q,anchorbase,scale=1]
	\draw[to-,thin] (0.3,-0.4)\botlabel{i} to (.3,.5);
	\multopendot{.3,.25}{west}{\del_{i}};
\filledanticlockwisebubble{1.1,0};	
\stringlabel{1.1,-.3}{\tau i};
\Circledbar{.3,-.05}{u};
\bubblelabel{1.1,0}{u};
\end{tikzpicture}\right]_{u:-1}
\!\!-
\zeta_{i}
\begin{tikzpicture}[Q,anchorbase,scale=1]
	\draw[to-,thin] (0.28,-0.4)\botlabel{i} to (.28,.5);
	\multopendot{.28,.2}{west}{\del_{i}};
\stringlabel{1.2,-.4}{\tau i};
\anticlockwisebubble{1.2,-.1};
\teleporter{.27,-.1}{1,-.1};
\end{tikzpicture}
,\\\label{internal4}
\begin{tikzpicture}[Q,anchorbase,scale=1]
\draw[to-,thin] (0,-.4)\botlabel{i} to (0,.5);
\clockwiseinternalbubble{0,.05};
\end{tikzpicture}
&:=
(-1)^{\del_{\tau i}+1}
\zeta_{\tau i}
\left[\begin{tikzpicture}[Q,anchorbase,scale=1]
	\draw[to-,thin] (-0.3,-0.4)\botlabel{i} to (-.3,.5);
	\multopendot{-.3,.3}{east}{\del_{\tau i}};
\filledclockwisebubble{-1.1,0};	
\bubblelabel{-1.1,0}{u};
\stringlabel{-1.1,-.3}{\tau i};
\Circledbar{-.3,-.05}{u};
\end{tikzpicture}\right]_{u:-1}
\!\!-
(-1)^{\del_{\tau i}+1}\zeta_{\tau i}\ 
\begin{tikzpicture}[Q,anchorbase,scale=1]
	\draw[to-,thin] (-0.28,-0.4)\botlabel{i} to (-.28,.5);
	\multopendot{-.28,.2}{east}{\del_{\tau i}};
\clockwisebubble{-1.2,-.1};
\stringlabel{-1.2,-.4}{\tau i};
\teleporter{-.27,-.1}{-1,-.1};
\end{tikzpicture}.
\end{align}
These definitions apply for all possible weights $\lambda$ labelling the rightmost 2-cells---the definitions are independent of the weight so we have omitted it, following our usual practice.
However, the degrees of internal bubbles depend on the choice of the
weight $\lambda$; they are listed in the bottom part of \cref{table3}.

It is easy to see that 
internal bubbles slide over cups and caps in the obvious way.
Also, they commute with dots and with 
other internal bubbles on the same string.
Moreover, the clockwise and counterclockwise internal bubbles
are mutual inverses. This is proved in \cref{Appendix-B}, together with various other technical identities
involving internal bubbles which are needed in the proof of the important \cref{boxingday} below.

\subsection{Construction of the 2-functor}\label{sinternalbubbles}

Let $\UU^\imath(\del,\zeta)$ be as in \cref{def2iqg}.
Now we would like to construct a 2-functor from 
$\UU^\imath(\del,\zeta)$ to $\UUloc(\del,\zeta)$. However, this does not quite make sense since the object set of $\UUloc(\del,\zeta)$ is $X$, whereas the object set of $\UU^\imath$ is $X^\imath$.
So first we must modify $\UUloc(\del,\zeta)$, contracting its object set 
along the quotient map $X \twoheadrightarrow X^\imath$.
We denote the result by $\widehat{\UUloc}(\del,\zeta)$.
By definition, this has object set
$X^\imath$, and a 1-morphism 
$G \one_\lambda = \one_\mu G:\lambda
\rightarrow \mu$ in $\widehat{\UUloc}(\del,\zeta)$ is a word
$G$ in the generators $E_i$ and $F_i$ whose total weight $\varpi \in X$
obtained by adding $\alpha_i$ 
for each $E_i$ and $-\alpha_i$ for each $F_i$ 
satisfies $\mu = \lambda + \varpi$.
For two such words $G \one_\lambda$ and $H \one_\lambda$, 
the 2-morphism space
$\Hom_{\widehat{\UUloc}(\del,\zeta)}(G \one_\lambda, H \one_\lambda)$ is the graded vector space
$\bigoplus_{n \in \Z}\Hom_{\widehat{\UUloc}(\del,\zeta)}(G \one_\lambda, H \one_\lambda)_n$ where
$$
\Hom_{\widehat{\UUloc}(\del,\zeta)}(G \one_\lambda, H \one_\lambda)_n :=
\prod_{\hat \lambda} \Hom_{\UUloc(\del,\zeta)}(G \one_{\hat\lambda}, H\one_{\hat\lambda})_n
$$ 
taking the product over all
pre-images $\hat\lambda \in X$ of $\lambda \in X^\imath$.
Any 2-morphism in this space can be represented as a tuple $f = (f_{\hat \lambda})$ with components indexed by the pre-images $\hat\lambda$ 
of $\lambda$. Horizontal and vertical composition in $\widehat{\UUloc}(\del,\zeta)$ are
induced by the ones in $\UUloc(\del,\zeta)$ in an obvious way.

Let $\Add(\widehat{\UUloc}(\del,\zeta))$ be the additive envelope of $\widehat{\UUloc}(\del,\zeta)$.
It has the same objects as $\widehat{\UUloc}(\del,\zeta)$,
but its 1-morphisms are finite direct sums of 1-morphisms in $\widehat{\UUloc}(\del,\zeta)$, and 2-morphisms are matrices of 2-morphisms in $\widehat{\UUloc}(\del,\zeta)$. Vertical composition of 2-morphisms is by matrix multiplication. In the sequel, we will suppress the matrices, denoting 2-morphisms in $\Add(\widehat{\UUloc}(\del,\zeta))$ simply as finite sums of 2-morphisms. The string labels at the top and bottom of each string diagram in such a sum determine the matrix entry that it belongs to in an unambiguous way.

Now we are ready to state and prove the main theorem. 
This generalizes \cite[Th.~4.2]{BWWbasis} which is essentially the same result for the special case of the nil-Brauer category. The monoidal functor constructed in 
\cite[Th.~4.2]{BWWbasis}
composed with the monoidal automorphism of the nil-Brauer category that negates crossings and dots
is equal to the 2-functor here for the split iquantum group of rank one.
We will refer to this often in the proof of the next 
theorem since we need all of the 
calculations made in \cite{BWWbasis} for strings labelled by $i = \tau i$.

For $i \in I$ with $i \neq \tau i$, 
we define $\sgn(i) \in \{\pm 1\}$ so that $i \in I_{\sgn(i)}$
in the decomposition \cref{usefuldecomposition}.
These signs are used in an essential way in the
next theorem in order to break some symmetry. 

\begin{theo}\label{boxingday}
There is a strict graded 2-functor
$\Xi^\imath:\UU^\imath(\del,\zeta) \rightarrow \Add(\widehat{\UUloc}(\del,\zeta))$
which is the identity on objects,
takes $B_i \one_\lambda \mapsto F_i \one_\lambda \oplus 
q_i^{\del_{i}-h_{i}(\lambda)}E_{\tau i} \one_\lambda$, and is defined on 2-morphisms by
\begin{align}\label{psi1}
\Xi^\imath\left(\ \begin{tikzpicture}[iQ,centerzero]
\draw[-] (0,-0.3) \botlabel{i} -- (0,0.3);
\closeddot{0,0};
\region{0.2,0}{\lambda};
\end{tikzpicture}\right)_{\hat \lambda}
&
:=
\begin{tikzpicture}[Q,centerzero]
\draw[to-] (0,-0.3) \botlabel{i} -- (0,0.3);
\opendot{0,0};
\region{0.2,0}{\hat\lambda};
\end{tikzpicture}-
\begin{tikzpicture}[Q,centerzero]
\draw[-to] (0,-0.3) \botlabel{\tau i} -- (0,0.3);
\opendot{0,0};
\region{0.2,0}{\hat\lambda};
\end{tikzpicture},
\\\label{psi2}
\Xi^\imath\left(\ \begin{tikzpicture}[iQ,baseline=0mm]
\draw[-] (-0.25,-0.15) \botlabel{\tau i} to [out=90,in=90,looseness=3](0.25,-0.15);
\region{0.45,0.1}{\lambda};
\end{tikzpicture}\right)_{\hat \lambda}\ &:=
\begin{tikzpicture}[Q,centerzero,scale=.8]
	\draw[-,thin] (-0.4,-0.3)\botlabel{i} to[out=90, in=180] (0,0.3);
	\draw[-to,thin] (-0,0.3) to[out = 0, in = 90] (0.4,-0.3);
  \region{.7,0}{\hat\lambda};
\end{tikzpicture}+
\begin{tikzpicture}[Q,centerzero,scale=.8]
	\draw[to-,thin] (-0.4,-0.3)\botlabel{\tau i} to [out=90, in=180] (0,0.3);
	\draw[-,thin] (-0,0.3)to[out = 0, in = 90] (0.4,-0.3);
    \clockwiseinternalbubble{.3,.1};
  \region{.75,0}{\hat\lambda};
\end{tikzpicture},\\\label{psi3}
\Xi^\imath\left(\ \begin{tikzpicture}[iQ,baseline=-1.5mm]
\draw[-] (-0.25,0.15) \toplabel{\tau i} to [out=-90,in=-90,looseness=3](0.25,0.15);
\region{0.45,0}{\lambda};
\end{tikzpicture}\right)_{\hat\lambda}\ &:=
\begin{tikzpicture}[Q,centerzero,scale=.8]
	\draw[to-,thin] (-0.4,0.3) \toplabel{i}to[out=-90, in=180] (0,-0.3);
	\draw[-,thin] (-0,-0.3) to[out = 0, in = -90] (0.4,0.3);
    \anticlockwiseinternalbubble{-.3,-.1};
  \region{.65,0}{\hat\lambda};
\end{tikzpicture}+\begin{tikzpicture}[Q,centerzero,scale=.8] 
	\draw[-,thin] (-0.4,0.3)\toplabel{\tau i} to[out=-90, in=180] (0,-0.3);
	\draw[-to,thin] (-0,-0.3) to[out = 0, in = -90] (0.4,0.3);
  \region{.65,0}{\hat\lambda};
\end{tikzpicture},\\
\label{whybeta}
\Xi^\imath\left(\!\begin{tikzpicture}[baseline=-.6mm,iQ]
\draw[-] (-0.25,0) arc(180:-180:0.25);
\node at (-0.42,0) {\strandlabel{\tau i}};
\region{.97,0}{\lambda};
\node at (.55,0) {$(u)$};
\end{tikzpicture}\!\right)_{\hat\lambda} \ 
&:=
(-1)^{h_{\tau i}(\hat\lambda)} u^{\del_i}\ 
\begin{tikzpicture}[Q,centerzero,scale=1]
\draw[to-] (-0.68,0) arc(180:-180:0.25);
\node[black] at (0.26,0) {$(-u)$};
\node at (-.43,-.4) {\strandlabel{\tau i}};
\region{.85,0}{\hat\lambda};
\end{tikzpicture} 
\begin{tikzpicture}[Q,centerzero,scale=1]
\draw[-to] (-.25,0) arc(180:-180:0.25);
\node[black] at (.54,0) {$(u)$};
\node at (0,-0.4) {\strandlabel{i}};
\end{tikzpicture}\ ,\\
\label{psi5a}
\Xi^\imath\left(\ \begin{tikzpicture}[iQ,centerzero]
\draw[-] (-0.3,-0.3) \botlabel{i} -- (0.3,0.3);
\draw[-] (0.3,-0.3) \botlabel{j} -- (-0.3,0.3);
\region{0.35,0}{\lambda};
\end{tikzpicture}\right)_{\hat\lambda}&:=
-r_{i,j}^{-1}\begin{tikzpicture}[Q,centerzero,scale=.7]
	\draw[to-] (0.6,-.6)\botlabel{j} to (-0.6,.6);
	\draw[to-] (-0.6,-.6)\botlabel{i} to (0.6,.6);
  \region{.5,0}{\hat\lambda};
\end{tikzpicture}
-r_{\tau j, \tau i}r_{\tau i, j} \begin{tikzpicture}[Q,centerzero,scale=.7]
	\draw[-to,thin] (0.6,-.6)\botlabel{\tau j} to (-0.6,.6);
	\draw[-to,thin] (-0.6,-.6)\botlabel{\tau i} to (0.6,.6);
  \region{.5,0}{\hat\lambda};
\end{tikzpicture}
-r_{\tau i, j}\begin{tikzpicture}[Q,centerzero,scale=.7]
	\draw[to-,thin] (0.6,-.6)\botlabel{j} to (-0.6,.6);
	\draw[-to,thin] (-0.6,-.6)\botlabel{\tau i} to (0.6,.6);
  \region{.5,0}{\hat\lambda};
\end{tikzpicture}
-r_{j,i}^{-1}r_{\tau j, i}^{-1}\begin{tikzpicture}[Q,centerzero,scale=.7]
	\draw[-to,thin] (0.6,-.6)\botlabel{\tau j} to (-0.6,.6);
	\draw[to-,thin] (-0.6,-.6)\botlabel{i} to (0.6,.6);
\anticlockwiseinternalbubble{-.3,.3};
\clockwiseinternalbubble{.3,-.3};
  \region{.7,0.1}{\hat\lambda};
\end{tikzpicture}\\\intertext{if $i \neq j$ and $i \neq \tau j$,}
\Xi^\imath\left(\ \begin{tikzpicture}[iQ,centerzero]
\draw[-] (-0.3,-0.3) \botlabel{i} -- (0.3,0.3);
\draw[-] (0.3,-0.3) \botlabel{i} -- (-0.3,0.3);
\region{0.35,0}{\lambda};
\end{tikzpicture}\right)_{\hat\lambda}&:=
-\begin{tikzpicture}[Q,centerzero,scale=.7]
	\draw[to-] (0.6,-.6)\botlabel{i} to (-0.6,.6);
	\draw[to-] (-0.6,-.6)\botlabel{i} to (0.6,.6);
  \region{.5,0}{\hat\lambda};
\end{tikzpicture}
- \begin{tikzpicture}[Q,centerzero,scale=.7]
	\draw[-to,thin] (0.6,-.6)\botlabel{\tau i} to (-0.6,.6);
	\draw[-to,thin] (-0.6,-.6)\botlabel{\tau i} to (0.6,.6);
  \region{.5,0}{\hat\lambda};
\end{tikzpicture}
-\sgn(i)\begin{tikzpicture}[Q,centerzero,scale=.7]
	\draw[to-,thin] (0.6,-.6)\botlabel{i} to (-0.6,.6);
	\draw[-to,thin] (-0.6,-.6)\botlabel{\tau i} to (0.6,.6);
  \region{.5,0}{\hat\lambda};
\end{tikzpicture}
+\sgn(i)\begin{tikzpicture}[Q,centerzero,scale=.7]
	\draw[-to,thin] (0.6,-.6)\botlabel{\tau i} to (-0.6,.6);
	\draw[to-,thin] (-0.6,-.6)\botlabel{i} to (0.6,.6);
\anticlockwiseinternalbubble{-.3,.3};
\clockwiseinternalbubble{.3,-.3};
  \region{.7,0.1}{\hat\lambda};
\end{tikzpicture}
-\begin{tikzpicture}[Q,centerzero,scale=.7]
	\draw[-to,thin] (-0.5,-.6)\botlabel{\tau i} to (-0.5,.6);
	\draw[to-,thin] (.1,-.6)\botlabel{i} to (.1,.6);
\teleporter{-.5,0}{.1,0};
  \region{.4,0}{\hat\lambda};
\end{tikzpicture}
+
\begin{tikzpicture}[Q,centerzero,scale=.7]
	\draw[to-,thin] (-0.5,-.6)\botlabel{i} to (-0.5,.6);
	\draw[-to,thin] (.1,-.6)\botlabel{\tau i} to (.1,.6);
\teleporter{-.5,0}{.1,0};
  \region{.4,0}{\hat\lambda};
\end{tikzpicture},\label{psi5b}\\
\Xi^\imath\left(\ \begin{tikzpicture}[iQ,centerzero]
\draw[-] (-0.3,-0.3) \botlabel{i} -- (0.3,0.3);
\draw[-] (0.3,-0.3) \botlabel{\tau i} -- (-0.3,0.3);
\region{0.35,0}{\lambda};
\end{tikzpicture}\right)_{\hat\lambda}&:=
\sgn(i)\begin{tikzpicture}[Q,centerzero,scale=.7]
	\draw[-to,thin] (0.6,-.6)\botlabel{i} to (-0.6,.6);
	\draw[-to,thin] (-0.6,-.6)\botlabel{\tau i} to (0.6,.6);
  \region{.5,0}{\hat\lambda};
\end{tikzpicture}
-\sgn(i)\begin{tikzpicture}[Q,centerzero,scale=.7]
	\draw[to-] (0.6,-.6)\botlabel{\tau i} to (-0.6,.6);
	\draw[to-] (-0.6,-.6)\botlabel{i} to (0.6,.6);
  \region{.5,0}{\hat\lambda};
\end{tikzpicture}
-\begin{tikzpicture}[Q,centerzero,scale=.7]
	\draw[to-,thin] (0.6,-.6)\botlabel{\tau i} to (-0.6,.6);
	\draw[-to,thin] (-0.6,-.6)\botlabel{\tau i} to (0.6,.6);
  \region{.5,0}{\hat\lambda};
\end{tikzpicture}
-\begin{tikzpicture}[Q,centerzero,scale=.7]
	\draw[-to,thin] (0.6,-.6)\botlabel{i} to (-0.6,.6);
	\draw[to-,thin] (-0.6,-.6)\botlabel{i} to (0.6,.6);
\anticlockwiseinternalbubble{-.3,.3};
\clockwiseinternalbubble{.3,-.3};
  \region{.7,0.1}{\hat\lambda};
\end{tikzpicture}
-\begin{tikzpicture}[Q,centerzero,scale=.7]
 	\draw[to-,thin] (-0.45,-.6)\botlabel{i} to[out=90,in=180] (-.1,-.1) to[out=0,in=90] (0.25,-.6);
 	\draw[-to,thin] (-0.45,.6) \toplabel{\tau i}to[out=-90,in=180] (-.1,.1) to[out=0,in=-90] (0.25,.6);
\teleporter{-.4,.3}{-.4,-.3};
     \clockwiseinternalbubble{.15,-.25};
  \region{.5,0.1}{\hat\lambda};
\end{tikzpicture}+\begin{tikzpicture}[Q,centerzero,scale=.7]
 	\draw[-to] (-0.45,-.6)\botlabel{\tau i} to[out=90,in=180] (-.1,-.1) to[out=0,in=90] (0.25,-.6);
 	\draw[to-] (-0.45,.6)\toplabel{i}to[out=-90,in=180] (-.1,.1) to[out=0,in=-90] (0.25,.6);
\teleporter{.19,.3}{.19,-.3};
\anticlockwiseinternalbubble{-.35,.25};
  \region{.55,0}{\hat\lambda};
\end{tikzpicture}\label{psi5c}\\\intertext{if $i \neq \tau i$, and}\label{psi5d}
\Xi^\imath\left(\ \begin{tikzpicture}[iQ,centerzero]
\draw[-] (-0.3,-0.3) \botlabel{i} -- (0.3,0.3);
\draw[-] (0.3,-0.3) \botlabel{i} -- (-0.3,0.3);
\region{0.35,0}{\lambda};
\end{tikzpicture}\right)_{\hat\lambda}&:=\!
-\begin{tikzpicture}[Q,centerzero,scale=.7]
	\draw[to-] (0.6,-.6)\botlabel{i} to (-0.6,.6);
	\draw[to-] (-0.6,-.6)\botlabel{i} to (0.6,.6);
  \region{.5,0}{\hat\lambda};
\end{tikzpicture}
- \!\!\begin{tikzpicture}[Q,centerzero,scale=.7]
	\draw[-to,thin] (0.6,-.6)\botlabel{i} to (-0.6,.6);
	\draw[-to,thin] (-0.6,-.6)\botlabel{i} to (0.6,.6);
  \region{.5,0}{\hat\lambda};
\end{tikzpicture}
-\!\!\begin{tikzpicture}[Q,centerzero,scale=.7]
	\draw[to-,thin] (0.6,-.6)\botlabel{i} to (-0.6,.6);
	\draw[-to,thin] (-0.6,-.6)\botlabel{i} to (0.6,.6);
  \region{.5,0}{\hat\lambda};
\end{tikzpicture}
-\!\!\begin{tikzpicture}[Q,centerzero,scale=.7]
	\draw[-to,thin] (0.6,-.6)\botlabel{i} to (-0.6,.6);
	\draw[to-,thin] (-0.6,-.6)\botlabel{i} to (0.6,.6);
\anticlockwiseinternalbubble{-.3,.3};
\clockwiseinternalbubble{.3,-.3};
  \region{.7,0.1}{\hat\lambda};
\end{tikzpicture}\!
-\begin{tikzpicture}[Q,centerzero,scale=.7]
	\draw[-to,thin] (-0.5,-.6)\botlabel{i} to (-0.5,.6);
	\draw[to-,thin] (.1,-.6)\botlabel{i} to (.1,.6);
\teleporter{-.5,0}{.1,0};
  \region{.4,0}{\hat\lambda};
\end{tikzpicture}
+
\begin{tikzpicture}[Q,centerzero,scale=.7]
	\draw[to-,thin] (-0.5,-.6)\botlabel{i} to (-0.5,.6);
	\draw[-to,thin] (.1,-.6)\botlabel{i} to (.1,.6);
\teleporter{-.5,0}{.1,0};
  \region{.4,0}{\hat\lambda};
\end{tikzpicture}\!
-\begin{tikzpicture}[Q,centerzero,scale=.7]
 	\draw[to-,thin] (-0.45,-.6)\botlabel{i} to[out=90,in=180] (-.1,-.1) to[out=0,in=90] (0.25,-.6);
 	\draw[-to,thin] (-0.45,.6) \toplabel{i}to[out=-90,in=180] (-.1,.1) to[out=0,in=-90] (0.25,.6);
\teleporter{-.4,.3}{-.4,-.3};
     \clockwiseinternalbubble{.15,-.25};
  \region{.5,0.1}{\hat\lambda};
\end{tikzpicture}
\!+\begin{tikzpicture}[Q,centerzero,scale=.7]
 	\draw[-to] (-0.45,-.6)\botlabel{i} to[out=90,in=180] (-.1,-.1) to[out=0,in=90] (0.25,-.6);
 	\draw[to-] (-0.45,.6)\toplabel{i}to[out=-90,in=180] (-.1,.1) to[out=0,in=-90] (0.25,.6);
\teleporter{.19,.3}{.19,-.3};
\anticlockwiseinternalbubble{-.35,.25};
  \region{.55,0}{\hat\lambda};
\end{tikzpicture}
\end{align}
if $i = \tau i$.
\end{theo}

\begin{proof}
First, one checks that the degrees are consistent. This follows by inspecting \cref{table2,table3} in \cref{new2cat,newgrading}.
To construct the functor, we define it on objects and generating 1-morphisms as in the statement of the theorem, and on generating 2-morphisms
according to \cref{psi1,psi2,psi3,psi5a,psi5b,psi5c,psi5d} together with the following for all $i \in I$ with $i \neq \tau i$,
which specifies $\Xi^\imath$ on the fake bubble generators:
\begin{equation}\label{fakebubblemap}
\Xi^\imath\left(\begin{tikzpicture}[baseline=-1mm,iQ]
\draw[fakebubble] (-0.25,0) arc(180:-180:0.25);
\node at (-0.42,0) {\strandlabel{\tau i}};
\region{0.4,0}{\lambda};
\bubblelabel{0,0}{u};
\end{tikzpicture}\right)_{\hat\lambda}:=
\zeta_i(\lambda)
\left[u^{\del_i}\ 
\begin{tikzpicture}[Q,centerzero,scale=1]
\draw[to-] (-0.68,0) arc(180:-180:0.25);
\node[black] at (0.26,0) {$(-u)$};
\node at (-.43,-.4) {\strandlabel{\tau i}};
\region{.9,0}{\hat\lambda};
\end{tikzpicture} 
\begin{tikzpicture}[Q,centerzero,scale=1]
\draw[-to] (-.25,0) arc(180:-180:0.25);
\node[black] at (0.54,0) {$(u)$};
\node at (0,-0.4) {\strandlabel{i}};
\end{tikzpicture}\right]_{u:\geq 0}
\end{equation}
To show that it is well defined, we need to check that the images of the defining relations of $\UU^\imath(\del,\zeta)$ are satisfied in $\widehat{\UUloc}(\del,\zeta)$.
This will take up most of the rest of the proof.

Before we verify any further relations, we check that \cref{whybeta} holds.
If $i = \tau i$ then it follows from 
\cite[(4.2)]{BWWbasis}; one needs to compose with the automorphism of the nil-Brauer category that negates dots and crossings 
because the monoidal functor there is defined in a slightly different way to the 2-functor $\Xi^\imath$ here.
Now suppose that $i \neq \tau i$.
Since \cref{fakebubblemap} is true by the definition of $\Xi^\imath$, and remembering
the definition \cref{ibubblegeneratingfunction},
to establish \cref{whybeta} it suffices to show that
\begin{align*}
\Xi^\imath\left(
\begin{tikzpicture}[baseline=-1mm,iQ,scale=1.4]
\draw[-] (-0.25,0) arc(180:-180:0.25);
\stringlabel{-.42,0}{\tau i};
\region{0.5,0}{\lambda};
\Circled{.25,0}{u};
\end{tikzpicture}\right)_{\hat \lambda}=
\zeta_i(\lambda) \left[u^{\del_i}\ 
\begin{tikzpicture}[Q,centerzero,scale=1.4]
\draw[to-] (-0.68,0) arc(180:-180:0.25);
\node[black] at (0.14,0) {$(-u)$};
\stringlabel{-.43,-.4}{\tau i};
\region{.6,0}{\hat\lambda};
\end{tikzpicture}\  
\begin{tikzpicture}[Q,centerzero,scale=1.4]
\draw[-to] (-.25,0) arc(180:-180:0.25);
\node[black] at (.46,0) {$(u)$};
\stringlabel{0,-.4}{i};
\end{tikzpicture}\right]_{u:<0}.
\end{align*}
Using \cref{wombles}, the left hand side is
\begin{align*}
\begin{tikzpicture}[Q,centerzero,scale=.95]
\draw[to-] (0.68,0) arc(0:360:0.4);
\node at (.28,-.5) {\strandlabel{i}};
\anticlockwiseinternalbubble{-0.02,-.22};
\Circled{-0.02,.25}{u};
\region{.9,0}{\hat\lambda};
\end{tikzpicture}
+\  \begin{tikzpicture}[Q,centerzero,scale=.95]
\draw[to-] (-0.68,0) arc(180:-180:0.4);
\node at (-.28,-.5) {\strandlabel{\tau i}};
\clockwiseinternalbubble{0.02,.22};
\Circledbar{0.02,-.25}{u};
\region{0.4,0}{\hat\lambda};
\end{tikzpicture}&=
\zeta_i \left[
u^{\del_i}\ \begin{tikzpicture}[Q,centerzero,scale=1.3]
\draw[to-] (0.68,0) arc(0:360:0.25);
\stringlabel{.43,-.35}{i};
\stringlabel{-.4,-.3}{\tau i};
\Circled{.18,0}{u};
\filledanticlockwisebubble{-.4,0};	
\bubblelabel{-.4,0}{-u};
\region{-0.05,0}{\hat\lambda};
\end{tikzpicture} 
-u^{\del_i}\ \begin{tikzpicture}[Q,centerzero,scale=1.3]
\draw[to-] (0.68,0) arc(0:360:0.25);
\stringlabel{.43,-.35}{i};
\stringlabel{-.4,-.3}{\tau i};
\Circled{.23,.15}{u};
\anticlockwisebubble{-.4,0};	
\teleporter{.2,-.11}{-.23,-.11};
\region{-0.05,0.15}{\hat\lambda};
\end{tikzpicture} 
-
u^{\del_i}\ \begin{tikzpicture}[Q,centerzero,scale=1.3]
\draw[to-] (-0.68,0) arc(180:-180:0.25);
\stringlabel{-.43,-.35}{\tau i};
\stringlabel{.4,-.3}{i};
\Circledbar{-.18,0}{u};
\filledclockwisebubble{.4,0};	
\bubblelabel{.4,0}{u};
\region{0.05,0}{\hat\lambda};
\end{tikzpicture} 
+u^{\del_i}\ \begin{tikzpicture}[Q,centerzero,scale=1.3]
\draw[to-] (-0.68,0) arc(180:-180:0.25);
\stringlabel{-.43,-.35}{\tau i};
\stringlabel{.4,-.3}{i};
\Circledbar{-.23,-.15}{u};
\clockwisebubble{.4,0};	
\teleporter{-.19,.11}{.24,.11};
\region{0.05,-0.1}{\hat\lambda};
\end{tikzpicture} 
\right]_{u:<0}\\
&\stackrel{\cref{partialfractions}}{=}
\zeta_i 
\left[u^{\del_i}\;
\begin{tikzpicture}[Q,centerzero,scale=1.3]
\filledanticlockwisebubble{-.68,0};
\bubblelabel{-.68,0}{-u};
\stringlabel{-.68,-.3}{\tau i};
\region{-.3,0}{\hat\lambda};
\end{tikzpicture} \,
\begin{tikzpicture}[Q,centerzero,scale=1.3]
\draw[-to] (-.25,0) arc(180:-180:0.25);
\Circled{.25,0}{u};
\stringlabel{0,-.4}{i};
\end{tikzpicture}-u^{\del_i}\ \begin{tikzpicture}[Q,centerzero,scale=1.3]
\draw[to-] (-0.68,0) arc(180:-180:0.25);
\stringlabel{-.43,-.35}{\tau i};
\Circledbar{-.18,0}{u};
\region{0.1,0}{\hat\lambda};
\end{tikzpicture} \,
\begin{tikzpicture}[Q,centerzero,scale=1.3]
\filledclockwisebubble{0.1,0};
\bubblelabel{.1,0}{u};
\stringlabel{0.1,-.3}{i};
\end{tikzpicture}
-u^{\del_i}\ \begin{tikzpicture}[Q,centerzero,scale=1.3]
\draw[to-] (-0.68,0) arc(180:-180:0.25);
\stringlabel{-.43,-.4}{\tau i};
\Circledbar{-.18,0}{u};
\region{.1,0}{\hat\lambda};
\end{tikzpicture} \;
\begin{tikzpicture}[Q,centerzero,scale=1.3]
\draw[-to] (-.25,0) arc(180:-180:0.25);
\Circled{.25,0}{u};
\stringlabel{0,-.4}{i};
\end{tikzpicture}
\right]_{u:<0}\!\\&\:=
\zeta_i\left[ u^{\del_i}\;
\left(
\begin{tikzpicture}[Q,centerzero,scale=1.3]
\filledanticlockwisebubble{-.68,0};
\bubblelabel{-.68,0}{-u};
\stringlabel{-.68,-.3}{\tau i};
\region{-.3,0}{\hat\lambda};
\end{tikzpicture} -
\begin{tikzpicture}[Q,centerzero,scale=1.3]
\draw[to-] (-0.68,0) arc(180:-180:0.25);
\stringlabel{-.43,-.35}{\tau i};
\Circledbar{-.18,0}{u};
\region{0.1,0}{\hat\lambda};
\end{tikzpicture}
\right)
\left(\begin{tikzpicture}[Q,centerzero,scale=1.3]
\filledclockwisebubble{0.1,0};
\bubblelabel{.1,0}{u};
\stringlabel{0.1,-.3}{i};
\region{.5,0}{\hat\lambda};
\end{tikzpicture}
+\begin{tikzpicture}[Q,centerzero,scale=1.3]
\draw[to-] (-0.68,0) arc(-180:180:0.25);
\stringlabel{-.43,-.4}{\tau i};
\Circledbar{-.18,0}{u};
\region{.1,0}{\hat\lambda};
\end{tikzpicture}
\right)
\right]_{u:<0}\\
&\stackrel{\cref{bubblegeneratingfunctions}}{=}
\zeta_i (-1)^{h_{\tau i}(\hat\lambda)} c_{\tau i}(\hat\lambda)^{-1}
c_i(\hat\lambda)
\left[ u^{\del_i}\ 
\begin{tikzpicture}[Q,centerzero,scale=1.4]
\draw[to-] (-0.68,0) arc(180:-180:0.25);
\node[black] at (0.14,0) {$(-u)$};
\stringlabel{-.43,-.4}{\tau i};
\region{.6,0}{\hat\lambda};
\end{tikzpicture}\  
\begin{tikzpicture}[Q,centerzero,scale=1.4]
\draw[-to] (-.25,0) arc(180:-180:0.25);
\node[black] at (.46,0) {$(u)$};
\stringlabel{0,-.4}{i};
\end{tikzpicture}
\right]_{u:<0},
\end{align*}
which is what we want due to the definitions
\cref{howthingswork,purer}.

To complete the proof, it remains to check the defining relations \cref{ibubblerel,iinfgrass,ibubslide,izigzag,ipivotal,ipitchfork,idotslide,icurl,iquadratic,ibraid}.
This is done in \cref{Appendix-C}. The calculation requires the other
properties of internal bubbles derived in \cref{Appendix-B}, and
also uses \cref{necessary} proved in \cref{Appendix-A}
in order to check the difficult
ibraid relation.
\end{proof}

\subsection{Special case: categorification of comultiplication}\label{sscm}

Let $\UU$ be any 2-quantum group
with parameters $(Q_{i,j}(x,y))_{i,j \in I}$
and normalization homomorphisms $c_i$,
leading to $r_{i,j}$ as in \cref{gollygosh}.
We also need
the 2-quantum group we denote by $\UU \odot \UU$,
which is defined from the doubled Cartan datum 
from \cref{iexisting}
using the parameters
\begin{align}
Q_{i^+,j^+}(x,y) = Q_{i^-,j^-}(x,y) &:= Q_{i,j}(x,y),&
Q_{i^+,j^-}(x,y)=Q_{i^-,j^+}(x,y) &:= 1,\\
c_{i^+}\big((\lambda^+,\lambda^-)\big) &:= c_i(\lambda^+),&c_{i^-}\big
((\lambda^+,\lambda^-)\big)
&:= c_i(\lambda^-),\\
r_{i^+,j^+} = r_{i^-,j^-}&:= r_{i,j},&
r_{i^+,j^-}=r_{i^-,j^+}&:= 1.
\end{align}
We think of $\UU\odot\UU$ as a commuting pair of copies of $\UU$,
with generating 1-morphisms
$E_{i^+} 
\one_{(\lambda^+,\lambda^-)}$,
$F_{i^+} \one_{(\lambda^+,\lambda^-)}$,
$E_{i^-} \one_{(\lambda^+,\lambda^-)}$
and $F_{i^-} \one_{(\lambda^+,\lambda^-)}$.
In fact, we prefer to work with a variant
$\UU\circledast\UU$, which is the wide and full sub-2-category of $(\UU\odot\UU)_q$ generated by
the 1-morphisms
\begin{align}
(E_i\otimes 1) \one_{(\lambda^+,\lambda^-)}
&:= E_{i^+} \one_{(\lambda^+,\lambda^-)},
&
(K_i \otimes E_i) \one_{(\lambda^+,\lambda^-)}
&:= q_i^{h_i(\lambda^+)} 
E_{i^-} \one_{(\lambda^+,\lambda^-)},\label{starb1}\\
(F_i\otimes K_i^{-1}) \one_{(\lambda^+,\lambda^-)}
&:= q_i^{-h_i(\lambda^-)} F_{i^+} \one_{(\lambda^+,\lambda^-)},
&
(1 \otimes F_i) \one_{(\lambda^+,\lambda^-)}
&:= 
F_{i^-} \one_{(\lambda^+,\lambda^-)},\label{starb2}
\end{align}
for all $i \in I$ and $\lambda^+,\lambda^- \in X$.
We will be drawing string diagrams to represent 2-morphisms in $\UU\circledast\UU$ rather than in $\UU\odot\UU$,
so that 
upward strings labelled $i^+$ and $i^-$
denote the identity 2-endomorphisms of 
$E_i \otimes 1$ and $K_i \otimes E_i$, and 
downward strings labelled $i^+$ and $i^-$
denote the identity 2-endomorphisms of $F_i \otimes K_i^{-1}$ and  $1 \otimes F_i$, respectively. This just means that we are using another grading convention for string diagrams.
In the ungraded setting,
$\UU\circledast \UU$ and $\UU \odot\UU$ are the same $\kk$-linear 2-category.
In string diagrams for 2-morphisms in $\UU\circledast\UU$, we use the convention of labelling a 2-cell by
$\lambda^{\pm}$ instead of by $(\lambda^+,\lambda^-)$ to save a little space.

Next, we pass to a localization $\UU\:\underline{\circledast}\:\UU$
of $\UU\circledast\UU$
by adjoining
two-sided inverses to all 2-morphisms given by a dot on an upward string labelled $i^+$ {\em minus} a dot on an upward string labelled $i^-$. We denote these two-sided inverses by
\begin{align}
\begin{tikzpicture}[Q,anchorbase,scale=1]
\draw[-to,thin] (-.5,-.3)\botlabel{i^+} to (-.5,.3);
\draw[-,ultra thick] (-.05,-.3) to (-.05,.3);
\draw[-to,thin] (0.4,-.3)\botlabel{i^-} to (0.4,.3);
\Teleporter{-.5,0}{0.4,0};
\region{.95,0}{\lambda^\pm};
\end{tikzpicture}\!\!
&:=
\left(\:\:
\begin{tikzpicture}[Q,anchorbase,scale=1]
\draw[-,ultra thick] (-.05,-.3) to (-.05,.3);
\draw[-to,thin] (-.5,-.3)\botlabel{i^+} to (-.5,.3);
\draw[-to,thin] (0.4,-.3)\botlabel{i^-} to (0.4,.3);
\Pinpin{.4,0}{-.5,0}{-1.2,0}{x-y};
\region{.95,0}{\lambda^\pm};
\end{tikzpicture}\!\right)^{-1}\!\!\!\!,
&
\begin{tikzpicture}[Q,anchorbase,scale=1]
\draw[-to,thin] (-.5,-.3)\botlabel{i^-} to (-.5,.3);
\draw[-,ultra thick] (-.05,-.3) to (-.05,.3);
\draw[-to,thin] (0.4,-.3)\botlabel{i^+} to (0.4,.3);
\Teleporter{-.5,0}{0.4,0};
\region{.95,0}{\lambda^\pm};
\end{tikzpicture}\!\!
&:=
\left(\:\:
\begin{tikzpicture}[Q,anchorbase,scale=1]
\draw[-,ultra thick] (-.05,-.3) to (-.05,.3);
\draw[-to,thin] (-.5,-.3)\botlabel{i^-} to (-.5,.3);
\draw[-to,thin] (0.4,-.3)\botlabel{i^+} to (0.4,.3);
\Pinpin{.4,0}{-.5,0}{-1.2,0}{y-x};
\region{.95,0}{\lambda^\pm};
\end{tikzpicture}\!\right)^{-1}\!\!\!\!.
\label{Teleporter}
\end{align}
They are similar to the teleporters defined earlier in the section but with a different sign (which is why we are using a new color).
As we did before, 
we extend the notation so that we can also attach teleporters to downward strings.
Then we define internal bubbles in $\UU\:\underline{\circledast}\:\UU$
by
\begin{align}\label{newinternal1}
\begin{tikzpicture}[Q,anchorbase,scale=1]
\draw[-to,thin] (0,-.4)\botlabel{i^+} to (0,.5);
\Clockwiseinternalbubble{0,.05};
\end{tikzpicture}
&:=
\left[\begin{tikzpicture}[Q,anchorbase,scale=1]
\draw[-to,thin] (0.3,-0.4)\botlabel{i^+} to (.3,.5);
\filledclockwisebubble{.9,0.05};	
\stringlabel{.9,-.28}{i^-};
\bubblelabel{.9,0.05}{u};
\Circled{.3,0.05}{u};
\end{tikzpicture}\right]_{u:-1}\!\!
+
\begin{tikzpicture}[Q,anchorbase,scale=1]
	\draw[-to,thin] (0.28,-0.4)\botlabel{i^+} to (.28,.5);
\stringlabel{1,-.28}{i^-};
\clockwisebubble{1,0.05};
\Teleporter{.27,0.05}{.8,0.05};
\end{tikzpicture},&
\begin{tikzpicture}[Q,anchorbase,scale=1]
\draw[-to,thin] (0,-.4)\botlabel{i^+} to (0,.5);
\Anticlockwiseinternalbubble{0,.05};
\end{tikzpicture}
&:=
\left[\begin{tikzpicture}[Q,anchorbase,scale=1]
	\draw[-to,thin] (-0.3,-0.4)\botlabel{i^+} to (-.3,.5);
\filledanticlockwisebubble{-.9,0.05};	
\bubblelabel{-.9,0.05}{u};
\stringlabel{-.9,-.28}{i^-};
\Circled{-.3,.05}{u};
\end{tikzpicture}\right]_{u:-1}\!\!
+
\begin{tikzpicture}[Q,anchorbase,scale=1]
	\draw[-to,thin] (-0.28,-0.4)\botlabel{i^+} to (-.28,.5);
\anticlockwisebubble{-1,.05};
\stringlabel{-1,-.28}{i^-};
\Teleporter{-.27,.05}{-.8,.05};
\end{tikzpicture},\\\label{newinternal2}
\begin{tikzpicture}[Q,anchorbase,scale=1]
\draw[to-,thin] (0,-.4)\botlabel{i^+} to (0,.5);
\Anticlockwiseinternalbubble{0,.05};
\end{tikzpicture}
&:= 
\left[\begin{tikzpicture}[Q,anchorbase,scale=1]
	\draw[to-,thin] (0.3,-0.4)\botlabel{i^+} to (.3,.5);
\filledanticlockwisebubble{.9,0.05};	
\stringlabel{.9,-.28}{i^-};
\Circled{.3,.05}{u};
\bubblelabel{.9,0.05}{u};
\end{tikzpicture}\right]_{u:-1}
\!\!+
\begin{tikzpicture}[Q,anchorbase,scale=1]
	\draw[to-,thin] (0.28,-0.4)\botlabel{i^+} to (.28,.5);
\stringlabel{1,-.28}{i^-};
\anticlockwisebubble{1,.05};
\Teleporter{.27,.05}{.8,.05};
\end{tikzpicture}
,&
\begin{tikzpicture}[Q,anchorbase,scale=1]
\draw[to-,thin] (0,-.4)\botlabel{i^+} to (0,.5);
\Clockwiseinternalbubble{0,.05};
\end{tikzpicture}
&:=
\left[\begin{tikzpicture}[Q,anchorbase,scale=1]
	\draw[to-,thin] (-0.3,-0.4)\botlabel{i^+} to (-.3,.5);
\filledclockwisebubble{-.9,0.05};	
\bubblelabel{-.9,0.05}{u};
\stringlabel{-.9,-.28}{i^-};
\Circled{-.3,.05}{u};
\end{tikzpicture}\right]_{u:-1}
\!\!+
\begin{tikzpicture}[Q,anchorbase,scale=1]
	\draw[to-,thin] (-0.28,-0.4)\botlabel{i^+} to (-.28,.5);
\clockwisebubble{-1,.05};
\stringlabel{-1,-.28}{i^-};
\Teleporter{-.27,.05}{-.8,.05};
\end{tikzpicture},\\\label{newinternal3}
\begin{tikzpicture}[Q,anchorbase,scale=1]
\draw[-to,thin] (0,-.4)\botlabel{i^-} to (0,.5);
\Clockwiseinternalbubble{0,.05};
\end{tikzpicture}
&:=
\left[\begin{tikzpicture}[Q,anchorbase,scale=1]
	\draw[-to,thin] (0.3,-0.4)\botlabel{i^-} to (.3,.5);
\filledclockwisebubble{.9,0.05};	
\stringlabel{.9,-.3}{i^+};
\Circled{.3,0.05}{u};
\bubblelabel{.9,0.05}{u};
\end{tikzpicture}\right]_{u:-1}\!\!
-
\begin{tikzpicture}[Q,anchorbase,scale=1]
	\draw[-to,thin] (0.28,-0.4)\botlabel{i^-} to (.28,.5);
\stringlabel{1,-.3}{i^+};
\clockwisebubble{1,0.05};
\Teleporter{.27,0.05}{.8,0.05};
\end{tikzpicture},&
\begin{tikzpicture}[Q,anchorbase,scale=1]
\draw[-to,thin] (0,-.4)\botlabel{i^-} to (0,.5);
\Anticlockwiseinternalbubble{0,.05};
\end{tikzpicture}
&:=
\left[\begin{tikzpicture}[Q,anchorbase,scale=1]
	\draw[-to,thin] (-0.3,-0.4)\botlabel{i^-} to (-.3,.5);
\filledanticlockwisebubble{-.9,0.05};	
\bubblelabel{-.9,0.05}{u};
\stringlabel{-.9,-.3}{i^+};
\Circled{-.3,.05}{u};
\end{tikzpicture}\right]_{u:-1}\!\!
-
\begin{tikzpicture}[Q,anchorbase,scale=1]
	\draw[-to,thin] (-0.28,-0.4)\botlabel{i^-} to (-.28,.5);
\anticlockwisebubble{-1,.05};
\stringlabel{-1,-.3}{i^+};
\Teleporter{-.27,.05}{-.8,.05};
\end{tikzpicture},\\\label{newinternal4}
\begin{tikzpicture}[Q,anchorbase,scale=1]
\draw[to-,thin] (0,-.4)\botlabel{i^-} to (0,.5);
\Anticlockwiseinternalbubble{0,.05};
\end{tikzpicture}
&:= 
\left[\begin{tikzpicture}[Q,anchorbase,scale=1]
	\draw[to-,thin] (0.3,-0.4)\botlabel{i^-} to (.3,.5);
\filledanticlockwisebubble{.9,0.05};	
\stringlabel{.9,-.3}{i^+};
\Circled{.3,.05}{u};
\bubblelabel{.9,0.05}{u};
\end{tikzpicture}\right]_{u:-1}
\!\!-
\begin{tikzpicture}[Q,anchorbase,scale=1]
	\draw[to-,thin] (0.28,-0.4)\botlabel{i^-} to (.28,.5);
\stringlabel{1,-.3}{i^+};
\anticlockwisebubble{1,.05};
\Teleporter{.27,.05}{.8,.05};
\end{tikzpicture}
,&
\begin{tikzpicture}[Q,anchorbase,scale=1]
\draw[to-,thin] (0,-.4)\botlabel{i^-} to (0,.5);
\Clockwiseinternalbubble{0,.05};
\end{tikzpicture}
&:=
\left[\begin{tikzpicture}[Q,anchorbase,scale=1]
	\draw[to-,thin] (-0.3,-0.4)\botlabel{i^-} to (-.3,.5);
\filledclockwisebubble{-.9,0.05};	
\bubblelabel{-.9,0.05}{u};
\stringlabel{-.9,-.3}{i^+};
\Circled{-.3,.05}{u};
\end{tikzpicture}\right]_{u:-1}
\!\!-
\begin{tikzpicture}[Q,anchorbase,scale=1]
	\draw[to-,thin] (-0.28,-0.4)\botlabel{i^-} to (-.28,.5);
\clockwisebubble{-1,.05};
\stringlabel{-1,-.3}{i^+};
\Teleporter{-.27,.05}{-.8,.05};
\end{tikzpicture}.
\end{align}
These have similar properties to the internal bubbles
defined before. For example, the clockwise and counterclockwise internal bubbles are two-sided inverses. Note also that the bubbles in the above definitions can be drawn on the other side of the string, e.g., we have that
$$
\left[\begin{tikzpicture}[Q,anchorbase,scale=1]
	\draw[-to,thin] (0.3,-0.4)\botlabel{i^+} to (.3,.5);
\filledclockwisebubble{.9,0.05};	
\stringlabel{.9,-.28}{i^-};
\Circled{.3,0.05}{u};
\bubblelabel{.9,0.05}{u};
\end{tikzpicture}\right]_{u:-1}\!\!
+
\begin{tikzpicture}[Q,anchorbase,scale=1]
	\draw[-to,thin] (0.28,-0.4)\botlabel{i^+} to (.28,.5);
\stringlabel{1,-.28}{i^-};
\clockwisebubble{1,0.05};
\Teleporter{.27,0.05}{.8,0.05};
\end{tikzpicture}
=
\left[\begin{tikzpicture}[Q,anchorbase,scale=1]
	\draw[-to,thin] (0.3,-0.4)\botlabel{i^+} to (.3,.5);
\filledclockwisebubble{-.3,0.05};	
\stringlabel{-.3,-.28}{i^-};
\Circled{.3,0.05}{u};
\bubblelabel{-.3,0.05}{u};
\end{tikzpicture}\right]_{u:-1}\!\!
+
\begin{tikzpicture}[Q,anchorbase,scale=1]
	\draw[-to,thin] (0.28,-0.4)\botlabel{i^+} to (.28,.5);
\stringlabel{-.5,-.28}{i^-};
\clockwisebubble{-.5,0.05};
\Teleporter{.28,0.05}{-.3,0.05};
\end{tikzpicture}\ .
$$
These statements follow from \cref{slushy,bubblyinverses} since the new internal bubbles are images of (a special case of) the old ones under an automorphism $\Omega\otimes J$; see the last sentence of the proof of the next theorem.
We also point out that the formulae defining these internal bubbles are the same as the formulae defining analogous morphisms in the Heisenberg category from \cite[(5.27)--(5.28)]{K0} and they appeared already for $\mathfrak{sl}_2$ in \cite[(6.3)--(6.4)]{unfurling}.

Finally, we contract the object set of $\UU\:\underline{\circledast}\:\UU$
along the quotient map ${\pmb X} \twoheadrightarrow X, (\lambda^+,\lambda^-) \mapsto \lambda^++\lambda^-$
to obtain a graded 2-category $\UU\:\underline{\widehat{\circledast}}\:\UU$. 
Its object set is $X$. A 1-morphism $G\one_\lambda = \one_\mu G:\lambda\rightarrow\mu$ in $\UU\:\underline{\widehat{\circledast}}\:\UU$
is a monomial $G$ in $E_i \otimes 1, K_i \otimes E_i, 
F_i \otimes K_i^{-1}, 
1 \otimes F_i\:(i \in I)$
whose total weight obtained by
adding $\alpha_i$ for each $E_i\otimes 1$ or $K_i \otimes E_{i}$
and $-\alpha_i$ for each $F_i \otimes K_i^{-1}$ or $1 \otimes F_i$
is equal to $\mu-\lambda$. Then
$$
\Hom_{\UU\:\underline{\widehat{\circledast}}\:\UU}(
G\one_\lambda, H\one_\lambda)_n := \prod_{\lambda^++\lambda^-=\lambda}
\Hom_{\UU\:\underline{\circledast}\:\UU}\big(G\one_{(\lambda^+,\lambda^-)},
H\one_{(\lambda^+,\lambda^-)}\big)_n.
$$

Now we can state a reformulation of \cref{boxingday} in this special case. The 2-functor $\Xi$ that it defines is a categorical version of the comultiplication $\Delta:\U \rightarrow \U \otimes_{\Q(q)} \U$.
The devoted reader could reprove this theorem from scratch by checking relations directly---it is easier 
than 
the proof of \cref{boxingday}
but not by much.

\begin{theo}\label{valentinesday}
There is a unique strict graded 2-functor
$\Xi:\UU \rightarrow \Add(\UU\:\underline{\widehat{\circledast}}\:\UU)$
which is the identity on objects,
takes $E_i \one_\lambda \mapsto (E_i \otimes 1) \one_\lambda \oplus (K_i \otimes E_i) \one_\lambda$
and $F_i \one_\lambda \mapsto (1 \otimes F_i) \one_\lambda
\oplus (F_i \otimes K_i^{-1}) \one_\lambda$, 
and is defined on 2-morphisms by
\begin{align*}
\Xi\left(\ \begin{tikzpicture}[Q,centerzero,scale=.9]
\draw[-to] (0,-0.3) \botlabel{i} -- (0,0.3);
\opendot{0,0};
\region{0.3,0}{\lambda};
\end{tikzpicture}\right)_{\!(\lambda^+,\lambda^-)}\!\!\!
&
:=
\begin{tikzpicture}[Q,centerzero,scale=.9]
\draw[-to] (0,-0.3) \botlabel{i^+} -- (0,0.3);
\opendot{0,0};
\region{0.5,0}{\lambda^\pm};
\end{tikzpicture}+
\begin{tikzpicture}[Q,centerzero,scale=.9]
\draw[-to] (0,-0.3) \botlabel{i^-} -- (0,0.3);
\opendot{0,0};
\region{0.5,0}{\lambda^\pm};
\end{tikzpicture}\ ,\\
\Xi\left(\ \begin{tikzpicture}[Q,centerzero,scale=.9]
\draw[-to] (-0.3,-0.3) \botlabel{i} -- (0.3,0.3);
\draw[-to] (0.3,-0.3) \botlabel{j} -- (-0.3,0.3);
\region{0.45,0}{\lambda};
\end{tikzpicture}\right)_{\!(\lambda^+,\lambda^-)}\!\!\!&:=\!
\begin{tikzpicture}[Q,centerzero,scale=.5]
	\draw[-to] (0.6,-.6)\botlabel{j^+} to (-0.6,.6);
	\draw[-to] (-0.6,-.6)\botlabel{i^+} to (0.6,.6);
  \region{1,0}{\lambda^\pm};
\end{tikzpicture}
\!+\!\! \begin{tikzpicture}[Q,centerzero,scale=.5]
	\draw[-to,thin] (0.6,-.6)\botlabel{j^-} to (-0.6,.6);
	\draw[-to,thin] (-0.6,-.6)\botlabel{i^-} to (0.6,.6);
  \region{1,0}{\lambda^\pm};
\end{tikzpicture}
\!\!+\!(-1)^{\delta_{i,j}}\!\!\!\begin{tikzpicture}[Q,centerzero,scale=.5]
	\draw[-to,thin] (0.6,-.6)\botlabel{j^+} to (-0.6,.6);
	\draw[-to,thin] (-0.6,-.6)\botlabel{i^-} to (0.6,.6);
  \region{1,0}{\lambda^\pm};
\end{tikzpicture}
\!\!+\!r_{j,i}^{-1}\!\begin{tikzpicture}[Q,centerzero,scale=.5]
	\draw[-to,thin] (0.6,-.6)\botlabel{j^-} to (-0.6,.6);
	\draw[-to,thin] (-0.6,-.6)\botlabel{i^+} to (0.6,.6);
\Babyclockwiseinternalbubble{-.3,.3};
\Babyanticlockwiseinternalbubble{.3,-.3};
  \region{1.1,0.1}{\lambda^\pm};
\end{tikzpicture}
\!\!-\!\delta_{i,j}\!\begin{tikzpicture}[Q,centerzero,scale=.5]
	\draw[-to,thin] (-0.5,-.6)\botlabel{i^-} to (-0.5,.6);
	\draw[-to,thin] (.3,-.6)\botlabel{i^+} to (.3,.6);
\Teleporter{-.5,0}{.3,0};
  \region{1.1,0}{\lambda^\pm};
\end{tikzpicture}
\!\!+\!\delta_{i,j}\!
\begin{tikzpicture}[Q,centerzero,scale=.5]
	\draw[-to,thin] (-0.5,-.6)\botlabel{i^+} to (-0.5,.6);
	\draw[-to,thin] (.3,-.6)\botlabel{i^-} to (.3,.6);
\Teleporter{-.5,0}{.3,0};
  \region{1.1,0}{\lambda^\pm};
\end{tikzpicture}\!,\\
\Xi\left(\ \begin{tikzpicture}[Q,centerzero]
\draw[-to] (-0.25,0.15) \toplabel{i} to [out=-90,in=-90,looseness=3](0.25,0.15);
\region{0.45,0}{\lambda};
\end{tikzpicture}\!\right)_{\!(\lambda^+,\lambda^-)}\!\!\! &:=
\begin{tikzpicture}[Q,anchorbase,scale=.7]
	\draw[-,thin] (-0.4,0.3) \toplabel{i^+}to[out=-90, in=180] (0,-0.3);
	\draw[-to,thin] (-0,-0.3) to[out = 0, in = -90] (0.4,0.3);
    \Babyanticlockwiseinternalbubble{.3,-.1};
  \region{1,0}{\lambda^\pm};
\end{tikzpicture}+\begin{tikzpicture}[Q,anchorbase,scale=.7] 
	\draw[-,thin] (-0.4,0.3)\toplabel{i^-} to[out=-90, in=180] (0,-0.3);
	\draw[-to,thin] (-0,-0.3) to[out = 0, in = -90] (0.4,0.3);
  \region{.9,0}{\lambda^\pm};
\end{tikzpicture},\qquad\qquad\quad\!
\Xi\left(\ \begin{tikzpicture}[Q,centerzero]
\draw[-to] (-0.25,-0.15) \botlabel{i} to [out=90,in=90,looseness=3](0.25,-0.15);
\region{0.45,0.1}{\lambda};
\end{tikzpicture}\!\right)_{\!(\lambda^+,\lambda^-)}\!\!\! :=
\begin{tikzpicture}[Q,anchorbase,scale=.7]
	\draw[-,thin] (-0.4,-0.3)\botlabel{i^+} to[out=90, in=180] (0,0.3);
	\draw[-to,thin] (-0,0.3) to[out = 0, in = 90] (0.4,-0.3);
    \Babyclockwiseinternalbubble{-.3,.1};
  \region{.9,0}{\lambda^\pm};
\end{tikzpicture}+
\begin{tikzpicture}[Q,anchorbase,scale=.7]
	\draw[-,thin] (-0.4,-0.3)\botlabel{i^-} to [out=90, in=180] (0,0.3);
	\draw[-to,thin] (-0,0.3)to[out = 0, in = 90] (0.4,-0.3);
  \region{.9,0}{\lambda^\pm};
\end{tikzpicture}
\end{align*}
for $i,j \in I$ and $\lambda^+,\lambda^-\in X$
with $\lambda^++\lambda^-=\lambda$.
We also have that
\begin{align*}
\Xi\left(\ \begin{tikzpicture}[Q,centerzero,scale=.9]
\draw[to-] (0,-0.3) \botlabel{i} -- (0,0.3);
\opendot{0,0};
\region{0.3,0}{\lambda};
\end{tikzpicture}\right)_{\!(\lambda^+,\lambda^-)}\!\!\!
&
=
\begin{tikzpicture}[Q,centerzero,scale=.9]
\draw[to-] (0,-0.3) \botlabel{i^+} -- (0,0.3);
\opendot{0,0};
\region{0.5,0}{\lambda^\pm};
\end{tikzpicture}+
\begin{tikzpicture}[Q,centerzero]
\draw[to-] (0,-0.3) \botlabel{i^-} -- (0,0.3);
\opendot{0,0};
\region{0.5,0}{\lambda^\pm};
\end{tikzpicture},\\
\Xi\left(\ \begin{tikzpicture}[Q,centerzero,scale=.9]
\draw[to-] (-0.3,-0.3) \botlabel{i} -- (0.3,0.3);
\draw[to-] (0.3,-0.3) \botlabel{j} -- (-0.3,0.3);
\region{0.45,0}{\lambda};
\end{tikzpicture}\right)_{\!(\lambda^+,\lambda^-)}\!\!\!&=\!
\begin{tikzpicture}[Q,centerzero,scale=.5]
	\draw[to-] (0.6,-.6)\botlabel{j^+} to (-0.6,.6);
	\draw[to-] (-0.6,-.6)\botlabel{i^+} to (0.6,.6);
  \region{1,0}{\lambda^\pm};
\end{tikzpicture}
\!\!+\!\!\begin{tikzpicture}[Q,centerzero,scale=.5]
	\draw[to-,thin] (0.6,-.6)\botlabel{j^-} to (-0.6,.6);
	\draw[to-,thin] (-0.6,-.6)\botlabel{i^-} to (0.6,.6);
  \region{1,0}{\lambda^\pm};
\end{tikzpicture}\!\!+\!(-1)^{\delta_{i,j}}r_{i,j}\!\begin{tikzpicture}[Q,centerzero,scale=.5]
	\draw[to-,thin] (0.6,-.6)\botlabel{j^+} to (-0.6,.6);
	\draw[to-,thin] (-0.6,-.6)\botlabel{i^-} to (0.6,.6);
\Babyclockwiseinternalbubble{-.3,.3};
\Babyanticlockwiseinternalbubble{.3,-.3};
  \region{1.1,0.1}{\lambda^\pm};
\end{tikzpicture}
\!\!+\!\!\!\begin{tikzpicture}[Q,centerzero,scale=.5]
	\draw[to-,thin] (0.6,-.6)\botlabel{j^-} to (-0.6,.6);
	\draw[to-,thin] (-0.6,-.6)\botlabel{i^+} to (0.6,.6);
  \region{1,0}{\lambda^\pm};
\end{tikzpicture}
\!\!+\!\delta_{i,j}\!
\begin{tikzpicture}[Q,centerzero,scale=.5]
	\draw[to-,thin] (-0.5,-.6)\botlabel{i^-} to (-0.5,.6);
	\draw[to-,thin] (.3,-.6)\botlabel{i^+} to (.3,.6);
\Teleporter{-.5,0}{.3,0};
  \region{1.1,0}{\lambda^\pm};
\end{tikzpicture}\!-\delta_{i,j}\!\begin{tikzpicture}[Q,centerzero,scale=.5]
	\draw[to-,thin] (-0.5,-.6)\botlabel{i^+} to (-0.5,.6);
	\draw[to-,thin] (.3,-.6)\botlabel{i^-} to (.3,.6);
\Teleporter{-.5,0}{.3,0};
  \region{1.1,0}{\lambda^\pm};
\end{tikzpicture}
\!,\\
\Xi\left(\ \begin{tikzpicture}[Q,centerzero,scale=.9]
\draw[-to] (-0.3,-0.3) \botlabel{i} -- (0.3,0.3);
\draw[to-] (0.3,-0.3) \botlabel{j} -- (-0.3,0.3);
\region{0.45,0}{\lambda};
\end{tikzpicture}\right)_{\!(\lambda^+,\lambda^-)}\!\!\!&=\!
\begin{tikzpicture}[Q,centerzero,scale=.5]
	\draw[to-] (0.6,-.6)\botlabel{j^+} to (-0.6,.6);
	\draw[-to] (-0.6,-.6)\botlabel{i^+} to (0.6,.6);
\Babyclockwiseinternalbubble{-.3,.3};
\Babyanticlockwiseinternalbubble{.3,-.3};
\region{1.1,0}{\lambda^\pm};
\end{tikzpicture}
\!\!+\!\! \begin{tikzpicture}[Q,centerzero,scale=.5]
	\draw[to-,thin] (0.6,-.6)\botlabel{j^-} to (-0.6,.6);
	\draw[-to,thin] (-0.6,-.6)\botlabel{i^-} to (0.6,.6);
  \region{1,0}{\lambda^\pm};
\end{tikzpicture}
\!\!+\!r_{j,i}^{-1}\!\begin{tikzpicture}[Q,centerzero,scale=.5]
	\draw[to-,thin] (0.6,-.6)\botlabel{j^+} to (-0.6,.6);
	\draw[-to,thin] (-0.6,-.6)\botlabel{i^-} to (0.6,.6);
  \region{1,0}{\lambda^\pm};
\end{tikzpicture}
\!\!+\!(-1)^{\delta_{i,j}}\!\!\!\begin{tikzpicture}[Q,centerzero,scale=.5]
	\draw[to-,thin] (0.6,-.6)\botlabel{j^-} to (-0.6,.6);
	\draw[-to,thin] (-0.6,-.6)\botlabel{i^+} to (0.6,.6);
  \region{1.1,0.1}{\lambda^\pm};
\end{tikzpicture} 
\!\!+\!\!\delta_{i,j}\!
\begin{tikzpicture}[Q,centerzero,scale=.5]
 	\draw[-to,thin] (-0.45,-.6)\botlabel{i^-} to[out=90,in=180] (-.1,-.1) to[out=0,in=90] (0.25,-.6);
 	\draw[-to,thin] (-0.45,.6) \toplabel{i^+}to[out=-90,in=180] (-.1,.1) to[out=0,in=-90] (0.25,.6);
\Teleporter{-.4,.3}{-.4,-.3};
     \Babyclockwiseinternalbubble{.15,.25};
  \region{1,0}{\lambda^\pm};
\end{tikzpicture}-\!\!\delta_{i,j}\!\begin{tikzpicture}[Q,centerzero,scale=.5]
 	\draw[-to] (-0.45,-.6)\botlabel{i^+} to[out=90,in=180] (-.1,-.1) to[out=0,in=90] (0.25,-.6);
 	\draw[-to] (-0.45,.6)\toplabel{i^-}to[out=-90,in=180] (-.1,.1) to[out=0,in=-90] (0.25,.6);
\Teleporter{.19,.3}{.19,-.3};
\Babyanticlockwiseinternalbubble{-.35,-.25};
  \region{1,0}{\lambda^\pm};
\end{tikzpicture}\!,\\
\Xi\left(\ \begin{tikzpicture}[Q,centerzero,scale=.9]
\draw[to-] (-0.3,-0.3) \botlabel{i} -- (0.3,0.3);
\draw[-to] (0.3,-0.3) \botlabel{j} -- (-0.3,0.3);
\region{0.45,0}{\lambda};
\end{tikzpicture}\right)_{\!(\lambda^+,\lambda^-)}\!\!\!&=\!
\begin{tikzpicture}[Q,centerzero,scale=.5]
	\draw[-to] (0.6,-.6)\botlabel{j^+} to (-0.6,.6);
	\draw[to-] (-0.6,-.6)\botlabel{i^+} to (0.6,.6);
\region{1,0}{\lambda^\pm};
\end{tikzpicture}
\!\!+\!\! \begin{tikzpicture}[Q,centerzero,scale=.5]
	\draw[-to,thin] (0.6,-.6)\botlabel{j^-} to (-0.6,.6);
	\draw[to-,thin] (-0.6,-.6)\botlabel{i^-} to (0.6,.6);
\Babyclockwiseinternalbubble{-.3,.3};
\Babyanticlockwiseinternalbubble{.3,-.3};
\region{1.1,0}{\lambda^\pm};
\end{tikzpicture}
\!\!+\!\!\begin{tikzpicture}[Q,centerzero,scale=.5]
	\draw[-to,thin] (0.6,-.6)\botlabel{j^+} to (-0.6,.6);
	\draw[to-,thin] (-0.6,-.6)\botlabel{i^-} to (0.6,.6);
  \region{1,0}{\lambda^\pm};
\end{tikzpicture}
\!\!+\!(-1)^{\delta_{i,j}} r_{i,j}\!\begin{tikzpicture}[Q,centerzero,scale=.5]
	\draw[-to,thin] (0.6,-.6)\botlabel{j^-} to (-0.6,.6);
	\draw[to-,thin] (-0.6,-.6)\botlabel{i^+} to (0.6,.6);
  \region{1.1,0.1}{\lambda^\pm};
\end{tikzpicture}
\!\!-\!\delta_{i,j}\!
\begin{tikzpicture}[Q,centerzero,scale=.5]
 	\draw[to-,thin] (-0.45,-.6)\botlabel{i^-} to[out=90,in=180] (-.1,-.1) to[out=0,in=90] (0.25,-.6);
 	\draw[to-,thin] (-0.45,.6) \toplabel{i^+}to[out=-90,in=180] (-.1,.1) to[out=0,in=-90] (0.25,.6);
\Teleporter{-.4,.3}{-.4,-.3};
     \Babyanticlockwiseinternalbubble{.15,-.25};
  \region{.95,0.1}{\lambda^\pm};
\end{tikzpicture}\!+\!\delta_{i,j}\!\begin{tikzpicture}[Q,centerzero,scale=.5]
 	\draw[to-] (-0.45,-.6)\botlabel{i^+} to[out=90,in=180] (-.1,-.1) to[out=0,in=90] (0.25,-.6);
 	\draw[to-] (-0.45,.6)\toplabel{i^-}to[out=-90,in=180] (-.1,.1) to[out=0,in=-90] (0.25,.6);
\Teleporter{.19,.3}{.19,-.3};
\Babyclockwiseinternalbubble{-.35,.25};
  \region{.95,0}{\lambda^\pm};
\end{tikzpicture}\!,\\
\Xi\left(\ \begin{tikzpicture}[Q,centerzero]
\draw[to-] (-0.25,-0.15) \botlabel{i} to [out=90,in=90,looseness=3](0.25,-0.15);
\region{0.45,0.1}{\lambda};
\end{tikzpicture}\!\right)_{(\lambda^+,\lambda^-)}\!\!\!&=
\begin{tikzpicture}[Q,anchorbase,scale=.7]
	\draw[to-,thin] (-0.4,-0.3)\botlabel{i^+} to[out=90, in=180] (0,0.3);
	\draw[-,thin] (-0,0.3) to[out = 0, in = 90] (0.4,-0.3);
  \region{.9,0}{\lambda^\pm};
\end{tikzpicture}+
\begin{tikzpicture}[Q,anchorbase,scale=.7]
	\draw[to-,thin] (-0.4,-0.3)\botlabel{i^-} to [out=90, in=180] (0,0.3);
	\draw[-,thin] (-0,0.3)to[out = 0, in = 90] (0.4,-0.3);
    \Babyanticlockwiseinternalbubble{.3,.1};
    \region{1,0}{\lambda^\pm};
\end{tikzpicture},\qquad\qquad\quad\!
\Xi\left(\ \begin{tikzpicture}[Q,centerzero]
\draw[to-] (-0.25,0.15) \toplabel{i} to [out=-90,in=-90,looseness=3](0.25,0.15);
\region{0.45,0}{\lambda};
\end{tikzpicture}\!\right)_{\!(\lambda^+,\lambda^-)}\!\!\! =
\begin{tikzpicture}[Q,anchorbase,scale=.7]
	\draw[to-,thin] (-0.4,0.3) \toplabel{i^+}to[out=-90, in=180] (0,-0.3);
	\draw[-,thin] (-0,-0.3) to[out = 0, in = -90] (0.4,0.3);
  \region{.9,0}{\lambda^\pm};
\end{tikzpicture}\!\!+\begin{tikzpicture}[Q,anchorbase,scale=.7] 
	\draw[to-,thin] (-0.4,0.3)\toplabel{i^-} to[out=-90, in=180] (0,-0.3);
	\draw[-,thin] (-0,-0.3) to[out = 0, in = -90] (0.4,0.3);
    \Babyclockwiseinternalbubble{-.3,-.1};
    \region{.9,0}{\lambda^\pm};
\end{tikzpicture}\!,\\
\Xi\left(\!\begin{tikzpicture}[baseline=-1mm,Q,scale=.9]
\draw[-to] (-0.25,0) arc(180:-180:0.25);
\node at (-0.42,0) {\strandlabel{i}};
\region{1,0}{\lambda};
\node at (.55,0) {$(u)$};
\end{tikzpicture}\!\right)_{\!(\lambda^+,\lambda^-)} \!\!\!
&=
\begin{tikzpicture}[Q,centerzero,scale=.9]
\draw[-to] (-0.68,0) arc(180:-180:0.25);
\node[black] at (0.13,0) {$(u)$};
\node at (-.43,-.4) {\strandlabel{i^+}};
\region{.92,0}{\lambda^\pm};
\end{tikzpicture} 
\begin{tikzpicture}[Q,centerzero,scale=.9]
\draw[-to] (-.25,0) arc(180:-180:0.25);
\node[black] at (.58,0) {$(u)$};
\node at (0,-0.4) {\strandlabel{i^-}};
\end{tikzpicture},\qquad\qquad\quad
\Xi\left(\!\begin{tikzpicture}[baseline=-1mm,Q,scale=.9]
\draw[to-] (-0.25,0) arc(180:-180:0.25);
\node at (-0.42,0) {\strandlabel{i}};
\region{1,0}{\lambda};
\node at (.55,0) {$(u)$};
\end{tikzpicture}\!\!\right)_{\!(\lambda^+,\lambda^-)}\!\!\!
=
\begin{tikzpicture}[Q,centerzero,scale=.9]
\draw[to-] (-0.68,0) arc(180:-180:0.25);
\node[black] at (0.13,0) {$(u)$};
\node at (-.43,-.4) {\strandlabel{i^+}};
\region{.92,0}{\lambda^\pm};
\end{tikzpicture} 
\begin{tikzpicture}[Q,centerzero,scale=.9]
\draw[to-] (-.25,0) arc(180:-180:0.25);
\node[black] at (.58,0) {$(u)$};
\node at (0,-0.4) {\strandlabel{i^-}};
\end{tikzpicture}.
\end{align*}
\end{theo}

\begin{proof}
We will deduce this from \cref{boxingday} by composing with some isomorphisms, paralleling the identity
$\Delta = (\omega \otimes \id) \circ d$ for the underlying quantum groups.
First, we need the isomorphism 
$$
D:\UU\stackrel{\sim}{\rightarrow} {\pmb\UU}^\imath
$$ 
from \cref{earphones}
from the 2-quantum group $\UU$
to the quasi-split 2-iquantum 
${\pmb\UU}^\imath$ of diagonal type.
Let ${\pmb\UU}$ be the 2-quantum group
associated to the doubled Cartan datum
with parameters as in \cref{this,that,theother,aargh,aaargh}.
\cref{boxingday} in this special case
gives us a strict graded 2-functor
$$
{\pmb\Xi}^\imath:{\pmb\UU}^\imath \rightarrow \widehat{\underline{\pmb\UU}}(\del,\zeta),
$$
where ${\pmb\UU}(\del,\zeta)$ is the version of
${\pmb\UU}$ with adapted grading from \cref{newgrading}, and $\del,\zeta$ are as in \cref{blippy,bloppy}.
Finally, we need the isomorphism of graded 2-categories
\begin{equation}\label{me}
\Omega\otimes J:{\pmb\UU}\rightarrow \UU\odot\UU
\end{equation}
defined on objects by $(\lambda^+,\lambda^-) \mapsto (-\lambda^+,\lambda^-)$, and on 1-morphisms 
by
\begin{align*}
E_{i^+} \one_{(\lambda^+,\lambda^-)}&\mapsto F_{i^+} \one_{(-\lambda^+,\lambda^-)},&
F_{i^+} \one_{(\lambda^+,\lambda^-)}&\mapsto E_{i^+} \one_{(-\lambda^+,\lambda^-)},\\
E_{i^-} \one_{(\lambda^+,\lambda^-)}&\mapsto E_{i^-} \one_{(-\lambda^+,\lambda^-)},&
F_{i^-} \one_{(\lambda^+,\lambda^-)}&\mapsto F_{i^-} \one_{(-\lambda^+,\lambda^-)}.
\end{align*}
On a string diagram representing a 2-morphism, the 2-functor $\Omega\otimes J$
applies the map 
$(\lambda^+,\lambda^-)\mapsto (-\lambda^+,\lambda^-)$ to 2-cell labels, 
reverses the orientation of all strings labelled by $I^+$,
and multiplies by
\begin{itemize}
\item
$-1$ for each
$\begin{tikzpicture}[Q,centerzero,scale=.7]
\draw[-to] (0,-0.3) \botlabel{i^-} -- (0,0.3);
\opendot{0,0};
\end{tikzpicture}$,
$\begin{tikzpicture}[Q,centerzero,scale=.7]
\draw[to-] (0,-0.3) \botlabel{i^-} -- (0,0.3);
\opendot{0,0};
\end{tikzpicture}$,
$\begin{tikzpicture}[Q,centerzero,scale=.8]
\draw[to-] (-0.25,-0.15) \botlabel{i^-} to [out=90,in=90,looseness=3](0.25,-0.15);
\end{tikzpicture}$ or
$\begin{tikzpicture}[Q,centerzero,scale=.8]
\draw[to-] (-0.25,0.15) \toplabel{i^-} to[out=-90,in=-90,looseness=3] (0.25,0.15);
\end{tikzpicture}\ \ (i \in I)$,
\item
$-r_{i,j}$ for each crossing of the form
$\begin{tikzpicture}[Q,centerzero,scale=.9]
\draw[to-] (-0.3,-0.3) \botlabel{i^+} -- (0.3,0.3);
\draw[to-] (0.3,-0.3) \botlabel{j^+} -- (-0.3,0.3);
\end{tikzpicture}$,
$\begin{tikzpicture}[Q,centerzero,scale=.9]
\draw[to-] (-0.3,-0.3) \botlabel{j^+} -- (0.3,0.3);
\draw[-to] (0.3,-0.3) \botlabel{i^+} -- (-0.3,0.3);
\end{tikzpicture}$,
$\begin{tikzpicture}[Q,centerzero,scale=.9]
\draw[-to] (-0.3,-0.3) \botlabel{i^-} -- (0.3,0.3);
\draw[-to] (0.3,-0.3) \botlabel{j^-} -- (-0.3,0.3);
\end{tikzpicture}$ or
$\begin{tikzpicture}[Q,centerzero,scale=.9]
\draw[-to] (-0.3,-0.3) \botlabel{j^-} -- (0.3,0.3);
\draw[to-] (0.3,-0.3) \botlabel{i^-} -- (-0.3,0.3);
\end{tikzpicture}\ (i,j \in I)$, and
\item
$-r_{i,j}^{-1}$ for each crossing of the form
$\begin{tikzpicture}[Q,centerzero,scale=.7]
\draw[-to] (-0.3,-0.3) \botlabel{i^+} -- (0.3,0.3);
\draw[to-] (0.3,-0.3) \botlabel{j^+} -- (-0.3,0.3);
\end{tikzpicture}
$,
$\begin{tikzpicture}[Q,centerzero,scale=.7]
\draw[-to] (-0.3,-0.3) \botlabel{j^+} -- (0.3,0.3);
\draw[-to] (0.3,-0.3) \botlabel{i^+} -- (-0.3,0.3);
\end{tikzpicture}$,
$\begin{tikzpicture}[Q,centerzero,scale=.7]
\draw[to-] (-0.3,-0.3) \botlabel{j^-} -- (0.3,0.3);
\draw[to-] (0.3,-0.3) \botlabel{i^-} -- (-0.3,0.3);
\end{tikzpicture}$ or
$\begin{tikzpicture}[Q,centerzero,scale=.7]
\draw[to-] (-0.3,-0.3) \botlabel{i^-} -- (0.3,0.3);
\draw[-to] (0.3,-0.3) \botlabel{j^-} -- (-0.3,0.3);
\end{tikzpicture}\ (i,j \in I)$
\end{itemize}
appearing in the original string diagram.
The fake bubble polynomials are transformed as follows:
\begin{align*}
\begin{tikzpicture}[Q,centerzero,scale=1.2]
\filledclockwisebubble{.9,0.05};	
\stringlabel{.9,-.28}{i^-};
\bubblelabel{.9,0.05}{u};
\end{tikzpicture}&\mapsto
\begin{tikzpicture}[Q,centerzero,scale=1.2]
\filledclockwisebubble{.9,0.05};	
\stringlabel{.9,-.28}{i^-};
\bubblelabel{.9,0.05}{-u};
\end{tikzpicture}\ ,&
\begin{tikzpicture}[Q,centerzero,scale=1.2]
\filledanticlockwisebubble{.9,0.05};	
\stringlabel{.9,-.28}{i^-};
\bubblelabel{.9,0.05}{u};
\end{tikzpicture}&\mapsto
\begin{tikzpicture}[Q,centerzero,scale=1.2]
\filledanticlockwisebubble{.9,0.05};	
\stringlabel{.9,-.28}{i^-};
\bubblelabel{.9,0.05}{-u};
\end{tikzpicture}\ ,&
\begin{tikzpicture}[Q,centerzero,scale=1.2]
\filledclockwisebubble{.9,0.05};	
\stringlabel{.9,-.28}{i^+};
\bubblelabel{.9,0.05}{u};
\end{tikzpicture}&\mapsto
\begin{tikzpicture}[Q,centerzero,scale=1.2]
\filledanticlockwisebubble{.9,0.05};	
\stringlabel{.9,-.28}{i^+};
\bubblelabel{.9,0.05}{u};
\end{tikzpicture}\ ,&
\begin{tikzpicture}[Q,centerzero,scale=1.2]
\filledanticlockwisebubble{.9,0.05};	
\stringlabel{.9,-.28}{i^+};
\bubblelabel{.9,0.05}{u};
\end{tikzpicture}&\mapsto
\begin{tikzpicture}[Q,centerzero,scale=1.2]
\filledclockwisebubble{.9,0.05};	
\stringlabel{.9,-.28}{i^+};
\bubblelabel{.9,0.05}{u};
\end{tikzpicture}\ .
\end{align*}

In fact, the 2-category ${\pmb\UU}$ can be viewed as combining copies of $\UU$ and of $\bar\UU$ which commute with each other; from this perspective, the automorphism $\Omega \otimes J$ agrees with  $\Omega$ from \cref{Omega}
on the copy of $\UU$, and it is the isomorphism $J$ from \cref{newyearseve} on the copy of $\bar\UU$.
The 2-category $\UU\circledast\UU$ is the image of
 ${\pmb\UU}(\del,\zeta)$ under 
the canonical extension
$\Omega\otimes J:{\pmb\UU}_q
\rightarrow (\UU\odot\UU)_q$ to the $q$-envelopes, with
the 1-morphisms \cref{starb1,starb2} being the images of
 $F_{i^+} \one_{(-\lambda^+,\lambda^-)}$,
$q_i^{h_i(\lambda^+)}E_{i^-} \one_{(-\lambda^+,\lambda^-)}$,
$q_i^{-h_i(\lambda^-)} E_{i^+} \one_{(-\lambda^+,\lambda^-)}$,
and $F_{i^-} \one_{(-\lambda^+,\lambda^-)}$, respectively.
This functor also extends further to additive envelopes.
Now everything is set up, and we can simply define
\begin{equation}
\Xi := (\Omega \otimes I) \circ {\pmb\Xi}^\imath \circ D:\UU\rightarrow 
\Add(\UU\:\underline{\widehat{\circledast}}\:\UU).
\end{equation}
It remains to calculate the images of the generating 2-morphisms to obtain the formulas in the statement of the theorem. In particular, the 
new internal bubbles \cref{newinternal1,newinternal2,newinternal3,newinternal4} are the
images under $\Omega\otimes J$ of the old internal bubbles
$\begin{tikzpicture}[Q,centerzero,scale=.7]
\draw[to-,thin] (0,-.4)\botlabel{i^+} to (0,.5);
\clockwiseinternalbubble{0,.05};
\end{tikzpicture}$,
$\begin{tikzpicture}[Q,centerzero,scale=.7]
\draw[to-,thin] (0,-.4)\botlabel{i^+} to (0,.5);
\anticlockwiseinternalbubble{0,.05};
\end{tikzpicture}$,
$\begin{tikzpicture}[Q,centerzero,scale=.7]
\draw[-to,thin] (0,-.4)\botlabel{i^+} to (0,.5);
\anticlockwiseinternalbubble{0,.05};
\end{tikzpicture}$, $\begin{tikzpicture}[Q,centerzero,scale=.7]
\draw[-to,thin] (0,-.4)\botlabel{i^+} to (0,.5);
\clockwiseinternalbubble{0,.05};
\end{tikzpicture}$,
$-\begin{tikzpicture}[Q,centerzero,scale=.7]
\draw[-to,thin] (0,-.4)\botlabel{i^-} to (0,.5);
\anticlockwiseinternalbubble{0,.05};
\end{tikzpicture}$, $-\begin{tikzpicture}[Q,centerzero,scale=.7]
\draw[-to,thin] (0,-.4)\botlabel{i^-} to (0,.5);
\clockwiseinternalbubble{0,.05};
\end{tikzpicture}$,
$-\begin{tikzpicture}[Q,centerzero,scale=.7]
\draw[to-,thin] (0,-.4)\botlabel{i^-} to (0,.5);
\clockwiseinternalbubble{0,.05};
\end{tikzpicture}$,
 $-\begin{tikzpicture}[Q,centerzero,scale=.7]
\draw[to-,thin] (0,-.4)\botlabel{i^-} to (0,.5);
\anticlockwiseinternalbubble{0,.05};
\end{tikzpicture}$, respectively; this depends on the observation made earlier that the internal bubbles can be drawn on the other side of the string due to \cref{slushy}.
\end{proof}

\begin{rem}\label{feeble}
The upward and downward crossings with internal bubbles in \cref{valentinesday}
can be written more simply using the observation that
\begin{align}
r_{j,i}^{-1}
\begin{tikzpicture}[Q,centerzero,scale=.6]
	\draw[-to,thin] (0.6,-.6)\botlabel{j^-} to (-0.6,.6);
	\draw[-to,thin] (-0.6,-.6)\botlabel{i^+} to (0.6,.6);
\Clockwiseinternalbubble{-.3,.3};
\Anticlockwiseinternalbubble{.3,-.3}
\end{tikzpicture}
&=
\begin{dcases}
\begin{tikzpicture}[Q,centerzero,scale=.6]
	\draw[-to,thin] (0.6,-.6)\botlabel{i^-} to (-0.6,.6);
	\draw[-to,thin] (-0.6,-.6)\botlabel{i^+} to (0.6,.6);
    \Teleporter{-.3,.3}{.3,.3};
    \Teleporter{-.3,-.3}{.3,-.3};
\end{tikzpicture}&\text{if $i = j$}\\
\begin{tikzpicture}[Q,centerzero,scale=.6]
	\draw[-to,thin] (0.6,-.6)\botlabel{j^-} to (-0.6,.6);
	\draw[-to,thin] (-0.6,-.6)\botlabel{i^+} to (0.6,.6);
    \Pinpin{.25,-.25}{-.25,-.25}{-2,-.25}{Q_{i,j}(x,y)};
    \end{tikzpicture}&\text{if $i\neq j$,}
\end{dcases}
&
r_{i,j} \begin{tikzpicture}[Q,centerzero,scale=.6]
	\draw[to-,thin] (0.6,-.6)\botlabel{j^+} to (-0.6,.6);
	\draw[to-,thin] (-0.6,-.6)\botlabel{i^-} to (0.6,.6);
\Clockwiseinternalbubble{-.3,.3};
\Anticlockwiseinternalbubble{.3,-.3}
\end{tikzpicture}
&=
\begin{dcases}
\begin{tikzpicture}[Q,centerzero,scale=.6]
	\draw[to-,thin] (0.6,-.6)\botlabel{i^+} to (-0.6,.6);
	\draw[to-,thin] (-0.6,-.6)\botlabel{i^-} to (0.6,.6);
    \Teleporter{-.3,.3}{.3,.3};
    \Teleporter{-.3,-.3}{.3,-.3};
\end{tikzpicture}&\text{if $i = j$}\\
\begin{tikzpicture}[Q,centerzero,scale=.6]
	\draw[to-,thin] (0.6,-.6)\botlabel{j^+} to (-0.6,.6);
	\draw[to-,thin] (-0.6,-.6)\botlabel{i^-} to (0.6,.6);
    \Pinpin{.25,.25}{-.25,.25}{-2,.25}{{'}Q_{i,j}(x,y)};
\end{tikzpicture}&\text{if $i\neq j$.}
\end{dcases}
\end{align}
This follows using \cref{monday} (actually, its image under $\Omega\otimes J$).
\end{rem}

\begin{rem}
In the $\mathfrak{sl}_2$ case, a very similar 2-functor to $\Xi$ was constructed in \cite[Th.~3.1]{unfurling}. Denoting the one in \cite{unfurling} by $\widetilde{\Xi}$ and using $+$ and $-$ instead of the red and blue strings there, it is defined by
\begin{align*}
\widetilde{\Xi}\left(\ \begin{tikzpicture}[Q,centerzero,scale=.9]
\draw[-to] (0,-0.3) -- (0,0.3);
\opendot{0,0};
\region{0.3,0}{\lambda};
\end{tikzpicture}\right)_{\!(\lambda^+,\lambda^-)}\!\!\!
&
:=
\begin{tikzpicture}[Q,centerzero,scale=.9]
\draw[-to] (0,-0.3) \botlabel{+} -- (0,0.3);
\opendot{0,0};
\region{0.5,0}{\lambda^\pm};
\end{tikzpicture}+
\begin{tikzpicture}[Q,centerzero,scale=.9]
\draw[-to] (0,-0.3) \botlabel{-} -- (0,0.3);
\opendot{0,0};
\region{0.5,0}{\lambda^\pm};
\end{tikzpicture}\ ,\\
\widetilde{\Xi}\left(\ \begin{tikzpicture}[Q,centerzero,scale=.9]
\draw[-to] (-0.3,-0.3) -- (0.3,0.3);
\draw[-to] (0.3,-0.3) -- (-0.3,0.3);
\region{0.45,0}{\lambda};
\end{tikzpicture}\right)_{\!(\lambda^+,\lambda^-)}\!\!\!&:=\!
\begin{tikzpicture}[Q,centerzero,scale=.5]
	\draw[-to] (0.6,-.6)\botlabel{+} to (-0.6,.6);
	\draw[-to] (-0.6,-.6)\botlabel{+} to (0.6,.6);
  \region{1,0}{\lambda^\pm};
\end{tikzpicture}
\!+\!\! \begin{tikzpicture}[Q,centerzero,scale=.5]
	\draw[-to,thin] (0.6,-.6)\botlabel{-} to (-0.6,.6);
	\draw[-to,thin] (-0.6,-.6)\botlabel{-} to (0.6,.6);
  \region{1,0}{\lambda^\pm};
\end{tikzpicture}
\!\!-\!\!\begin{tikzpicture}[Q,centerzero,scale=.5]
	\draw[-to,thin] (0.6,-.6)\botlabel{+} to (-0.6,.6);
	\draw[-to,thin] (-0.6,-.6)\botlabel{-} to (0.6,.6);
\Teleporter{-.3,-.3}{.3,-.3};
\region{1.1,0.1}{\lambda^\pm};
\end{tikzpicture}
\!\!+\!\!\begin{tikzpicture}[Q,centerzero,scale=.5]
	\draw[-to,thin] (0.6,-.6)\botlabel{-} to (-0.6,.6);
	\draw[-to,thin] (-0.6,-.6)\botlabel{+} to (0.6,.6);
\Teleporter{-.3,-.3}{.3,-.3};
  \region{1.1,0.1}{\lambda^\pm};
\end{tikzpicture}
\!\!-\!\!\begin{tikzpicture}[Q,centerzero,scale=.5]
	\draw[-to,thin] (-0.5,-.6)\botlabel{-} to (-0.5,.6);
	\draw[-to,thin] (.3,-.6)\botlabel{+} to (.3,.6);
\Teleporter{-.5,0}{.3,0};
  \region{1.1,0}{\lambda^\pm};
\end{tikzpicture}
\!\!+\!\!
\begin{tikzpicture}[Q,centerzero,scale=.5]
	\draw[-to,thin] (-0.5,-.6)\botlabel{+} to (-0.5,.6);
	\draw[-to,thin] (.3,-.6)\botlabel{-} to (.3,.6);
\Teleporter{-.5,0}{.3,0};
  \region{1.1,0}{\lambda^\pm};
\end{tikzpicture}\!,\\
\widetilde{\Xi}\left(\ \begin{tikzpicture}[Q,centerzero]
\draw[-to] (-0.25,0.15)  to [out=-90,in=-90,looseness=3](0.25,0.15);
\region{0.45,0}{\lambda};
\end{tikzpicture}\!\right)_{\!(\lambda^+,\lambda^-)}\!\!\! &:=
\begin{tikzpicture}[Q,anchorbase,scale=.7]
	\draw[-,thin] (-0.4,0.3) \toplabel{+}to[out=-90, in=180] (0,-0.3);
	\draw[-to,thin] (-0,-0.3) to[out = 0, in = -90] (0.4,0.3);
  \region{.9,0}{\lambda^\pm};
\end{tikzpicture}+\begin{tikzpicture}[Q,anchorbase,scale=.7] 
	\draw[-,thin] (-0.4,0.3)\toplabel{-} to[out=-90, in=180] (0,-0.3);
	\draw[-to,thin] (-0,-0.3) to[out = 0, in = -90] (0.4,0.3);
  \region{.9,0}{\lambda^\pm};
\end{tikzpicture},\qquad\qquad\quad
\widetilde{\Xi}\left(\ \begin{tikzpicture}[Q,centerzero]
\draw[-to] (-0.25,-0.15)  to [out=90,in=90,looseness=3](0.25,-0.15);
\region{0.45,0.1}{\lambda};
\end{tikzpicture}\!\right)_{\!(\lambda^+,\lambda^-)}\!\!\! :=
\begin{tikzpicture}[Q,anchorbase,scale=.7]
	\draw[-,thin] (-0.4,-0.3)\botlabel{+} to[out=90, in=180] (0,0.3);
	\draw[-to,thin] (-0,0.3) to[out = 0, in = 90] (0.4,-0.3);
  \region{.9,0}{\lambda^\pm};
\end{tikzpicture}+
\begin{tikzpicture}[Q,anchorbase,scale=.7]
	\draw[-,thin] (-0.4,-0.3)\botlabel{-} to [out=90, in=180] (0,0.3);
	\draw[-to,thin] (-0,0.3)to[out = 0, in = 90] (0.4,-0.3);
  \region{.9,0}{\lambda^\pm};
\end{tikzpicture}\ .\\\intertext{However, unlike $\Xi$, $\widetilde{\Xi}$ does not preserve degrees of 2-morphisms. Although rightward cups and caps are simpler, the leftward ones are more complicated:}
\Xi\left(\ \begin{tikzpicture}[Q,centerzero]
\draw[to-] (-0.25,-0.15) to [out=90,in=90,looseness=3](0.25,-0.15);
\region{0.45,0.1}{\lambda};
\end{tikzpicture}\!\right)_{(\lambda^+,\lambda^-)}\!\!\!&=
\begin{tikzpicture}[Q,anchorbase,scale=.7]
	\draw[to-,thin] (-0.4,-0.3)\botlabel{+} to[out=90, in=180] (0,0.3);
	\draw[-,thin] (-0,0.3) to[out = 0, in = 90] (0.4,-0.3);
    \Babyanticlockwiseinternalbubble{.3,.1};
    \region{1,0}{\lambda^\pm};
\end{tikzpicture}+
\begin{tikzpicture}[Q,anchorbase,scale=.7]
	\draw[to-,thin] (-0.4,-0.3)\botlabel{-} to [out=90, in=180] (0,0.3);
	\draw[-,thin] (-0,0.3)to[out = 0, in = 90] (0.4,-0.3);
    \Babyanticlockwiseinternalbubble{.3,.1};
    \region{1,0}{\lambda^\pm};
\end{tikzpicture},\qquad\qquad\quad
\Xi\left(\ \begin{tikzpicture}[Q,centerzero]
\draw[to-] (-0.25,0.15) to [out=-90,in=-90,looseness=3](0.25,0.15);
\region{0.45,0}{\lambda};
\end{tikzpicture}\!\right)_{\!(\lambda^+,\lambda^-)}\!\!\! =
\begin{tikzpicture}[Q,anchorbase,scale=.7]
	\draw[to-,thin] (-0.4,0.3) \toplabel{+}to[out=-90, in=180] (0,-0.3);
	\draw[-,thin] (-0,-0.3) to[out = 0, in = -90] (0.4,0.3);
    \Babyclockwiseinternalbubble{-.3,-.1};
    \region{.9,0}{\lambda^\pm};
\end{tikzpicture}+\begin{tikzpicture}[Q,anchorbase,scale=.7] 
	\draw[to-,thin] (-0.4,0.3)\toplabel{-} to[out=-90, in=180] (0,-0.3);
	\draw[-,thin] (-0,-0.3) to[out = 0, in = -90] (0.4,0.3);
    \Babyclockwiseinternalbubble{-.3,-.1};
    \region{.9,0}{\lambda^\pm};
\end{tikzpicture}.
\end{align*}
For general $\mathfrak{g}$, there is also another 2-functor of a
similar nature constructed in \cite[Lem.~3.5]{unfurling} whose codomain is
another 2-quantum group obtained by ``unfurling'' the Cartan datum.
\end{rem}

%% file: s5-nondegeneracy.tex
\setcounter{section}{4}

\section{\texorpdfstring{Non-degeneracy and the bilinear form on $\U^\imath(\del)$}{}}\label{s5-nondegeneracy}

For any symmetrizable Cartan
datum and any choice of parameters, the 2-quantum group $\UU$ is non-degenerate, that is,
the 2-morphism spaces
in $\UU$ have explicit diagrammatic bases.
This was conjectured originally in \cite{KL3}, where it was proved for
$\mathfrak{sl}_n$.
The conjecture was subsequently proved first 
for all finite types, and then in general in \cite[Th.~3.6]{unfurling}.
By \cite[Th.~2.7, Prop.~3.12]{KL3},
the validity of the Non-Degeneracy Conjecture implies that graded ranks of 2-morphism spaces in 
$\UU$ can be computed using 
Lusztig's symmetric bilinear form
\begin{equation}\label{lf}
(\cdot,\cdot): \dot\U_\Z \times \dot\U_\Z \rightarrow \Z[q,q^{-1}]
\end{equation}
from \cite[Th.~26.1.2]{Lubook} (we recall its definition in \cref{recallithere}).
In this section, we are going to use the 2-functor $\Xi^\imath$ from
\cref{boxingday} to deduce analogous results for
2-iquantum groups with geometric parameters.

\subsection{Notation for words and monomials}\label{stupidsubsection}

Let $\langle I \rangle$ be the set of words
$\bi = i_1 \cdots i_l$ for $l \geq 0$ and $i_1,\dots,i_l \in I$.
We denote the length $l$ of the word by $l(\bi)$ and let
$\wt^\imath(\bi)$ be the image of 
$\alpha_{i_1} + \cdots + \alpha_{i_l}$ in $X^\imath$.
Words of this form will be used to index various monomials:
\begin{align}\label{wordnotation1}
b_\bi &:= b_{i_1} \cdots b_{i_l},&
B_\bi &:= B_{i_1} \cdots B_{i_l},\\\label{wordnotation1b}
\theta_\bi &:= \theta_{i_1}\cdots \theta_{i_l},&
\Theta_\bi &:= \Theta_{i_1}\otimes\cdots\otimes\Theta_{i_l}.
\end{align}
The first of these is an element of $\U^\imath$.
The second has no meaning by itself but, given an iweight $\lambda \in X^\imath$, the notation
$B_\bi\one_\lambda$ denotes the 1-morphism in the 2-iquantum group $\UU^\imath$ that is the horizontal composition of the
generating 1-morphisms according to the sequence.
The third is a monomial in Lusztig's algebra $\f$, which will be introduced in \cref{sgraphical}.
The fourth denotes an object in a certain
quiver Hecke category which is TBD (``to be defined''---see \cref{iQHA}) but
it seems convenient to include 
the notation $\Theta_\bi$ here already.

More generally, we will eventually need {\em divided power words}
of the form $\bi = i_1^{(n_1)} \cdots i_l^{(n_l)}$
for $l \geq 0, i_1,\dots, i_l \in I$ and $n_1,\dots,n_l \geq 1$.
We denote the set of all such by $\llangle I \rrangle$.
The set $\langle I \rangle$ is a subset of $\llangle I \rrangle$ in the obvious way.
For such a divided power word, we
define $\wt^\imath(\bi)$ to be the image of $n_1 \alpha_{i_1}+\cdots+n_l \alpha_{i_l}$
in $X^\imath$, and let
\begin{align}\label{wordnotation2}
b_\bi &:= b_{i_1}^{(n_1)} \cdots b_{i_l}^{(n_l)},&
B_\bi &:= B_{i_1}^{(n_1)} \cdots B_{i_l}^{(n_l)},\\\label{wordnotation2b}
\theta_\bi &:= \theta_{i_1}^{(n_1)}\cdots \theta_{i_l}^{(n_l)},&
\Theta_\bi &:= \Theta_{i_1}^{(n_1)}\otimes\cdots\otimes\Theta_{i_l}^{(n_l)}.
\end{align}
Given $\bi \in \llangle I \rrangle$ and an iweight $\lambda$, the notation $b_\bi 1_\lambda$ denotes the corresponding product of divided/idivided powers in $\dot\U^\imath$
from \cref{idividedpowerrelation}. The other monomial introduced in  \cref{wordnotation2} 
involves $B_i^{(n)}$ which is TBD (see \cref{ssidempotents}). The monomials \cref{wordnotation2b}
involve
$\theta_i^{(n)} := \theta_i^n / [n]^!_{q_i}$
and $\Theta_i^{(n)}$ which is TBD (see \cref{ssgtb}).

There is also the set $\langle I^+, I^-\rangle$
of words in the alphabet ${\pmb I} = I^+\sqcup I^- =
\{i^+, i^-\:|\:i \in I\}$.
We define the weight 
$\wt(\bi) \in X$ of such a word to be the sum of the weights of its letters, with the understanding that 
$\wt(i^+) = \alpha_i$ and $\wt(i^-) = -\alpha_i$.
We let $g_{i^+} := e_i, g_{i^-} := f_i,
G_{i^+} := E_i$ and $G_{i^-} := F_i$, then for
$\bi = i_1 \cdots i_l \in \langle I^+, I^-\rangle$ we define
\begin{align}\label{wordnotation3}
g_{\bi} &:= g_{i_1} \cdots g_{i_l},&
G_{\bi} &:= G_{i_1} \cdots G_{i_l}.
\end{align}
The monomial $g_\bi$ is an element of $\U$, and $G_\bi \one_\lambda$ is a 1-morphism in $\UU$ for any $\lambda \in X$.
One could also introduce the set
$\llangle I^+, I^-\rrangle$ of divided power words
 which index monomials in divided powers of $e_i,f_i$ or $E_i,F_i$, but actually we will not need these subsequently.

\subsection{Non-degeneracy of 2-quantum groups}\label{nondeg}

We begin by giving a self-contained proof of the non-degeneracy of 2-quantum
groups for all symmetrizable Kac-Moody types. Although similar in spirit to the proof given in
\cite{unfurling}, we will prove it using the comultiplication 2-functor
$\Xi$ 
from \cref{valentinesday}
instead of the unfurling 2-functor from
\cite[Lem.~3.5]{unfurling}.
Let $\UU$ be as in \cref{existing}.
We will forget the grading on $\UU$, viewing it as a $\kk$-linear 2-category.

A {\em 2-representation} $\catR$ of $\UU$ is 
the data of a strict $\kk$-linear 2-functor
from $\UU$ to the 2-category of $\kk$-linear categories. It means that
we are given $\kk$-linear categories $\one_\lambda \catR$ for each $\lambda \in X$
and $\kk$-linear functors
$E_i: \one_\lambda\catR \rightarrow \one_{\lambda+\alpha_i}\catR$
and $F_i:\one_\lambda\catR\rightarrow \one_{\lambda-\alpha_i}\catR$
for each $\lambda \in X$ and $i \in I$.
Moreover,
every 2-morphism in $\UU$ induces a natural transformation between the appropriate 
compositions of these functors in a way that respects horizontal and vertical composition.
We abuse notation by writing $\catR$ also for the
$\kk$-linear category $\coprod_{\lambda \in X} \one_\lambda \catR$
or, if we are in a setting in which each $\one_\lambda \catR$ is additive, for the additive $\kk$-linear category $\bigoplus_{\lambda \in X} \one_\lambda \catR$. Then we refer to $\one_\lambda \catR$ as the {\em $\lambda$-weight subcategory}.

There is a natural notion of morphism $M:\catR \rightarrow \catS$ between two 
2-representations. It is the data of a $\kk$-linear functor which restricts to
$\kk$-linear functors $M_\lambda:\one_\lambda \catR \rightarrow \one_\lambda \catS$ for each $\lambda \in X$,
plus natural isomorphisms
$E_i^{\catS} \circ M_\lambda \cong M_{\lambda+\alpha_i} \circ E_i^{\catR}$ and $F_i^{\catS} \circ M_\lambda \cong M_{\lambda-\alpha_i} \circ F_i^{\catR}$ which are compatible with the natural transformations arising 
from 2-morphisms in $\UU$; see 
\cite[Def.~4.6]{BD} for more details.
We use the term {\em equivariant functor}
rather than ``morphism of 2-representations'', 
and call it an {\em equivariant equivalence} if each $M_\lambda$ is an equivalence of categories.

One can also take some commutative $\kk$-algebra $\KK$ and consider {\em $\KK$-linear 2-representations} and {\em $\KK$-linear equivariant 
functors} between them. These are just 2-representations and equivariant functors as defined above for the base change $\UU\otimes_\kk \KK$.
We say that $\catR$ is a {\em $\KK$-finite 2-representation} to indicate that the $\kk$-algebra 
$\KK$ is an algebraically closed field, $\catR$ is a $\KK$-linear 2-representation, and all morphism spaces of
$\catR$ are finite-dimensional as vector spaces over $\KK$.

Assume from now on that $\catR$ is a $\KK$-finite 2-representation.
There is a canonically induced structure of $\KK$-finite 2-representation on the additive Karoubi envelope ($=$the idempotent completion of the additive envelope) of 
$\catR$.
The Yoneda embedding induces
a contravariant equivariant between this and another 
$\KK$-finite 2-representation we denote by $\proj{\catR}$, whose $\lambda$-weight subcategory
is the category $\proj{\one_\lambda\catR}$ of finitely generated
projective left $\one_\lambda\catR$-modules. 
We will work here with the latter construction, the advantage 
being that $\proj{\one_\lambda\catR}$ is a full subcategory of the {\em Abelian} category $\mod{\one_\lambda\catR}$ of locally finite-dimensional left
$\one_\lambda\catR$-modules, that is, 
the category of $\KK$-linear functors 
from $\one_\lambda\catR$ to the category of
finite-dimensional vector spaces over $\KK$.
A disadvantage\footnote{As usual, one could avoid this by working with right modules, but that creates its own problems in a different place---one has to reverse the tensor product of bimodules.} is that 
the Yoneda equivalence being used is contravariant, so that an application of the anti-involution $\Psi$ from \cref{Psiinv} is needed in the construction of the required natural transformations.

One can think about $\mod{\one_\lambda\catR}$ in the more traditional language of
modules by passing to the path algebra of $\one_\lambda\catR$. This
is spelled out in \cite{BD}, where Abelian categories of this sort are called
{\em Schurian categories}. In a Schurian category, the
endomorphism algebra of any finitely generated object is
finite-dimensional, and the endomorphism algebra of an irreducible
object is one-dimensional. Also, it makes sense to talk about composition
multiplicities of any object, which are finite. However, objects (even finitely generated
ones) are not necessarily of finite length,
and morphism spaces between objects that are not finitely generated can be infinite-dimensional.

Let $V$ be a
finitely generated 
$\one_\lambda\catR$-module 
for $\lambda \in X$. Since
$\End_{\one_\lambda \catR}(V)$ is finite-dimensional, we can talk about
the (monic) minimal
polynomial $m_f(x) \in \KK[x]$
of any morphism $f:V \rightarrow V$.
By adjunction properties, both of $E_i V$ and $F_i V$ are finitely
generated, so it makes sense to define
$m_{V,i}(x)$ and $n_{V,i}(x)$ be the minimal polynomials of
\begin{align}
\begin{tikzpicture}[Q,centerzero,scale=.8]
\draw[-to] (0,-0.3) \botlabel{i} -- (0,0.3);
\opendot{0,0};
\draw[-,darkg,thick] (0.4,-.3) to (0.4,.3);
\catlabel{.65,0}{V};
\end{tikzpicture}&:E_i V \rightarrow E_i V,
&
\begin{tikzpicture}[Q,centerzero,scale=.8]
\draw[to-] (0,-0.3) \botlabel{i} -- (0,0.3);
\opendot{0,0};
\draw[-,darkg,thick] (0.4,-.3) to (0.4,.3);
\catlabel{.65,0}{V};
\end{tikzpicture}&:F_i V \rightarrow F_i V,
\end{align}
respectively.
The {\em spectrum} of $\catR$ is the set of roots
of the minimal polynomials $m_{V,i}(x)$ 
for all $i \in I$, all finitely generated $\one_\lambda \catR$-modules $V$ and all $\lambda \in X$.
Equivalently, by an argument with adjunctions, it is the set of roots of the minimal polynomials
$n_{V,i}(x)$ for all $i$ and $V$.

\begin{lem}\label{overtheoceanstill}
Let $V$ be a
finitely generated left $\one_\lambda\catR$-module and
$i,j \in I$. 
All roots of the minimal polynomials $m_{E_j V,i}(x), m_{F_j V, i}(x)$, $n_{E_j V,i}(x)$
and $n_{F_j V,i}(x)$
are roots of $m_{V,i}(x)$ or $n_{V,i}(x)$
or, when $j \neq i$, $Q_{i,j}(x,\beta)$
for roots $\beta$ of $m_{V,j}(x)$ or $n_{V,j}(x)$.
\end{lem}

\begin{proof}
First, we look at $m_{F_j V,i}(x)$.
From \cref{Altquadratic} (actually, the image of this relation under \cref{Omegainv}), we have that
$$
\begin{tikzpicture}[Q,centerzero,scale=1.2]
\draw[-to] (-0.15,-0.4) \botlabel{i} -- (-0.15,0.4);
\draw[to-] (0.15,-0.4) \botlabel{j} -- (0.15,0.4);
\Circled{-.15,-.1}{u};
\draw[-,darkg,thick] (.45,-.4) to (.45,.4);
\catlabel{.6,0}{V};
\Pin{-.15,.2}{-.9,.2}{n_{V,i}(x)};
\end{tikzpicture}
=
(-1)^{\delta_{i,j}} \begin{tikzpicture}[Q,centerzero,scale=1.2]
\draw[-to] (-0.2,-0.4) \botlabel{i} to[out=45,in=down] (0.15,0) to[out=up,in=-45] (-0.2,0.4);
\draw[to-] (0.2,-0.4) \botlabel{j} to[out=135,in=down] (-0.15,0) to[out=up,in=225] (0.2,0.4);
\Circled{.15,0}{u};
\draw[-,darkg,thick] (.5,-.4) to (.5,.4);
\catlabel{.65,0}{V};
\Pin{-.09,.3}{-.9,.3}{n_{V,i}(x)};
\end{tikzpicture}.
$$
Multiplying by $m_{V,i}(u)$ and taking the $u^{-1}$-coefficient using
\cref{trick} gives that
$$
\begin{tikzpicture}[Q,centerzero,scale=1.2]
\draw[-to] (-0.15,-0.4) \botlabel{i} -- (-0.15,0.4);
\draw[to-] (0.15,-0.4) \botlabel{j} -- (0.15,0.4);
\Pin{-.15,0}{-1.35,0}{m_{V,i}(x)n_{V,i}(x)};
\draw[-,darkg,thick] (.45,-.4) to (.45,.4);
\catlabel{.6,0}{V};
\end{tikzpicture}
=0.
$$
It follows that $m_{F_j V,i}(x)$ divides $m_{V,i}(x) n_{V,i}(x)$.
This shows that roots of $m_{F_j V,i}(x)$
are roots of $m_{V,i}(x)$ or $n_{V,i}(x)$.
A similar argument shows that roots of
$n_{E_j V,i}(x)$ are roots of $m_{V,i}(x)$ or $n_{V,i}(x)$.

Next suppose that $\alpha$ is a root of 
$m_{E_j V, i}(x)$ for $j \neq i$.
For some root $\beta$ of $m_{V,j}(x)$,
we can find a simultaneous eigenvector 
$v \in E_i E_j V$
for the commuting endomorphisms
$\begin{tikzpicture}[Q,centerzero,scale=.7]
\draw[-to] (-.4,-.3) \botlabel{i}--(-.4,.3);
 \draw[-to] (0,-0.3) \botlabel{j} -- (0,0.3);
\opendot{-.4,0};
\draw[-,darkg,thick] (0.4,-.3) to (0.4,.3);
\catlabel{.65,0}{V};
\end{tikzpicture}$
 and $\begin{tikzpicture}[Q,centerzero,scale=.7]
\draw[-to] (-.4,-.3) \botlabel{i}--(-.4,.3);
 \draw[-to] (0,-0.3) \botlabel{j} -- (0,0.3);
\opendot{0,0};
\draw[-,darkg,thick] (0.4,-.3) to (0.4,.3);
\catlabel{.65,0}{V};
\end{tikzpicture}$ with eigenvalues
$\alpha$ and $\beta$, respectively. 
We have that
$$
\begin{tikzpicture}[Q,centerzero,scale=1.2]
\draw[-to] (-0.15,-0.4) \botlabel{i} to(-0.15,0.4);
\draw[-to] (0.15,-0.4) \botlabel{j} to (0.15,0.4);
\Pin{-.15,.23}{-.95,.23}{m_{V,i}(x)};
\Pinpin{.15,-.23}{-.15,-.23}{-.9,-.22}{Q_{i,j}(x,y)};
\draw[-,darkg,thick] (.45,-.4) to (.45,.4);
\catlabel{.6,0}{V};
\end{tikzpicture}
=
\begin{tikzpicture}[Q,centerzero,scale=1.2]
\draw[-to] (-0.2,-0.4) \botlabel{i} to[out=45,in=down] (0.15,0) to[out=up,in=-45] (-0.2,0.4);
\draw[-to] (0.2,-0.4) \botlabel{j} to[out=135,in=down] (-0.15,0)
to[out=up,in=225] (0.2,0.4);
\Pin{-.09,.3}{-.95,.3}{m_{V,i}(x)};
\draw[-,darkg,thick] (.5,-.4) to (.5,.4);
\catlabel{.65,0}{V};
\end{tikzpicture}
=
\begin{tikzpicture}[Q,centerzero,scale=1.2]
\draw[-to] (-0.2,-0.4) \botlabel{i} to[out=45,in=down] (0.15,0) to[out=up,in=-45] (-0.2,0.4);
\draw[-to] (0.2,-0.4) \botlabel{j} to[out=135,in=down] (-0.15,0)
to[out=up,in=225] (0.2,0.4);
\Pin{.17,0}{.95,0}{m_{V,i}(x)};
\draw[-,darkg,thick] (1.65,-.4) to (1.65,.4);
\catlabel{1.8,0}{V};
\end{tikzpicture} = 0.
$$
Acting on $v$, we deduce that $m_{V,i}(\alpha) Q_{i,j}(\alpha,\beta) =
0$.
So $\alpha$ is a root of $m_{V,i}(x) Q_{i,j}(x,\beta)$. This shows for $j \neq i$ 
that roots of $m_{E_j V,i}(x)$ are roots of $m_{V,i}(x)$ or $Q_{i,j}(x,\beta)$ for roots $\beta$ of $m_{V,j}(x)$.
A similar argument shows that 
roots of $n_{F_j V,i}(x)$ are roots of $n_{V,i}(x)$ or $Q_{i,j}(x,\beta)$ for roots $\beta$ of $n_{V,j}(x)$, again assuming $j \neq i$.

Finally, we show that any root of $m_{E_i V,i}(x)$ is a root of $m_{V,i}(x)$. A similar argument shows that any root of $n_{F_i V,i}(x)$ is a root of $n_{V,i}(x)$ to finish the proof of the lemma.
Let $\alpha$ be a root of $m_{E_i V,i}(x)$. 
For some $\beta \in \kk$, we can find
a simultaneous eigenvector
$v \in E_i^2 V$ for 
 the commuting endomorphisms
$x := \begin{tikzpicture}[Q,centerzero,scale=.7]
\draw[-to] (-.4,-.3) \botlabel{i}--(-.4,.3);
 \draw[-to] (0,-0.3) \botlabel{i} -- (0,0.3);
\opendot{-.4,0};
\draw[-,darkg,thick] (0.4,-.3) to (0.4,.3);
\catlabel{.65,0}{V};
\end{tikzpicture}$
 and $y := \begin{tikzpicture}[Q,centerzero,scale=.7]
\draw[-to] (-.4,-.3) \botlabel{i}--(-.4,.3);
 \draw[-to] (0,-0.3) \botlabel{i} -- (0,0.3);
\opendot{0,0};
\draw[-,darkg,thick] (0.4,-.3) to (0.4,.3);
\catlabel{.65,0}{V};
\end{tikzpicture}$ with eigenvalues
$\alpha$ and $\beta$, respectively.
Also let $s := \begin{tikzpicture}[Q,centerzero,scale=.7]
\draw[-to] (-.4,-.3) \botlabel{i}--(0,.3);
 \draw[-to] (0,-0.3) \botlabel{i} -- (-.4,0.3);
\draw[-,darkg,thick] (0.4,-.3) to (0.4,.3);
\catlabel{.65,0}{V};
\end{tikzpicture}$.
We have that $y s v = s x v - v = \alpha sv - v$.
If $sv$ is a multiple of $v$ then $sv = 0$ since $s^2 = 0$, hence, the identity in the previous sentence gives that $v = 0$. This is not the case. 
So $v$ and $sv$ are linearly independent.
On the subspace of $E_i^2 V$ spanned by
$v$ and $sv$, the matrix of $y$ is $\left(\begin{smallmatrix}\beta&-1\\0&\alpha\end{smallmatrix}\right)$. This shows that $\alpha$ is an eigenvalue of $y$, that is, it is a root of $m_{V,i}(x)$.
\end{proof}

Consider the $\kk$-linear 2-category $\UU\circledast\UU$ from \cref{sscm} with weight lattice 
$X \oplus X$. We are ignoring the grading, so it is 
just the same as the 2-category $\UU\odot\UU$ generated by a commuting pair of copies of $\UU$.
Also recall the 2-categories $\UU\:\underline{\circledast}\:\UU$ and 
$\UU\:\widehat{\underline{\circledast}}\:\UU$.
If $\catR$ and $\catS$ are
$\KK$-finite 2-representations of $\UU$,
there is a $\KK$-finite 2-representation
$\catR \boxtimes \catS$ of $\UU \circledast \UU$
defined by letting
$\one_{(\lambda^+,\lambda^-)} (\catR \boxtimes \catS)$
be the $\KK$-linearized
Cartesian product
$(\one_{\lambda^+} \catR) \boxtimes (\one_{\lambda^-} \catS)$.
The functors $E_{i^+}$ and $F_{i^+}$ are
$E_i \boxtimes \id$ and $F_i \boxtimes \id$, while
$E_{j^-}$ and $F_{j^-}$ are $\id \boxtimes E_j$ and $\id \boxtimes F_j$, which commute strictly with $E_{i^+}$ and $F_{i^+}$. 
Crossings with one string labelled $i^+\:(i \in I)$ and the other labelled $j^-\:(j \in I)$ act as identity natural transformations.
The natural transformations defining the actions of the other generating 2-morphisms of $\UU\circledast\UU$ are the obvious ones from the actions of $\UU$ on $\catR$ and $\catS$.

Suppose in addition that $\catR$ and $\catS$ have disjoint spectra. Then the action of $\UU\circledast\UU$ on $\catR\boxtimes \catS$ extends to an action of the localization
$\UU\:\underline{\circledast}\:\UU$. Using the 2-functor $\Xi$ from \cref{valentinesday}, we obtain from this a $\KK$-linear 2-representation $\catR \otimes \catS$ of $\UU$ itself
with
$$
\one_\lambda (\catR\otimes \catS)
:= \bigoplus_{\lambda^++\lambda^-=\lambda}
(\one_{\lambda^+}\catR)\boxtimes (\one_{\lambda^-}\catS).
$$
It might not be $\KK$-finite.

For this construction to be useful, we need a supply of $\KK$-finite 2-representations.
For each $\lambda \in X$, there is a corresponding
{\em left regular 2-representation} $\UU \one_\lambda$
of $\UU$.
This is the 
2-representation with
$\one_\kappa \UU \one_\lambda := \HOM_{\UU}(\lambda,\kappa)$. 
The functors $E_i$ and $F_i$, and the natural transformations associated to 2-morphisms,
are defined by composing horizontally on the left.
The {\em generalized cyclotomic quotients} (GCQ for short) in the next definition first appeared in \cite[Prop.~5.6]{canonical}.
The {\em cyclotomic quotients}, which are a special case, have a much longer history.

\begin{defin}[Generalized cyclotomic quotients]\label{gcq}
Suppose that we are given a pair $\mu,\nu \in X^+$ of dominant weights, a commutative $\kk$-algebra $\KK$, and monic polynomials $\mu_i(x), \nu_i(x) \in \KK[x]$
for each $i \in I$
with $\deg(\mu_i(x)) = h_i(\mu)$ and $\deg(\nu_i(x)) = h_i(\nu)$. 
The {\em generalized cyclotomic quotient}
of $\UU$ associated to this data is the $\KK$-linear 
2-representation
\begin{equation*}
\catH(\mu|\nu) := (\UU \one_{\nu-\mu} \otimes_\kk \KK) / \catI(\mu|\nu)
\end{equation*}
where $\catI(\mu|\nu)$ is the $\KK$-linear sub-2-representation of $\UU \one_{\nu-\mu}\otimes_\kk \KK$ generated by 
\begin{equation*}
\left\{\begin{tikzpicture}[Q,centerzero]
\draw[-to] (0,-0.3) \botlabel{i} -- (0,0.3);
\Pin{0,0}{-.75,0}{\mu_i(x)};
\region{0.4,0}{\nu-\mu};
\end{tikzpicture},\ 
c_i(\nu-\mu)\ \begin{tikzpicture}[baseline=-1mm,Q]
\draw[to-] (-0.25,0) arc(180:-180:0.25);
\node at (0,-.38) {\strandlabel{i}};
\region{0.65,0}{\nu-\mu};
\dottybubblelabel{0,0}{n};
\end{tikzpicture}
-\left[\nu_i(u)/\mu_i(u)\right]_{u:h_i(\mu-\nu)-n} \id_{\one_{\nu-\mu}}\:\bigg|\:
i \in I, 0 < n \leq h_i(\nu)
\right\}.
\end{equation*}
Equivalently, by \cite[Lem.~4.14]{BD}, $\catI(\mu|\nu)$
is generated by 
\begin{equation*}
\left\{\begin{tikzpicture}[Q,centerzero]
\draw[to-] (0,-0.3) \botlabel{i} -- (0,0.3);
\Pin{0,0}{-.75,0}{\nu_i(x)};
\region{0.4,0}{\nu-\mu};
\end{tikzpicture},\ 
c_i(\nu-\mu)^{-1}\ \begin{tikzpicture}[baseline=-1mm,Q]
\draw[to-] (-0.25,0) arc(-180:180:0.25);
\node at (0,-.38) {\strandlabel{i}};
\region{0.65,0}{\nu-\mu};
\dottybubblelabel{0,0}{n};
\end{tikzpicture}
-\left[\mu_i(u)/\nu_i(u)\right]_{u:h_i(\nu-\mu)-n} \id_{\one_{\nu-\mu}}\:\bigg|\:
i \in I, 0 < n \leq h_i(\mu)
\right\},
\end{equation*}
and 
we have that
\begin{align}\label{hot}
\begin{tikzpicture}[baseline=-1mm,Q]
\draw[-to] (-0.25,0) arc(-180:180:0.25);
\node at (0,-.38) {\strandlabel{i}};
\region{1.25,0}{\nu-\mu};
\node at (.53,0) {$(u)$};
\end{tikzpicture}
&= 
\frac{\nu_i(u)}{\mu_i(u)} 
\id_{\one_{\nu-\mu}},
&
\begin{tikzpicture}[baseline=-1mm,Q]
\draw[to-] (-0.25,0) arc(-180:180:0.25);
\node at (0,-.38) {\strandlabel{i}};
\region{1.25,0}{\nu-\mu};
\node at (.53,0) {$(u)$};
\end{tikzpicture}
&= \frac{\mu_i(u)}{\nu_i(u)} 
\id_{\one_{\nu-\mu}}
\end{align}
in $\catH(\mu|\nu)$.
We denote the object $\one_{\nu-\mu}$ of $\one_{\nu-\mu} \catH(\mu|\nu)$ by $V(\mu|\nu)$.
It is a generating object in the sense that any other object can be obtained from $V(\mu|\nu)$ by applying a sequence of the functors $E_i, F_i\:(i \in I)$.
If either $\mu = 0$ or $\nu = 0$, we call $\catH(\mu|\nu)$ simply a {\em cyclotomic quotient}.
\end{defin}

We will soon need a theorem
proved in the graded setting in \cite[Cor.~3.20]{Web}
giving a Morita equivalent realization of cyclotomic quotients in terms of cyclotomic quiver Hecke algebras. Another proof of the result valid also when the grading is forgotten 
is given in \cite[Th.~4.25]{Rou2}; this depends on \cite[Th.~4.24]{Rou2}, for which Rouquier cites \cite{KK,Web} noting that the arguments from \cite{KK} are also valid in the ungraded setting. 
We state the theorem shortly, after some more definitions.

\begin{defin}[Quiver Hecke categories]\label{qhc}
The {\em quiver Hecke category} $\catQH$ is the strict graded monoidal category with generating objects 
$E_i\:(i \in I)$ and generating morphisms 
\begin{align*}
\begin{tikzpicture}[Q,centerzero,scale=.7]
\draw[-to] (0,-0.3) \botlabel{i} -- (0,0.3);
\opendot{0,0};
\end{tikzpicture}&:E_i \rightarrow E_i,&
\begin{tikzpicture}[Q,centerzero,scale=.7]
\draw[-to] (-0.3,-0.3) \botlabel{i} -- (0.3,0.3);
\draw[-to] (0.3,-0.3) \botlabel{j} -- (-0.3,0.3);
\end{tikzpicture}&: E_i \otimes E_j \rightarrow E_j \otimes E_i
\end{align*}
of degrees $2d_i$ and $-d_i a_{i,j}$, respectively,
subject to the relations \cref{dotslide,quadratic,braid}. 
For $\mu$ and $\mu_i(x) \in \KK[x]$ as in \cref{gcq},
let $\catI(\mu)$
be the left tensor ideal of $\catQH\otimes_\kk \KK$
generated by
$\Big\{\begin{tikzpicture}[Q,centerzero]
\draw[-to] (0,-0.2) \botlabel{i} -- (0,0.2);
\Pin{0,0}{-.75,0}{\mu_i(x)};
\end{tikzpicture}\:\Big|\:i \in I\Big\}$.
Then we pass to the quotient category
\begin{align*}
\catQH(\mu) &:= (\catQH \otimes_\kk \KK) / \catI(\mu).
\end{align*}
This $\KK$-linear category 
is a {\em cyclotomic quiver Hecke category}. If $\KK$ is a graded $\kk$-algebra and each $\mu_i(x)$ is a homogeneous polynomial then $\catQH(\mu)$ is also a graded category.
\end{defin}

For $l \geq 0$, the locally unital endomorphism algebra
\begin{equation}\label{qha}
\QH_l :=
\bigoplus_{\substack{i_1,\dots,i_l \in I\\j_1,\dots,j_l \in I}}
\Hom_{\catQH}(E_{j_1}\otimes\cdots\otimes E_{j_l},
E_{i_1}\otimes\cdots\otimes E_{i_l})
\end{equation}
is the {\em quiver Hecke algebra}
introduced in \cite{Rou, KL1}.
If $I$ is finite then $\QH_l$ is a unital algebra, but in general it is merely locally unital.
It has a well-known basis as a free $\kk$-module
whose elements are the composites of a monomial
$x_1^{n_1}\cdots x_l^{n_l}\:(n_1,\dots,n_l \geq 0)$,
a string diagram representing a reduced expression for a permutation $w \in S_l$, and an idempotent indexed by
a sequence
$(i_1,\dots,i_l) \in I^l$.
Since there are no non-zero morphisms 
$E_{j_1}\otimes\cdots\otimes E_{j_l} \rightarrow
E_{i_1}\otimes\cdots\otimes E_{i_{l'}}$ for $l \neq l'$, we obtain from this a basis for each morphism space in $\catQH$.

Switching attention to the quotient category $\catQH(\mu)$, the locally unital endomorphism algebra
\begin{equation}\label{cqha}
\QH_l(\mu) :=
\bigoplus_{\substack{i_1,\dots,i_l\in I\\j_1,\dots,j_l\in I}}
\Hom_{\catQH(\mu)}(E_{j_1}\otimes\cdots\otimes E_{j_l},
E_{i_1}\otimes\cdots\otimes E_{i_l})
\end{equation}
is a {\em cyclotomic quiver Hecke algebra}.
By \cite[Cor.~3.26]{Web}, it is free of finite rank as a $\KK$-module; in particular, it is a unital algebra.
It follows that each morphism space in $\catQH(\mu)$ is free of finite rank as a $\KK$-module.
There is a $\KK$-linear functor
\begin{align}\label{obviousone}
F:\catQH(\mu) &\rightarrow \catH(\mu|\nu)
\end{align}
defined by acting on the generating object $V(\mu|\nu)$.

\begin{theo}[Kang-Kashiwara, Rouquier, Webster]\label{cyclotomic}
Assuming that $\nu = 0$, 
the functor $F$  
is fully faithful, and it
induces an equivalence 
$\proj{\catQH(\mu)}\stackrel{\sim}{\rightarrow}\proj{\catH(\mu|0)}$.
\end{theo}

The next lemma is
a partial extension of \cref{cyclotomic} to GCQs.
This can also be deduced from results in \cite{canonical,unfurling}, but the point is to give an independent proof of it using \cref{valentinesday}.
See also \cref{best} for a stronger version.

\begin{lem}\label{better}
In the setup of \cref{gcq}, suppose that the ground ring $\kk$ is an integral domain and that $\KK$ is an algebraically closed field.
Assume that the polynomials $\mu_i(x),\nu_i(x)$ for $i \in I$ are chosen so that all roots of the polynomials $\mu_i(x)$ (resp., $\nu_i(x)$) are algebraic (resp., transcendental) over the field of fractions of $\kk$.
Then the functor $F:
\catQH(\mu) \rightarrow \catH(\mu|\nu)$ 
from \cref{obviousone}
is also fully faithful when the dominant weight $\nu$ is non-zero.
\end{lem}

\begin{proof}
We claim that the spectrum of 
$\catH(\mu|0)$ consists of elements of $\KK$ that are algebraic over $\Frac(\kk)$, and the spectrum of $\catH(0|\nu)$ consists of elements of $\KK$
that are transcendental over $\Frac(\kk)$.
To justify this, the minimal polynomials
$m_{V(\mu|0),i}(x)$ divide $\mu_i(x)$ and $n_{V(\mu|0),i}(x) = 1$, so the roots of 
these polynomials are algebraic over $\Frac(\kk)$ by assumption. The minimal polynomials
$n_{V(0|\nu),i}(x)$ divide $\nu_i(x)$ and
$m_{V(0|\nu),i}(x) = 1$, so the roots of these polynomials are transcendental over $\Frac(\kk)$ by assumption.
The other objects of $\catH(\mu|0)$ and $\catH(0|\nu)$ are of the form $G_\bi V(\mu|0)$ or $G_\bi V(0|\nu)$
for some word $\bi \in \langle I^+, I^- \rangle$.
For $\alpha,\beta \in \KK$ with $Q_{i,j}(\alpha,\beta) = 0$ for some $i \neq j$, 
$\alpha$ is algebraic over $\Frac(\kk)$ if and only if $\beta$ is algebraic over $\Frac(\kk)$.
Using this, \cref{overtheoceanstill} and induction on the length of $\bi$, it follows that the roots of $m_{G_\bi V(\mu|0),i}(x)$ and $n_{G_\bi V(\mu|0),i}(x)$
are algebraic over $\Frac(\kk)$, and the roots of 
$m_{G_\bi V(0|\nu),i}(x)$ and $n_{G_\bi V(0|\nu),i}(x)$ are transcendental over $\Frac(\kk)$.
The claim follows.

The claim implies that the spectra of $\catH(\mu|0)$ and $\catH(0|\nu)$ are disjoint. So we obtain a
$\KK$-linear 2-representation
$\catH(\mu|0) \otimes \catH(0|\nu)$
of $\UU$ via the construction with
\cref{valentinesday}
explained before \cref{gcq}.
Denoting its object $(V(\mu|0),V(0|\nu))$ 
simply by $V$,
we have that
\begin{align*}
\begin{tikzpicture}[Q,centerzero]
\draw[-to] (0,-0.3) \botlabel{i^+} -- (0,0.3);
\Pin{0,0}{-.75,0}{\mu_i(x)};
\draw[-,darkg,thick] (0.4,-.3) to (0.4,.3);
\catlabel{.6,0}{V};\end{tikzpicture}
=
\begin{tikzpicture}[Q,centerzero]
\draw[-to] (0,-0.3) \botlabel{i^-} -- (0,0.3);
\draw[-,darkg,thick] (0.4,-.3) to (0.4,.3);
\catlabel{.6,0}{V};
\end{tikzpicture}&=0,&
\begin{tikzpicture}[Q,centerzero]
\draw[to-] (0,-0.3) \botlabel{i^+} -- (0,0.3);
\draw[-,darkg,thick] (0.4,-.3) to (0.4,.3);
\catlabel{.6,0}{V};
\end{tikzpicture}=
\begin{tikzpicture}[Q,centerzero]
\draw[to-] (0,-0.3) \botlabel{i^-} -- (0,0.3);
\Pin{0,0}{-.75,0}{\nu_i(x)};
\draw[-,darkg,thick] (0.4,-.3) to (0.4,.3);
\catlabel{.6,0}{V};
\end{tikzpicture}
&=0,\\
\begin{tikzpicture}[baseline=-1mm,Q]
\draw[-to] (-0.25,0) arc(-180:180:0.25);
\node at (0,-.38) {\strandlabel{i^+}};
\node at (.53,0) {$(u)$};
\draw[-,darkg,thick] (1,-.3) to (1,.3);
\catlabel{1.2,0}{V};
\end{tikzpicture}
&= 
\frac{1}{\mu_i(u)} 
\id_{V},&
\begin{tikzpicture}[baseline=-1mm,Q]
\draw[-to] (-0.25,0) arc(-180:180:0.25);
\node at (0,-.38) {\strandlabel{i^-}};
\draw[-,darkg,thick] (1,-.3) to (1,.3);
\catlabel{1.2,0}{V};
\node at (.53,0) {$(u)$};
\end{tikzpicture}
&=\nu_i(u)
\id_{V}.
\end{align*}
The faithfulness of the 
 functor in \cref{cyclotomic} implies that there is a faithful functor
$$
\widetilde{F}:\catQH(\mu) \rightarrow \catH(\mu|0) \otimes \catH(0|\nu)
$$
defined by acting on $V$.
There is a canonical equivariant functor $C:\UU \one_{\nu-\mu} \rightarrow
\catH(\mu|0) \otimes \catH(0|\nu)$ taking $\one_{\nu-\mu}$ to $V$. The image of
$\begin{tikzpicture}[Q,centerzero]
\draw[-to] (0,-0.3) \botlabel{i} -- (0,0.3);
\Pin{0,0}{-.75,0}{\mu_i(x)};
\region{0.45,0}{\nu-\mu};
\end{tikzpicture}$
is
$\begin{tikzpicture}[Q,centerzero]
\draw[-to] (0,-0.3) \botlabel{i^+} -- (0,0.3);
\Pin{0,0}{-.75,0}{\mu_i(x)};
\draw[-,darkg,thick] (.3,-.3) to (.3,.3);
\catlabel{.5,0}{V};
\end{tikzpicture}+\begin{tikzpicture}[Q,centerzero]
\draw[-to] (0,-0.3) \botlabel{i^-} -- (0,0.3);
\Pin{0,0}{-.75,0}{\mu_i(x)};
\draw[-,darkg,thick] (.3,-.3) to (.3,.3);
\catlabel{.5,0}{V};
\end{tikzpicture} = 0$,
and the image of 
$\begin{tikzpicture}[baseline=-1mm,Q]
\draw[-to] (-0.25,0) arc(-180:180:0.25);
\node at (0,-.38) {\strandlabel{i}};
\region{1.15,0}{\nu-\mu};
\node at (.53,0) {$(u)$};
\end{tikzpicture}$
is 
$\begin{tikzpicture}[baseline=-1mm,Q]
\draw[-to] (-0.25,0) arc(-180:180:0.25);
\node at (0,-.38) {\strandlabel{i^+}};
\node at (.53,0) {$(u)$};
\end{tikzpicture}
\begin{tikzpicture}[baseline=-1mm,Q]
\draw[-to] (-0.3,0) arc(-180:180:0.25);
\node at (0,-.38) {\strandlabel{i^-}};
\node at (.51,0) {$(u)$};
\draw[-,darkg,thick] (.9,-.3) to (.9,.3);
\catlabel{1.1,0}{V};
\end{tikzpicture}=\nu_i(u) / \mu_i(u)\id_V$.
So $C$ factors through the quotient to induce  
an equivariant functor $\bar C:\catH(\mu|\nu)
\rightarrow \catH(\mu|0) \otimes \catH(0|\nu)$.
Finally we observe that the following diagram of functors commutes:
$$
\begin{tikzcd}
\catH(\mu|\nu)\arrow[rr,"\bar C"]&&\catH(\mu|0) \otimes \catH(0|\nu)\\
&\catQH(\mu)\arrow[ul,"F"]\arrow[ur,"\widetilde{F}" below right]
\end{tikzcd}
$$
We deduce that $F$ is faithful because $\widetilde{F}$ is faithful. It is obvious that $F$ is full.
\end{proof}

Now the tools are in place, and we can state and prove Khovanov-Lauda's
Non-Degeneracy Conjecture for 2-quantum groups.
Using the conventions from
\cref{stupidsubsection}, for $\bi, \bj \in \langle I^+,I^-\rangle$ with $\wt(\bi) = \wt(\bj)$, 
$G_\bi \one_\lambda$ and $G_\bj \one_\lambda$ are 1-morphisms belonging the same morphism category of $\UU$.
By an $\bi \times \bj$ 
{\em shape}, we mean an oriented string diagram $\vec{D}$ 
with no dots or closed components (``bubbles'')
such that $\vec{D}\ \text{\regionlabel{\lambda}}$
represents a 2-morphism $G_\bj \one_\lambda \Rightarrow G_\bi\one_\lambda$ in $\UU$. Such a shape $\vec{D}$ is {\em reduced} if no string has more than one critical point (i.e., no zig-zags), no pair of strings cross each other twice (i.e., no bigons), and there are no self-intersections of a string with itself (i.e., no loops).
Note in a reduced shape that the boundary points of any cup are on the top edge and the boundary points of any cap are on the bottom edge; there are also {\em propagating strings} which have one boundary point on the top edge and one on the bottom edge.
We say that two shapes are {\em equivalent} if they define the same matching between the boundary points.
Then we fix a set $\oshape(\bi\times\bj)$ of representatives for the equivalence classes of reduced $\bi \times \bj$ shapes.
Assuming this set is non-empty, 
we have that $l(\bi) \equiv l(\bj) \pmod{2}$, so it makes sense to consider $l := (|l(\bi)|+|l(\bj)|)/2$; this is the number of connected components in any $\bi \times \bj$
shape.
Then, for each $\vec{D} \in \oshape(\bi\times\bj)$, we fix a choice of points $p_1,\dots,p_l$ away from crossings and critical points, one on each connected component of the shape, ordering them according to the lexicographic ordering of their Cartesian coordinates.
For $f \in \kk[x_1,\dots,x_l]$, we let
$\vec{D}(f)\ \text{\regionlabel{\lambda}}:G_\bj\one_\lambda\Rightarrow G_\bi \one_\lambda$
be the 2-morphism obtained from $\vec{D}$ by pinning  $f$ to the points $p_1,\dots,p_l$.

\begin{theo}\label{nd}
For $\lambda \in X$, the $\kk$-algebra $\End_{\UU}(\one_\lambda)$ is freely generated by either of the sets
\begin{equation}\label{bar3}
\bigg\{
\begin{tikzpicture}[baseline=-1mm,Q]
\draw[to-] (-0.25,0) arc(180:-180:0.25);
\node at (0,-.36) {\strandlabel{i}};
\region{0.5,0}{\lambda};
\dottybubblelabel{0,0}{n};
\end{tikzpicture}\:\bigg|\: i \in I, n \geq 1\bigg\}
\quad\text{or}\quad
\bigg\{
\begin{tikzpicture}[baseline=-1mm,Q]
\draw[-to] (-0.25,0) arc(180:-180:0.25);
\node at (0,-.36) {\strandlabel{i}};
\region{0.5,0}{\lambda};
\dottybubblelabel{0,0}{n};
\end{tikzpicture}\:\bigg|\: i \in I, n \geq 1\bigg\}.
\end{equation}
Moreover, for words $\bi,\bj \in \langle I^+,I^- \rangle$
with $\wt(\bi)=\wt(\bj), l(\bi) \equiv l(\bj) \pmod{2}$
and $l := (l(\bi)+l(\bj))/2$,
the 2-morphism space $\Hom_{\UU}(G_\bj\one_\lambda, G_\bi\one_\lambda)$ is free as a right $\End_\UU(\one_\lambda)$-module with basis given by the 
morphisms 
\begin{equation}\label{basis}
\left\{\vec{D}(x_1^{n_1} \cdots x_l^{n_l})\ \text{\regionlabel{\lambda}}\:\Big|\:\vec{D} \in \oshape(\bi\times\bj),
n_1,\dots,n_l \geq 0
\right\}.
\end{equation}
\end{theo}

\begin{proof}
We assume for the proof that the normalization functions
$c_i$ are all 1. This just amounts to rescaling diagrams by units (recall \cref{snow}), so it is permitted.
Note also that $\End_\UU(\one_\lambda)$
is freely generated by the first set in \cref{bar3} if and only if it is freely generated by the second set.
This follows from the infinite Grassmannian relation \cref{infgrass}. We will prove the result for the first set.

Let $\R := \kk\big[\beta_i^{(n)}\:|\:i \in I, n \geq 1\big]$
and $f_\lambda:\R \rightarrow \End_\UU(\one_\lambda)$
be the homomorphism mapping $$
\beta_i^{(n)}\mapsto 
\begin{tikzpicture}[baseline=-1mm,Q]
\draw[to-] (-0.25,0) arc(180:-180:0.25);
\node at (0,-.36) {\strandlabel{i}};
\region{0.5,0}{\lambda};
\dottybubblelabel{0,0}{n};
\end{tikzpicture}.
$$
We view
$\Hom_{\UU}(G_\bj\one_\lambda, G_\bi \one_\lambda)$ as a right $\R$-module with action of $\beta \in \R$
defined by horizontally composing with $f_\lambda(\beta)$ on the right of a string diagram.
We need to show that the set \cref{basis} generates 
$\Hom_{\UU}(G_\bj\one_\lambda, G_\bi\one_\lambda)$
as a right $\R$-module, and that this set is $\R$-linearly independent. The spanning part follows because there is
a straightening algorithm to take any string diagram representing a 2-morphism $G_\bj\one_\lambda \Rightarrow G_\bi\one_\lambda$ to an $\R$-linear combination of elements of \cref{basis}. This is explained in \cite[Sec.~3.2.3]{KL3} and \cite[Sec.~8]{Lauda}; the algorithm proceeds by induction on the number of crossings in the diagram. 

It remains to prove that the set \cref{basis} is $\R$-linearly independent.
For this, we 
make two more reductions:
\begin{itemize}
\item
We can assume that the ground ring $\kk$ is the quotient of
the polynomial algebra
$$
\Z[t_{i,j}, t_{i,j;r,s}\:|\:i \neq j\text{ in }I, 0 \leq r < -a_{i,j}, 0 \leq s < -a_{j,i}]
$$
by relations asserting 
that each $t_{i,j}$ is invertible
and $t_{i,j;r,s} = t_{j,i;s,r}$,
taking $Q_{i,j}(x,y)$ to be the generic polynomial defined by \cref{generic}.
The grading on $\kk$ is defined so that $t_{i,j}$ is of degree 0 and $t_{i,j;r,s}$ is of degree $-2 d_i a_{i,j} - 2d_i r - 2d_j s$.
The 2-category $\UU$ over any other ground ring 
for any choice of parameters can be obtained from this by base change. So if we can prove the linear independence in this special case, it is true in general.
\item
As noted in \cite[Rem.~3.16]{KL3}, 
the proof of linear independence reduces further to the case
that $\bi = i_1^+ \cdots i_l^+$  and $\bj = j_1^+\cdots i_l^+$ 
for $l \geq 0$ and 
$(i_1,\dots,i_l), (j_1,\dots,j_l) \in I^l$ in the same $S_l$-orbit.
So now we have that $G_\bi = E_{i_1} \cdots E_{i_l}$
and $G_\bj = E_{j_1} \cdots E_{j_l}$.

\end{itemize}
Having made these reductions, consider a linear 
relation
\begin{equation}\label{quebec}
\sum_{\substack{\vec{D} \in \oshape(\bi\times\bj)\\n_1,\dots,n_l \geq 0}}
\vec{D}(x_1^{n_1} \cdots x_l^{n_l}) \ \text{\regionlabel{\lambda}}\ 
f_\lambda(\beta(\vec{D};n_1,\dots,n_l))=0
\end{equation}
in $\Hom_\UU(G_\bj\one_\lambda, G_\bi\one_\lambda)$
for $\beta(\vec{D};n_1,\dots,n_l) \in \R$.
Now we pick sufficiently large dominant weights $\mu,\nu \in X^+$
so that
\begin{enumerate}
\item  $\lambda = \nu-\mu$;
\item
every $\beta(\vec{D};n_1,\dots,n_l)$ lies in the subalgebra 
$\kk\big[\beta_i^{(n)}\:|\:i \in I, 1\leq n \leq h_i(\nu)\big]$ of $\R$;
\item
$\beta(\vec{D};n_1,\dots,n_l) \neq 0\ \Rightarrow$
$\deg\big(\vec{D}(x_1^{n_1} \cdots x_l^{n_l})\big) < \min(d_{i_1} h_{i_1}(\mu),\dots,d_{i_l} h_{i_l})$.
\end{enumerate}
Consider the GCQ $\catH(\mu|\nu)$ of $\UU$ 
defined over the algebraic closure $\KK$ of the field extension of $\Frac(\kk)$ obtained by
adjoining indeterminates
$\{z_{i,n}\:|\:i \in I, 1 \leq n \leq h_i(\nu)\}$,
taking the monic polynomials
$\mu_i(x) := x^{h_i(\mu)}$ and
$\nu_i(x) := 
(x-z_{i,1})
\cdots (x-z_{i,h_i(\nu)})$.
Acting with \cref{quebec} on the generating object, we deduce using \cref{hot} that
\begin{equation}\label{manitoba}
\sum_{\substack{\vec{D} \in \oshape(\bi\times\bj)\\n_1,\dots,n_l \geq 0}}
f(\beta(\vec{D};n_1,\dots,n_l))
\vec{D}(x_1^{n_1} \cdots x_l^{n_l}) \ 
\begin{tikzpicture}[Q,centerzero]
\draw[-,darkg,thick] (0.4,-.3) to (0.4,.3);
\catlabel{.9,0}{V(\mu|\nu)};
\end{tikzpicture}
=0
\end{equation}
where $f:\R\rightarrow \KK$ is the $\kk$-algebra homomorphism mapping $\beta_i^{(n)}$ to the $x^{h_i(\nu)-n}$-coefficient of $\nu_i(x)$ if $1 \leq n \leq h_i(\nu)$ and to $0$ if $n > h_i(\nu)$. Up to some signs, $f\big(\beta_i^{(1)}\big),
\dots, f\big(\beta_i^{(h_i(\nu))}\big)$ are elementary symmetric polynomials in $z_{i,1},\dots,z_{i,h_i(\nu)}$.
So the elements $f\big(\beta_i^{(n)}\big)$ for $i \in I$ and $1 \leq n \leq h_i(\nu)$ are algebraically independent over $\Frac(\kk)$.

The objects $G_\bi V(\mu|\nu)$ and $G_\bj V(\mu|\nu)$ are in the image of the functor
$F$ from \cref{better}.
As $F$ is
fully faithful by that lemma, we can rewrite \cref{manitoba} 
equivalently as
\begin{equation}\label{alberta}
\sum_{\substack{\vec{D} \in \oshape(\bi\times\bj)\\n_1,\dots,n_l \geq 0}}
f(\beta(\vec{D};n_1,\dots,n_l))
\vec{D}(x_1^{n_1} \cdots x_l^{n_l}) \ 
=0
\end{equation}
in the morphism space $\Hom_{\catQH(\mu)}(E_{j_1}\otimes\cdots\otimes E_{j_l},E_{i_1}\otimes\cdots\otimes E_{i_l})$ of the cyclotomic quiver Hecke category. 
This morphism space
is naturally identified with a 
morphism space in the $\KK$-linear category
$(\catQH^0 / \catJ) \otimes_\kk \KK$,
where $\catQH^0$ is the quiver Hecke category that
is the full monoidal 
subcategory of $\catQH$ generated by the morphisms $E_{i_1},\dots,E_{i_l}$ (cf. \cite[Cor.~3.8]{Rou}), and $\catJ$ is the 
left tensor ideal of $\catQH^0$ generated
by the homogeneous morphisms
$\begin{tikzpicture}[Q,centerzero]
\draw[-to] (0,-0.2) \botlabel{i_1} -- (0,0.2);
\multopendot{0,0}{west}{h_{i_1}(\mu)};
\end{tikzpicture},
\dots,\begin{tikzpicture}[Q,centerzero]
\draw[-to] (0,-0.2) \botlabel{i_l} -- (0,0.2);
\multopendot{0,0}{west}{h_{i_l}(\mu)};
\end{tikzpicture}$.
Note that $\catQH^0 / \catJ$ is a {\em graded} $\kk$-linear category.
The hypothesis (3) means that all $\vec{D}(x_1^{n_1}\cdots 
x_l^{n_l})$ such that $\beta(\vec{D};n_1,\dots,n_l) \neq 0$
are of degree strictly smaller than the degrees of any of the generating morphisms of $\mathring\catI(\mu)$. Morphism spaces of
$\catQH^0 / \catJ$, hence, of $(\catQH^0 / \catJ) \otimes_\kk \KK$ in such small degrees have the same bases as the corresponding morphism spaces of $\catQH^0$. 
From this, we deduce that all of these $\vec{D}(x_1^{n_1}\cdots x_l^{n_l})$ are linearly independent.
Hence, $f(\beta(\vec{D};n_1,\dots,n_l)) = 0$ for all $\vec{D} \in \oshape(\bi\times\bj)$
and $n_1,\dots,n_l \geq 0$. 
The restriction of $f:\R \rightarrow \KK$ to 
$\kk\big[\beta_i^{(n)}\:|\:i \in I, 1\leq n \leq h_i(\nu)\big]$ is injective
since the images of its generators are algebraically independent over $\kk$. Each $\beta(\vec{D};n_1,\dots,n_l)$
lies in this subalgebra by the hypothesis (2).
So we have that $\beta(\vec{D};n_1,\dots,n_l) = 0$
for all $\vec{D}$ and $n_1,\dots,n_l$.
\end{proof}

We have already discussed some of the well-known 
consequences of \cref{nd}. Here is one more application, which strengthens \cref{better}.

\begin{cor}\label{best}
The functor $F:
\catQH(\mu) \rightarrow \catH(\mu|\nu)$ 
from \cref{obviousone}
is fully faithful for all $\mu, \nu \in X^+$ and all choices of the commutative $\kk$-algebra $\KK$ and the polynomials
$\mu_i(x), \nu_i(x) \in \KK[x]$.
\end{cor}

\begin{proof}
This follows from \cref{nd} by the argument from the proof of \cite[Lem.~5.6]{HKM}.
\end{proof}

\subsection{Non-degeneracy of 2-iquantum groups with geometric parameters}\label{geometricnondeg}

Now we would like to 
upgrade \cref{nd} from 2-quantum groups to 2-iquantum groups. 
Let notation be as in \cref{new2cat}.
We remind the reader that the existence of an admissible choice of parameters implicitly puts a restriction on the Cartan matrix, as discussed in \cref{unfortunately}.

Take words $\bi, \bj \in \langle I \rangle$
with $\wt^\imath(\bi) = \wt^\imath(\bj)$.
For any $\lambda \in X^\imath$, 
$B_\bi\one_\lambda$ and $B_\bj\one_\lambda$
are objects of the same 1-morphism category of $\UU^\imath$.
An {\em $\bi \times \bj$-shape} now means an unoriented
string diagram $D$ with
no dots and no bubbles so that
$D\ \text{\regionlabel{\lambda}}: B_\bj\one_\lambda \Rightarrow B_\bi\one_\lambda$ is a 2-morphism in $\UU^\imath$ for any $\lambda \in X^\imath$. As before, such a shape is {\em reduced} if there are no zig-zags, bigons or loops, and
two shapes are {\em equivalent} if they define the same matching between the boundary points. 
Then we fix a set $\shape(\bi\times\bj)$ of representatives for the equivalence classes of reduced $\bi \times \bj$ shapes
and, for $D \in \shape(\bi\times\bj)$ with $l$ connected components, we fix a normally-ordered choice of points $p_1,\dots,p_l$ away from crossings and critical points, one on each connected component of the shape.
For 
$f \in \kk[x_1,\dots,x_l]$, we let
$D(f)\ \text{\regionlabel{\lambda}}:B_\bj\one_\lambda\Rightarrow B_\bi\one_\lambda$
be the 2-morphism in $\UU^\imath$ 
obtained from $D$ by pinning $f$ to the points $p_1,\dots,p_l$.
Also recall the decomposition $I = I_1 \sqcup I_0 \sqcup I_{-1}$
from \cref{usefuldecomposition}.

\begin{defin}[Non-degeneracy]\label{defnondeg}
Let
$\R^\imath := \kk\big[\beta_i^{(m)}, \beta_j^{(n)}\:\big|\:i \in I_1, j \in I_0, m,n \geq 1, n\text{ odd}\big]$
viewed as a graded algebra with $\beta_i^{(m)}$ 
and $\beta_j^{(n)}$ in
degrees $2d_i m$ and $2d_j n$, respectively.
Given $\lambda \in X^\imath$, let
$f_\lambda^\imath:\R^\imath \rightarrow \End_{\UU^\imath}(\one_\lambda)$
be the graded algebra homomorphism such that
\begin{align}\label{extraordinary}
f_\lambda^\imath\big(\beta_i^{(m)}\big) &=
\zeta_i(\lambda)^{-1}
\begin{tikzpicture}[baseline=-1mm,iQ]
\draw[-] (-0.25,0) arc(180:-180:0.25);
\node at (-.47,0) {\strandlabel{\tau i}};
\region{0.5,0}{\lambda};
\dottybubblelabel{0,0}{m};
\end{tikzpicture}
,& f_\lambda^\imath\big(\beta_j^{(n)}\big)&
=
\zeta_j(\lambda)^{-1}
\begin{tikzpicture}[baseline=-1mm,iQ]
\draw[-] (-0.25,0) arc(180:-180:0.25);
\node at (-.42,0) {\strandlabel{j}};
\region{0.5,0}{\lambda};
\dottybubblelabel{0,0}{n};
\end{tikzpicture}
\end{align}
for $i \in I_1$, $j \in I_0$ and $m,n \geq 1$ with $n$ odd.
We use it to make $\Hom_{\UU^\imath}(B_\bj\one_\lambda, B_\bi\one_\lambda)$ into a graded right $\R^\imath$-module for any $\bi,\bj \in \langle I \rangle$ with $\wt^\imath(\bi)=\wt^\imath(\bj)$, so $\beta \in \R^\imath$
acts by horizontal composition on the right with
$f^\imath_\lambda(\beta)$.
We say that the 2-iquantum group $\UU^\imath$ is {\em non-degenerate} if
the 2-morphism space $\Hom_{\UU^\imath}(B_\bj\one_\lambda, B_\bi\one_\lambda)$ is free as a right $\R^\imath$-module with basis 
\begin{equation}\label{ibasis}
\left\{D(x_1^{n_1} \cdots x_l^{n_l})\ \text{\regionlabel{\lambda}}\:\big|\:D \in \shape(\bi\times\bj), 
n_1,\dots,n_l \geq 0
\right\}
\end{equation}
for each $\lambda \in X^\imath$, $\bi, \bj \in \langle I \rangle$
with $\wt^\imath(\bi)=\wt^\imath(\bj), 
l(\bi) \equiv l(\bj)\pmod{2}$,
and $l := (l(\bi)+l(\bj))/2$.
Taking $\bi = \bj = \varnothing$, non-degeneracy implies that
$f^\imath_\lambda:\R^\imath \rightarrow \End_{\UU^\imath}(\one_\lambda)$ is an isomorphism.
\end{defin}

The spanning part of non-degeneracy is relatively straightforward:

\begin{lem}\label{easyspanning}
The 2-morphism space
$\Hom_{\UU^\imath}(B_\bj\one_\lambda, B_\bi\one_\lambda)$ is generated as a right $\R^\imath$-module by the set \cref{ibasis}
for each $\lambda \in X^\imath$, $\bi, \bj \in \langle I \rangle$
with $\wt^\imath(\bi)=\wt^\imath(\bj), 
l(\bi) \equiv l(\bj)\pmod{2}$,
and $l := (l(\bi)+l(\bj))/2$.
\end{lem}

\begin{proof}
Any dotted bubble representing a 2-endomorphism of $\one_\lambda$
can be expressed as a polynomial in the elements of the set
$$
\left\{\begin{tikzpicture}[baseline=-1mm,iQ]
\draw[-] (-0.25,0) arc(180:-180:0.25);
\node at (-.47,0) {\strandlabel{\tau i}};
\region{0.5,0}{\lambda};
\dottybubblelabel{0,0}{m};
\end{tikzpicture}\ ,\begin{tikzpicture}[baseline=-1mm,iQ]
\draw[-] (-0.25,0) arc(180:-180:0.25);
\node at (-.38,0) {\strandlabel{j}};
\region{0.5,0}{\lambda};
\dottybubblelabel{0,0}{n};
\end{tikzpicture}\:\Big|\:i \in I_1, j \in I_0, m,n \geq 1, n\text{ odd}
\right\}.
$$
This is clear from
the relation \cref{iinfgrass} for bubbles
 $\begin{tikzpicture}[baseline=-1mm,iQ]
\draw[-] (-0.25,0) arc(180:-180:0.25);
\node at (-.47,0) {\strandlabel{\tau i}};
\region{0.5,0}{\lambda};
\dottybubblelabel{0,0}{m};
\end{tikzpicture}$ for $i \in I_{-1}$.
For bubbles with $i \in I_0$, it is not such an obvious consequence of \cref{iinfgrass}, but using this relation any
$\begin{tikzpicture}[baseline=-1mm,iQ]
\draw[-] (-0.25,0) arc(180:-180:0.25);
\node at (-.37,0) {\strandlabel{i}};
\region{0.5,0}{\lambda};
\dottybubblelabel{0,0}{n};
\end{tikzpicture}$ with $i \in I_0$ and $n$ even can be written as a quadratic expression involving these bubbles for smaller odd values of $n$; e.g., see the paragraph before \cite[Cor.~2.6]{BWWbasis} where this is deduced from an identity involving Schur's $q$-functions.
In view of this, the 2-morphisms \cref{ibasis} generate
$\Hom_{\UU^\imath}(B_\bj\one_\lambda, B_\bi\one_\lambda)$ as a right $\End_{\UU^\imath}(\one_\lambda)$-module
because there is a straightening algorithm similar to the Khovanov-Lauda straightening algorithm for 2-quantum groups mentioned in the proof of \cref{nd}. This
proceeds by induction on the number of crossings, using all of the defining relations of $\UU^\imath$.
\end{proof}

The difficult part of establishing non-degeneracy of 2-iquantum groups is the linear independence.
The following treats the case of geometric parameters as in \cref{tojapan}.

\begin{theo}\label{ind}
Assuming that the Cartan matrix is symmetric and that the parameters
satisfy $Q_{i,j}(x,y) = t_{i,j}(x-y)^{-a_{i,j}}$
for all $i \neq j$ (with $t_{i,j} \in \kk^\times$ as usual), the 2-iquantum group $\UU^\imath$ is non-degenerate.
\end{theo}

\begin{proof}
In view of \cref{easyspanning},
we must show that the 2-morphisms \cref{ibasis} are $\R^\imath$-linearly independent.
Let $\UU(\del,\zeta)$, $\UUloc(\del,\zeta)$
and $\widehat{\UUloc}(\del,\zeta)$ be as in \cref{newgrading,sinternalbubbles}.
We are going to work with
the $\kk$-linear 2-category 
$\UU_c$ obtained 
by completing $\UU(\del,\zeta)$ with respect to the grading. It has the same objects and 1-morphisms as $\UU(\del,\zeta)$ and,
for words $\bk,\bl \in \langle I^+,I^- \rangle$ 
of the same weight,
its 2-morphism space
$\Hom_{\UU_c}(G_\bl\one_\lambda, G_\bk\one_\lambda)$ is the product $\prod_{n \in \Z} \Hom_{\UU(\del,\zeta)}(G_\bl\one_\lambda, G_\bk\one_\lambda)_n$. The horizontal and vertical compositions of 2-morphisms in $\UU(\del,\zeta)$ extend to the grading completion because $\Hom_{\UU(\del,\zeta)}(G_\bl\one_\lambda, G_\bk\one_\lambda)_n = \{0\}$ for $n \ll 0$.
The
2-morphism spaces in $\UU_c$ have topological bases arising from \cref{nd}: given $\lambda \in X$,
$\End_{\UU_c}(\one_\lambda)$ is the grading completion of
$\End_{\UU(\del,\zeta)}(\one_\lambda)$
and $\bk,\bl$ as before with
$l(\bk)\equiv l(\bl)\pmod{2}$
and $l := (l(\bk)+l(\bl))/2$,
the space
$\Hom_{\UU_c}(G_\bl\one_\lambda, G_\bk\one_\lambda)$ is free as a topological $\End_{\UU_c}(\one_\lambda)$-module 
with basis
$$
\big\{\vec{D}(x_1^{n_1},\dots,x_l^{n_l})\:\big|\:\vec{D} \in \oshape(\bk\times\bl), n_1,\dots,n_l \geq 0\big\}.
$$
This makes $\UU_c$ better for the present purpose than
the localization $\UUloc(\del,\zeta)$, where we do not know any reasonable basis.

Fix $\hbar \in \kk^\times$.
There is an automorphism (not graded)
\begin{align*}
\eta:\UU(\del,\zeta)&\rightarrow \UU(\del,\zeta),\\
\begin{tikzpicture}[Q,centerzero]
\draw[-to] (0,-0.3) \botlabel{i} -- (0,0.3);
\opendot{0,0};
\region{0.2,0}{\lambda};
\end{tikzpicture}&\mapsto 
\begin{tikzpicture}[Q,centerzero]
\draw[-to] (0,-0.3) \botlabel{i} -- (0,0.3);
\Pin{0,0}{-.6,0}{x+\hbar};
\region{0.2,0}{\lambda};
\end{tikzpicture},\quad
\begin{tikzpicture}[Q,centerzero,scale=.9]
\draw[-to] (-0.3,-0.3) \botlabel{i} -- (0.3,0.3);
\draw[-to] (0.3,-0.3) \botlabel{j} -- (-0.3,0.3);
\region{0.35,0}{\lambda};
\end{tikzpicture}\mapsto \begin{tikzpicture}[Q,centerzero,scale=.9]
\draw[-to] (-0.3,-0.3) \botlabel{i} -- (0.3,0.3);
\draw[-to] (0.3,-0.3) \botlabel{j} -- (-0.3,0.3);
\region{0.35,0}{\lambda};
\end{tikzpicture},\quad
\begin{tikzpicture}[Q,centerzero]
\draw[-to] (-0.25,-0.15) \botlabel{i} to [out=90,in=90,looseness=3](0.25,-0.15);
\region{0.45,0.1}{\lambda};
\end{tikzpicture}
\mapsto \begin{tikzpicture}[Q,centerzero]
\draw[-to] (-0.25,-0.15) \botlabel{i} to [out=90,in=90,looseness=3](0.25,-0.15);
\region{0.45,0.1}{\lambda};
\end{tikzpicture},\quad
\begin{tikzpicture}[Q,centerzero]
\draw[-to] (-0.25,0.15) \toplabel{i} to[out=-90,in=-90,looseness=3] (0.25,0.15);
\region{0.45,-.1}{\lambda};
\end{tikzpicture}\mapsto 
\begin{tikzpicture}[Q,centerzero]
\draw[-to] (-0.25,0.15) \toplabel{i} to[out=-90,in=-90,looseness=3] (0.25,0.15);
\region{0.45,-.1}{\lambda};
\end{tikzpicture}.
\end{align*}
This follows by an easy relations check;
the key point is that $Q_{i,j}(x+\hbar,y+\hbar) = Q_{i,j}(x,y)$.
We also have that
\begin{align*}
\eta\left(\begin{tikzpicture}[baseline=-1mm,Q]
\draw[to-] (-0.25,0) arc(180:-180:0.25);
\node at (0,-.4) {\strandlabel{i}};
\region{0.95,0}{\lambda};
\node at (.55,0) {$(u)$};
\end{tikzpicture}\!\!\right)
&=\begin{tikzpicture}[baseline=-1mm,Q]
\draw[to-] (-0.25,0) arc(180:-180:0.25);
\node at (0,-.4) {\strandlabel{i}};
\region{1.6,0}{\lambda};
\node at (.86,0) {$(u-\hbar)$};
\end{tikzpicture},&
\eta\left(\begin{tikzpicture}[baseline=-1mm,Q]
\draw[to-] (-0.25,0) arc(-180:180:0.25);
\node at (0,-.4) {\strandlabel{i}};
\region{0.95,0}{\lambda};
\node at (.55,0) {$(u)$};
\end{tikzpicture}\!\!\right)
&=\begin{tikzpicture}[baseline=-1mm,Q]
\draw[to-] (-0.25,0) arc(-180:180:0.25);
\node at (0,-.4) {\strandlabel{i}};
\region{1.6,0}{\lambda};
\node at (.86,0) {$(u-\hbar)$};
\end{tikzpicture}\ .
\end{align*}
The proof of this is explained in \cite[Lem.~5.3]{BWWbasis}.
The 2-morphisms obtained by applying $\eta$ to \cref{bensidea,teletubby} have non-zero constant terms (for \cref{bensidea} this uses the assumption 
about the parameters again).
So they are invertible in $\UU_c$.
Hence, $\eta:\UU(\del,\zeta)
\rightarrow \UU_c$ induces a strict $\kk$-linear 
2-functor
$\underline{\eta}:\UUloc(\del,\zeta)\rightarrow \UU_c$.
Finally, we collapse the object set of $\UU_c$ along the fibers of the quotient map $X \twoheadrightarrow X^\imath$ in the same way that
$\widehat{\UUloc}(\del,\zeta)$
was constructed from $\UUloc(\del,\zeta)$,
to obtain a $\kk$-linear 2-category $\widehat{\UU}_c$ with object set $X^\imath$, and a strict $\kk$-linear 2-functor
$\underline{\widehat{\eta}}:\underline{\widehat{\UU}}(\del,\zeta)
\rightarrow \widehat{\UU}_c$.
We will use the composition
$\underline{\widehat{\eta}} \circ \Xi^\imath:
\UU^\imath\rightarrow \widehat{\UU}_c$ for $\Xi^\imath$ from \cref{boxingday} to complete the proof of the theorem.

Suppose for a contradiction 
that there is a non-trivial linear relation
\begin{equation}\label{britishcolumbia}
\sum_{\substack{D \in \shape(\bi\times\bj)\\n_1,\dots,n_l \geq 0}}
D(x_1^{n_1} \cdots x_l^{n_l})\ \text{\regionlabel{\lambda}}\ 
f^\imath_\lambda(\beta(D;n_1,\dots,n_l)) = 0
\end{equation}
for $\beta(D;n_1,\dots,n_l) \in \R^\imath$.
Choose a lift $\hat\lambda \in X$ of $\lambda$ and 
a shape $D \in \shape(\bi\times\bj)$
which has a maximal number of crossings amongst
all $D'$ such that
$\beta(D';n_1,\dots,n_l) \neq 0$ for some 
$n_1,\dots,n_l$. Let $\vec{D}$ be the string diagram obtained from $D$ by orienting strings so that each cap and each cup is directed from left to right
and each propagating string is directed from bottom to top,
labelling endpoints by the same label $i$ as in $D$ if downward and by $\tau i$ if upward. This produces words $\bk,\bl \in \langle I^+,I^-\rangle$
such that
$\vec{D}$ is an $\bk \times \bl$ shape.
Since $D$ is reduced, $\vec{D}$ is reduced, and we may assume that $\vec{D} \in \oshape(\bk\times\bl)$, taking the same choice of distinguished points for $\vec{D}$ as for $D$.
We observe since it is reduced 
that $\vec{D}$ only involves upward, downward and rightward crossings and rightward cups and caps.
Let $f:\R^\imath \rightarrow \End_{\UU_c}(\one_{\hat\lambda})$ be the homomorphism mapping
\begin{align*}
\beta_i^{(m)} &\mapsto
(-1)^{h_{\tau i}(\hat\lambda)}
\left[ 
\begin{tikzpicture}[Q,centerzero,scale=1]
\draw[to-] (-0.68,0) arc(180:-180:0.25);
\node[black] at (0.26,0) {$(-u)$};
\node at (-.43,-.4) {\strandlabel{\tau i}};
\region{.85,0}{\hat\lambda};
\end{tikzpicture} 
\begin{tikzpicture}[Q,centerzero,scale=1]
\draw[-to] (-.25,0) arc(180:-180:0.25);
\node[black] at (.54,0) {$(u)$};
\node at (0,-0.4) {\strandlabel{i}};
\end{tikzpicture}\right]_{u:-\lambda_i-m}\!\!\!, \quad
& \beta_j^{(n)}&\mapsto 
(-1)^{h_{\tau j}(\hat\lambda)}
\left[ 
\begin{tikzpicture}[Q,centerzero,scale=1]
\draw[to-] (-0.68,0) arc(180:-180:0.25);
\node[black] at (0.26,0) {$(-u)$};
\node at (-.43,-.4) {\strandlabel{j}};
\region{.85,0}{\hat\lambda};
\end{tikzpicture} 
\begin{tikzpicture}[Q,centerzero,scale=1]
\draw[-to] (-.25,0) arc(180:-180:0.25);
\node[black] at (.54,0) {$(u)$};
\node at (0,-0.4) {\strandlabel{j}};
\end{tikzpicture}\right]_{u:-n}\!\!\!,
\end{align*}
for $i \in I_1$, $j \in I_0$ and $m,n \geq 1$ with $n$ odd.
These are algebraically independent elements of
$\End_{\UU_c}(\one_{\hat\lambda})$ (see \cite[(5.9)]{BWWbasis} which explains this for $i \in I_0$). So $f$ is injective.
Now we apply the 2-functor 
$\underline{\widehat{\eta}} \circ \Xi^\imath$ to \cref{britishcolumbia} then project onto the $G_\bk \one_{\hat\lambda} \times G_\bl \one_{\hat\lambda}$-matrix entry of the result
to obtain the equality
\begin{equation}\label{saskatchewan}
\sum_{n_1,\dots,n_l \geq 0}
c_D D(x_1^{n_1} \cdots x_l^{n_l})\ \text{\regionlabel{\hat\lambda}}\ 
f(\beta(D;n_1,\dots,n_l)) +
(*) = 0
\end{equation}
in $\Hom_{\UU_c}\big(G_\bl\one_{\hat\lambda},
G_\bk\one_{\hat\lambda}\big)$,
where $c_D$ is a non-zero scalar depending only on $D$ arising as a product of the coefficients in \cref{psi1,psi2,psi3,psi5a,psi5b,psi5c,psi5d},
and $(*)$ is an $\End_{\UU_c}(\one_{\hat\lambda})$-linear combination of diagrams of shapes that are not equivalent to $\vec{D}$ with the same number as or fewer crossings than $D$. The terms in $(*)$ may involve teleporters, and leftward cups, caps and crossings decorated with internal bubbles.
All terms in $(*)$ can be ``straightened'' to rewrite them as 
$\End_{\UU_c}(\one_{\hat\lambda})$-linear combination of 
$\vec{D}'(x_1^{n_1} \cdots x_l^{n_l})$ for $\vec{D}' \in 
\oshape(\bk\times\bl)$ different from $\vec{D}$.
Using what we know about the topological basis of
$\Hom_{\UU_c}(G_\bl\one_{\hat\lambda}, G_\bk\one_{\hat\lambda})$, we deduce that
from \cref{saskatchewan}
that $f(\beta(D;n_1,\dots,n_l)) = 0$ for all $n_1,\dots,n_l \geq 0$. As $f$ is injective, this shows that
all $\beta(D;n_1,\dots,n_l)$ are zero, which is a contradiction.
\end{proof}

We have some hope that the general approach in the proofs of \cref{nd,ind} can also be used to prove non-degeneracy for non-geometric values of parameters. However, it seems to be hard to understand the spectrum when one makes some
sufficiently generic deformation of the polynomials $Q_{i,j}(x,y)$, so we have been unsuccessful in our attempts to prove the following:

\begin{conj}\label{conj1}
$\UU^\imath$ is non-degenerate in all cases.
\end{conj}

\subsection{Graphical interpretation of bilinear forms}\label{sgraphical}

In this subsection, we go back to the general quasi-split iquantum group setup of \cref{iqg}.
We will not be using 2-iquantum groups so we do not need any choice of parameters $Q_{i,j}(x,y)$ or normalization functions
$c_i$ to have been made. 
This means that the results in this subsection (except for \cref{bcor})
are valid for {\em all} quasi-split iquantum groups, 
even the ones excluded by \cref{unfortunately}.
When there is some underlying 
quasi-split iquantum group $\UU^\imath$ which is non-degenerate in the sense of \cref{defnondeg},
the bilinear form on $\dot\U^\imath$ described combinatorially here 
computes graded dimensions of 2-morphism spaces in $\UU^\imath$
(this is the content of \cref{bcor}).

We need to make a few more reminders about quantum groups.
Let $\f$ be Lusztig's $\Q(q)$-algebra 
generated by $\theta_i\:(i \in I)$. 
It is isomorphic both to the subalgebra $\U^+$
of $\U$ generated by all $e_i\:(i \in I)$ and to the subalgebra
$\U^-$ of $\U$ generated by all $f_i\:(i \in I)$ via the $\Q(q)$-algebra homomorphisms
\begin{align}
\f\stackrel{\sim}{\rightarrow} \U^+, \quad 
x &\mapsto x^+,&
\f\stackrel{\sim}{\rightarrow} \U^-, \quad 
x &\mapsto x^-,&
\end{align}
defined by setting $\theta_i^+ := e_i$ and $\theta_i^- := f_i$.
The algebra $\f$ is naturally graded by
\begin{equation}\label{poset}
\Lambda := \sum_{i \in I} \N \alpha_i
\end{equation} 
with $\theta_i$ being of degree $\alpha_i$ (Lusztig denotes this monoid by $\N[I]$).
Let ${_i}R$ and $R_i$ be the $\Q(q)$-linear endomorphisms of $\f$ denoted 
${_i}r / (1-q_i^{-2})$ and $r_i / (1-q_i^{-2})$ in \cite[$\S$1.2.13]{Lubook}.
They are uniquely determined by the properties 
that ${_i}R(1) = R_i(1) = 0$, 
${_i}R(\theta_j) = R_i(\theta_j) = \delta_{i,j} / (1-q_i^{-2})$, and
\begin{align}\label{moaaz}
{_i}R(xy) &= {_i}R(x) y + q_i^{h_i(\alpha)} x\ {_i}R(y),
&
R_i(xy) &= q_i^{h_i(\beta)} R_i(x) y + x R_i(y)
\end{align}
for $x \in \f_{\alpha}$ and $y \in \f_{\beta}$.
By \cite[Prop.~3.1.6(b)]{Lubook}, for $i\in I$ and $x\in \f$, we have that
\begin{align} \label{Eiy}
[e_i,x^-] =q_i^{-1} \big(k_i\, {_i}R(x)^- -R_i(x)^- k_i^{-1}\big).
\end{align}
Lusztig's non-degenerate symmetric bilinear form 
$(\cdot, \cdot):\f\times\f\rightarrow \Q(q)$ 
is characterized by the properties $(1,1)=1$ and 
\begin{align}\label{adj}
    (\theta_i x, y) &= (x, {_i}R(y)), 
    &
    (x \theta_i, y) &= (x, R_i(y))
\end{align}
for $x,y \in \f$ and $i \in I$; see \cite[1.2.13(a)]{Lubook}.

Recall from \cite[$\S$19.1.1]{Lubook} the linear anti-involution $\rho$ of the algebra $\U$ such that 
\begin{align}\label{flop}
\rho(e_i)&=q_i k_i f_i,& \rho(f_i)&=q_i^{-1} e_i k_i^{-1},&
\rho(q^h)=q^h
\end{align}
for all $i\in I$ and $h \in Y$.
It induces an anti-automorphism of $\dot\U$
which fixes each $1_\lambda\:(\lambda \in X)$ and satisfies
\begin{align}\label{poor}
\rho(e_i 1_\lambda) &=q_i^{1+h_i(\lambda)} 1_\lambda f_i,
&
\rho(f_i 1_\lambda) &= q_i^{1-h_i(\lambda)} 1_\lambda e_i.
\end{align}
By \cite[Prop.
4.6]{BW18QSP}, $\rho$ preserves the subalgebra $\U^\imath$, and
\begin{align}
    \rho(b_i) =q_i^{-1-\del_{\tau i}} b_{\tau i}k_{\tau i}k_i^{-1}.
\end{align}
This induces an anti-automorphism of $\dot\U^\imath$
which fixes each $1_\lambda\:(\lambda \in X^\imath)$ and satisfies
\begin{equation}\label{paw}
\rho(b_i 1_\lambda) = q_i^{1+\del_i-\lambda_i} 1_\lambda b_{\tau i}.
\end{equation}
The powers of $q$ in \cref{poor,paw} are degrees of cups and caps in $\UU$ and $\UU^\imath$; cf.
\cref{table1,table2}.

Denote by $L(\lambda)$ the irreducible $\U$-module of highest weight $\lambda \in X^+$, with $\eta_\lambda$ being a highest weight vector.
There is a unique symmetric bilinear form 
$(\cdot, \cdot)_{\lambda}$
on $L(\lambda)$
such that $(\eta_{\lambda}, \eta_{\lambda})_\lambda =1$
and
\begin{equation}\label{greenfish}
(uv,w)_\lambda = (v, \rho(u) w)_\lambda
\end{equation}
for all $v,w \in L(\lambda)$ and $u \in \U$.
Moreover, $(v,w)_\lambda = 0$ for $v$ and $w$ from different weight spaces of $L(\lambda)$.
By \cite[Prop.~19.3.7]{Lubook},
the form $(\cdot,\cdot)$ on $\f$ 
is related to these forms by
\begin{equation}\label{redfish}
(x,y) =\lim_{\lambda\rightarrow\infty} (x^-\eta_{\lambda}, y^-\eta_{\lambda})_{\lambda}
\end{equation}
for $x,y\in \f$,
where the limit is taken over 
$\lambda \in X$ with
$h_i(\lambda) \rightarrow \infty$ for all $i \in I$, and the convergence is in  
$\Q\lround q^{-1}\rround$.

From now on, $\lambda$ denotes an iweight in $X^\imath$, and $\hat\lambda \in X$ is a pre-image.
According to \cite[Def. 6.25]{BW18QSP},
there is a symmetric bilinear form $(\cdot, \cdot)^\imath:\dot\U^\imath\times\dot\U^\imath \rightarrow \Q(q)$ such that 
$(x1_\lambda,y1_{\lambda'})^\imath=0$ for $x,y \in \U^\imath$ and $\lambda \neq \lambda'$ in $X^\imath$ and
\begin{equation}\label{bluefish}
(x 1_\lambda ,y 1_\lambda)^\imath =\lim_{\hat\lambda\rightarrow\infty} (x\eta_{\hat\lambda}, y\eta_{\hat\lambda})_{\hat\lambda}
\end{equation}
for $x,y\in \U^\imath$, where the limit is taken over pre-images $\hat\lambda \in X$ of $\lambda$ with
$h_i(\hat\lambda) \rightarrow \infty$ for all $i \in I$, and convergence is in  
$\Q\lround q^{-1}\rround$.
By \cite[Cor.~6.26]{BW18QSP},
we have that
\begin{equation}\label{iadj}
(ux1_\lambda,y1_\lambda)^\imath =(x 1_\lambda, \rho(u)y1_\lambda)^\imath
\end{equation}
for $u,x,y\in \U^\imath$ and $\lambda \in X^\imath$.
Also $(\cdot,\cdot)^\imath$ is non-degenerate by \cite[Th.~6.27]{BW18QSP} and \cite[Th.~7.6]{BW21iCB}.
The following theorem 
shows that $\U^\imath 1_\lambda$
equipped with this bilinear form is isometric
to $\f$ with its bilinear form
$(\cdot,\cdot)$. 

\begin{theo}\label{jbversion}
Fix $\lambda \in X^\imath$.
There is a unique $\Q(q)$-linear isomorphism
$\jmath:\U^\imath 1_\lambda \stackrel{\sim}{\rightarrow} \f$ such that
\begin{equation}\label{holding}
\lim_{\hat\lambda\rightarrow\infty}\big(x \eta_{\hat\lambda}, y^- \eta_{\hat\lambda}\big)_{\hat \lambda} 
=
\big(\jmath(x1_\lambda),y\big)
\end{equation}
for $x \in \U^\imath$ and $y \in \f$, where the limit is taken over pre-images $\hat\lambda$ of $\lambda$ with $h_i(\hat\lambda)\rightarrow \infty$ for all $i \in I$, and convergence is in $\Q\lround u^{-1}\rround$.
Moreover:
\begin{enumerate}
\item $\jmath(1_\lambda) = 1$;
\item $\jmath(b_ix 1_\lambda) =\theta_i \, \jmath(x 1_\lambda) + q_i^{1+\del_i-\kappa_i} {_{\tau i}}R\big(\jmath(x 1_\lambda)\big)$
for any $x\in \U^\imath$ which is homogeneous in the sense that $x 1_\lambda = 1_\kappa x$ for $\kappa \in X^\imath$;
\item
$(x 1_\lambda,y1_\lambda)^\imath =\big(\jmath(x 1_\lambda), \jmath(y 1_\lambda) \big)$
for all $x,y \in \U^\imath$. 
\end{enumerate}
\end{theo}

\begin{proof}
Suppose that $x$ is some element of $\U^\imath$ with
$1_\kappa x = x 1_\lambda$ for $\kappa \in X^\imath$,
and that we are given an element $\jmath(x1_\lambda) \in \f$
such that \cref{holding} is true for all $y \in \f$. Take some $i \in I$. We claim that there is a unique element $\jmath(b_i x1_\lambda)$ of $\f$ such that
\begin{equation}\label{holding2}
\lim_{\hat\lambda\rightarrow\infty}\big(b_i x \eta_{\hat\lambda}, y^- \eta_{\hat\lambda}\big)_{\hat \lambda} 
=
\big(\jmath(b_i x1_\lambda), y\big)
\end{equation}
for all $y \in \f$, namely, the element
$\jmath(b_i x1_\lambda)$ defined by the formula in (2).
The uniqueness of $\jmath(b_i x 1_\lambda)$ is clear from the non-degeneracy of the form $(\cdot,\cdot)$.
For existence, we have for $y \in \f$
and $\jmath(b_i x1_\lambda)$ defined by this formula that
\begin{align*}
\big(\jmath(b_i x1_\lambda), y\big)
&=\left(\theta_i \jmath(x 1_\lambda) + q_i^{1+\del_i-\kappa_i} {_{\tau i}}R\big(\jmath(x 1_\lambda)\big), y \right)\\
&\ \substack{\cref{adj}\\{\textstyle=}}\ \left(\jmath(x 1_\lambda),{_i}R(y)  + q_i^{1+\del_i-\kappa_i} \theta_{\tau i} y \right)\\
&\ \substack{\cref{holding}\\{\textstyle=}}\ 
\lim_{\hat\lambda\rightarrow\infty}\left(x \eta_{\hat\lambda},\big({_i}R(y)^- + q_i^{1+\del_i-\kappa_i} f_{\tau i} y^-\big) \eta_{\hat\lambda}\right)_{\hat\lambda}\\
&\ \ \;\substack{{\textstyle=}}\ \ \  \lim_{\hat\lambda\rightarrow\infty}\left(x \eta_{\hat\lambda},\big({_i}R(y)^- -k_i^{-1}R_i(y)^- k_i^{-1} + q_i^{1+\del_i-\kappa_i} f_{\tau i} y^-\big) \eta_{\hat\lambda}\right)_{\hat\lambda}\\
&\ \substack{\cref{Eiy}\\{\textstyle=}}\ \lim_{\hat\lambda\rightarrow\infty}\left(x \eta_{\hat\lambda},\big(q_i k_i^{-1}e_i y^- + q_i^{1+\del_i-\kappa_i} f_{\tau i} y^-\big) \eta_{\hat\lambda}\right)_{\hat\lambda}
\\&
\ \substack{\cref{adj}\\{\textstyle=}\\\cref{flop}}\ \lim_{\hat\lambda\rightarrow\infty}\left(\big(f_i+q_i^{2+\del_i-\kappa_i} k_{\tau i}^{-1}e_{\tau i}\big)x \eta_{\hat\lambda}, y^- \eta_{\hat\lambda}\right)_{\hat\lambda}\\
&\ \ =\ \ \lim_{\hat\lambda\rightarrow\infty}\left(
(f_i + q_i^{\del_i}
e_{\tau i}k_i^{-1})x\eta_{\hat\lambda},y^- \eta_{\hat\lambda}\right)_{\hat\lambda}
= \lim_{\hat\lambda\rightarrow\infty}\left(
b_i x\eta_{\hat\lambda},y^- \eta_{\hat\lambda}\right)_{\hat\lambda},
\end{align*}
which checks \cref{holding2}.

Now we can prove existence of a linear map $\jmath$ satisfying \cref{holding}. We set $\jmath(1_\lambda) := 1$.
Then \cref{holding} holds for $x=1$ and all $y \in \f$ by the definition \cref{redfish}.
For any monomial $x$
in $\{b_i \:|\: i \in I\}$, we can use the claim in the previous paragraph plus induction on the length of the monomial
to construct a vector $\jmath(x1_\lambda) \in \f$
such that \cref{holding} holds for all $y \in \f$.
Then we pick a basis for
$\U^\imath 1_\lambda$ 
consisting
of some such monomials applied to $1_\lambda$, and define $\jmath$ on these basis elements as just explained, extending linearly to obtain the desired $\Q(q)$-linear map 
$\jmath:\U^\imath 1_\lambda \rightarrow \f$
such that \cref{holding} is true for all $x \in \U^\imath$.
The uniqueness of $\jmath$ follows from the non-degeneracy of the form $(\cdot,\cdot)$.

Next, we establish the properties (1)--(3).
The first is clear from \cref{redfish}. The second follows from the first paragraph of the proof.
For the third, in view of \cref{holding}, it suffices to show that
\begin{equation}
\lim_{\hat\lambda\rightarrow\infty}
\left(x \eta_{\hat\lambda}, \jmath(y1_\lambda)^-\eta_{\hat\lambda}\right)_{\hat\lambda}=(x 1_\lambda, y 1_\lambda)^\imath
\label{holding3}
\end{equation}
for $x,y \in \U^\imath$.
To see this, we may assume that $y$ is a monomial in $\{b_i\:|\:i \in I\}$, and proceed by induction on its length, the case $y=1$ being clear from \cref{bluefish}. For the induction step, we assume that \cref{holding3}
is true for a monomial $y$ with $y 1_\lambda = 1_\kappa y$
and all $x \in \U^\imath$, and prove it for $b_i y$, using (2) for the first equality:
\begin{align*}\lim_{\hat\lambda\rightarrow\infty}&\left(x \eta_{\hat\lambda}, \jmath(b_i y1_\lambda)^- \eta_{\hat\lambda}\right)_{\hat\lambda}
=\lim_{\hat\lambda\rightarrow\infty}
\left(x \eta_{\hat\lambda}, \big(f_i \jmath(y 1_\lambda)^- + q_i^{1+\del_i-\kappa_i} {_{\tau i}}R(\jmath(y 1_\lambda))^-\big)\eta_{\hat\lambda}\right)_{\hat\lambda}\\
&\ \ =\ \ 
\lim_{\hat\lambda\rightarrow\infty}
\left(x \eta_{\hat\lambda}, \big(f_i \jmath(y 1_\lambda)^- + q_i^{-1-\del_{\tau i}}k_i^{-1} 
(k_{\tau i}\; {_{\tau i}}R(\jmath(y 1_\lambda))^- - 
R_{\tau i}(\jmath(y 1_\lambda))^-k_{\tau i}^{-1} 
)\big)\eta_{\hat\lambda}\right)_{\hat\lambda}\\
&\stackrel{\cref{Eiy}}{=}
\lim_{\hat\lambda\rightarrow\infty}
\left(x \eta_{\hat\lambda}, \big(f_i \jmath(y 1_\lambda)^- + q_i^{-\del_{\tau i}}k_i^{-1} 
e_{\tau i}\jmath(y 1_\lambda)^-\big)\eta_{\hat\lambda}\right)_{\hat\lambda}
=\lim_{\hat\lambda\rightarrow\infty}
\left(x \eta_{\hat\lambda}, b_i\, \jmath(y 1_\lambda)^-\eta_{\hat\lambda}\right)_{\hat\lambda}\\
&\stackrel{\cref{greenfish}}{=}\lim_{\hat\lambda\rightarrow\infty}
\left(\rho(b_i) x \eta_{\hat\lambda}, \jmath(y 1_\lambda)^-\eta_{\hat\lambda}\right)_{\hat\lambda}
\stackrel{\cref{holding3}}{=}
\left(\rho(b_i) x 1_\lambda, y 1_\lambda\right)^\imath
\stackrel{\cref{iadj}}{=}
\left(x 1_\lambda, b_i y 1_\lambda\right)^\imath.
\end{align*}
Now (3) is proved.

Finally, we must show that $\jmath$ is an isomorphism. 
Using 
(1)--(2) and induction on length, it follows that any
monomial in $\{\theta_i\:|\:i \in I\}$ lies in the image of $\jmath$. Hence, $\jmath$ is surjective. It is injective
by (3) and the non-degeneracy of the form $(\cdot,\cdot)^\imath$.
\end{proof}

For $\bi,\bj \in \langle I \rangle$, we now need the set $\shape(\bi\times\bj)$ of reduced $\bi \times \bj$ shapes
introduced before \cref{defnondeg}; this set is empty unless $l(\bi) \equiv l(\bj)\pmod{2}$.
We defined $\shape(\bi\times\bj)$ 
in terms of string diagrams representing morphisms in $\UU^\imath$, but the definition can also be formulated in purely diagrammatic terms---it is just a diagram representing 
a matching between the words $\bi$ and $\bj$ which is reduced in the
sense that the number of crossings is as small as possible. Recall the
definitions of $b_\bi$ and $\theta_\bi$ from \cref{wordnotation1,wordnotation1b}.

\begin{theo}\label{bthm}
For $\bi, \bj \in \langle I \rangle$ and $\lambda \in X^\imath$,
we have that
\begin{align}\label{wet0}
(\theta_\bi, \theta_\bj) &= \sum_{\substack{D \in \shape(\bi\times\bj)\\\text{cup-cap-free}}}
q^{-\deg(D\ \text{\regionlabel{\lambda}})}\ \big/\ (1-q_{k_1}^{-2}) \cdots (1-q_{k_l}^{-2}),\\
\label{wet1}
\left( \jmath(b_\bi 1_\lambda), \theta_\bj \right)
&=
\sum_{\substack{D \in \shape(\bi\times\bj)\\\text{\em cap-free}}}
q^{-\deg(D \ \text{\regionlabel{\lambda}})} \ \big/\  (1-q_{k_1}^{-2})\cdots (1-q_{k_l}^{-2}),\\
\big(b_\bi 1_\lambda, b_\bj 1_\lambda\big)^\imath &= 
\sum_{D \in \shape(\bi\times\bj)}
q^{-\deg(D \ \text{\regionlabel{\lambda}})} \ \big/\  (1-q_{k_1}^{-2})\cdots (1-q_{k_l}^{-2}),\label{wet2}
\end{align}
where 
$\deg\big(D \ \text{\regionlabel{\lambda}}\big)$ is computed according to \cref{table2}, and $k_1,\dots,k_l \in I$ are the labels at the distinguished points $p_1,\dots,p_l$ of the strings in $D$.
\end{theo}

\begin{proof}
The first formula \cref{wet0} is well known, e.g., it 
was exploited already in \cite{KL1}.

Next we prove \cref{wet1}. 
We proceed by induction on $l(\bi)$.
If $l(\bi) = 0$ then both sides are 0 unless $l(\bj) = 0$, when they are both 1. So the induction base is true.
Now we assume \cref{wet1} is true for some $\bi$ and all $\bj$, 
take $h \in I$,
and prove the result 
for the slightly longer word $h \bi$ and all $\bj$. 
Let $\kappa := \lambda - \wt^\imath(\bi)$.
By \cref{jbversion}(2) and \cref{adj}, we have that
\begin{equation}\label{saltyseaweed}
\left( \jmath(b_h b_\bi 1_\lambda), \theta_\bj\right)
=
\left( \theta_h \jmath(b_\bi 1_\lambda), \theta_\bj\right)
+ 
q_h^{1+\del_h-\kappa_h}\left(\jmath(b_\bi 1_\lambda), \theta_{\tau h}\theta_\bj\right).
\end{equation}
Now consider a cup-free shape $D \in \shape(h \bi\times\bj)$.
Its vertices along the top edge are labelled according to the letters of the word $h \bi$, and the vertices along the bottom edge are labelled according to $\bj$.
Let $k_1,\dots,k_l$ be the labels of its strings at the distinguished points.
We say $D$ is of {type I} if the string with one endpoint at the top left vertex labelled $h$
is propagating, that is, its other endpoint is at the bottom of the diagram. We say it is of type II otherwise, in which case this string is a cup connecting the top left vertex to another vertex at the top of the diagram.
The argument is completed using \cref{saltyseaweed} on 
adding the equations established in the following two claims:

\vspace{2mm}
\noindent
\underline{Claim 1}:
$\displaystyle\left( \theta_h \jmath(\bi 1_\lambda), \theta_\bj \right)
= \sum_{\substack{D \in \shape(h \bi\times\bj)\\\text{cap-free of type I}}}\!\!
q^{-\deg(D\ \text{\regionlabel{\lambda}})} \ \big/\  (1-q_{k_1}^{-2})\cdots (1-q_{k_l}^{-2})
$.

\vspace{1mm}
\noindent
\underline{Proof}.
Suppose that $\bj = j_1 \cdots j_n$.
For $1 \leq r \leq n$, 
let $\bj[r] := j_1 \cdots j_{r-1} j_{r+1}\cdots j_n$.
By \cref{moaaz}, we have that
$$
{_h}R(\theta_\bj) = 
\sum_{\substack{1 \leq r \leq n \\ j_r = h}}
q_h^{a_{h,j_1}+\cdots + a_{h,j_{r-1}}}
\theta_{\bj[r]}\ \big/ \ (1-q_h^{-2}).
$$
Also using \cref{adj}, we deduce that
\begin{equation}\label{dark}
\left( \theta_h \jmath(b_\bi 1_\lambda), \theta_\bj \right)
=
\left(\jmath(b_\bi 1_\lambda), {_h}R(\theta_\bj) \right)
=
\sum_{\substack{1 \leq r \leq n \\ j_r = h}}
q_h^{a_{h,j_1}+\cdots + q_{h,j_{r-1}}}
\left(\jmath(b_\bi 1_\lambda), \theta_{\bj[r]}\right)\ \big/ \ (1-q_h^{-2}).
\end{equation}
Now let $D \in \shape(h \bi\times\bj)$ be a cup-free shape of type I.
Consider the string in $D$ which has 
one endpoint at the top left vertex.
Define $r(D)$ to be the vertex number of its other endpoint, indexing vertices at the bottom of the diagram 
by $1,\dots,n$ from left to right.
Let $\widetilde{D} \in \shape(\bi\times\bj[r(D)])$ be the cap-free shape obtained by removing this string from $D$. The function 
\begin{align*}
D = \begin{tikzpicture}[centerzero,iQ]
\draw (-.7,-.7) to (-.7,.7);
\draw (.7,-.7) to (.7,.7);
\node at (-0.3,-.63) {$\cdots$};
\node at (0.34,-.63) {$\cdots$};
\node at (0,.63) {$\cdot$};
\node at (-0.32,.63) {$\cdot$};
\node at (0.32,.63) {$\cdot$};
\draw [thick,decoration={brace,raise=0.4cm,mirror},decorate] (-.7,-.4) -- (.7,-.4);
\draw [thick,decoration={brace,raise=0.4cm},decorate] (-.7,.4) -- (.7,.4);
\node at (0,1.1) {$\scriptstyle \bi$};
\node at (0,-1.1) {$\scriptstyle \bj$};
\draw (-1.1,.7) \toplabel{h} to[out=-90,in=90,looseness=1] (0,-.7);
\draw[fill=lightgray] (-.8,.5) to (.8,.5) to (.8,-.5) to (-.8,-.5) to (-.8,.5);
\draw[dotted] (-1.1,.7)  to[out=-90,in=90,looseness=1] (0,-.7);
\node at (0.2,-.3) {$\scriptscriptstyle r(D)$};
\end{tikzpicture}
\quad&\mapsto\quad
\widetilde{D} =
\begin{tikzpicture}[centerzero,iQ]
\draw (-.6,-.7) to (-.6,.7);
\draw (.6,-.7) to (.6,.7);
\node at (0,-.63) {$\cdot$};
\node at (0,.63) {$\cdot$};
\node at (0.27,-.63) {$\cdot$};
\node at (0.27,.63) {$\cdot$};
\node at (-0.27,-.63) {$\cdot$};
\node at (-0.27,.63) {$\cdot$};
\draw [thick,decoration={brace,raise=0.4cm,mirror},decorate] (-.6,-.4) -- (.6,-.4);
\draw [thick,decoration={brace,raise=0.4cm},decorate] (-.6,.4) -- (.6,.4);
\node at (0,1.1) {$\scriptstyle \bi$};
\node at (0,-1.15) {$\scriptstyle \bj[r(D)]$};
\draw[fill=lightgray] (-.7,.5) to (.7,.5) to (.7,-.5) to (-.7,-.5) to (-.7,.5);
\end{tikzpicture}\ 
\end{align*}
is a bijection
from the set of cup-free shapes $D \in \shape(h \bi\times\bj)$
of type I to the disjoint union 
of the sets of cup-free shapes in $\shape(\bi\times\bj[r])$
for $1 \leq r \leq n$ with $j_r = h$;
the inverse bijection inserts a string of color $h$ in the obvious way.
In passing from $D$ to $\widetilde{D}$ we have removed crossings of a string labelled $h$ with strings labelled $j_1,\dots,j_{r(D)-1}$, so 
$\deg(\widetilde{D}\ \text{\regionlabel{\lambda}}) = \deg(D\ \text{\regionlabel{\lambda}})
+ d_h a_{h,j_1}+\cdots+d_h a_{h,j_{r(D)-1}}$.
Putting these things together shows that
\begin{multline*}
\sum_{\substack{D \in \shape(h \bi\times \bj)\\\text{cap-free of type I}}}\!\!
q^{-\deg(D\ \text{\regionlabel{\lambda}})} \ \big/\  (1-q_{k_1}^{-2})\cdots (1-q_{k_l}^{-2})
\\=
\sum_{\substack{D \in \shape(h \bi\times\bj)\\\text{cap-free of type I}}}\!\!
q_h^{a_{h,j_1}+\cdots+a_{h,j_{r(D)-1}}}q^{-\deg(\widetilde{D}\ \text{\regionlabel{\lambda}})} \ \big/\  (1-q_{k_1}^{-2})\cdots (1-q_{k_l}^{-2})
\\=
\sum_{\substack{1 \leq r \leq n\\j_r = h}}
q_h^{a_{h,j_1}+\cdots+a_{h,j_{r-1}}}
\sum_{\substack{\widetilde{D} \in \shape(\bi\times \bj[r])\\\text{cap-free}}}
q^{-\deg(\widetilde{D}\ \text{\regionlabel{\lambda}})} 
\ \big/\  (1-q_{k_1}^{-2})\cdots (1-q_{k_l}^{-2}).
\end{multline*}
By the induction hypothesis, the final expression here is equal to
$$
\sum_{\substack{1 \leq r \leq n\\j_r = h}}
q_h^{a_{h,j_1}+\cdots+a_{h,j_{r-1}}} 
\left(\jmath(b_\bi 1_\lambda), \theta_{\bj[r]}\right)
\ \big/\  (1-q_h^{-2}),
$$
which is the desired
$\left( \theta_h \jmath(b_\bi 1_\lambda), \theta_\bj \right)$
by \cref{dark}.

\vspace{2mm}
\noindent
\underline{Claim 2}:
$q_h^{1+\del_h-\kappa_h}\displaystyle\left(\jmath(b_\bi 1_\lambda), \theta_{\tau h}\theta_\bj \right)
= 
\sum_{\substack{D \in \shape(h \bi,\bj)\\\text{cap-free of type II}}}\!\!
q^{-\deg(D\ \text{\regionlabel{\lambda}})} \ \big/\  (1-q_{k_1}^{-2})\cdots (1-q_{k_l}^{-2})$.

\vspace{1mm}
\noindent
\underline{Proof}.
Consider the bijection
from the set of cap-free shapes $D \in \shape(h \bi\times \bj)$
of type II to the set of cap-free shapes $\widetilde{D} \in \shape(\bi\times (\tau h)\bj)$ defined so that 
\begin{align*}
D  = \begin{tikzpicture}[centerzero,iQ]
\draw (-.4,.6) to (-.4,-.6);
\draw (.4,.6) to (.4,-.6);
\node at (0.02,-.53) {$\cdots$};
\node at (0.02,.53) {$\cdots$};
\draw [thick,decoration={brace,mirror,raise=0.4cm},decorate] (-.4,-.3) -- (.4,-.3);
\draw [thick,decoration={brace,raise=0.4cm},decorate] (-.4,.3) -- (.4,.3);
\node at (0,1) {$\scriptstyle \bi$};
\node at (0,-1) {$\scriptstyle \bj$};
\draw (-.8,.6) \toplabel{h} to[out=-90,in=-135,looseness=1] (-.5,-.1);
\draw[fill=lightgray] (-.5,.4) to (.5,.4) to (.5,-.4) to (-.5,-.4) to (-.5,.4);
\end{tikzpicture}
\quad&\mapsto\quad
\widetilde{D} =
\begin{tikzpicture}[centerzero,iQ]
\draw (-.4,.6) to (-.4,-.6);
\draw (.4,.6) to (.4,-.6);
\node at (0.02,-.53) {$\cdots$};
\node at (0.02,.53) {$\cdots$};
\draw [thick,decoration={brace,mirror,raise=0.4cm},decorate] (-.4,-.3) -- (.4,-.3);
\draw [thick,decoration={brace,raise=0.4cm},decorate] (-.4,.3) -- (.4,.3);
\node at (0,1) {$\scriptstyle \bi$};
\node at (0,-1) {$\scriptstyle \bj$};
\draw (-.8,-.6) \botlabel{\tau h} to[out=90,in=-135,looseness=1] (-.5,-.1);
\draw[fill=lightgray] (-.5,.4) to (.5,.4) to (.5,-.4) to (-.5,-.4) to (-.5,.4);
\end{tikzpicture}\ .
\end{align*}
This removes a cup of degree
$d_h(1+\del_{\tau h} - (\kappa+\alpha_h)_{\tau h})
=-d_h(1+\del_h-\kappa_h)$.
so we have that 
$q^{-\deg(D\ \text{\regionlabel{\lambda}})} = q_h^{1+\del_h-\kappa_h} 
q^{-\deg(\widetilde{D}\ \text{\regionlabel{\lambda}})}$.
Using the induction hypothesis for the first equality, we deduce:
\begin{align*}
q_h^{1+\del_h-\kappa_h}\!\displaystyle\left(\jmath(b_\bi 1_\lambda), \theta_{\tau h}\theta_\bj \right)
&=
q_h^{1+\del_h-\kappa_h}\!
\sum_{\substack{\widetilde{D} \in \shape(\bi\times(\tau h) \bj)\\\text{cap-free}}}\!\!
q^{-\deg(\widetilde{D}\ \text{\regionlabel{\lambda}})} \ \big/\  (1-q_{k_1}^{-2})\cdots (1-q_{k_l}^{-2})\\
&=
\!\!\sum_{\substack{D \in \shape(h \bi\times\bj)\\\text{cap-free of type II}}}\!\!\!
q^{-\deg(D\ \text{\regionlabel{\lambda}})} \ \big/\  (1-q_{k_1}^{-2})\cdots (1-q_{k_l}^{-2}).
\end{align*}

\vspace{2mm}
Now \cref{wet1} is proved. We move on to \cref{wet2}, which we prove by induction on the length of the word $\bj$.
If $\bj$ is of length 0 then $\jmath\left(b_\bj 1_\lambda\right) = 1$ and
$\left(b_\bi 1_\lambda, b_\bj 1_\lambda\right)^\imath
= \left(\jmath(b_\bi 1_\lambda), 1\right)$ 
for all $\bi \in \langle I \rangle$ by \cref{jbversion}.
This is computed by the stated sum over shapes thanks to \cref{wet1}.
This checks the induction base. 
For the induction step, we assume that \cref{wet2} is true for some 
word $\bj \in \langle I \rangle$  and all words $\bi$, and prove it for $h \bj$ for some $h \in I$.
Let $\kappa$ be the iweight with
$b_\bj \one_\lambda = \one_\kappa b_\bj$.
We have that
\begin{align*}
\big(b_\bi 1_\lambda, b_h b_\bj 1_\lambda\big)^\imath 
&\ \substack{\cref{iadj}\\{\textstyle=}\\\cref{paw}}\ 
q_h^{1+\del_h-\kappa_h}\big(b_{\tau h} b_\bi 1_\lambda, b_\bj 1_\lambda\big)^\imath\\&\ \  =\  
q_h^{1+\del_h-\kappa_h}
\sum_{D \in \shape((\tau h) \bi\times\bj)}
q^{-\deg(D \text{\regionlabel{\lambda}})} \ \big/\  (1-q_{k_1}^{-2}\cdots (1-q_{k_l}^{-2})\\
&\ \ =\ \sum_{\widetilde{D} \in \shape(\bi\times h \bj)}
q^{-\deg(\widetilde{D} \text{\regionlabel{\lambda}})} \ \big/\  (1-q_{k_1}^{-2})\cdots (1-q_{k_l}^{-2}).
\end{align*}
The second equality here follows by the induction hypothesis. The final one follows by a similar argument
to the proof of Claim 2, using the bijection
$\shape((\tau h) \bi\times \bj)
\stackrel{\sim}{\rightarrow} \shape(\bi\times h \bj), D \mapsto \widetilde{D}$
defined by
\begin{align*}
D = \begin{tikzpicture}[centerzero,iQ]
\draw (-.4,.6) to (-.4,-.6);
\draw (.4,.6) to (.4,-.6);
\node at (0.02,-.53) {$\cdots$};
\node at (0.02,.53) {$\cdots$};
\draw [thick,decoration={brace,mirror,raise=0.4cm},decorate] (-.4,-.3) -- (.4,-.3);
\draw [thick,decoration={brace,raise=0.4cm},decorate] (-.4,.3) -- (.4,.3);
\node at (0,1) {$\scriptstyle \bi$};
\node at (0,-1) {$\scriptstyle \bj$};
\draw (-.8,.6) \toplabel{\tau h} to[out=-90,in=135,looseness=1] (-.5,-.1);
\draw[fill=lightgray] (-.5,.4) to (.5,.4) to (.5,-.4) to (-.5,-.4) to (-.5,.4);
\end{tikzpicture}\quad&\mapsto\quad
\widetilde{D} =
\begin{tikzpicture}[centerzero,iQ]
\draw (-.4,.6) to (-.4,-.6);
\draw (.4,.6) to (.4,-.6);
\node at (0.02,-.53) {$\cdots$};
\node at (0.02,.53) {$\cdots$};
\draw [thick,decoration={brace,mirror,raise=0.4cm},decorate] (-.4,-.3) -- (.4,-.3);
\draw [thick,decoration={brace,raise=0.4cm},decorate] (-.4,.3) -- (.4,.3);
\node at (0,1) {$\scriptstyle \bi$};
\node at (0,-1) {$\scriptstyle \bj$};
\draw (-.8,-.6) \botlabel{h} to[out=90,in=-135,looseness=1] (-.5,-.1);
\draw[fill=lightgray] (-.5,.4) to (.5,.4) to (.5,-.4) to (-.5,-.4) to (-.5,.4);
\end{tikzpicture}\ .
\end{align*}
This completes the proof of \cref{wet2} and hence the theorem. 
\end{proof}

\begin{rem}\label{recallithere}
As well as the symmetric bilinear forms $(\cdot,\cdot)$ on $\f$ and
$(\cdot,\cdot)^\imath$ on $\dot\U^\imath$, there is Lusztig's non-degenerate 
symmetric form
$(\cdot,\cdot)$ on $\dot\U$, which was mentioned already in \cref{lf}. This is uniquely determined by the following properties:
\begin{itemize}
\item
$(1_\kappa x 1_\lambda, 1_{\kappa'} y 1_{\lambda'}) = 0$ for
any $x,y \in \U$ and weights with $\kappa \neq \kappa'$ or $\lambda \neq \lambda'$;
\item
$(x^+ 1_\lambda, y^+ 1_\lambda)
=(x^- 1_\lambda, y^- 1_\lambda) = (x, y)$ for any $\lambda \in X$ and $x,y \in \f$;
\item
$(ux 1_\lambda, y 1_\lambda) = (x 1_\lambda, \rho(u) y 1_\lambda)$
for any $\lambda \in X$ and $u,x,y \in \U$.
\end{itemize}
In \cite[Th.~2.7]{KL3}, Khovanov and Lauda gave the following graphical interpretation of this bilinear form:
for $\bi,\bj \in \langle I^+,I^-\rangle$ and $\lambda \in X$
we have that
\begin{equation}\label{kldiagrams}
\big(g_\bi 1_\lambda, g_\bj 1_\lambda\big)
= \sum_{\vec{D} \in \oshape(\bi\times\bj)}
q^{-\deg(\vec{D}\ \text{\regionlabel{\lambda}})} \ \big/\  (1-q_{k_1}^{-2})\cdots (1-q_{k_l}^{-2})
\end{equation}
where $\deg(\vec{D}\ \text{\regionlabel{\lambda}})$
is the degree of this 2-morphism $G_\bj \one_\lambda \rightarrow G_\bi \one_\lambda$  computed according to \cref{table1},
and $k_1,\dots,k_l \in I$ are the labels of the strings in $\vec{D}$. 
In view of \cref{earphones}, 
this identity is a special case of the identity \cref{wet2} in \cref{bthm}.
\end{rem}

For the purposes of categorification, we prefer to replace the bilinear forms $(\cdot,\cdot)$ and $(\cdot,\cdot)^\imath$
with sesquilinear 
forms (anti-linear with respect to the bar involution in 
the first argument).
To define these, let $\psi:\f \rightarrow \f$ be the usual bar involution,
that is, the anti-linear algebra involution with $\psi(\theta_i) = \theta_i$ for all $i \in I$. Also
$\psi^\imath:\dot\U^\imath \rightarrow \dot\U^\imath$ is the 
anti-linear bar involution on the iquantum group, with
$\psi^\imath(b_i 1_\lambda) = b_i 1_\lambda$ for all $i \in I$
and $\lambda \in X^\imath$.
We define
\begin{align}\label{newby1}
\langle\cdot,\cdot\rangle:\f\times\f&\rightarrow \Q(q),&
\langle x, y \rangle &:= (\psi(x), y),\\\label{newby2}
\langle\cdot,\cdot\rangle^\imath:\dot\U^\imath \times \dot\U^\imath &\rightarrow \Q(q),&
\left\langle x 1_\lambda, y 1_{\lambda'} \right\rangle^\imath
&:= \left(\psi^\imath(x 1_\lambda), y 1_{\lambda'}\right)^\imath.
\end{align}
To write down the analog of \cref{jbversion} for these sesquilinear forms, 
we also need the linear maps 
\begin{align}\label{newpropsofj}
{_i}\widetilde{R}&:= \psi \circ 
{_i}R \circ \psi :\f \rightarrow \f,
&\tilde\jmath := \psi \circ \jmath\circ \psi^\imath:\dot\U^\imath 1_\lambda \stackrel{\sim}{\rightarrow} \f.
\end{align}
From \cref{moaaz,adj}, we get that
\begin{align}\label{newmoaaz1}
{_i}\widetilde{R}(1) &= 0,&
{_i}\widetilde{R}(\theta_j)&= \frac{\delta_{i,j}}{1-q_i^2},&
{_i}\widetilde{R}(xy) &= {_i}\widetilde{R}(x) y + q_i^{-h_i(\alpha)}x\, {_i}\widetilde{R}(y),\end{align}\begin{align}\label{newmoaaz2}
\big\langle \theta_i x,y \big\rangle &= 
\big\langle x, {_i}R(y) \big\rangle,&
\big\langle x,\theta_i y \big\rangle &= 
\big\langle {_i}\widetilde{R}(x), y \big\rangle,
\end{align}
for $x \in \f_\alpha, y \in \f_\beta$. From \cref{jbversion}, we get that
\begin{align}
\label{propsofjbar1}
\tilde\jmath(1) &= 1,\\\label{propsofjbar1b}
\tilde\jmath(b_i x 1_\lambda) 
&= \theta_i \tilde\jmath(x 1_\lambda) + 
q_i^{\kappa_i-\del_i-1} {_{\tau i}}\widetilde{R}\left(\tilde\jmath(x 1_\lambda)\right),\\
\big\langle x 1_\lambda, y 1_\lambda \big\rangle^\imath
&= \big\langle \tilde\jmath(x 1_\lambda), \jmath(y 1_\lambda) \big\rangle\label{propsofjbar2}
\end{align}
for $x,y \in \U^\imath$ such that $x 1_\lambda = 1_\kappa x$.
For $\bi \in \llangle I \rrangle$, we define\footnote{
The letters $\delta$ and $\atled$ are the lower case forms of $\Delta$ and $\nabla$.
This notation will be justified in the next section when we relate
$\delta_\bi 1_\lambda$ and $\atled_\bi 1_\lambda$ to certain standard and costandard modules $\Delta(\bi)$ and $\nabla(\bi)$.} $\delta_\bi 1_\lambda$
and $\atled_\bi 1_\lambda$
to be the unique elements of $\dot\U^\imath 1_\lambda$ such that
\begin{align}\label{abitwacky}
\tilde\jmath(\delta_\bi 1_\lambda) &= \theta_\bi,&
\jmath(\atled_\bi 1_\lambda) &= \theta_\bi.
\end{align}
Note that $\delta_\bi 1_\lambda = \psi^\imath(\atled_\bi 1_\lambda)$.
Now we state some corollaries to \cref{bthm}.
The first simply restates the theorem using the new notation.

\begin{cor} 
For $\bi, \bj \in \langle I \rangle$, we have that
\begin{align}
\big\langle \delta_\bi 1_\lambda, \atled_\bj 1_\lambda \big\rangle^\imath
&=
\sum_{\substack{D \in \shape(\bi\times \bj)\\\text{cup-cap-free}}}
q^{-\deg(D \ \text{\regionlabel{\lambda}})} \ \big/\  (1-q_{k_1}^{-2})\cdots (1-q_{k_l}^{-2}),\label{wet5}\\
\big\langle b_\bi 1_\lambda, \atled_\bj 1_\lambda \big\rangle^\imath
&=
\sum_{\substack{D \in \shape(\bi\times \bj)\\\text{cap-free}}}
q^{-\deg(D \ \text{\regionlabel{\lambda}})} \ \big/\  (1-q_{k_1}^{-2})\cdots (1-q_{k_l}^{-2}),\label{wet4}\\
\big\langle \delta_\bi 1_\lambda, b_\bj 1_\lambda \big\rangle^\imath
&=
\sum_{\substack{D \in \shape(\bi\times \bj)\\\text{cup-free}}}
q^{-\deg(D \ \text{\regionlabel{\lambda}})} \ \big/\  (1-q_{k_1}^{-2})\cdots (1-q_{k_l}^{-2}),\label{wet5b}\\
\big\langle b_\bi 1_\lambda, b_\bj 1_\lambda\big\rangle^\imath &= 
\sum_{D \in \shape(\bi\times \bj)}
q^{-\deg(D \ \text{\regionlabel{\lambda}})} \ \big/\  (1-q_{k_1}^{-2})\cdots (1-q_{k_l}^{-2}),\label{wet3}
\end{align}
$k_1,\dots,k_l$ being the string labels of $D$ as before.
\end{cor}

\begin{proof}
This follows from \cref{wet0,wet1,wet2} and the definitions \cref{newby1,newby2} because
$\theta_\bi$ and $b_\bi \one_\lambda$ are invariant under $\psi$ and $\psi^\imath$, respectively.
For \cref{wet5b}, we also used
$\langle x 1_\lambda, y 1_{\lambda'} \rangle^\imath
= \left\langle \psi^\imath(y 1_{\lambda'}), \psi^\imath(x 1_\lambda)\right\rangle^\imath$.
\end{proof}

In the next corollary, the involution $\tau:I \rightarrow I$ enters into the combinatorics
via the introduction of a lower finite partial order $\leq$ on the monoid $\Lambda$ from \cref{poset}
defined by
\begin{equation}\label{ordering}
\alpha \leq \beta\Leftrightarrow (\beta - \alpha) \in \sum_{i \in I} \N (\alpha_i + \alpha_{\tau i}).
\end{equation}
Note that $\alpha \leq \beta \Rightarrow \wt^\imath(\alpha)=\wt^\imath(\beta)$.
We also define
a function
\begin{align}\label{thenewfunc}
|\cdot|:\llangle I \rrangle &\rightarrow \Lambda,&
|\bi| := n_1\alpha_{i_1}+\cdots+n_l\alpha_{i_l}
\end{align}
for $\bi = i_1^{(n_1)} \cdots i_l^{(n_l)}$.

\begin{cor}\label{singapore}
Let $\bi, \bj \in \llangle I \rrangle$ and $\lambda \in X^\imath$.
If $\langle b_\bi 1_\lambda, \atled_\bj 1_\lambda\rangle^\imath \neq 0$
then we have that $|\bj| \leq |\bi|$.
\end{cor}

\begin{proof}
Since it just amounts to scaling by a non-zero scalar, we may as well
assume that $\bi, \bj \in \langle I \rangle$.
Suppose that $\langle b_\bi 1_\lambda, \atled_\bj 1_\lambda \rangle^\imath \neq 0$.
Then \cref{wet4} implies that
there exists a cap-free shape $D \in \shape(\bi\times\bj)$.
It follows that the word
$\bj$ can be obtained from $\bi$
by removing pairs of letters of the form $(h, \tau h)$ (one such pair for each cup in $D$),
then permuting the remaining letters.
Each such removal makes 
$|\bi|$ smaller in the partial order $\leq$.
\end{proof}

Finally in this subsection, we assume that we are once again given the additional data needed to define the 2-iquantum group $\UU^\imath$.

\begin{theo}\label{bcor}
Assume that $\UU^\imath$ is non-degenerate in the sense of \cref{defnondeg}.
For $\bi,\bj \in \langle I \rangle$ with $\wt^\imath(\bi)=\wt^\imath(\bj)$, we have that
\begin{equation}
\Hom_{\UU^\imath}(B_\bj \one_\lambda, B_\bi \one_\lambda)
\cong 
\left(\R^\imath\right)^{\oplus 
\overline{\langle b_\bi
1_\lambda, b_\bj 1_\lambda\rangle^\imath}}
\end{equation}
as graded right $\R^\imath$-modules
(see {\em Conventions}).
\end{theo}

\begin{proof}
Both sides are 0 unless $\wt^\imath(\bi) = \wt^\imath(\bj)$. 
When $\wt^\imath(\bi)=\wt^\imath(\bj)$, the assumed non-degeneracy means that the hom space on the left hand side
is a free right $\R^\imath$-module with basis \cref{ibasis}
obtained by pinning dots to elements of $\shape(\bi\times \bj)$.
Hence, its graded rank is the image under the bar involution of
the expression \cref{wet3}.
\end{proof}

\subsection{Application: the Balagović-Kolb-Letzter relation}\label{another}
We end the section by using the isometry $\jmath$ from \cref{jbversion} to give a conceptual explanation of the Balagović-Kolb-Letzter case of the iSerre relation.

\begin{lem}\label{nahacurrypre}
Suppose that $i \in I$ and $\lambda \in X^\imath$ with pre-image
$\hat\lambda \in X$.
For $n \geq 0$,
we have that
$$
b_{i^{(n)}} 1_\lambda = 
\begin{dcases}
\delta_{i^{(n)}} 1_\lambda&\text{if $i \neq \tau i$}\\
\sum_{m=0}^{\lfloor\frac{n}{2}\rfloor}
\frac{q_i^{m(2m-1)}}{(1-q_i^{4})
  (1-q_i^{8})\cdots(1-q_i^{4m})} \delta_{i^{(n-2m)}} 1_\lambda&\text{if $i = \tau
  i$ and $n \equiv h_i(\hat\lambda)\pmod{2}$}\\
\sum_{m=0}^{\lfloor\frac{n}{2}\rfloor}
\frac{q_i^{m(2m+1)}}{(1-q_i^{4})
  (1-q_i^{8})\cdots(1-q_i^{4m})} \delta_{i^{(n-2m)}} 1_\lambda&\text{if $i = \tau i$
  and $n \not\equiv h_i(\hat\lambda)\pmod{2}$.}
\end{dcases}
$$
\end{lem}

\begin{proof}
In the case $i = \tau i$, this follows from
\cite[Th.~2.7]{BWWa1} (and we will not use this formula subsequently). 
To prove it in the (much easier!) case that $i \neq \tau i$,
it suffices to show that
$b_{i^n} 1_\lambda = \delta_{i^n} 1_\lambda$ 
(we simply multiplied both sides by $[n]^!_{q_i}$). To prove this, we show that both sides
pair in the same way with the spanning set $\{\atled_\bj 1_\lambda\:|\:\bj \in \langle I \rangle\}$ 
for $\dot\U^\imath 1_\lambda$.
Since $i \neq \tau i$, the weight $n \alpha_i$ is minimal in the poset $\Lambda$.
From \cref{singapore}, 
it follows that the only $\bj \in \langle I \rangle$
such that 
$\left\langle b_{i^n} 1_\lambda, \atled_\bj 1_\lambda\right\rangle^\imath \neq 0$ is $\bj = i^n$. This is clearly also the case for
$\delta_{i^n} 1_\lambda$ by \cref{wet5}.
Also by \cref{wet5,wet4} we have that
$\left\langle b_{i^n} 1_\lambda, \atled_{i^n} 1_\lambda\right\rangle^\imath
= \left\langle \delta_{i^n} 1_\lambda, \atled_{i^n} 1_\lambda\right\rangle^\imath$.
\end{proof}

\begin{lem}\label{nahacurry}
Suppose that $i,j \in I$ with $i \neq \tau i$ and $i \neq j$. 
For $m \geq 1$ and $0 \leq n \leq m$, we have that
\begin{equation}
b_{i^{(n)} j i^{(m-n)}} 1_\lambda = 
\begin{dcases}
\delta_{i^{(n)} j i^{(m-n)}} 1_\lambda&\text{if $i \neq \tau j$}\\
\delta_{i^{(n)} j i^{(m-n)}} 1_\lambda
+
f^\lambda_{n,m;i}(q)\, \delta_{i^{(m-1)}} 1_\lambda&\text{if $i = \tau j$}
\end{dcases}
\end{equation}
where 
\begin{align}\label{veday}
f^\lambda_{n,m;i}(q) := 
\frac{q_i^{1+\lambda_i - \del_i-(m-n-1)(1-a_{i,\tau i})-m}\qbinom{m-1}{n-1}_{q_i}}
{1-q_i^2}
+
\frac{
q_i^{1+(m-n-1)(1-a_{i,\tau i})+\del_i-\lambda_i}\qbinom{m-1}{n}_{q_i}}{1-q_i^2}.
\end{align}
\end{lem}

\begin{proof}
The strategy is the same as the proof of \cref{nahacurrypre} just explained.
The weight $m \alpha_i +\alpha_j$ is minimal in the poset $\Lambda$ if $i \neq \tau j$, and if $i = \tau j$ then 
the only $\beta \in \Lambda$ with
$\beta < m+ \alpha_i + \alpha_j$ is $\beta = (m-1) \alpha_i$.
Also \cref{wet5,wet4} imply that
$\left\langle b_{i^n j i^{m-n}} 1_\lambda, \atled_\bj 1_\lambda \right\rangle
= \left\langle \delta_{i^n j i^{m-n}} 1_\lambda, \atled_\bj 1_\lambda\right\rangle$
for any $\bj \in \langle I \rangle$ with $|\bj| = m \alpha_i + \alpha_j$.
It follows that
$$
b_{i^n j i^{m-n}} 1_\lambda
= 
\begin{cases}
\delta_{i^n j i^{m-n}} 1_\lambda &\text{if $i \neq \tau j$}\\
\delta_{i^n j i^{m-n}} + f(q) \delta_{i^{m-1}} 1_\lambda&\text{if $i = \tau j$}
\end{cases}
$$
for some $f(q) \in \Q(q)$.
This already proves the lemma in the case $i \neq \tau j$.
When $i = \tau j$, we still need to compute $f(q)$. 
Applying $\langle\cdot,\atled_{i^{m-1}} 1_\lambda\rangle^\imath$ to the equation gives that
$$
f(q^{-1})=
\frac{\left\langle
b_{i^n j i^{m-n}} 1_\lambda, \atled_{i^{m-1}} 1_\lambda\right\rangle^\imath}{\left\langle\delta_{i^{m-1}} 1_\lambda,\atled_{i^{m-1}} 1_\lambda\right\rangle^\imath}.
$$
In view of \cref{wet5,wet4}, it follows that $f(q)=\sum_D q^{\deg(D \ \text{\regionlabel{\lambda}})}$ summing over 
\begin{align}\label{somediagrams}
D \in \Biggl\{\:\:
\begin{tikzpicture}[iQ,baseline=-6mm]
\draw (0,0)\toplabel{j} to [out=-90,in=-90,looseness=2] (.8,0) \toplabel{i};
\draw (-.2,0) \toplabel{i} to (-.1,-1);
\draw (-.6,0) \toplabel{i} to (-.5,-1);
\draw (.2,0) \toplabel{i} to (0.2,-1);
\draw (.6,0) \toplabel{i} to (.6,-1);
\draw (1,0) \toplabel{i} to (.9,-1);
\draw (1.4,0) \toplabel{i} to (1.3,-1);
\node at (1.05,-.6) {$\scriptstyle\cdot$};
\node at (1.15,-.6) {$\scriptstyle\cdot$};
\node at (1.25,-.6) {$\scriptstyle\cdot$};
\node at (-.25,-.6) {$\scriptstyle\cdot$};
\node at (-.35,-.6) {$\scriptstyle\cdot$};
\node at (-.45,-.6) {$\scriptstyle\cdot$};
\node at (.3,-.6) {$\scriptstyle\cdot$};
\node at (.4,-.6) {$\scriptstyle\cdot$};
\node at (.5,-.6) {$\scriptstyle\cdot$};
\draw [thick,decoration={brace,raise=0.4cm},decorate] (.95,-0.13) -- (1.45,-0.13);
\node at (1.35,.5) {$\scriptscriptstyle m\!-\!n\!-\!1\!-\!k$};
\draw [thick,decoration={brace,raise=0.4cm},decorate] (-.7,-0.13) -- (-.2,-0.13);
\node at (-.45,.5) {$\scriptscriptstyle n$};
\draw [thick,decoration={brace,raise=0.4cm},decorate] (0.2,-0.13) -- (0.7,-0.13);
\node at (.45,.5) {$\scriptscriptstyle k$};
\multcloseddot{.1,-.32}{east}{r\!};
\end{tikzpicture},
\begin{tikzpicture}[iQ,baseline=-6mm]
\draw (0,0)\toplabel{j} to [out=-90,in=-90,looseness=2] (-.8,0) \toplabel{i};
\draw (.2,0) \toplabel{i} to (.1,-1);
\draw (.6,0) \toplabel{i} to (.5,-1);
\draw (-.2,0) \toplabel{i} to (-.2,-1);
\draw (-.6,0) \toplabel{i} to (-.6,-1);
\draw (-1,0) \toplabel{i} to (-.9,-1);
\draw (-1.4,0) \toplabel{i} to (-1.3,-1);
\node at (-1.05,-.6) {$\scriptstyle\cdot$};
\node at (-1.15,-.6) {$\scriptstyle\cdot$};
\node at (-1.25,-.6) {$\scriptstyle\cdot$};
\node at (.24,-.6) {$\scriptstyle\cdot$};
\node at (.34,-.6) {$\scriptstyle\cdot$};
\node at (.44,-.6) {$\scriptstyle\cdot$};
\node at (-.3,-.6) {$\scriptstyle\cdot$};
\node at (-.4,-.6) {$\scriptstyle\cdot$};
\node at (-.5,-.6) {$\scriptstyle\cdot$};
\draw [thick,decoration={brace,mirror,raise=0.4cm},decorate] (-.95,-0.13) -- (-1.45,-0.13);
\node at (-1.35,.5) {$\scriptscriptstyle n\!-\!1\!-\!l$};
\draw [thick,decoration={brace,mirror,raise=0.4cm},decorate] (.7,-0.13) -- (.2,-0.13);
\node at (.45,.5) {$\scriptscriptstyle m\!-\!n$};
\draw [thick,decoration={brace,mirror,raise=0.4cm},decorate] (-0.2,-0.13) -- (-0.7,-0.13);
\node at (-.45,.5) {$\scriptscriptstyle l$};
\multcloseddot{-.1,-.32}{west}{\!r};
\end{tikzpicture}\:\:\Bigg|\:\:
\text{$0 \leq k \leq m-n-1$, $0 \leq l \leq n-1$, $r \geq 0$}
\Biggr\}.
\end{align}
We deduce that
\begin{align}
f(q) &= 
\sum_{k=0}^{m-n-1}
\frac{q_i^{1+\del_i - (\lambda-(m-n-1)\alpha_i)_i-2k}}{1-q_i^{2}}
+\sum_{l=0}^{n-1}
\frac{q_i^{1+\del_j - (\lambda-(m-n)\alpha_i)_j-2l}}{1-q_i^{2}}\notag\\
\label{Zset}
&=
\frac{q_i^{1+(m-n-1)(1-a_{i,j})+\del_i-\lambda_i}[m-n]_{q_i}}
{1-q_i^2}
+
\frac{q_i^{1+\lambda_i-\del_i -(m-n-1)(1-a_{i,j})-m}[n]_{q_i}}{1-q_i^2}.
\end{align}
To deduce the formula for $f^\lambda_{n,m;i}(q)$
in the statement,
multiply by $[m-1]^!_{q_i} / [m-n]^!_{q_i} [n]^!_{q_i}$.
\end{proof}

Now we can give a new proof of
the relation \cref{secretary} in the difficult case $i = \tau j$.
Assume that $i \neq j=\tau i$ and set $m := 1-a_{i,j}$.
Applying the inverse of the isomorphism $\widetilde\jmath$ to
the Serre relation in $\f$ gives that
$\sum_{n=0}^{m} (-1)^n \delta_{i^{(n)} j i^{(m-n)}} 1_\lambda =
0$. 
Since $\delta_{i^{(n)} j i^{(m-n)}} 1_\lambda =
b_{i^{(n)} j i^{(m-n)}} 1_\lambda -
f^\lambda_{n,m;i}(q) b_{i^{(m-1)}} 1_\lambda$ by \cref{nahacurrypre,nahacurry},
we deduce that
\begin{equation*}
\sum_{n=0}^{m} (-1)^n b_i^{(n)} b_j b_i^{(m-n)}1_\lambda=
\sum_{n=0}^{m} (-1)^n f^\lambda_{n,m;i}(q) b_i^{(m-1)}1_\lambda.
\end{equation*}
The iSerre relation \cref{secretary} for $i = \tau j$, which was proved originally in \cite[Th. 3.6]{BK}, follows from this using the next lemma to simplify the right hand side.

\begin{lem}\label{nextlem}
Let $m := 1-a_{i,\tau i}$ for $i \in I$ with $i \neq \tau i$.
We have that
$$
\sum_{n=0}^{m} (-1)^n f^\lambda_{n,m;i}(q) = 
\prod_{r=1}^{m-1}
(q_i^r-q_i^{-r})\times
\frac{(-1)^{m-1} q_i^{\lambda_i-\del_i-\binom{m}{2}}-q_i^{
\binom{m}{2}+\del_i-\lambda_i}}{q_i-q_i^{-1}}.
$$
\end{lem}

\begin{proof}
We replace each $f^\lambda_{n,m;i}(q)$ by the sum of two fractions that is its definition \cref{veday} to rewrite the left hand side of the identity to be proved as the sum of two summations. Then we reindex the first summation to see that it is equal to
\begin{align*}
\frac{q_i^{\lambda_i-\del_i-\binom{m}{2}}}{q_i-q_i^{-1}}
\left(\sum_{n=0}^{m-1} (-1)^n
q_i^{nm-\binom{m}{2}}\qbinom{m-1}{n}_{q_i}\right)
-\frac{q_i^{\binom{m}{2}+\del_i-\lambda_i}}{q_i-q_i^{-1}}
\left(\sum_{n=0}^{m-1}(-1)^n
q_i^{\binom{m}{2}-nm}\qbinom{m-1}{n}_{q_i}\right).
\end{align*}
To see that this expression is equal to the desired right hand side, it remains to apply 
the elementary identity
\begin{equation}
\displaystyle\prod_{r=1}^{m-1} (q^r-q^{-r})=
\sum_{n=0}^{m-1} (-1)^n q^{\binom{m}{2}-nm}\qbinom{m-1}{n}_q
=(-1)^{m-1}
\sum_{n=0}^{m-1} (-1)^n q^{nm-\binom{m}{2}}\qbinom{m-1}{n}_q,
\end{equation}
which follows from \cite[1.3.1(c)]{Lubook} taking $z = -v^2$.
\end{proof}

%% file: s6-categorification.tex
\setcounter{section}{5}

\section{Identification of the Grothendieck ring}\label{s6-categorification}

Assuming the non-degeneracy of the 2-quantum group $\UU$ which we (re-)proved in \cref{nd},
the arguments in
\cite[Sec.~3.8, Sec.~3.9]{KL3} show that 
the Grothendieck ring of the appropriate completion of $\UU$ is isomorphic to $\dot\U_\Z$, i.e.,
$\UU$ categorifies $\dot\U_\Z$. The main goal in this section is to
prove the analogous statement for all quasi-split iquantum
groups which are non-degenerate 
 in the sense of \cref{defnondeg}. This includes all of the cases with symmetric Cartan matrix and geometric parameters thanks to \cref{ind}.
We assume throughout the section that $\kk_0$ is a field, necessarily
of characteristic $\neq 2$ if $a_{i,\tau i} \neq 0$ for any $i \in I$,
and that the 2-iquantum group
$\UU^\imath$ is non-degenerate.
We let
\begin{equation}\label{alone}
\R^\imath := \kk\big[\beta_i^{(m)}, \beta_j^{(n)}\:\big|\:i \in I_1, j \in I_0, m,n \geq 1, n\text{ odd}\big]
\end{equation}
as in \cref{defnondeg}.
Due to the assumed non-degeneracy of $\UU^\imath$, we have the isomorphism
$f^\imath_\lambda:\R^\imath
\stackrel{\sim}{\rightarrow}
\End_{\UU^\imath}(\one_\lambda)$ for any $\lambda \in X^\imath$. 

\subsection{Idempotents and divided/idivided powers}\label{ssidempotents}

We begin by constructing a Grothendieck ring from $\UU^\imath$
which, like $\dot\U^\imath_\Z$, is
a locally unital $\Z[q,q^{-1}]$-algebra.
Starting from $\UU^\imath$, 
we pass first to the $q$-envelope
$\UU^\imath_q$ defined at the start of \cref{newgrading}. Then we pass from there to the underlying $\kk$-linear 
2-category consisting of the same objects and 1-morphisms as $\UU^\imath_q$ but taking only the 2-morphisms that are homogeneous of degree 0. Finally we take the additive Karoubi envelope ($=$ idempotent completion of additive envelope) of each of the morphism categories of that 
to obtain a $\kk$-linear 2-category denoted $\Kar(\UU^\imath_q)$.
Then we define
\begin{equation}\label{groth}
K_0\left(\Kar(\UU^\imath_q)\right) := \bigoplus_{\kappa,\lambda \in X^\imath}
K_0\left(\one_\kappa\Kar(\UU^\imath_q)\one_\lambda\right)
\end{equation}
with
$K_0\left(\one_\kappa \Kar(\UU^\imath_q)\one_\lambda\right)$
being the usual split Grothendieck group of
$\HOM_{\Kar(\UU^\imath_q)}(\lambda,\kappa)$.
We make
$K_0\left(\Kar(\UU^\imath_q)\right)$ into a $\Z[q,q^{-1}]$-algebra
with multiplication induced by horizontal composition, and the action of $q$ defined by $q [B_\bi \one_\lambda] := 
[q B_\bi \one_\lambda]$ for $\bi \in \langle I \rangle$.
It is a locally unital algebra
with distinguished idempotents given by the isomorphism classes
$1_\lambda := [\one_\lambda]$ for $\lambda \in X^\imath$.

Fixing $i \in I$ and $\lambda \in X^\imath$ for a while, 
let $b_i 1_\lambda = 1_{\lambda-\alpha_i} b_i \in K_0\left(\Kar(\UU^\imath_q)\right)$ be 
the 2-isomorphism class of the 1-morphism $B_i \one_\lambda$.
There are also {\em divided powers} $B_i^{(n)} 1_\lambda = 1_{\lambda-n\alpha_i} B_i^{(n)}$ if $i \neq \tau i$ or {\em
idivided powers} $B_i^{(n)} 1_\lambda = 1_{\lambda-n\alpha_i} B_i^{(n)}$
if $i = \tau i$. 
The distinction according to whether $i \neq \tau i$ or $i = \tau i$ is important---in the former case the defining relations imply that there is a homomorphism from the usual
nil-Hecke algebra $\NH_n$ (graded so that the polynomial generators are all of degree $2 d_i$)
to $\End_{\UU^\imath}(B_i^n \one_\lambda)$
whereas in the latter case one needs to work with the nil-Brauer category as in \cref{aging}.
In either case, we let $1_{i^{(n)}} \in \End_{\UU^\imath}(B_i^n \one_\lambda)$ be the 2-endomorphism defined by following the pattern:
\begin{align*}
1_{\varnothing} &:= \id_{\one_\lambda},
&1_{i}
&:= 
\begin{tikzpicture}[anchorbase,scale=1.4,iQ]
\draw (0,-.4)\botlabel{i}  to (0,.4);
\end{tikzpicture}\:,&
1_{i^{(2)}} &:=
\begin{tikzpicture}[anchorbase,scale=1.4,iQ]
\draw (0.3,-.4)\botlabel{i}  to (-.3,.4);
\draw (-.3,-.4) \botlabel{i} to (0.3,.4);
\closeddot{-.15,.2};
\end{tikzpicture}\:,&
1_{i^{(3)}} &:=
\begin{tikzpicture}[anchorbase,scale=1.4,iQ]
\draw (0.3,-.4)\botlabel{i}  to (-.3,.4);
\draw (-.3,-.4) \botlabel{i} to (0.3,.4);
\draw (0,.4) to[out=-45,in=90] (.3,0) to[out=-90,in=45] (0,-.4)\botlabel{i} ;
\closeddot{.08,.33};
\closeddot{-.25,.33};
\closeddot{-.14,.19};
\end{tikzpicture}\:\;,&
1_{i^{(4)}} &:=
\begin{tikzpicture}[anchorbase,scale=1.4,iQ]
\draw (0.3,-.4)\botlabel{i} to (-.3,.4);
\draw (-.3,-.4)\botlabel{i}  to (0.3,.4);
\draw (.1,.4) to[out=-50,in=90,looseness=.75] (.3,.15) to[out=-90,in=35,looseness=.75] (-.1,-.4)\botlabel{i} ;
\draw (-.1,.4) to[out=-35,in=90,looseness=.75] (.3,-.15) to[out=-90,in=50,looseness=.75] (.1,-.4)\botlabel{i} ;
\closeddot{-.26,.345};
\closeddot{-0.03,.345};
\closeddot{.15,.345};
\closeddot{.07,.25};
\closeddot{-.19,.25};
\closeddot{-.12,.16};
\end{tikzpicture}\:\;,\qquad\cdots
\end{align*}
We will use
the ``thick calculus'' notation developed in \cite[Sec.~4]{BWWa1} to denote this in general:
\begin{equation}\label{endef2}
1_{i^{(n)}} =
\begin{tikzpicture}[iQ,centerzero,scale=1.1]
\draw[ultra thick] (0,.4) to (0,-.4)\botlabel{i^n};
\multcloseddot{0,0.15}{east}{\rho};
\cross{0,-.1};
\end{tikzpicture}
\end{equation}
The thick string labelled $i^n$ here is a shorthand for $n$ parallel thin strings each of color $i$.
If $i \neq \tau i$, it is well known that $1_{i^{(n)}}$ is the image of a primitive homogeneous idempotent in $\NH_n$. When $i = \tau i$, it is the image of 
a primitive homogeneous idempotent in the nil-Brauer category; see \cite[Cor.~4.24]{BWWa1}.
The essential property which implies that this is an idempotent is that
\begin{equation}\label{essence}
\begin{tikzpicture}[iQ,centerzero]
\draw[ultra thick] (0,.4) to (0,-.4)\botlabel{i^n};
\multcloseddot{0,0}{east}{\rho};
\cross{0,-.22};
\cross{0,.22};
\end{tikzpicture}
=\begin{tikzpicture}[iQ,centerzero]
\draw[ultra thick] (0,.4) to (0,-.4)\botlabel{i^n};
\cross{0,0};
\end{tikzpicture}\ .
\end{equation}
Applying the symmetry $\Psi^\imath$ from \cref{eastchinasea}, we obtain another idempotent
\begin{equation}\label{florence}
{'}1_{i^{(n)}} := \Psi^\imath\left(1_{i^{(n)}}\right) = \begin{tikzpicture}[iQ,centerzero,scale=1.1]
\draw[ultra thick] (0,.4) to (0,-.4)\botlabel{i^n};
\multcloseddot{0,-0.15}{east}{\rho};
\cross{0,.1};
\end{tikzpicture}\ .
\end{equation}
The idempotents $1_{i^{(n)}}$ and ${'}1_{i^{(n)}}$ are conjugate:
we have that $1_{i^{(n)}} = u_{i^{(n)}} \circ v_{i^{(n)}}$
and ${'}1_{i^{(n)}} = v_{i^{(n)}} \circ u_{i^{(n)}}$ where $u_{i^{(n)}}$ consists of the dots at the top and $v_{i^{(n)}}$ consists of the crossings at the bottom of \cref{endef2}.
We need both since $\Psi^\imath$ will be used in some subsequent constructions.
Let 
\begin{align}\label{colderthanitshouldbe}
B_i^{(n)} \one_\lambda &:= 
\left(q_i^{-\binom{n}{2}} B_i^n \one_\lambda, 1_{i^{(n)}}\right)\cong
\left(q_i^{\binom{n}{2}} B_i^n \one_\lambda, {'}1_{i^{(n)}}\right).
\end{align}
With this definition, the notation \cref{wordnotation2} now makes sense.
If $\bi = i_1^{(n_1)} \cdots i_l^{(n_l)} \in \llangle I \rrangle$ 
then
\begin{align}\label{diving}
B_\bi \one_\lambda =
(q^{-\deg(\bi)} B_{i_1}^{n_1} \cdots
B_{i_l}^{n_l} \one_\lambda, 1_\bi)\cong
(q^{\deg(\bi)} B_{i_1}^{n_1} \cdots
B_{i_l}^{n_l}\one_\lambda, {'}1_\bi)
\end{align}
where
\begin{align}\label{snorkeling}
1_\bi &:= 1_{i_1^{(n_1)}} \cdots 1_{i_l^{(n_l)}},&
{'}1_\bi &:= {'}1_{i_1^{(n_1)}} \cdots {'}1_{i_l^{(n_l)}},&
\deg(\bi) &:= d_{i_1} \binom{n_1}{2}+\cdots+d_{i_l}\binom{n_l}{2}.
\end{align}

\begin{lem}\label{rentalcar}
The elements $[B_i^{(n)} \one_\lambda]$ of $K_0\left(\Kar(\UU^\imath_q)\right)$
satisfy 
\begin{equation}\label{iDividedpowerrelation}
[B_i^{(n)} \one_\lambda] =
\begin{dcases}
\frac{1}{[n]_{q_i}^!}
\left[B_i^n \one_\lambda\right]&\text{if $i \neq \tau i$}\\
\frac{1}{[n]_{q_i}^!}
\left[\prod_{\substack{m=0\\m\equiv h_i(\hat\lambda)\!\!\!\!\!\pmod{2}}}^{n-1}
\!\!\!\left(B_i^2 - [m]_{q_i}^2\right)\one_\lambda\right] 
&\text{if $i = \tau i$ and $n$ is even}\\
\frac{1}{[n]_{q_i}^!}
\left[B_i\!\!\!\!\!\!\!
\prod_{\substack{m=1\\m\equiv h_i(\hat\lambda)\!\!\!\!\!\pmod{2}}}^{n-1}
\!\!\!\left(B_i^2 - [m]_{q_i}^2\right)\one_\lambda\right]
&\text{if $i = \tau i$ and $n$ is odd.}\\
\end{dcases}
\end{equation}
This is the same as the identity satisfied by the elements 
$b_i^{(n)} 1_\lambda$ of $\dot\U^\imath_\Z$ from \cref{idividedpowerrelation}.
\end{lem}

\begin{proof}
When $i \neq \tau i$, this follows in the same way as the analogous statement for 2-quantum groups; e.g., see \cite[Lem.~4.1]{Rou} (which ignores the grading) or \cite[(3.54)--(3.55)]{KL3}. Ultimately, it 
depends on a well-known result about primitive idempotents in the nil-Hecke algebra $\NH_n$.
When $i = \tau i$, the result follows from \cite[Th.~4.21]{BWWa1},
which shows that $B_i B^{(n)}_i \one_\lambda$
is isomorphic to a direct sum of graded shifts of $B^{(n+1)}_i \one_\lambda$ and $B^{(n-1)}_i \one_\lambda$ 
in a way that agrees with the recursive formula for $b_i b_i^{(n)} 1_\lambda$ from \cite{BeW18}.
Although \cite{BWWa1} works over an algebraically closed field,
the decomposition of $\id_{B_i^n \one_\lambda}$ into mutually
orthogonal idempotents and the conjugacy of these idempotents
established there is valid in our current setup.
\end{proof}

\subsection{Classification of indecomposables and standard modules}\label{ssgtb}

We next explain how to
classify indecomposable objects in
$\Kar(\UU^\imath_q)$. For this, we are going to switch to using 
the language of modules. 
For $\lambda \in X^\imath$, we let 
\begin{equation}
\UU^\imath \one_\lambda := \coprod_{\kappa \in X^\imath} \one_\kappa \UU^\imath \one_\lambda
\end{equation}
with $q$-envelope $\UU^\imath_q \one_\lambda$.
Recall that ${'}\UU^\imath$ is the 2-iquantum group defined using the parameters ${'}Q_{i,j}(x,y) = r_{i,j}r_{j,i}Q_{i,j}(x,y)$.
The existence of the isomorphism $\Psi^\imath:(\UU^\imath)^\op
\stackrel{\sim}{\rightarrow} {'}\UU^\imath$ from
\cref{eastchinasea}
implies that ${'}\UU^\imath$ is non-degenerate.
Let
\begin{align}\label{pathalgebra}
\HH
:= \bigoplus_{\bi,\bj \in \langle I \rangle} 
\Hom_{{'}\UU^\imath}(B_\bj
  \one_\lambda, B_\bi \one_\lambda),
\end{align}
which is the path algebra of the graded category
${'}\UU^\imath \one_\lambda$.
We keep $\lambda$ fixed for the remainder of the subsection, so will drop the superscript, denoting $\HH$
simply by $\H$. This means that most of our subsequent
notation related to $\H$ also depends implicitly on $\lambda$ (one can always add a superscript $\lambda$ when necessary to avoid ambiguity).

Note that $\H$ is naturally a graded $\R^\imath$-algebra,
with $\beta \in \R^\imath$ acting by horizontal composition on the right with $f_\lambda^\imath(\beta) \in \End_{{^\prime}\UU^\imath}(\one_\lambda)$ (recall \cref{extraordinary}).
For $\bi =i_1\cdots i_l \in \langle I \rangle$, the idempotent $1_\bi \in \H$
is the identity endomorphism of $B_\bi \one_\lambda$,
and $\H$ is a locally unital algebra
with these as its distinguished idempotents, that is, we have that
$$
\H = \bigoplus_{\bi, \bj \in \langle I \rangle} 1_\bi
\H 1_\bj.
$$
An important point is that
each of the spaces
$1_\bi \H 1_\bj$ is locally finite-dimensional and
bounded below as a vector space over $\kk_0$, that is, 
its graded components are all finite-dimensional and they are zero in
sufficiently negative degrees. This follows from \cref{easyspanning}.
Thus, $\H$ is
a graded algebra which is locally finite-dimensional and bounded below.

Let $\gmod{\H}$ be the category of graded left $\H$-modules
$V = \bigoplus_{\bi \in \langle I \rangle} 1_\bi V$.
As a matter of notation, we denote the morphism space
between two graded left $\H$-modules 
in the category $\gmod{\H}$
by $\Hom_\H(U,V)_0$, reserving
$\Hom_\H(U,V)$ for the graded hom $\bigoplus_{n \in \Z} \Hom_\H(U,V)_n$
where $\Hom_\H(U,V)_n = \Hom_\H(U, q^{-n} V)_0$ consists of the
$\H$-module homomorphisms which map $U_i$ into $V_{i+n}$ for all $i \in \Z$.
There are also the categories $\gproj{\H}$ and $\ginj{\H}$
of finitely generated graded projective left $\H$-modules 
and finitely cogenerated graded injective 
left $\H$-modules, respectively.

In the definition of $\H$, we used
the 2-iquantum group
${'}\UU^\imath$ rather than the usual
$\UU^\imath$. The reason for this peculiarity 
is so that we can use
the isomorphism of graded 2-categories
$\Psi^\imath:\big(\UU^\imath\big)^\op\stackrel{\sim}{\rightarrow}
{'}\UU^\imath$ from \cref{eastchinasea}. Recalling \cref{awkwardinverse}, it induces an isomorphism
$\Kar(\UU^\imath_q \one_\lambda)\stackrel{\sim}{\rightarrow} \Kar\big({'}\UU^\imath_{q^{-1}}\one_\lambda\big)^\op$.
But also Yoneda gives an equivalence
$\contra:\Kar\big({'}\UU^\imath_{q^{-1}} \one_\lambda\big)^\op \rightarrow \gproj{\H}$.
Composing this with the isomorphism induced by $\Psi^\imath$, we obtain an equivalence
we denote simply by
\begin{equation}\label{yonedaop}
\cov:
\Kar\left(\UU^\imath_q \one_\lambda\right)
\rightarrow \gproj{\H}
\end{equation}
such that $\cov\circ q = q \circ \cov$. 
Using \cref{diving,snorkeling},
for $\bi \in \llangle I \rrangle$,
we have that
\begin{equation}\label{projy}
P(\bi) := q^{\deg(\bi)} \H\, 1_\bi
\cong q^{-\deg(\bi)} \H\, {'}1_\bi =
\cov(B_\bi \one_\lambda).
\end{equation}
There has been an application of $\Psi^\imath$ here to switch $1_\bi$ and ${'}1_\bi$.

For a graded right $\H$-module,
let $V^\circledast$ be $\bigoplus_{\bi \in \langle I \rangle} 
(V 1_\bi)^*$ (linear dual over the field $\kk_0$), which is 
a graded left $\H$-module. If $V$ is locally finite-dimensional and bounded below (resp., finitely generated and projective)
then $V^\circledast$ is locally finite-dimensional and bounded above
(resp., finitely cogenerated and injective). This duality functor 
is needed to define the {\em Nakayama functor}
\begin{equation}\label{nakayama}
\Nak := \Hom_\H(-,\H)^\circledast:\gproj{\H} \rightarrow \ginj{\H},
\end{equation}
which is an equivalence of categories commuting with the grading shift functor.
For $\bi \in \llangle I \rrangle$, it takes $P(\bi)$ to 
\begin{equation}\label{injy}
I(\bi) := q^{\deg(\bi)} (1_\bi \H)^\circledast,
\end{equation}
this being a finitely cogenerated injective graded left $\H$-module.

Now we are going to study $\gmod{\H}$ by appealing to
the general theory of algebras with {\em graded
triangular bases} from \cite{GTB}, taking the ground field to be
$\kk_0$. 
We need to equip $\H$ with a graded triangular basis
consisting of the products $xhy$ for
$x \in X(\bi,\bj), h \in H(\bj,\bk)$ and $y \in Y(\bk,\bl)$
and $\bi,\bj,\bk,\bl \in \langle I \rangle$.
The sets $X(\bi,\bj), H(\bj,\bk)$ and $Y(\bk,\bl)$ are defined 
in the next paragraph.
Referring to \cite[Def.~1.1]{GTB} for the other language being used,
all of the distinguished idempotents $1_\bi\:(\bi \in \langle I \rangle)$ are special.
The weight poset required in the definition of graded triangular basis
is the poset $\Lambda$ from \cref{poset} ordered by the
partial order $\leq$ from \cref{ordering}.
The required function $\langle I \rangle\rightarrow \Lambda$
is the function $\bi \mapsto |\bi|$ defined by \cref{thenewfunc}.
The fibers 
$\langle I \rangle_\alpha := \left\{\bi \in \langle I \rangle\:\big|\:|\bi|=\alpha\right\}$
of this function are finite, 
so it makes sense to define
\begin{equation}
e_\alpha := \sum_{\bi \in \langle I \rangle_\alpha}
1_\bi
\end{equation}
for $\alpha \in \Lambda$.
We refer to these as {\em weight idempotents}.
Any left $\H$-module $V$ decomposes as $V = \bigoplus_{\alpha \in \Lambda} e_\alpha V$. We refer to the subspace $e_\alpha V$ as the {\em $\alpha$-weight space} of $V$.

Here, we define $X(\bi,\bj), H(\bj,\bk)$ and $Y(\bk,\bl)$.
For $\bi,\bj \in \langle I \rangle$,
we choose the set $\shape(\bi\times\bj)$ of representatives for isotopy classes of reduced $\bi\times \bj$ shapes so that all cups in such a shape 
are in the top third of the diagram, all crossings of
propagating strings are in the middle third of the diagram, and all caps
are in the bottom third. 
We also fix a choice of basis for the algebra 
$\R^\imath$ as a vector space over $\kk_0$.
Then:
\begin{itemize}
\item The set $X(\bi,\bj)$ consists of all $D \in \shape(\bi,\bj)$
which only
involve cups (no caps, no crossings of propagating strings or
bubbles) with some number of dots added at the distinguished
points on each cup.
\item The set $H(\bj,\bk)$ consists of all $D \in \shape(\bj,\bk)$
which only involve propagating strings (no cups or caps)
with some number of dots added at the distinguished points on
these strings and a bubble labelled by one of the basis vectors for $\R^\imath$
 on the right hand boundary.
\item The set $Y(\bk,\bl)$ consists of all $D \in \shape(\bk,\bl)$
which only
involve caps (no cups, no crossings of propagating strings or
bubbles) with some number of dots added at the distinguished
points on each cap.
\end{itemize}
The assumption of non-degeneracy implies that the monomials $xhy$ indexed by these sets do
indeed give a basis for $\H$ as a vector space over $\kk_0$.
The following picture illustrates a typical element of the resulting graded triangular basis,
for strings colored by $\bi$ at the top boundary and by $\bl$
at the bottom boundary, some non-negative multicities labelling the
dots, and a basis vector $f$ of $\R^\imath$:
$$
\begin{tikzpicture}[anchorbase,scale=.66]
\draw (-1.5,-3)  to [out=90,in=-90] (0.5,3);
\draw (0,-3) to [out=90,in=-90] (-1.5,3);
\draw (1.5,-3) to [out=90,in=-70]
(0,0.5) to [out=110,in=-90] (-0.5,3);
\draw (-2,-3) to [out=90,in=-90] (-2.5,3);
\draw (-3,3)  to [out=-90,in=180] (0,1.25) to [out=0,in=-90] (3.5,3);
\draw (-2,3) to [out=-90,in=180] (0.5,2) to [out=0,in=-90] (2,3);
\draw (-1,3) to [out=-90,in=180] (-.5,2.5) to [out=0,in=-90] (0,3);
\draw (1.55,3) to [out=-90,in=180] (2,2.5) to [out=0,in=-90] (2.5,3);
\draw (.95,3) to [out=-90,in=180] (2,2) to [out=0,in=-90] (3,3);
\draw (1,-3) to [out=90,in=180] (1.5,-2.5) to [out=0,in=90] (2,-3);
\draw (.5,-3) to [out=90,in=180] (1.75,-2) to [out=0,in=90] (3,-3);
\draw (-2.5,-3) to [out=90,in=180] (.5,-1.5) to [out=0,in=90] (2.5,-3);
\draw (-1,-3) to [out=90,in=180] (-.75,-2.6) to [out=0,in=90] (-0.5,-3);
\draw[dashed,red] (-3.2,1) to (3.5,1);
\draw[dashed,red] (-3.2,-1) to (3.5,-1);
\draw[-] (2.85,.1) arc(180:-180:0.25);
\dottybubblelabel{3.1,.1}{f};
\closeddot{-2.98,2.75};
\closeddot{-.94,2.75};
\closeddot{2.45,2.75};
\closeddot{1,2.05};
\closeddot{-2.25,0};
\closeddot{-.96,.5};
\closeddot{-.88,-.7};
\closeddot{-.97,-2.8};
\closeddot{2.99,-2.8};
\closeddot{1.97,-2.8};
\closeddot{2.86,2.5};
\closeddot{.3,-0.15};
\closeddot{2.38,-2.45};
\node[rectangle,rounded corners,draw,fill=red!15!white,inner sep=4pt] at (-5.5,0) {$h$};
\node[rectangle,rounded corners,draw,fill=red!15!white,inner sep=4.3pt] at (-5.5,2.1) {$x$};
\node[rectangle,rounded corners,draw,fill=red!15!white,inner sep=4pt] at (-5.5,-2.1) {$y$};
\end{tikzpicture}
$$

The next step in the development from \cite{GTB}
is to introduce
the {\em Cartan algebras}, which are defined to be
the subquotients
\begin{equation}
\H_\alpha := \bar e_\alpha \H_{\geq \alpha} \bar e_\alpha
\end{equation}
for $\alpha \in \Lambda$,
where
$\H_{\geq \alpha}$ is quotient of $\H$ by the 2-sided ideal generated
by $\{e_\beta\:|\:\beta \in \Lambda \,\text{with}\,\beta
\not\geq \alpha\}$, and $\bar x$ denotes the
canonical image of an element $x \in \H$ in $\H_{\geq \alpha}$.
Note this is a unital graded algebra, unlike $\H$ itself which is merely locally unital.
We have that
$$
\H_\alpha = \bigoplus_{\bi,\bj \in \langle I \rangle_\alpha}
\bar 1_\bi \H_\alpha \bar 1_\bj
$$
with $\bar 1_\bi \H_\alpha \bar 1_\bj$ having basis $\left\{\bar h\:|\:h \in
H(\bi,\bj)\right\}$.
In fact, as proved in the subsequent lemma, the Cartan algebra $\H_\alpha$ is a quiver Hecke algebra over the ground ring
$\R^\imath$:

\begin{defin}[iQuiver Hecke algebras]\label{iQHA}
We define the {\em iquiver Hecke category}
$\catQH^\imath$ to be the strict $\R^\imath$-linear graded monoidal category
with generating objects $\Theta_i\:(i \in I)$
and generating morphisms
\begin{align}
\begin{tikzpicture}[iQ,centerzero,scale=.7]
\draw[-] (0,-0.3) \botlabel{i} -- (0,0.3);
\closeddot{0,0};
\end{tikzpicture}&:\Theta_i \rightarrow \Theta_i,&
\begin{tikzpicture}[iQ,centerzero,scale=.7]
\draw[-] (-0.3,-0.3) \botlabel{i} -- (0.3,0.3);
\draw[-] (0.3,-0.3) \botlabel{j} -- (-0.3,0.3);
\end{tikzpicture}&: \Theta_i \otimes \Theta_j \rightarrow \Theta_j \otimes \Theta_i
\end{align}
of degrees $2 d_i$ and $-d_i a_{i,j}$, respectively,
subject to the relations
\begin{align}
\begin{tikzpicture}[iQ,centerzero]
\draw[-] (-0.3,-0.3) \botlabel{i} -- (0.3,0.3);
\draw[-] (0.3,-0.3) \botlabel{j} -- (-0.3,0.3);
\closeddot{-0.15,-0.15};
\end{tikzpicture}
-
\begin{tikzpicture}[iQ,centerzero]
\draw[-] (-0.3,-0.3) \botlabel{i} -- (0.3,0.3);
\draw[-] (0.3,-0.3) \botlabel{j} -- (-0.3,0.3);
\closeddot{0.15,0.15};
\end{tikzpicture}
&= \delta_{i,j}  \ 
\begin{tikzpicture}[iQ,centerzero]
\draw[-] (-0.2,-0.3) \botlabel{i} -- (-0.2,0.3);
\draw[-] (0.2,-0.3) \botlabel{i} -- (0.2,0.3);
\end{tikzpicture} = \begin{tikzpicture}[iQ,centerzero]
\draw[-] (-0.3,-0.3) \botlabel{i} -- (0.3,0.3);
\draw[-] (0.3,-0.3) \botlabel{j} -- (-0.3,0.3);
\closeddot{-0.15,0.15};
\end{tikzpicture}
-
\begin{tikzpicture}[iQ,centerzero]
\draw[-] (-0.3,-0.3) \botlabel{i} -- (0.3,0.3);
\draw[-] (0.3,-0.3) \botlabel{j} -- (-0.3,0.3);
\closeddot{0.15,-0.15};
\end{tikzpicture} \ ,\label{dotslideiQH}\\\label{quadraticiQH}
\begin{tikzpicture}[iQ,centerzero,scale=1.1]
\draw[-] (-0.2,-0.4) \botlabel{i} to[out=45,in=down] (0.15,0) to[out=up,in=-45] (-0.2,0.4);
\draw[-] (0.2,-0.4) \botlabel{j} to[out=135,in=down] (-0.15,0) to[out=up,in=225] (0.2,0.4);
\end{tikzpicture}
&=
\begin{tikzpicture}[iQ,centerzero,scale=1.1]
\draw[-] (-0.2,-0.4) \botlabel{i} -- (-0.2,0.4);
\draw[-] (0.2,-0.4) \botlabel{j} -- (0.2,0.4);
\Pinpin{0.2,0}{-0.2,0}{-1.1,0}{{'}Q^\imath_{i,j}(x,y)};
\end{tikzpicture}\ ,\\\label{braidiQH}
\begin{tikzpicture}[iQ,centerzero,scale=1.1]
\draw[-] (-0.4,-0.4) \botlabel{i} -- (0.4,0.4);
\draw[-] (0,-0.4) \botlabel{j} to[out=135,in=down] (-0.32,0) to[out=up,in=225] (0,0.4);
\draw[-] (0.4,-0.4) \botlabel{k} -- (-0.4,0.4);
\end{tikzpicture}
\ -\
\begin{tikzpicture}[iQ,centerzero,scale=1.1]
\draw[-] (-0.4,-0.4) \botlabel{i} -- (0.4,0.4);
\draw[-] (0,-0.4) \botlabel{j} to[out=45,in=down] (0.32,0) to[out=up,in=-45] (0,0.4);
\draw[-] (0.4,-0.4) \botlabel{k} -- (-0.4,0.4);
\end{tikzpicture}
&=\delta_{i,k}\ 
\begin{tikzpicture}[iQ,centerzero,scale=1.1]
\draw[-] (-0.3,-0.4) \botlabel{i} -- (-0.3,0.4);
\draw[-] (0,-0.4) \botlabel{j} -- (0,0.4);
\draw[-] (0.3,-0.4) \botlabel{i} -- (0.3,0.4);
\Pinpinpin{.3,0}{0,0}{-.3,0}{-1.7,0}{
\frac{{'}Q^\imath_{i,j}(x,y)-{'}Q^\imath_{i,j}(z,y)}{x-z}};
\end{tikzpicture}\ .
\end{align}
Then, for $\alpha\in\Lambda$ of height $l$, we let
\begin{equation}
\QH^\imath_\alpha := \bigoplus_{\bi,\bj \in \langle I \rangle_\alpha}
\Hom_{\catQH^\imath}(\Theta_\bj,
\Theta_\bi),
\end{equation}
notation as in \cref{wordnotation1b}.
\end{defin}

\begin{lem}\label{sycamore}
We have that 
$\H_\alpha \cong \QH^\imath_\alpha$ as graded $\R^\imath$-algebras.
\end{lem}

\begin{proof}
We remind that $\H$ was defined using ${'}\UU^\imath$ in \cref{pathalgebra} so that its parameters 
are ${'}Q_{i,j}(x,y) = r_{i,j}r_{j,i} Q_{i,j}(x,y)$,
and ${'}Q_{i,j}^\imath(x,y) = (-1)^{\delta_{i, \tau j}}
{'} Q_{i,j}(x,y)$.
With this in mind, it is obvious on comparing
the defining relations \cref{dotslideiQH,quadraticiQH,braidiQH} with
\cref{idotslide,iquadratic,ibraid} that there is a
homomorphism 
$\QH^\imath_\alpha \rightarrow \H_\alpha$ taking diagrams to cosets of 2-morphisms represented by the
same diagrams with $\lambda$ labelling the rightmost 2-cell.
It is an isomorphism because it maps a basis
to a basis.
\end{proof}

We will identify $\QH^\imath_\alpha$ and $\H_\alpha$ 
via the isomorphism from the proof of \cref{sycamore}.
The general theory developed by Khovanov and Lauda in \cite{KL1} shows that $\catQH^\imath$ categorifies 
the $\Z[q,q^{-1}]$-form $\f_\Z$ for the algebra $\f$ generated by the divided powers
$\theta_i^{(n)} := \theta_i^n / [n]^!_{q_i}$. By this, we mean that 
$K_0\left(\Kar(\catQH^\imath_q)\right) \cong \f_\Z$ as $\Z[q,q^{-1}]$-algebras. 
The canonical isomorphism 
\begin{equation}\label{KLequiv}
\KL:
\f_\Z\stackrel{\sim}{\rightarrow}K_0\left(\Kar(\catQH^\imath_q)\right)
\end{equation}
maps $\theta_i \mapsto [\Theta_i]$ for each $i \in I$.
There are also divided powers \begin{equation}
\Theta_i^{(n)} := \Big(q_i^{-\binom{n}{2}}\Theta_i^{\otimes n},
  1_{i^{(n)}}\Big)
\end{equation}
where $1_{i^{(n)}} \in \QH^\imath_{n \alpha_i}$ is the idempotent defined by \cref{endef2}.
Then we can define $\theta_\bi$ and $\Theta_\bi$ for 
$\bi \in \llangle I \rrangle$ as in \cref{wordnotation2b},
and have that $\KL(\theta_\bi)=[\Theta_\bi]$.

There is one more essential piece of data:
for each $\alpha \in \Lambda$, we need to be given a set
$\B_\alpha$
and irreducible graded left $\H_\alpha$-modules
$L_\alpha(\bb)$ for $\bb \in \B_\alpha$
such that
\begin{equation}
\left\{L_\alpha(\bb)\:\big|\:\bb \in \B_\alpha\right\}
\end{equation} 
is a full set of irreducible graded left $\H_\alpha$-modules up to grading shift.
As at the end of \cite[Sec.~3.2]{KL1}, we assume the grading shift on each $L_\alpha(\bb)$ is chosen so that it is graded self-dual.
We will not concern ourselves with the much-studied problem of determining an {\em explicit parameterization} of irreducibles here---all that is important for us is that the {\em set} of isomorphism classes of the modules 
$L_\alpha(\bb)\:(\bb \in \B_\alpha)$ does not depend on any choices.
Then we let $P_\alpha(\bb)$ be a projective cover 
and $I_\alpha(\bb)$ be an injective hull of
$L_\alpha(\bb)$ in $\gmod{\H_\alpha}$. 

Similar to \cref{yonedaop}, there is an equivalence
\begin{equation}\label{otheryonedaop}
\cov:\Kar(\catQH^\imath_q) \rightarrow \bigoplus_{\alpha \in \Lambda}
\gproj{\H_\alpha}
\end{equation}
defined using the contravariant Yoneda equivalence together with the anti-automorphism that reflects string diagrams in a horizontal axis.
It maps $\Theta_\bi\:(\bi \in \llangle I \rrangle_\alpha)$ to a finitely generated graded projective module that is isomorphic to
\begin{equation}
P_\alpha(\bi) := q^{\deg(\bi)} \H_\alpha\; 1_\bi.
\end{equation}
For $\bb\in \B_\alpha$, let $\Theta_\bb$ be an indecomposable object of
$\Kar(\QH^\imath_q)$ such that
\begin{equation}
\cov(\Theta_\bb) \cong P_\alpha(\bb).
\end{equation}
Similar to \cref{nakayama}, there is an equivalence
$\Nak:\gproj{\H_\alpha} \rightarrow \ginj{\H_\alpha}$,
hence, $\Kar(\catQH^\imath_q)$ is also equivalent
to $\bigoplus_{\alpha \in \Lambda} \ginj{\H_\alpha}$.
For $\bi \in \llangle I \rrangle_\alpha$, we let
\begin{equation}
I_\alpha(\bi) := q^{\deg(\bi)} (1_\bi\;\H_\alpha)^\circledast
\cong \Nak(P_\alpha(\bi)).
\end{equation}

With these preliminaries about Cartan algebras out of the way, 
we can return to the study of $\gmod{H}$.
For $\alpha \in \Lambda$, the functor
$\jmath^\alpha:\gmod{\H_{\geq \alpha}} \rightarrow \gmod{\H_\alpha}$
defined by truncation with respect to the idempotent $\bar e_\alpha$
has a left adjoint $\jmath^\alpha_!$ and a right adjoint $\jmath^\alpha_*$ defined by
\begin{align}\label{stfun}
\jmath^\alpha_! &:= \H_{\geq \alpha} \bar e_\alpha \otimes_{\H_\alpha} -,&
\jmath^\alpha_* &:= \bigoplus_{\beta \in \Lambda} \Hom_{\H_\alpha}(\bar e_\alpha \H_{\geq \alpha} \bar e_\beta, -).
\end{align}
Composing with inflation from $\H_{\geq\alpha}$ to $\H$,
we can view them as functors
$\jmath^\alpha_!:\gmod{\H_\alpha} \rightarrow \gmod{\H}$
and $\jmath^\alpha_*:\gmod{\H_\alpha}\rightarrow\gmod{\H}$,
which we call 
{\em standardization} and {\em costandardization}, respectively.
An important point is that both of these functors are exact. This
follows from the next lemma, which shows that 
$\H_{\geq \alpha} \bar e_\alpha$ is a projective right
$\H_\alpha$-module and
$\bar e_\alpha \H_{\geq \alpha}$ is a projective left $\H_\alpha$-module.

\begin{lem}\label{standardize}
For $\bi \in \langle I \rangle$, we have that 
$$
\bar 1_\bi \H_{\geq \alpha} \bar e_\alpha
= \bigoplus_{\substack{\bj \in \langle I \rangle_\alpha\\x \in
  X(\bi,\bj)}} \bar x \H_\alpha
\qquad\text{and}\qquad
\bar e_\alpha \H_{\geq \alpha} \bar 1_\bi
= \bigoplus_{\substack{\bj \in \langle I \rangle_\alpha\\y \in
    Y(\bj,\bi)}} \H_\alpha \bar y
$$
with $\bar x \H_\alpha \cong q^{\deg(x)} \bar 1_\bj \H_\alpha$ as a
graded right $\H_\alpha$-module, and
$\H_\alpha \bar y \cong q^{\deg(y)} \H_\alpha \bar 1_\bj$ as a
graded left $\H_\alpha$-module, respectively.
\end{lem}

\begin{proof}
This follows from \cite[(4.4)--(4.5)]{GTB}.
\end{proof}

Now we introduce
{\em standard} and {\em costandard modules}, denoted using $\Delta$ and $
\nabla$, respectively. 
Let
\begin{equation}\label{bigB}
\B := \bigsqcup_{\alpha \in \Lambda} \B_\alpha,
\end{equation}
assuming also that this set
is disjoint from $\llangle I \rrangle$.
For $\alpha \in \Lambda$, $\bb \in \B_\alpha$ and $\bi \in \llangle I \rrangle_\alpha$, we let
\begin{align}
\Delta(\bb) &:= \jmath^\alpha_! P_\alpha(\bb),&
\nabla(\bb) &:= \jmath^\alpha_* I_\alpha(\bb),\\
\Delta(\bi) &:= \jmath^\alpha_! P_\alpha(\bi),&
\nabla(\bi) &:= \jmath^\alpha_* I_\alpha(\bi).
\intertext{
Since $P_\alpha(\bi)$ is a finite direct sum of grading shifts of the indecomposable projectives $P_\alpha(\bb)\:(\bb \in \B_\alpha)$, 
the standard module $\Delta(\bi)$ is a finite direct sum of grading shifts of the indecomposable standard modules $\Delta(\bb)\:(\bb \in \B_\alpha)$.
Similarly, the costandard module 
$\nabla(\bi)$ is a finite direct sum of grading shifts of
the indecomposable costandard modules $\nabla(\bb)\:(\bb \in \B_\alpha)$.
Also, for $\bb \in \B_\alpha$ again, there are the {\em proper standard} and
{\em proper costandard modules}}
\bar\Delta(\bb) &:= \jmath^\alpha_! L_\alpha(\bb),&
\bar\nabla(\bb) &:= \jmath^\alpha_* L_\alpha(\bb).
\end{align}
All of these are graded left $\H$-modules.
Fundamental to the entire theory is that
\begin{equation}\label{CPSfactoid}
\operatorname{Ext}^i_\H(\Delta(\ba), \bar\nabla(\bb))
\cong \operatorname{Ext}^i_\H(\bar\Delta(\ba), \nabla(\bb))
\cong \begin{cases}
\kk&\text{if $i = 0$ and $\ba = \bb$}\\
\{0\}&\text{otherwise.}
\end{cases}
\end{equation}
See \cite[Cor.~4.6, Th.~7.5]{GTB} for the proof of this.
Also, $\Delta(\ba)$ has a $\bar\Delta$-flag in the sense of 
\cite[Def.~6.3]{GTB}
and 
$\nabla(\ba)$ has a $\bar\nabla$-flag
in the sense of \cite[Def.~6.4]{GTB}, with
\begin{align}
(\Delta(\ba):\bar\Delta(\bb))_q 
= \begin{cases}
[P_\alpha(\ba):L_\alpha(\bb)]_q
&\text{if $\ba,\bb \in \B_\alpha$ for some $\alpha \in \Lambda$}\\
0&\text{otherwise,}
\end{cases}\\
(\nabla(\ba):\bar\nabla(\bb))_q 
=
\begin{cases}
[I_\alpha(\ba):L_\alpha(\bb)]_q
&\text{if $\ba,\bb \in \B_\alpha$ for some $\alpha \in \Lambda$}\\
0&\text{otherwise.}
\end{cases}
\end{align}
This is obvious from the definitions and exactness of the standardization/costandardization functors; see also \cite[(8.9), (8.15)]{GTB}.

Now everything is in place to be able to apply the general theory of graded triangular bases. It shows that the set $\B$ from \cref{bigB}
indexes a full set 
\begin{equation}
\{L(\bb)\:|\:\bb \in \B\}
\end{equation}
of irreducible graded left $\H$-modules $L(\bb)$ up to grading shift. 
By definition, $L(\bb)$ is the irreducible head of
$\bar\Delta(\bb)$ and the irreducible socle of $\bar\nabla(\bb)$. 
Assuming that $\bb \in \B_\alpha$,
the projective cover $P(\bb)$
of $\bar\Delta(\bb)$ in $\gmod{\H}$ has a $\Delta$-flag in the sense of
\cite[Def.~6.3]{GTB} with top section $\Delta(\bb)$ and finitely many other sections indexed 
by the weights $\beta \in \Lambda$ with $\beta < \alpha$.
By the version of BGG reciprocity from \cite[Cor.~8.4]{GTB},
the $\beta$-section of this $\Delta$-flag
is isomorphic to $\bigoplus_{\ba \in \B_\beta} \Delta(\ba)^{\oplus (P(\bb):\Delta(\ba))_q}$
with
\begin{equation}\label{BGG}
(P(\bb):\Delta(\ba))_q := 
\grdim_{q^{-1}} \Hom_{\H}(P(\bb), \bar\nabla(\ba)) = [\bar\nabla(\ba):L(\bb)]_{q^{-1}} \in \N\lround q \rround.
\end{equation}
Also $[\bar\nabla(\ba):L(\bb)]_{q}$ is the usual graded composition multiplicity and
$[\bar\nabla(\ba):L(\bb)]_{q^{-1}}$ is its image under the bar involution.
There is a dual 
statement describing the injective hull $I(\bb)$ of $\bar\nabla(\bb)$, which will not be needed subsequently; see \cite[Cor.~8.9]{GTB}.
The following theorem is just a restatement of part of this discussion:

\begin{theo}\label{TBA}
For each $\bb \in \B$, let $B_\bb \one_\lambda$ be the unique (up to
isomorphism) indecomposable object of $\Kar(\UU^\imath_q)\one_\lambda$ 
such that $\cov(B_\bb \one_\lambda) \cong P(\bb)$. Then
the 1-morphisms $B_\bb \one_\lambda\:(\bb \in \B)$
give a full set of indecomposable objects in 
$\Kar(\UU^\imath_q \one_\lambda)$ up to grading shift.
\end{theo}

There are two more useful categories of $\H$-modules,
$\standard{\H}$ 
and $\costandard{\H}$, which are the full subcategories of $\gmod{\H}$ consisting of modules with a $\Delta$-flag in the sense of \cite[Def.~6.3]{GTB} or a $\nabla$-flag in the sense of \cite[Def.~6.4]{GTB}, respectively. These are exact categories (but not Abelian).
If $V$ is a module with a $\Delta$-flag, it has a filtration with sections indexed by 
$\Lambda$, all but finitely many of which are zero,
with $\beta$-section isomorphic to $\bigoplus_{\bb\in \B_\beta}
\Delta(\bb)^{\oplus (V:\Delta(\bb))_q}$ where
\begin{equation}\label{deltaflags}
(V:\Delta(\bb))_q = \grdim_{q^{-1}} \Hom_\H(V, \bar\nabla(\bb))
\in \N\lround q\rround.
\end{equation}
The Grothendieck group
$K_0(\standard{\H})$ is the free $\Z\lround q\rround$-module
with basis $[\Delta(\bb)]\:(\bb \in \B)$
and, for $V$ as in \cref{deltaflags}, we have that
$[V] = \sum_{\bb \in \B} (V:\Delta(\bb))_q [\Delta(\bb)]$.
By the BGG reciprocity discussed before \cref{BGG}, finitely generated projectives have $\Delta$-flags, and the inclusion functor
$\gproj{\H} \rightarrow \standard{\H}$ induces an embedding
\begin{equation}\label{projintostandard}
K_0\big(\gproj{\H}\big) \hookrightarrow K_0\big(\standard{\H}\big)
\end{equation}
of $\Z[q,q^{-1}]$-modules.
There are analogous dual statements about modules with $\nabla$-flags.
In particular, $\ginj{\H}$ embeds into $\costandard{\H}$,
and $K_0(\costandard{\H})$ is the free $\Z\lround q^{-1}\rround$-module with basis $[\nabla(\bb)]\:(\bb \in \B)$.
We refer to \cite{GTB} for further background.

\begin{lem}\label{totaipei}
For $\alpha \in \Lambda$,
there is a partial order $\leq_\alpha$ on $\B_\alpha$ 
and an injective function
$\bi_\alpha:\B_\alpha \hookrightarrow \llangle I \rrangle_\alpha$
such that 
$$
P_\alpha(\bi_\alpha(\bb)) \cong P_\alpha(\bb) \oplus\!
\bigoplus_{\ba \in \B_\alpha\,\text{with}\,\ba <_\alpha \bb} \!\!\! P_\alpha(\ba)^{\oplus m_{\bb,\ba}(q)}
$$
for some $m_{\bb,\ba}(q) \in \N[q,q^{-1}]$.
\end{lem}

\begin{proof}
This is a well-known general property of quiver Hecke
algebras. 
Its proof is explained in \cite[Sec.~3.2]{KL1} (the strategy of the
proof can be seen already in \cite[Sec.~11]{Groj}).
\end{proof}

\begin{lem}\label{halfwaytotaipei}
For $\bi \in \llangle I \rrangle_\alpha$, $P(\bi)$ is isomorphic to the
direct sum of the
projective cover of
$\Delta(\bi)$ and grading shifts of finitely many other indecomposable projectives
$P(\ba)$ for $\ba \in \bigcup_{\beta < \alpha} \B_\beta$.
\end{lem}

\begin{proof}
Take any $\beta\in \Lambda$ and $\bb \in \B_\beta$. Since $P(\bi) = q^{\deg(\bi)} \H 1_\bi$, we have that
$$
\Hom_\H(P(\bi), L(\bb)) \cong 
q^{-\deg(\bi)} 1_\bi L(\bb).
$$
This is $\{0\}$ unless the $\alpha$-weight space of $L(\bb)$ is
non-zero. By the general theory of graded triangular bases, $\beta$ is the lowest
weight of $L(\bb)$, so this implies that
$\beta \leq \alpha$.
Consequently, $P(\bi)$ is a direct sum of finitely many indecomposable
projectives that are grading shifts of
$P(\bb)$ for $\bb \in \bigcup_{\beta \leq\alpha}\B_\beta$.
It remains to observe for $\bb \in \B_\alpha$ that
$$
\Hom_\H(\Delta(\bi), L(\bb)) =
\Hom_{\H_{\geq \alpha}}(\Delta(\bi), L(\bb))
\cong \Hom_{\H_\alpha}(P_\alpha(\bi), e_\alpha L(\bb))
\cong
q^{-\deg(\bi)} 1_\bi L(\bb)).
$$
Here, we used adjunction plus the definition $P_\alpha(\bi) =
q^{\deg(\bi)} \H_\alpha 1_\bi$.
So 
$$
\Hom_\H(\Delta(\bi), L(\bb))\cong \Hom_\H(P(\bi), L(\bb)).
$$
We deduce for $\bb \in \B_\alpha$ that the graded multiplicity of the summand $P(\bb)$ of $P(\bi)$ is 
the same as the graded multiplicity of $P(\bb)$
in the projective cover of $\Delta(\bi)$.
The lemma follows.
\end{proof}

\begin{theo}\label{taipeiflight}
There is a partial order $\leq$ on $\B$ 
and an injective function
$\bi:\B \hookrightarrow \llangle I \rrangle$
such that the following hold:
\begin{enumerate}
\item
$\bi$ maps $\B_\alpha$ into $\llangle I \rrangle_\alpha$ for each $\alpha \in \Lambda$;
\item
if $\ba \in \B_\alpha$, $\bb \in \B_\beta$ and $\ba \leq \bb$
then $\alpha \leq \beta$;
\item
for any $\bb \in \B$,
$$
\displaystyle P(\bi(\bb)) \cong P(\bb) \oplus\!
\bigoplus_{\ba \in \B\text{ with }\ba < \bb} \!\!\!P(\ba)^{\oplus
  m_{\bb,\ba}(q)}
$$
for some $m_{\bb,\ba}(q) \in \N[q,q^{-1}]$.
\end{enumerate}
\end{theo}

\begin{proof}
We define the required partial order $\leq$ on $\B$ by 
declaring for $\ba \in \B_\beta$ and $\bb \in \B_\alpha$
that $\ba \leq \bb$ if {\em either}
$\beta < \alpha$ in the partial order on $\Lambda$, 
{\em or} $\beta = \alpha$ and
$\ba \leq_\alpha \bb$ in the partial order from \cref{totaipei}. Of course this has the property (1).
We define the function $\bi:\B \hookrightarrow \llangle I \rrangle$
so that $\bi|_{\B_\alpha}$ is the function $\bi_\alpha$ from \cref{totaipei}.
Of course this has the property (2).
It remains to prove the property (3).
Suppose that $\bb \in \B_\alpha$. Applying the exact standardization
functor
$\jmath^\alpha_!$ to the decomposition in \cref{totaipei}, we get that
$$
\Delta(\bi(\bb)) \cong \Delta(\bb) \oplus\!
\bigoplus_{\ba \in \B_\alpha\text{ with }\ba <_\alpha \bb} \!\!\!\Delta(\ba)^{\oplus
  m_{\bb,\ba}(q)}.
$$
By \cref{halfwaytotaipei}, $P(\bi(\bb))$ is the direct sum of the
projective cover of this module, which is
$$
P(\bb) \oplus\!
\bigoplus_{\ba \in \B_\alpha\text{ with }\ba <_\alpha \bb} \!\!\!P(\ba)^{\oplus
  m_{\bb,\ba}(q)},
$$
and grading shifts of finitely many other indecomposable projectives
$P(\ba)$ for $\ba \in \B_\beta$ with $\beta < \alpha$.
This proves the theorem.
\end{proof}

\begin{cor}\label{internationalhouse}
The classes $[B_\bi \one_\lambda]$ for $\bi \in \llangle I \rrangle$
span $K_0\left(\Kar(\UU^\imath_q)\one_\lambda\right)$ as a $\Z[q,q^{-1}]$-module.
\end{cor}

\begin{proof}
Using the equivalence $\cov:\Kar(\UU^\imath_q)\one_\lambda \rightarrow
\gproj{\H}$, the uni-triangularity property from \cref{taipeiflight}(3) implies that each
$[B_\bb]\:(\bb \in \B)$ can be written as a $\Z[q,q^{-1}]$-linear combination of $[B_\bi]$ for $\bi
\in \llangle I \rrangle$.
\end{proof}

For the final result in this subsection,
which generalizes \cite[Th.~5.13]{BWWa1},
we use a well-known adjoint pair $(I_{\alpha_i,\alpha}, R_{\alpha_i,\alpha})$ of functors
\begin{align}
I_{\alpha_i,\alpha}:\gmod{\H_\alpha} &\rightarrow \gmod{\H_{\alpha_i+\alpha}},&
R_{\alpha_i,\alpha}:\gmod{\H_{\alpha_i+\alpha}} &\rightarrow \gmod{\H_\alpha}.
\end{align}
To define them, let $\bar e_{\alpha_i,\alpha}$ be the idempotent in $\H_{\alpha_i+\alpha}$ that is the image of $e_{\alpha_i,\alpha} := 
1_i \otimes e_\alpha$.
There is an algebra homomorphism 
$\iota_{\alpha_i,\alpha}:\H_\alpha \hookrightarrow \H_{\alpha_i+\alpha}$
defined in terms of string diagrams by tensoring on the left with a vertical string of color $i$. This maps $\bar e_\alpha$ to $\bar e_{\alpha_i,\alpha}$.
Using it, we view $\H_\alpha$ as a unital subalgebra of
$\bar e_{\alpha_i,\alpha} \H_{\alpha_i+\alpha} \bar e_{\alpha_i,\alpha}$.
Then $\H_{\alpha_i+\alpha} \bar e_{\alpha_i,\alpha}$ is a graded  $(\H_{\alpha_i+\alpha}, \H_\alpha)$-bimodule and we can define
\begin{align}
I_{\alpha_i, \alpha} &:= \H_{\alpha_i+\alpha} \bar e_{\alpha_i,\alpha}
\otimes_{\H_{\alpha}} -,&
R_{\alpha_i,\alpha} &:= 
\bar e_{\alpha_i,\alpha} \H_{\alpha+\alpha_i}\otimes_{\H_{\alpha+\alpha_i}} -
\cong \Hom_{\H_{\alpha+\alpha_i}}(\H_{\alpha_i+\alpha} \bar e_{\alpha_i,\alpha},-).
\end{align}
It is clear that $(I_{\alpha_i,\alpha}, R_{\alpha_i,\alpha})$ is an
adjoint pair (with adjunction that is homogeneous of degree zero).
Both functors are additive graded functors which take finitely generated projectives to finitely generated projectives, so they induce $\Z[q,q^{-1}]$-linear maps
between the split Grothendieck groups $K_0\big(\gproj{\H_\alpha}\big)$ and $K_0\big(\gproj{\H_{\alpha+\alpha_i}})$.
The following lemma is well known; e.g., see \cite{KK}.

\begin{lem}\label{nofoodleft}
Under the identification of split Grothendieck groups of $\gproj{\H_\alpha}$ and $\gproj{\H_{\alpha+\alpha_i}}$
with the $\Z[q,q^{-1}]$-forms for $\f_\alpha$ and $\f_{\alpha+\alpha_i}$,
the map induced by $I_{\alpha_i,\alpha}$ corresponds to left multiplication by $\theta_i$, and the map induced by $R_{\alpha_i,\alpha}$ corresponds to the map ${_i}\widetilde{R}$ from \cref{newmoaaz1}.
\end{lem}

Acting on the left with $B_i$ defines an endofunctor of 
$\UU^\imath \one_\lambda$. It induces a graded functor we also denote by
\begin{equation}\label{bim1}
B_i: \gmod{\H} \rightarrow\gmod{\H}.
\end{equation}
This is defined explicitly by tensoring with the $(\H,\H)$-bimodule 
$\bigoplus_{\bj \in \langle I \rangle} \H 1_{i \bj}$
viewed as a left $\H$-module in the natural way and as a right
$\H$-module with action defined by pulling-back the natural right
action along the homomorphism $\H \rightarrow
\H$ which maps a string diagram to the string diagram obtained by
adding a vertical string labelled $i$ on the left boundary.

\begin{lem}\label{vargaslemma}
The functor $B_i:\gmod{\H}\rightarrow\gmod{\H}$ is isomorphic to the functor
defined by tensoring with the graded $(\H,\H)$-bimodule
$$
\bigoplus_{\beta \in \Lambda} \bigoplus_{\bi \in \langle I \rangle_\beta}
q_i^{(\lambda-\beta)_i - \del_i-1} 
1_{(\tau i) \bi} \H
$$
with the natural right action of $\H$ and the left action defined
by pull-back
along the homomorphism which adds a vertical string labelled $\tau i$
on the left boundary.
\end{lem}

\begin{proof}
This is similar to \cite[Lem.~2.5]{BV}. The idea is to define
 mutually inverse bimodule homomorphisms by attaching a suitable cup or cap to the string diagrams.
\end{proof}

The short exact sequence in the following theorem 
is important because it categorifies
the equality
\begin{equation}\label{stench}
b_i \widetilde{\jmath}^{\:-1}(x) = \widetilde{\jmath}^{\:-1}(\theta_i x)
+ q_i^{(\lambda-\alpha)_i-\del_i-1} \widetilde{\jmath}^{\:-1}\big({_{\tau i}}\widetilde{R}(x)\big)
\end{equation} 
for $x \in \f_\alpha$ from \cref{propsofjbar1b},
where $\widetilde{\jmath}^{\:-1}$ is the inverse of the isomorphism
$\widetilde{\jmath}:\f\rightarrow \dot\U^\imath 1_\lambda$.
This follows using also \cref{nofoodleft}.

\begin{theo}\label{lastgasp}
Take $\alpha \in \Lambda$ and $i \in I$.
There is a short exact
sequence of functors
$$
0\longrightarrow
q_i^{(\lambda-\alpha)_i-\del_i-1}
\jmath^{\alpha-\alpha_{\tau i}}_! \circ R_{\alpha_{\tau i},\alpha-\alpha_{\tau i}}
\longrightarrow B_i \circ \jmath^\alpha_! \longrightarrow
\jmath^{\alpha+\alpha_i}_! \circ I_{\alpha_i,\alpha} \longrightarrow 0
$$
from $\gmod{\H_\alpha}$ to $\gmod{\H}$.
The first term should be interpreted as the zero functor if
$\alpha-\alpha_{\tau i} \notin \Lambda$.
\end{theo}

\begin{proof}
All three functors appearing in the short exact sequence are isomorphic to the functors defined by tensoring with certain graded 
$(\H, \H_\alpha)$-bimodules, as follows:
\begin{itemize}
\item
$j^{\alpha-\alpha_{\tau i}}_! \circ R_{\alpha_{\tau i},\alpha-\alpha_{\tau i}}$
is tensoring with 
$q_i^{(\lambda-\alpha)_i-\del_i-1}\H_{\geq (\alpha-\alpha_{\tau i})} \bar e_{\alpha-\alpha_{\tau i}} \otimes_{\H_{\alpha-\alpha_{\tau i}}} 
\bar e_{\alpha_{\tau i},\alpha-\alpha_{\tau i}} \H_\alpha$;
\item
$B_i \circ j^\alpha_!$ is tensoring with
$\bigoplus_{\beta \in \Lambda} \bigoplus_{\bi \in \langle I \rangle_\beta}
q_i^{(\lambda-\beta)_i - \del_i-1} 
\bar e_{(\tau i) \bi} \H_{\geq \alpha} \bar e_{\alpha}$ (see \cref{vargaslemma});
\item
$j^{\alpha+\alpha_i}_! \circ I_{\alpha_i,\alpha}$ is tensoring with
$\H_{\geq (\alpha+\alpha_i)} \bar e_{\alpha+\alpha_i} \otimes_{\H_{\alpha+\alpha_i}} \H_{\alpha+\alpha_i} \bar e_{\alpha_i,\alpha}$.
\end{itemize}
We claim that there is
a short exact sequence of graded bimodules and degree-preserving bimodule homomorphisms:
\begin{multline*}
0 \longrightarrow 
q_i^{(\lambda-\alpha)_i-\del_i-1}
\H_{\geq (\alpha-\alpha_{\tau i})} \bar e_{\alpha-\alpha_{\tau i}} \otimes_{\H_{\alpha-\alpha_{\tau i}}} 
\bar e_{\alpha_{\tau i},\alpha-\alpha_{\tau i}} \H_\alpha\stackrel{f}{\longrightarrow}\\
\bigoplus_{\beta \in \Lambda} \bigoplus_{\bi \in \langle I \rangle_\beta}
q_i^{(\lambda-\beta)_i - \del_i-1} 
\bar 1_{(\tau i) \bi} \H_{\geq \alpha} \bar e_{\alpha}\stackrel{g}{\longrightarrow}
\H_{\geq (\alpha+\alpha_i)} \bar e_{\alpha+\alpha_i} \otimes_{\H_{\alpha+\alpha_i}} \H_{\alpha+\alpha_i} \bar e_{\alpha_i,\alpha}
\longrightarrow 0.
\end{multline*}
The graded bimodule homomorphism $f$ is defined
on basis vectors so that
\begin{align*}
f:\begin{tikzpicture}[anchorbase,scale=1.25]
\draw[-] (0.4,-0.05)--(0.4,-0.4);
\draw[-] (-0.2,0.15)--(-0.2,0.55);
\draw[-] (-0.2,-0.05)--(-0.2,-0.4);
\draw[-] (0.4,0.15)--(0.4,0.55);
\node at (-0.05,0.42) {$\cdot$};
\node at (0.1,0.42) {$\cdot$};
\node at (0.25,0.42) {$\cdot$};
\node at (-0.05,-0.27) {$\cdot$};
\node at (0.1,-0.27) {$\cdot$};
\node at (0.25,-0.27) {$\cdot$};
\node[rectangle,rounded corners,draw,fill=blue!15!white,inner sep=4.2pt] at (0.1,.1) {$\hspace{3mm}\scriptstyle x\hspace{3mm}$};
\draw [thick,decoration={brace,raise=0.4cm},decorate] (-.2,.3) -- (.4,.3);
\stringlabel{.1,.85}{\bi};
\draw [thick,decoration={brace,mirror,raise=0.4cm},decorate] (-.2,-.15) -- (.4,-.15);
\stringlabel{.1,-.7}{\bj};
\region{.7,0.1}{\lambda};
\end{tikzpicture}
\:\otimes\:
\begin{tikzpicture}[anchorbase,scale=1.25]
\draw[-] (-0.45,-1.1)--(-0.45,-0.7);
\draw[-] (0.45,-1.1)--(0.45,-0.7);
\draw[-] (-0.45,-0.1)\toplabel{\tau i}--(-0.45,-0.4);
\draw[-] (-0.25,-0.1)--(-0.25,-0.4);
\draw[-] (0.45,-0.1)--(0.45,-0.4);
\node at (-0.1,-0.25) {$\cdot$};
\node at (0.1,-0.25) {$\cdot$};
\node at (0.3,-0.25) {$\cdot$};
\node at (-0.25,-.96) {$\cdot$};
\node at (0,-.96) {$\cdot$};
\node at (0.25,-.96) {$\cdot$};
\node[rectangle,rounded corners,draw,fill=blue!15!white,inner sep=4.2pt] at (0,-.6) {$\hspace{5mm}\scriptstyle h\hspace{5mm}$};
\draw [thick,decoration={brace,raise=0.4cm},decorate] (-.25,-.35) -- (.45,-.35);
\stringlabel{.1,.2}{\bj};
\draw [thick,decoration={brace,mirror,raise=0.4cm},decorate] (-.45,-.85) -- (.45,-.85);
\stringlabel{0,-1.38}{\bk};
\region{.75,-.6}{\lambda};
\end{tikzpicture}
&\mapsto
 \begin{tikzpicture}[anchorbase,scale=1.25]
\draw[-] (-0.45,-1.05)--(-0.45,-0.7);
\draw[-] (0.45,-1.05)--(0.45,-0.7);
\draw[-] (-0.45,0.5)\toplabel{\tau i}--(-0.45,-0.4);
\draw[-] (0.4,-0.05)--(0.4,-0.4);
\draw[-] (-0.2,0.15)--(-0.2,0.5);
\draw[-] (-0.2,-0.05)--(-0.2,-0.4);
\draw[-] (0.4,0.15)--(0.4,0.5);
\node at (-0.05,0.4) {$\cdot$};
\node at (0.1,0.4) {$\cdot$};
\node at (0.25,0.4) {$\cdot$};
\node at (-0.05,-0.25) {$\cdot$};
\node at (0.1,-0.25) {$\cdot$};
\node at (0.25,-0.25) {$\cdot$};
\node at (-0.25,-0.95) {$\cdot$};
\node at (0,-0.95) {$\cdot$};
\node at (0.25,-0.95) {$\cdot$};
\node[rectangle,rounded corners,draw,fill=blue!15!white,inner sep=4.2pt] at (0.1,.1) {$\hspace{3mm}\scriptstyle x\hspace{3mm}$};
\node[rectangle,rounded corners,draw,fill=blue!15!white,inner sep=4.2pt] at (0,-.6) {$\hspace{5mm}\scriptstyle h\hspace{5mm}$};
\draw [thick,decoration={brace,raise=0.4cm},decorate] (-.23,.23) -- (.43,.23);
\stringlabel{.1,.78}{\bi};
\draw [thick,decoration={brace,mirror,raise=0.4cm},decorate] (-.45,-.78) -- (.45,-.78);
\stringlabel{0,-1.32}{\bk};
\region{.7,0.1}{\lambda};
\end{tikzpicture}
\end{align*}
for $\beta\geq\alpha, \bi \in \langle I \rangle_\beta, \bj \in \langle I\rangle_{\alpha-\alpha_{\tau i}},
\bk \in \langle I\rangle_\alpha$ and $x \in X(\bi,\bj), h \in H((\tau i) \bj, \bk)$.
The graded bimodule homomorphism $g$ is defined on basis vectors so that
the vectors in the image of $f$ are in the kernel of $g$ and 
$$
g:
\begin{tikzpicture}[anchorbase,scale=1.25]
\draw[-] (-0.45,0.5)\toplabel{\tau i}--(-0.45,-0.2)
to [out=-90,in=180] (-.2,-.5)
to [out=0,in=-90] (0,-.05);
\draw[ultra thick] (0.4,-0.05) to [out=-90,in=90] (0.6,-0.7);
\draw[-] (-0.15,0.2)--(-0.15,0.5);
\draw[ultra thick] (-0.2,-0.05) to [out=-90,in=90] (0.05,-0.4) to [out=-90,in=90] (0,-.75);
\draw[-] (0.4,0.2)--(0.4,0.5);
\draw[-] (0,-1.05)--(0,-1.35);
\draw[-] (0.6,-1.05)--(0.6,-1.35);
\node at (-0.05,0.4) {$\cdot$};
\node at (0.1,0.4) {$\cdot$};
\node at (0.25,0.4) {$\cdot$};
\node at (.14,-1.27) {$\cdot$};
\node at (0.32,-1.27) {$\cdot$};
\node at (0.5,-1.27) {$\cdot$};
\stringlabel{0.23,-.45}{\bj^-};
\stringlabel{0.1,-.15}{i};
\stringlabel{0.74,-.3}{\bj^+};
\multcloseddot{-.05,-.4}{east}{r};
\region{.8,0.1}{\lambda};
\node[rectangle,rounded corners,draw,fill=blue!15!white,inner sep=4.2pt] at (0.11,0.1) {$\hspace{3mm}\scriptstyle x\hspace{3mm}$};
\node[rectangle,rounded corners,draw,fill=blue!15!white,inner sep=4.2pt] at (.3,-.9) {$\hspace{3mm}\scriptstyle h\hspace{3mm}$};
\draw [thick,decoration={brace,raise=0.4cm},decorate] (-.2,.23) -- (.42,.23);
\stringlabel{.11,.78}{\bi};
\draw [thick,decoration={brace,mirror,raise=0.4cm},decorate] (-0.02,-1.1) -- (.62,-1.1);
\stringlabel{0.3,-1.6}{\bk};
\end{tikzpicture}\:\:
\mapsto
\:\:
\begin{tikzpicture}[anchorbase,scale=1.25]
\draw[-] (0.4,-0.05)--(0.4,-0.4);
\draw[-] (-0.2,0.15)--(-0.2,0.55);
\draw[-] (-0.2,-0.05)--(-0.2,-0.4);
\draw[-] (0.4,0.15)--(0.4,0.55);
\node at (-0.05,0.42) {$\cdot$};
\node at (0.1,0.42) {$\cdot$};
\node at (0.25,0.42) {$\cdot$};
\node at (-0.05,-0.27) {$\cdot$};
\node at (0.1,-0.27) {$\cdot$};
\node at (0.25,-0.27) {$\cdot$};
\node[rectangle,rounded corners,draw,fill=blue!15!white,inner sep=4.2pt] at (0.1,.1) {$\hspace{3mm}\scriptstyle x\hspace{3mm}$};
\draw [thick,decoration={brace,raise=0.4cm},decorate] (-.2,.3) -- (.4,.3);
\stringlabel{.1,.85}{\bi};
\draw [thick,decoration={brace,mirror,raise=0.4cm},decorate] (-.2,-.15) -- (.4,-.15);
\stringlabel{.1,-.7}{\bj};
\region{.7,0.1}{\lambda};
\end{tikzpicture}
\:\:\otimes
\begin{tikzpicture}[baseline=-1mm,scale=1.25]
\draw[-] (-0.45,-0.5) \botlabel{i} to [out=90,in=-90] (-.45,0) to [out=90,in=-120] (0,0.5);
\draw[-] (0.4,-0.15)--(0.4,-0.5);
\draw[ultra thick] (-0.15,0.1) to [out=90,in=-90] (-.45,0.5)\toplabel{\bj^-};
\draw[-] (-0.2,-0.15)--(-0.2,-0.5);
\draw[ultra thick] (0.35,0.1)--(0.35,0.5)\toplabel{\bj^+};
\region{.7,0}{\lambda};
\multcloseddot{-.45,-.2}{east}{r};
\node at (-0.05,-0.4) {$\cdot$};
\node at (0.1,-0.4) {$\cdot$};
\node at (0.25,-0.4) {$\cdot$};
\node[rectangle,rounded corners,draw,fill=blue!15!white,inner sep=4.2pt] at (0.1,-.07) {$\hspace{3mm}\scriptstyle h\hspace{3mm}$};
\draw [thick,decoration={brace,mirror,raise=0.4cm},decorate] (-0.22,-.25) -- (.42,-.25);
\stringlabel{0.1,-.75}{\bk};
\end{tikzpicture}
$$
for $\beta \geq \alpha+\alpha_i$, $\bi \in \langle I \rangle_\beta, 
\bj = \bj^- i \bj^+ \in \langle I \rangle_{\alpha+\alpha_i},
\bk \in \langle I \rangle_\alpha$ and $r \geq 0$.
The proof that these are well-defined bimodule homomorphisms
and that the sequence is exact is essentially the same as the prooof of
\cite[Th.~5.13]{BWWa1} so we omit any further details here.
To deduce the theorem from the claim, 
it remains to observe by
\cref{standardize} that the bimodule in the third position in the short exact sequence is projective
as a right $\H_\alpha$-module. Hence, this is a split short exact sequence of right $\H_\alpha$-modules, so it remains exact on applying $-\otimes_{\H_\alpha} V$ for any graded left $\H_\alpha$-module $V$. Thus, it defines an exact sequence of functors.
\end{proof}

\subsection{Categorification of the iSerre relation}\label{serrecat}

We continue to assume that $\UU^\imath$ is non-degenerate and $\kk_0$ is a field.
The two theorems in this subsection show that the analog of \cref{secretary} holds in $K_0\left(\Kar(\UU^\imath_q)\right)$.

First, we treat the case $i \neq \tau j$.
Explicit formulae for the differential $d$ and the splitting $s$ in the following theorem can be found in the proof.

\begin{theo}\label{wheniserre}
Suppose that $i \neq j$ and $i \neq \tau j$,
and let $m := 1-a_{i,j}$.
For any $\lambda \in X^\imath$, 
there is a split exact complex
\begin{equation}\label{complex}
\begin{tikzcd}
0\arrow[r,"d_{m+1}"]&\arrow[l,dotted,bend left,"s_m"]
B_i^{(m)} B_j \one_\lambda \arrow[r,"d_{m}"]&\arrow[l,dotted,bend left,"s_{m-1}"]\phantom{B_i}\cdots\phantom{B_i}\arrow[r,"d_2"]& \arrow[l,dotted,bend left,"s_{1}"]B_i B_j B_i^{(m-1)}\one_\lambda\arrow[r,"d_1"]&
\arrow[l,dotted,bend left,"s_{0}"]
B_j B_i^{(m)}\one_\lambda\arrow[r,"d_0"]&\arrow[l,dotted,bend left,"s_{-1}"]0
 \end{tikzcd}
\end{equation}
in
$\HOM_{\Kar(\UU_q^\imath)}(\lambda,\lambda-m\alpha_i - \alpha_j)$.
\end{theo}

\begin{proof}
When $i \neq \tau i$, this is proved
in exactly the same way as the analogous result for 2-quantum groups in \cite[Prop.~4.2]{Rou} (generalizing \cite[Cor.~7]{KL2}).

Now assume that $i = \tau i$. The relation looks like the usual Serre relation, but it is not the same since the idivided power is more complicated. However, the same proof strategy also works in the case $i=\tau i$. This is quite subtle since the nil-Brauer category is fussier than the nil-Hecke algebras, so we go through the details.
We define 2-morphisms in $\Kar(\UU^\imath_q)$ by
\begin{align*}
d_n := 
\!\begin{tikzpicture}[iQ,centerzero]
\draw[ultra thick] (-.4,-.8)\botlabel{i^n\:} to (-.4,.8)\toplabel{i^{n\!-\!1}\:\:};
\draw[ultra thick] (.4,-.8) \botlabel{\quad i^{m\!-\!n}}to (.4,.8)\toplabel{\quad\:\: i^{m\!-\!n\!+\!1}};
\draw (-.425,-.8) to (-.425,-.3) to[looseness=1.5,out=135,in=-135] (.37,.3)
to (.37,.8);
\draw (0,-.8)\botlabel{j} to (0,.8);
\cross{-.4,.5};
\multcloseddot{-.4,.7}{east}{\rho};
\cross{-.405,-.7};
\multcloseddot{-.405,-.5}{east}{\rho};
\cross{.395,.5};
\multcloseddot{.395,.7}{west}{\rho};
\cross{.4,-.7};
\multcloseddot{.4,-.5}{west}{\rho};
\region{.7,0}{\lambda};
\end{tikzpicture}
\!\stackrel{\cref{essence}}{=}\!\begin{tikzpicture}[iQ,centerzero]
\draw[ultra thick] (-.4,-.8)\botlabel{i^n\:} to (-.4,.8)\toplabel{i^{n\!-\!1}\:\:};
\draw[ultra thick] (.4,-.8) \botlabel{\quad i^{m\!-\!n}}to (.4,.8)\toplabel{\quad\:\: i^{m\!-\!n\!+\!1}};
\draw (-.425,-.8) to (-.425,-.3) to[looseness=1.5,out=135,in=-135] (.37,.3)
to (.37,.8);
\draw (0,-.8)\botlabel{j} to (0,.8);
\cross{-.4,.5};
\multcloseddot{-.4,.7}{east}{\rho};
\cross{.395,.5};
\multcloseddot{.395,.7}{west}{\rho};
\region{.7,0}{\lambda};
\end{tikzpicture}
\!\!: B_i^{(n)}  B_j B_i^{(m-n)}\one_\lambda &\Rightarrow
B_i^{(n-1)} B_j B_i^{(m-n+1)}\one_\lambda,\\
s_n := 
(-1)^n t_{i,j}^{-1}\!
\begin{tikzpicture}[centerzero,iQ]
\draw[ultra thick] (-.4,-.8)\botlabel{i^n\:\:} to (-.4,.8)\toplabel{\quad i^{n\!+\!1}};
\draw[ultra thick] (.4,-.8)\botlabel{\quad i^{m\!-\!n}} to (.4,.8)\toplabel{\quad\:\: i^{m\!-\!n\!-\!1}};
\draw (.43,-.8) to (.43,-.3) to[looseness=1.5,out=45,in=-45] (-.375,.3)
to (-.375,.8);
\draw (0,-.8) \botlabel{j} to (0,.8);
\cross{-.395,.5};
\multcloseddot{-.395,.7}{east}{\rho};
\cross{-.4,-.7};
\multcloseddot{-.4,-.5}{east}{\rho};
\cross{.4,.5};
\multcloseddot{.405,.7}{west}{\rho};
\cross{.408,-.7};
\multcloseddot{.408,-.5}{west}{\rho};
\region{.7,0}{\lambda};
\end{tikzpicture}\!\!
\stackrel{\cref{essence}}{=}\!
(-1)^n t_{i,j}^{-1}\!
\begin{tikzpicture}[centerzero,iQ]
\draw[ultra thick] (-.4,-.8)\botlabel{i^n\:\:} to (-.4,.8)\toplabel{\quad i^{n\!+\!1}};
\draw[ultra thick] (.4,-.8)\botlabel{\quad i^{m\!-\!n}} to (.4,.8)\toplabel{\quad\:\: i^{m\!-\!n\!-\!1}};
\draw (.432,-.8) to (.432,-.3) to[looseness=1.5,out=45,in=-45] (-.375,.3)
to (-.375,.8);
\draw (0,-.8)\botlabel{j} to (0,.8);
\cross{-.395,.5};
\multcloseddot{-.395,.7}{east}{\rho};
\cross{.4,.5};
\multcloseddot{.4,.7}{west}{\rho};
\region{.7,0}{\lambda};
\end{tikzpicture}
\!\!: B_i^{(n)}  B_j B_i^{(m-n)}\one_\lambda &\Rightarrow
B_i^{(n+1)} B_j B_i^{(m-n-1)}\one_\lambda,
\end{align*}
interpreting
$B_i^{(n)} B_j B_i^{(m-n)} \one_\lambda$
(hence, 2-morphisms with this as their domain or range)
as $0$ if $n < 0$ or $n > m$.

The following calculation checks that $d$ defines a differential making \cref{complex} into a complex:
\begin{align*}
d_n \circ d_{n+1} &= 
\begin{tikzpicture}[centerzero,iQ,scale=1.2]
\draw (-.42,-1.4) to (-.42,.3) to[looseness=1.5,out=135,in=-135] (.365,.9)
to (.365,1.4);
\cross{-.4,1.1};
\multcloseddot{-.4,1.3}{east}{\rho};
\cross{.3875,1.1};
\closeddot{.375,1.3};
\multcloseddot{.4,1.3}{west}{\rho};
\draw[ultra thick] (-.4,-1.4) \botlabel{i^{n\!+\!1}\:\:}to (-.4,1.4) \toplabel{i^{n\!-\!1}\:\:};
\draw[ultra thick] (.4,-1.4) \botlabel{\quad\:\: i^{m\!-\!n\!-\!1}}to (.4,1.4)\toplabel{\quad\:\:i^{m\!-\!n\!+\!1}};
\draw (-.435,-1.4) to (-.435,-1.1) to[looseness=1.5,out=135,in=-135] (.38,-.5)
to (.38,1.4);
\draw (0,-1.4)\botlabel{j} to (0,1.4);
\cross{-.4,-.2};
\multcloseddot{-.4,0}{east}{\rho};
\cross{.4,-.2};
\multcloseddot{.4,0}{west}{\rho};
\region{.9,0}{\lambda};
\end{tikzpicture}
\stackrel{\cref{essence}}{=}
\begin{tikzpicture}[centerzero,iQ,scale=1.2]
\draw (-.42,-1.4) to (-.42,.3) to[looseness=1.5,out=135,in=-135] (.365,.9)
to (.365,1.4);
\cross{-.4,1.1};
\multcloseddot{-.4,1.3}{east}{\rho};
\cross{.3875,1.1};
\closeddot{.375,1.3};
\multcloseddot{.4,1.3}{west}{\rho};
\draw[ultra thick] (-.4,-1.4) \botlabel{i^{n\!+\!1}\:\:}to (-.4,1.4) \toplabel{i^{n\!-\!1}\:\:};
\draw[ultra thick] (.4,-1.4) \botlabel{\quad\:\: i^{m\!-\!n\!-\!1}}to (.4,1.4)\toplabel{\quad\:\:i^{m\!-\!n\!+\!1}};
\draw (-.435,-1.4) to (-.435,-1.1) to[looseness=1.5,out=135,in=-135] (.38,-.5)
to (.38,1.4);
\draw (0,-1.4)\botlabel{j} to (0,1.4);
\region{.7,0}{\lambda};
\end{tikzpicture}
=\begin{tikzpicture}[centerzero,iQ,scale=1.2]
\multcloseddot{-.4,1.3}{east}{\rho};
\cross{.3875,1.1};
\closeddot{.375,1.3};
\multcloseddot{.4,1.3}{west}{\rho};
\draw[ultra thick] (-.4,-1.4) \botlabel{i^{n\!+\!1}\:\:}to (-.4,1.4) \toplabel{i^{n\!-\!1}\:\:};
\draw[ultra thick] (.4,-1.4) \botlabel{\quad\:\: i^{m\!-\!n\!-\!1}}to (.4,1.4)\toplabel{\quad\:\:i^{m\!-\!n\!+\!1}};
\draw (-.42,-1.4) to (-.42,-.2) to[looseness=1.5,out=135,in=-135] (.365,.5)
to (.365,1.4);
\cross{-.4,1.1};
\draw (-.435,-1.4) to (-.435,-.4) to[looseness=1.7,out=135,in=-135] (.38,.4)
to (.38,1.4);
\draw (0,-1.4)\botlabel{j} to (0,1.4);
\region{.7,0}{\lambda};
\end{tikzpicture}
\stackrel{\cref{ibraid}}{=}\begin{tikzpicture}[centerzero,iQ,scale=1.2]
\multcloseddot{-.4,1.3}{east}{\rho};
\cross{.3875,1.1};
\closeddot{.375,1.3};
\multcloseddot{.4,1.3}{west}{\rho};
\draw[ultra thick] (-.4,-1.4) \botlabel{i^{n\!+\!1}\:\:}to (-.4,1.4) \toplabel{i^{n\!-\!1}\:\:};
\draw[ultra thick] (.4,-1.4) \botlabel{\quad\:\: i^{m\!-\!n\!-\!1}}to (.4,1.4)\toplabel{\quad\:\:i^{m\!-\!n\!+\!1}};
\draw (-.42,-1.4) to (-.42,-.6) to[looseness=1.1,out=135,in=-135] (.365,.5)
to (.365,1.4);
\cross{-.4,1.1};
\draw (-.435,-1.4) to (-.435,-.75) to[looseness=1.9,out=135,in=-155] (.38,.3)
to (.38,1.4);
\draw (0,-1.4)\botlabel{j} to (0,1.4);
\region{.7,0}{\lambda};
\end{tikzpicture}\\&
\stackrel{\cref{ibraid}}{=}\begin{tikzpicture}[centerzero,iQ,scale=1.2]
\multcloseddot{-.4,1.3}{east}{\rho};
\cross{.3875,1.1};
\closeddot{.375,1.3};
\multcloseddot{.4,1.3}{west}{\rho};
\draw[ultra thick] (-.4,-1.4) \botlabel{i^{n\!+\!1}\:\:}to (-.4,1.4) \toplabel{i^{n\!-\!1}\:\:};
\draw[ultra thick] (.4,-1.4) \botlabel{\quad\:\: i^{m\!-\!n\!-\!1}}to (.4,1.4)\toplabel{\quad\:\:i^{m\!-\!n\!+\!1}};
\draw (-.42,-1.4) to (-.42,-.6) to[looseness=1.1,out=135,in=-145] (0,0) to[looseness=1,out=35,in=-110] (.365,.7)
to (.365,1.4);
\cross{-.4,1.1};
\draw (-.435,-1.4) to (-.435,-.8) to[looseness=1.8,out=135,in=-145] (.38,.4)
to (.38,1.4);
\draw (0,-1.4)\botlabel{j} to (0,1.4);
\region{.7,0}{\lambda};
\end{tikzpicture}
-\delta_{i,\tau i} 
r_{i,j}\!\!\!\!\!\!\!
\begin{tikzpicture}[centerzero,iQ,scale=1.2]
\multcloseddot{-.4,1.3}{east}{\rho};
\cross{.3875,1.1};
\closeddot{.375,1.3};
\multcloseddot{.4,1.3}{west}{\rho};
\cross{-.4,1.1};
\draw[ultra thick] (-.4,-1.4) \botlabel{i^{n\!+\!1}\:\:}to (-.4,1.4) \toplabel{i^{n\!-\!1}\:\:};
\draw[ultra thick] (.4,-1.4) \botlabel{\quad\:\: i^{m\!-\!n\!-\!1}}to (.4,1.4)\toplabel{\quad\:\:i^{m\!-\!n\!+\!1}};
\draw (-.42,-1.4) to (-.42,-.8) to[looseness=1,out=150,in=-30] (-.16,0) to[looseness=1,out=150,in=150] (-.435,-1) to (-.435,-1.4);
\draw (.375,1.3) to (.375,.5) to [out=-135,in=10,looseness=1] (.16,0) to [out=190,in=-135,looseness=1] (.36,.7) to (.36,1.3);
\draw (0.1,-1.4)\botlabel{j} to[out=90,in=-60,looseness=1] (0,-0.08) to [out=120,in=-90,looseness=1](-0.1,1.4);
\region{.7,0}{\lambda};
\pinpinpinmid{.1,0}{-0.02,0}{-.15,0}{-1.7,.8}{
\frac{Q_{i,j}(x,y)-Q_{i,j}(z,y)}{x-z}};
\end{tikzpicture}
=0.
\end{align*}
The final equality needs some explanation. The first term is zero because it involves the square of a crossing of strings of color $i$, which is 0 by \cref{iquadratic}.
The second term is zero because
$\frac{Q_{i,j}(x,y)-Q_{i,j}(z,y)}{x-z}$
is a linear combination of monomials of the form
$x^a y^b z^c$ with $a+c \leq m-2$, so either $a < n-1$
or $c \leq m-n-1$. The term arising
from the monomial $x^a y^b z^c$ is zero because
\begin{align*}
\begin{tikzpicture}[centerzero,scale=1.2]
\draw[ultra thick] (.45,.4) to[looseness=1.88,out=-160,in=160] (.45,-.4)\botlabel{\quad i^{n\!-\!1}};
\draw (-0.05,.4) to[looseness=1.88,out=-20,in=20] (-0.05,-.4)\botlabel{i};
\cross{0.05,0};
\multcloseddot{.37,0}{west}{a};
\end{tikzpicture}\!\!&=0\text{ if $a < n-1$,}\qquad
\begin{tikzpicture}[centerzero,scale=1.2]
\draw[-] (0.25,-.3) \botlabel{i} to[looseness=1,out=120,in=-45] (0.05,-.1) to [looseness=1,out=135,in=0] (-.3,.3) to [looseness=1,out=180,in=90] (-.5,.1) to[looseness=1,out=-90,in=180] (-.3,-.1) to [looseness=1,out=0,in=-135] (0.05,.3) to [looseness=1,out=45,in=-120](0.25,.5);
\multcloseddot{-.5,.1}{east}{c};
\draw[ultra thick] (-0.5,.5) to [out=-50,in=50,looseness=1] (-0.5,-.3)\botlabel{i^{m\!-\!n\!-\!1}\:\:\:};
\cross{-.35,.1};
\end{tikzpicture}\!\! =0\quad\text{if $c \leq m-n-1$}.
\end{align*}
These identities are \cite[Lem.~4.2 and Cor.~4.3]{BWWa1}.

Finally, we 
prove that $s$ is a splitting. We need to show that
\begin{align*}
s_{n-1} \circ d_n + d_{n+1} \circ s_n &=
(-1)^{n-1} t_{i,j}^{-1}
\left(\begin{tikzpicture}[iQ,centerzero,scale=1.2]
\draw[ultra thick] (-.4,0) to (-.4,1.6)\toplabel{\: i^{n}};
\draw[ultra thick] (.387,0) to (.387,1.6)\toplabel{\quad i^{m\!-\!n}};
\draw (.416,0) to (.416,.5) to[looseness=1.5,out=45,in=-45] (-.375,1.1)
to (-.375,1.6);
\draw (0,0) to (0,1.6);
\cross{-.395,1.3};
\multcloseddot{-.395,1.5}{east}{\rho};
\cross{.387,1.3};
\multcloseddot{.387,1.5}{west}{\rho};
\draw[ultra thick] (-.4,-1.4)\botlabel{i^n\:} to (-.4,0);
\draw[ultra thick] (.4,-1.4) \botlabel{\quad i^{m\!-\!n}}to (.4,0);
\draw (-.425,-1.4) to (-.425,-1.1) to[looseness=1.5,out=135,in=-135] (.375,-.5)
to (.375,0);
\draw (0,-1.4)\botlabel{j} to (0,0);
\cross{-.4,-.2};
\multcloseddot{-.4,0}{east}{\rho};
\cross{.395,-.2};
\multcloseddot{.395,0}{west}{\rho};
\region{.9,0}{\lambda};
\end{tikzpicture}
-
\begin{tikzpicture}[iQ,centerzero,scale=1.2]
\draw[ultra thick] (-.4,0) to (-.4,1.6)\toplabel{i^{n}\:};
\draw[ultra thick] (.4,0) to (.4,1.6)\toplabel{\quad i^{m\!-\!n}};
\draw (-.425,0) to (-.425,.5) to[looseness=1.5,out=135,in=-135] (.375,1.1)
to (.375,1.6);
\draw (0,0) to (0,1.6);
\cross{-.4,1.3};
\multcloseddot{-.4,1.5}{east}{\rho};
\cross{.395,1.3};
\multcloseddot{.395,1.5}{west}{\rho};
\region{.9,0}{\lambda};
\draw[ultra thick] (-.415,-1.4)\botlabel{i^n\:\:} to (-.415,0);
\draw[ultra thick] (.4,-1.4)\botlabel{\quad i^{m\!-\!n}} to (.4,0);
\draw (.425,-1.4) to (.425,-1.1) to[looseness=1.5,out=45,in=-45] (-.387,-.5)
to (-.387,0);
\draw (0,-1.4)\botlabel{j} to (0,0);
\cross{-.405,-.2};
\multcloseddot{-.405,0}{east}{\rho};
\cross{.4,-.2};
\multcloseddot{.4,0}{west}{\rho};
\end{tikzpicture}\right)
= \id_{B_i^{(n)} B_j B_i^{(m-n)}\one_\lambda}.
\end{align*}
To see this, we rewrite the expression inside the parentheses: it equals
\begin{align*}
\begin{tikzpicture}[iQ,centerzero,scale=1.2]
\draw[ultra thick] (-.4,-.8)\botlabel{i^n} to (-.4,.8);
\draw[ultra thick] (.4,-.8) \botlabel{\quad i^{m\!-\!n}}to (.4,.8);
\draw (-.425,-.8) to (-.425,-.5) to[looseness=1.5,out=135,in=-135] (.37,.15)
to (.37,.8);
\draw (.432,-.8) to (.432,-.3) to[looseness=1.5,out=45,in=-45] (-.375,.3)
to (-.375,.8);
\draw (0,-.8)\botlabel{j} to (0,.8);
\cross{-.395,.5};
\multcloseddot{-.395,.7}{east}{\rho};
\cross{.395,.5};
\multcloseddot{.395,.7}{west}{\rho};
\region{.7,0}{\lambda};
\end{tikzpicture}
-
\begin{tikzpicture}[centerzero,iQ,scale=1.2]
\draw[ultra thick] (-.4,-.8)\botlabel{i^n\:\:} to (-.4,.8);
\draw[ultra thick] (.4,-.8)\botlabel{\quad i^{m\!-\!n}} to (.4,.8);
\draw (-.425,-.8) to (-.425,-.3) to[looseness=1.5,out=135,in=-135] (.37,.3)
to (.37,.8);
\draw (.432,-.8) to (.432,-.5) to[looseness=1.5,out=45,in=-45] (-.375,.15)
to (-.375,.8);
\draw (0,-.8)\botlabel{j} to (0,.8);
\cross{-.395,.5};
\multcloseddot{-.395,.7}{east}{\rho};
\cross{.4,.5};
\multcloseddot{.4,.7}{west}{\rho};
\region{.7,0}{\lambda};
\end{tikzpicture}
&\stackrel{\cref{ibraid}}{=}
\begin{tikzpicture}[iQ,centerzero,scale=1.2]
\draw[ultra thick] (-.4,-.8)\botlabel{i^n} to (-.4,.8);
\draw[ultra thick] (.4,-.8) \botlabel{\quad i^{m\!-\!n}}to (.4,.8);
\draw (-.425,-.8) to (-.425,-.55) to[looseness=2.5,out=145,in=-35] (-.375,.25)
to (-.375,.8);
\draw (.432,-.8) to (.432,-.55) to[looseness=2.5,out=35,in=-145] (.375,.25)
to (.375,.8);
\draw (0,-.8)\botlabel{j} to (0,.8);
\cross{-.395,.5};
\multcloseddot{-.395,.7}{east}{\rho};
\cross{.395,.5};
\multcloseddot{.395,.7}{west}{\rho};
\pinpinpin{.2,0.05}{0,.05}{-.2,0.05}{-1.8,0.05}{\frac{Q_{i,j}(x,y)-Q_{i,j}(z,y)}{x-z}};
\region{.7,0}{\lambda};
\end{tikzpicture}\ .
\end{align*}
Like in the previous paragraph, $\frac{Q_{i,j}(x,y)-Q_{i,j}(z,y)}{x-z}$
is a linear combination of monomials of the form $x^a y^b z^c$
with $a+c \leq m-2$, and the monomials with $a+c=m-2$
come from the expansion 
$t_{i,j} (x^{m-1}-z^{m-1})/(x-z) = \sum_{a+c=m-2} t_{i,j} x^a y^c$.
Now we use the identities
\begin{align*}
\begin{tikzpicture}[centerzero,scale=1.2]
\draw[ultra thick] (.45,.4) to[looseness=1.88,out=-160,in=160] (.45,-.4)\botlabel{\quad i^{n\!-\!1}};
\draw (-0.05,.4) to[looseness=1.88,out=-20,in=20] (-0.05,-.4)\botlabel{i};
\cross{0.05,0};
\multcloseddot{.37,0}{west}{a};
\end{tikzpicture}\!\!&=\delta_{a,n-1} (-1)^{n-1}
\begin{tikzpicture}[centerzero,scale=1.2]
\draw[ultra thick] (.45,.4) to (.45,-.4)\botlabel{i^{n}};
\cross{.45,0};
\end{tikzpicture}
\text{ if $a \leq n-1$,}\qquad
\begin{tikzpicture}[centerzero,scale=1.2]
\draw[ultra thick] (-.45,.4) to[looseness=1.88,out=-20,in=20] (-.45,-.4)\botlabel{i^{m\!-\!n\!-\!1}\quad};
\draw (0.05,.4) to[looseness=1.88,out=-160,in=160] (0.05,-.4)\botlabel{i};
\cross{-0.05,0};
\multcloseddot{-.37,0}{east}{c};
\end{tikzpicture}\!\!\!\!\!\!\!\!&=\delta_{c,m-n-1}
\begin{tikzpicture}[centerzero,scale=1.2]
\draw[ultra thick] (.45,.4) to (.45,-.4)\botlabel{i^{m\!-\!n}};
\cross{.45,0};
\end{tikzpicture}
\text{ if $c \leq m-n-1$}
\end{align*}
from \cite[Lem.~4.2]{BWWa1} (this treats the new case when $i = \tau i$
but these identities are also well known when $i \neq \tau i$).
We deduce that the terms coming from all of the monomials except 
for $t_{i,j} x^{n-1} z^{m-n-1}$ evaluate to 0, and the remaining term contributes $(-1)^{n-1} t_{i,j} \id_{B_i^{(n)} B_j B_i^{(m-n)}\one_\lambda}$, as required to finish the argument.
\end{proof}

\begin{cor}\label{iserrepart1}
If $i \neq j$, $i \neq \tau j$ and $m := 1-a_{i,j}$, we have that
$$
\bigoplus^{m}_{\substack{n=0 \\ n\text{ even}}}
B_i^{(n)} B_j B_i^{(m-n)} \one_\lambda \cong
\bigoplus^{m}_{\substack{n=0\\ n\text{ odd}}}
B_i^{(n)} B_j B_i^{(m-n)} \one_\lambda.
$$
\end{cor}

For the categorical iSerre relation in the case $i\neq j = \tau j$,
we will work in terms of modules
over the path algebra $\H = \HH$ from \cref{pathalgebra}, adopting all of the module-theoretic notation
from the previous subsection.

\begin{lem}\label{landingintaipei}
Let $\H = \HH$ for $\lambda \in X^\imath$.
Suppose that $i \neq j$ and $i = \tau j$.
For $m \geq 1$ and $0 \leq n \leq m$, there is a short exact sequence of graded $\H$-modules
$$
0\longrightarrow P(i^{(m-1)})^{\oplus f^\lambda_{n,m;i}(q)} \longrightarrow P(i^{(n)} j i^{(m-n)}) \longrightarrow \Delta(i^{(n)} j i^{(m-n)})
\longrightarrow 0
$$
where $f^\lambda_{n,m;i}(q) \in \N\lround q \rround$ is defined by \cref{veday}.
\end{lem}

\begin{proof}
There is just one smaller weight that $m\alpha_i + \alpha_j$ in
the poset $\Lambda$, namely, $(m-1) \alpha_i$.
Since $\H_{(m-1) \alpha_i}$ is a nil-Hecke algebra, the set
$\B_{(m-1)\alpha_i}$ contains just one element which we denote by $\bb$,
and $\Delta(\bb) = \Delta(i^{(m-1)})$.
Also let $\bi := i^{(n)} j i^{(m-n)}$ for short.
By \cref{halfwaytotaipei}, $P(\bi)$ is isomorphic to the projective cover of $\Delta(\bi)$ plus
a finite direct sum of grading shifts of $P(\bb)$.
It follows that $P(\bi)$ has a two-step $\Delta$-flag with
top section isomorphic to $\Delta(\bi)$ and bottom
section isomorphic to $\Delta(i^{(m-1)})^{\oplus g(q)}$ for some $g(q) \in \N\lround q\rround$. Since $\Delta(i^{(m-1)}) \cong P(i^{(m-1)})$ by \cref{halfwaytotaipei} again, it just remains to prove that $g(q) = f^\lambda_{n,m;i}(q^{-1})$.
Using \cref{deltaflags} for the first equality and $P(i^n j i^{m-n}) \cong P(\bi)^{\oplus [n]^!_{q_i} [m-n]^!_{q_i}}$ for the second,
we have that
\begin{align*}
g(q) &= \grdim_{q^{-1}} \Hom_{\H}(P(\bi),
\bar\nabla(\bb))= 
\frac{1}{[n]^!_{q_i} [m-n]^!_{q_i}}
\grdim_{q^{-1}} 1_{i^n j i^{m-n}} \bar\nabla(\bb).
\end{align*}
By \cref{standardize}, $\bar\nabla(\bb)$
has an explicit basis consisting of the
functions $\varphi_{y,v}$ 
for $y \in 
Y(i^{m-1},i^n j i^{m-n})$ and $v$ running over a basis for
$L_{(m-1)\alpha_i}(\bb)$, this being the function that maps
$\bar y \mapsto v$ and all other basis vectors
of $\bar e_{(m-1)\alpha_i} \H_{\geq (m-1)\alpha_i}$ to $0$.
We have that $\deg(\varphi_{y,v}) = \deg(v)-\deg(y)$,
the graded dimension of $L{(m-1) \alpha_i}(\bb)$ is $[m-1]^!_{q_i}$,
and the set
$Y(i^{m-1},i^n j i^{m-n})$ may be chosen so that its elements are the mirror images of the diagrams in \cref{somediagrams}.
It follows that $g(q)=
[m-1]^!_{q_i} \times \sum_D q^{\deg(D)}$ summing over these diagrams,
that is, $g(q) = [m-1]_{q_i}^! f(q) / [n]^!_{q_i} [m-n]^!_{q_i}$
where $f(q)$ was computed in \cref{Zset}.
This equals $f^\lambda_{n,m;i}(q)$.
\end{proof}

Now we establish the categorical iSerre relation in the case $i = \tau j$.
Unlike \cref{wheniserre}, we do not know formulae for the differential or an explicit splitting.
 
\begin{theo}\label{wheniserre2}
Suppose that $i \neq j$ and $i = \tau j$. Fix $\lambda \in X^\imath$ and
let $m := 1 - a_{i,j}$. 
Let $f^\lambda_{n,m;i}(q) \in \N\lround q \rround$ be as in \cref{veday}.
There is a split exact complex
\begin{multline*}
\qquad0\longrightarrow
B_i^{(m-1)} \one_\lambda^{\oplus f^\lambda_{m,m;i}(q)}
\longrightarrow
B_i^{(m)} B_j \one_\lambda\oplus 
B_i^{(m-1)} \one_\lambda^{\oplus f^\lambda_{m-1,m;i}(q)}
\longrightarrow \\
\cdots\longrightarrow
B_i^{(n)} B_j B_i^{(m-n)} \one_\lambda\oplus
B_i^{(m-1)} \one_\lambda^{\oplus f^\lambda_{n-1,m;i}(q)}
\longrightarrow\cdots\\
\longrightarrow  B_i B_j B_i^{(m-1)} \one_\lambda
\oplus B_i^{(m-1)} \one_\lambda^{\oplus f^\lambda_{0,m;i}(q)}
\longrightarrow
B_j B_i^{(m)} \one_\lambda
\longrightarrow 0\qquad
\end{multline*}
in 
$\HOM_{\Kar(\UU_q^\imath)}(\lambda,\lambda-m\alpha_i - \alpha_j)$.
\end{theo}

\begin{proof}
Let 
$\alpha := m \alpha_i + \alpha_j \in \Lambda$.
Applying the standardization functor $\jmath^\alpha_!$
to the 
split exact complex in $\gmod{\H_\alpha}$ from \cite[Cor.~7]{KL2} and 
\cite[Prop.~4.2]{Rou} (which is proved in the same way as the $i \neq \tau i$ case
of \cref{wheniserre}) gives a split exact sequence which is the bottom row of the following diagram:
$$
\begin{tikzcd}
&0\arrow[d]&&0\arrow[d]&0\arrow[d]\\
0\arrow[r]&
P(i^{(m-1)})^{\oplus f_{m,m;i}^\lambda(q)}  \arrow[r]\arrow[d]&
\!\cdots\!\arrow[r]&
P(i^{(m-1)})^{\oplus f_{1,m;i}^\lambda(q)} 
\arrow[r]\arrow[d]&
P(i^{(m-1)})^{\oplus f_{0,m;i}^\lambda(q)}  \arrow[r]\arrow[d]&
0\\
0\arrow[r]&
P(i^{(m)} j) \arrow[r]\arrow[d]&
\!\cdots\!\arrow[r]&
P(iji^{(m-1)})\arrow[d]
\arrow[r]&P(ji^{(m)}) \arrow[d]\arrow[r]&
0\\
0\arrow[r]&
\Delta(i^{(m)} j) \arrow[r]\arrow[d]&
\!\cdots\!\arrow[r]&
\Delta(iji^{(m-1)})\arrow[d]
\arrow[r]&\Delta(ji^{(m)}) \arrow[r]\arrow[d]&
0\\
&0&&0&0
\end{tikzcd}
$$
The columns are the short exact sequences arising from \cref{landingintaipei}.
The maps along the bottom row lift to give horizontal maps (not necessarily differentials) along the middle and top rows making the diagram commute.
While this is not a double complex, we can still define a total complex for it of the form
\begin{multline*}
\qquad0\longrightarrow
P(i^{(m-1)})^{\oplus f^\lambda_{m,m;i}(q)}
\longrightarrow
P(i^{(m)} j)\oplus 
P(i^{(m-1)})^{\oplus f^\lambda_{m-1,m;i}(q)}
\longrightarrow \\
\cdots\longrightarrow
P(i^{(n)} j i^{(m-n)})\oplus
P(i^{(m-1)})^{\oplus f^\lambda_{n-1,m;i}(q)}
\longrightarrow\cdots\\
\longrightarrow  P(iji^{(m-1)})
\oplus P(i^{(m-1)})^{\oplus f^\lambda_{0,m;i}(q)}
\longrightarrow
P(j i^{(m)})
\longrightarrow 0\qquad
\end{multline*}
by adding maps $P(i^{(n)}ji^{(m-n)})\to P(i^{(m-1)})^{\oplus f^\lambda_{n-2,m;i}(q)}$ in the diagram above to make the total differential square to 0; these maps can be obtained by noting that the composition of two horizontal maps must lie in the image of the vertical map, since the bottom row is a complex and the columns are exact.
The obvious projection to the bottom row is an isomorphism, so the total complex is exact. Since all of the modules involved are finitely generated and projective, the complex is again split exact.
The complex in the statement of the theorem is obtained from this using the equivalence $\cov$.
\end{proof}

\begin{cor}\label{iserrepart2}
For $i \neq j$, $i = \tau j$ and $m := 1-a_{i,j}$,
we have that
$$
\bigoplus_{\substack{n=0 \\ n\text{ even}}}^{m}
B_i^{(n)} B_j B_i^{(m-n)} \one_\lambda
\oplus
\bigoplus_{\substack{n=0 \\ n\text{ odd}}}^{m}
B_i^{(m-1)} \one_\lambda^{\oplus f^\lambda_{n,m;i}(q)}
\cong
\bigoplus_{\substack{n=0\\ n\text{ odd}}}^{m}
B_i^{(n)} B_j B_i^{(m-n)} \one_\lambda
\oplus
\bigoplus_{\substack{n=0\\ n\text{ even}}}^{m}
B_i^{(m-1)} \one_\lambda^{\oplus f^\lambda_{n,m;i}(q)}\!.
$$
\end{cor}

\subsection{Main theorem}\label{mainthm}

Now everything is in place in order to be able to state and 
prove the main theorem. 
Recall the isomorphism $\KL$ from \cref{KLequiv} and the equivalences
$\cov$ from \cref{yonedaop,otheryonedaop}.
Setting 
$\theta_\bb := \KL^{-1}\big([\Theta_\bb]\big)$,
we obtain a basis $\big\{\theta_\bb\:\big|\:\bb \in \B\big\}$ for $\f_\Z$
as a free $\Z[q,q^{-1}]$-module
which we call the {\em orthodox basis}\footnote{This terminology was introduced in \cite{canonical}.}.
Let $\dot\U^\imath_{\Z\lround q\rround} := \dot\U^\imath_{\Z} \otimes_{\Z[q,q^{-1}]} \Z\lround q \rround$
and $\dot\U^\imath_{\Q\lround q\rround} := \dot\U^\imath_{\Z\lround q \rround} \otimes_\Z \Q$.
There is a natural inclusion
$\dot\U^\imath_\Z \hookrightarrow \dot\U^\imath_{\Z\lround q \rround}$,
and both $\dot\U^\imath_{\Z\lround q \rround}$ 
and $\dot\U^\imath$ embed into
$\dot\U^\imath_{\Q\lround q \rround}$.
Also let $\jmath_\lambda, \widetilde{\jmath}_\lambda:\dot\U^\imath 1_\lambda\stackrel{\sim}{\rightarrow}\f$ be the $\Q(q)$-linear isomorphisms from \cref{jbversion} and \cref{newpropsofj}; we have added the subscript $\lambda \in X^\imath$ to avoid any confusion as $\lambda$ varies. 

\begin{theo}\label{cthm}
Assume that $\kk_0$ is a field and the 2-iquantum group
$\UU^\imath$ is non-degenerate in the sense of \cref{defnondeg}.
\begin{enumerate}
\item
There is a unique isomorphism of locally unital $\Z[q,q^{-1}]$-algebras
\begin{align*}
\BWWI:\dot\U^\imath_\Z &\stackrel{\sim}{\rightarrow}
K_0\left(\Kar(\UU_q^\imath)\right)
\end{align*}
mapping
$b_i 1_\lambda \mapsto \left[B_i \one_\lambda\right]$
for each $i \in I$ and $\lambda \in X^\imath$.
It maps $b_\bi 1_\lambda \mapsto B_\bi \one_\lambda$ for any
$\bi \in \llangle I \rrangle$.
Setting
$b_\bb 1_\lambda := \BWWI^{-1}\big([B_\bb \one_\lambda]\big)$,
we obtain a basis $\big\{b_\bb 1_\lambda\:\big|\:\bb \in \B, \lambda \in X^\imath\big\}$ for $\dot\U^\imath_\Z$ as a free $\Z[q,q^{-1}]$-module,
which we call the {\em iorthodox basis}.
\item For $\lambda \in X^\imath$ and $\H = \HH$ as in \cref{pathalgebra},
there is a unique isomorphism of $\Z\lround q \rround$-modules
\begin{align*}
\widetilde{\BWWJ}:\dot \U^\imath_{\Z\lround q\rround} 1_\lambda
&\stackrel{\sim}{\rightarrow}
K_0\left(\standard{\H}\right)
\end{align*}
such that the following diagram commutes:
\begin{equation}\label{bigdiagram}
\begin{tikzcd}
\begin{smallmatrix}
\phantom{\textstyle X}\\{\textstyle \dot\U^\imath_\Z 1_\lambda}\\
{\color{purple}b_\bi 1_\lambda}\quad {\color{darkg}b_\bb 1_\lambda}
\end{smallmatrix}
\arrow[d,hookrightarrow]\arrow[r,"\BWWI" below,"\sim" above]&
\begin{smallmatrix}\phantom{\textstyle X}\\{\textstyle K_0\left(\Kar(\UU^\imath_q)\one_\lambda \right)}\\{\color{purple} [B_\bi \one_\lambda]}\quad {\color{darkg} [B_\bb \one_\lambda]}
\end{smallmatrix}
\arrow[r,"\sim" above,"{[\cov]}" below]&
\begin{smallmatrix}
\phantom{\textstyle{X}}\\{\textstyle K_0\left(\gproj{\HH}\right)}\\{\color{purple} [P(\bi)]}\quad {\color{darkg} [P(\bb)]}
\end{smallmatrix}
\arrow[d,hookrightarrow]\\
\begin{smallmatrix}{\color{purple}b_\bi 1_\lambda}\quad 
 {\color{darkg}b_\bb 1_\lambda}\\
{\textstyle\dot\U^\imath_{\Z\lround q\rround} 1_\lambda}\\
{\color{blue}\delta_\bi 1_\lambda}\quad
{\color{red}\delta_\bb 1_\lambda}
\end{smallmatrix}
\arrow[rr,"\sim" above,"\widetilde{\BWWJ}" below]&&
\begin{smallmatrix}
{\color{purple} [P(\bi)]}\quad {\color{darkg} [P(\bb)]}\\
{\textstyle K_0\left(\standard{\HH}\right)}\\{\color{blue}[\Delta(\bi)]}\quad 
{\color{red}[\Delta(\bb)]}
\end{smallmatrix}\\
\begin{smallmatrix}
{\color{blue}\theta_\bi}\quad {\color{red}\theta_\bb}
\\{\textstyle\f_\Z}\\
\phantom{\textstyle X}
\end{smallmatrix}
\arrow[r,"\sim" above,"\KL" below]
\arrow[u,hookrightarrow,"\widetilde{\jmath}^{\:-1}_\lambda" right]&
\begin{smallmatrix}
{\color{blue}[\Theta_\bi]}\quad
{\color{red}[\Theta_\bb]}\\{\textstyle K_0\left(\Kar(\catQH^\imath_q)\right)}\\
\phantom{\textstyle X}
\end{smallmatrix}
\arrow[r,"\sim" above,"{[\cov]}" below]&
\begin{smallmatrix}
{\color{blue}[P_\alpha(\bi)]}\quad 
{\color{red}[P_\alpha(\bb)]}\\
{\textstyle\bigoplus_{\alpha \in \Lambda}K_0\left(\gproj{\QH^\imath_\alpha}\right)}\\
\phantom{X}
\end{smallmatrix}
\arrow[u,hookrightarrow,"{[\jmath_!^\alpha]}" right]
\end{tikzcd}
\end{equation}
For $\bi \in \llangle I \rrangle$, the element
$\widetilde{\BWWJ}^{-1}\big([\Delta(\bi)]\big) \in \dot\U^\imath_{\Z\lround q\rround}$ is equal to $\delta_\bi = \widetilde{\jmath}_\lambda^{-1}(\theta_\bi)$
defined in \cref{abitwacky}.
Setting
$\delta_\bb 1_\lambda := \widetilde{\BWWJ}^{-1}\big([\Delta_\bb \one_\lambda]\big) = \widetilde{\jmath}^{\:-1}_\lambda(\theta_\bb)$,
we obtain a basis $\big\{\delta_\bb 1_\lambda\:\big|\:\bb \in \B, \lambda \in X^\imath\big\}$ for $\dot\U^\imath_{\Z\lround q\rround}$ as a free 
$\Z\lround q\rround$-module,
which we call the {\em standardized orthodox basis}.
\end{enumerate}
\end{theo}

\begin{proof}
(1) This is essentially the same as the proof of the analogous result for 2-quantum groups in \cite{KL3}. First, we construct $\BWWI$.
Morphism spaces in
$\HOM_{\UU^\imath_q}(\lambda,\kappa)$ 
are graded vector spaces over $\kk_0$ which are locally finite-dimensional (i.e., each graded component is finite-dimensional) 
and bounded below (i.e., they are zero in all but finitely many negative degrees).
This follows from \cref{easyspanning}.
So $\Kar(\UU^\imath_q)$ has the
Krull-Schmidt property.
It follows that $K_0\left(\Kar(\UU^\imath_q)\right)$
is free as a $\Z[q,q^{-1}]$-module with basis given by 
isomorphism classes of
indecomposable objects, up to grading shift.
This means that we can identify 
$K_0\left(\Kar(\UU^\imath_q)\right)$ with a $\Z[q,q^{-1}]$-subalgebra of
$K_0\left(\Kar(\UU^\imath_q)\right) \otimes_{\Z[q,q^{-1}]} \Q(q)$.
Also, by its definition, 
$\dot\U^\imath_\Z$ is a $\Z[q,q^{-1}]$-subalgebra of $\dot\U^\imath$.
Now we define a $\Q(q)$-algebra homomorphism 
$$
\BWWI_\Q:\dot\U^\imath \rightarrow K_0\left(\Kar(\UU^\imath_q)\right) \otimes_{\Z[q,q^{-1}]} \Q(q)
$$ 
on generators by mapping $1_\lambda \mapsto 1_\lambda$ and $b_i 1_\lambda \mapsto b_i 1_\lambda$
for $i \in I$ and $\lambda \in X^\imath$. 
\cref{rentalcar} implies that $b_i^{(n)} 1_\lambda \mapsto b_i^{(n)} 1_\lambda$. 
\cref{iserrepart1,iserrepart2} show
that the defining relations of $\dot\U^\imath$ from \cref{secretary} hold in $K_0\left(\Kar(\UU^\imath_q)\right) \otimes_{\Z[q,q^{-1}]} \Q(q)$ (to see this in the $i = \tau j$ case also needs \cref{nextlem}).
Hence, $\BWWI_\Q$ is well defined.
Then we define 
$\BWWI:\dot\U^\imath_\Z \rightarrow K_0\left(\Kar(\UU^\imath_q)\right)$ to be the restriction of $\BWWI_\Q$. This produces the desired $\Z[q,q^{-1}]$-algebra homomorphism.
The homomorphism $\BWWI$ is surjective thanks to
\cref{internationalhouse}.
Finally we show that it is injective:
by \cref{bcor}, if $y \in \ker \BWWI$ then 
 $\langle x,y\rangle^\imath = 0$ for all $x \in \dot\U^\imath_\Z$, hence, $y = 0$ since the form 
 $\langle\cdot,\cdot\rangle^\imath$ is non-degenerate.

\vspace{2mm}
\noindent
(2)
It is always the case that $K_0(\standard{\H}) \cong 
K_0(\gproj{\H}) \otimes_{\Z[q,q^{-1}]} \Z\lround q\rround$
for any algebra $\H$ with a graded triangular basis.
The existence of a unique isomorphism $\widetilde{\BWWJ}$ making the top square commute follows immediately.
The bottom square in the diagram commutes by \cref{lastgasp}, as noted already in \cref{stench}.
\end{proof}

\begin{rem}\label{pr}
The orthodox basis for $\f_\Z$ is closely related (but not necessarily equal) to the {\em canonical basis} of $\f_\Z$, which we denote by
$\big\{c_\bb\:\big|\:\bb \in \B\big\}$. 
We are using the same labelling set $\B$
both for the orthodox basis and for the canonical basis. This is
justified since both bases are {\em dual perfect bases} whose 
{\em associated crystal} is Kashiwara's crystal $\B(\infty)$, and via
this theory all labelling sets are canonically isomorphic.
By \cite{VV,Rou2}, for symmetric Cartan matrices with geometric parameters as in \cref{tojapan}, and assuming $\operatorname{char}(\kk_0) = 0$,
the orthodox basis and the canonical basis coincide: we have that
$[\Theta_\bb] = c_\bb$ for all $\bb \in \B$.
There is also another basis for $\dot\U^\imath_{\Z\lround q\rround}$, 
\end{rem}

\begin{rem}\label{tediousdiscussion}
The iorthodox basis $\big\{b_\bb 1_\lambda\:\big|\:\bb \in \B, \lambda \in X^\imath\big\}$
defined by \cref{cthm}(1)
should be closely related (but not 
necessarily equal) to the {\em icanonical basis}
for $\dot\U^\imath_\Z$ constructed in \cite{BW18QSP}.
Both bases are related to the standardized orthodox basis from \cref{cthm}(2) in a 
uni-triangular way, so they are naturally indexed by the same set $\B$.
There is also the
{\em standardized canonical basis}
$\big\{\widetilde{\jmath}^{\:-1}_\lambda(c_\bb)\:\big|\:\bb \in \B, \lambda \in X^\imath\big\}$, with similarly uni-triangular transition matrices to the iorthodox and icanonical bases. 
Entries in the transition matrix from the iorthodox to the standardized orthodox basis lie in $\N\lround q\rround$, and off-diagonal 
entries in the transition matrix from the icanonical basis to the standardized canonical basis lie in $q \Z\llbracket q\rrbracket$.
\end{rem}

\begin{rem}
In this section, we have generally preferred to work with finitely generated projectives and standard modules. But it is just as natural to work in terms of finitely cogenerated injectives and 
costandard modules. These fit into a similar commutative 
diagram to \cref{bigdiagram}:
\begin{equation}\label{bigdiagram2}
\begin{tikzcd}
\begin{smallmatrix}
 \phantom{\textstyle X}\\
 {\textstyle \dot\U^\imath_\Z 1_\lambda}\\
{\color{purple}b_\bi 1_\lambda}\quad 
 {\color{darkg}b_\bb 1_\lambda}
\end{smallmatrix}
\arrow[d,hookrightarrow]\arrow[r,"\BWWI" below,"\sim" above]&
\begin{smallmatrix}
\phantom{\textstyle X}\\
{\textstyle K_0\left(\Kar(\UU^\imath_q)\one_\lambda \right)}\\{\color{purple} [B_\bi \one_\lambda]}\quad {\color{darkg} [B_\bb \one_\lambda]}
\end{smallmatrix}
\arrow[rr,"\sim" above,"{[\Nak\circ \cov]}" below]&&
\begin{smallmatrix}
\phantom{\textstyle{X}}\\
{\textstyle K_0\left(\ginj{\HH}\right)}\\
\arrow[d,hookrightarrow]\\
{\color{purple} [I(\bi)]}\quad {\color{darkg} [I(\bb)]}
\end{smallmatrix}
\arrow[d,hookrightarrow]
\\
\begin{smallmatrix}
{\color{purple}b_\bi 1_\lambda}\quad 
{\color{darkg}b_\bb 1_\lambda}\\
{\textstyle\dot\U^\imath_{\Z\lround q^{-1}\rround} 1_\lambda}\\
{\color{blue}\atled_\bi 1_\lambda}\quad
{\color{red}\atled_\bb 1_\lambda}
\end{smallmatrix}
\arrow[rrr,"\sim" above,"\BWWJ" below]&&&
\begin{smallmatrix}
{\color{purple} [I(\bi)]}\quad {\color{darkg} [I(\bb)]}
\\
{\textstyle K_0\left(\costandard{\HH}\right)}\\
{\color{blue}[\nabla(\bi)]}\quad 
{\color{red}[\nabla(\bb)]}
\end{smallmatrix}\\
\arrow[u,"\jmath^{-1}_\lambda" right,hookrightarrow]
\begin{smallmatrix}
{\color{blue}\theta_\bi}\quad {\color{red}\theta_\bb}
\\{\textstyle\f_\Z}\\
\phantom{\textstyle X}
\end{smallmatrix}
\arrow[r,"\sim" above,"\KL" below]&
\begin{smallmatrix}
{\color{blue}[\Theta_\bi]}\quad
{\color{red}[\Theta_\bb]}\\{\textstyle K_0\left(\Kar(\catQH^\imath_q)\right)}\\
\phantom{\textstyle X}
\end{smallmatrix}
\arrow[rr,"\sim" above,"{[\Nak\circ\cov]}" below]&&
\begin{smallmatrix}
{\color{blue}[I_\alpha(\bi)]}\quad 
{\color{red}[I_\alpha(\bb)]}
\\
{\textstyle\bigoplus_{\alpha \in \Lambda}K_0\left(\ginj{\QH^\imath_\alpha}\right)}\\
\phantom{X}
\end{smallmatrix}
\arrow[u,hookrightarrow,"{[\jmath_*^\alpha]}" right]
\end{tikzcd}
\end{equation}
The proof is similar to the above and will be omitted.
We call the basis $\atled_\bb\:(\bb \in \B)$ 
for
$\dot\U^\imath_{\Z\lround q^{-1}\rround}$ as a free $\Z\lround q^{-1}\rround$-module obtained by applying $\jmath^{-1}_\lambda$
to $\theta_\bb\:(\bb \in \B)$ 
the {\em costandardized orthodox basis}; we have simply
that $\atled_\bb = \psi^\imath(\delta_\bb)$.
Applying $\jmath^{-1}_\lambda$ instead to the
canonical basis $c_\bb\:(\bb \in \B)$ for $\f$
produces the {\em costandardized canonical basis}.
We mention this in order to point out that
the {\em fused canonical basis} for modified forms of quantum groups defined in \cite{PBW} is the costandardized canonical basis 
in the special case of iquantum groups of diagonal type.
\end{rem}

\begin{rem}[Homomorphism formulae]
The sesquilinear form $\langle \cdot,\cdot\rangle^\imath$
from \cref{newby2} also has a natural interpretation in the categorification framework.
For finitely generated graded projective left 
$\H$-modules $V$ and $W$
corresponding to
$v,w \in \dot\U^\imath_\Z$ under the isomorphism $\cov \circ \BWWI$,
we have that
\begin{equation}\label{homformula1}
\Hom_\H(W,V) \cong \left(\R^\imath\right)^{\oplus \overline{\langle v,w\rangle^\imath}}
\end{equation}
as a graded $\R^\imath$-module.
To prove this, we may assume by \cref{taipeiflight} 
that $V = P(\bi)$ and $W = P(\bj)$
for $\bi,\bj \in \langle I\rangle$, and then it follows by \cref{bcor}.
Also, for $\H$-modules $V$ and $W$ with a $\Delta$-flag and a $\nabla$-flag, respectively, corresponding to $v \in \dot\U^\imath_{\Z\lround q\rround}$
and $w \in \dot\U^\imath_{\Z\lround q^{-1}\rround}$
under the isomorphisms $\widetilde{\BWWJ}$ and $\BWWJ$, we have that
\begin{equation}\label{homformula2}
\Hom_\H(V,W) \cong \left(\R^\imath\right)^{\oplus \langle v,w\rangle^\imath}
\end{equation}
as a graded $\R^\imath$-module.
To prove this, one can reduce to the case that $V = P(\bi)$ and $W = I(\bj)$,
and in this case the formula can be deduced from \cref{homformula1} using some duality.
\end{rem}

%% file: appendix-a.tex
\appendix
\setcounter{section}{0}

\section{Proof of \texorpdfstring{\cref{necessary}}{}}\label{Appendix-A}

For this appendix, we work in the weak 2-iquantum group
$\xUUi$ from \cref{extended},
assuming that $i \neq \tau i$ throughout.
We will need the following
for $f(x,y) \in \kk[x,y]$:
\begin{align}
\left[\frac{f(u,v)}{(u-x)(v-y)(u-z)}\right]_{\substack{u:-1\\v:-1}}=
\left[\frac{f(u,y)}{(u-x)(u-z)}\right]_{u:-1} &= \frac{f(x,y)-f(z,y)}{x-z}.\label{differentworld}
\end{align}
This follows from \cref{trick} and partial fractions.
We proceed with several lemmas, before proving \cref{necessary}.

\begin{lem}\label{coince}
For $i \neq \tau i$, the following holds in $\xUUi$:
\begin{align}\label{coincidence}
.
\end{align*}
Then we expand the final expression thus obtained using
\cref{iquadratic,ipitchfork,icurl} to see that it equals
\begin{multline*}
- %
\right]_{w:-1}\!\!\!\!.
\end{multline*}
It remains to observe that 
$$
%
\,.
$$
Indeed, if one denotes a dot on the left hand string by $x$
and dots on the bottom right or top right strings by $y$ and $z$, this can be rewritten as the elementary algebraic expression
\begin{multline*}
\frac{1}{(w+x)(v+y)(w+y)(u-z)(v+z)}
-\frac{1}{(w+x)(v+y)(u-z)(v+z)(w+z)}\\
-\frac{1}{(w+x)(u-x)(v+y)(w+y)(u-z)}
+\frac{1}{(u-x)(w+y)(u-z)(v+z)(w+z)}\\
+\frac{1}{(u-x)(v+y)(w+y)(u-z)(v+z)}
-\frac{1}{(u-x)(v+y)(u-z)(v+z)(w+z)}\\
=\frac{1}{(u-x)(w+x)(v+y)(w+y)(w+z)}
-
\frac{1}{(w+x)(w+y)(u-z)(v+z)(w+z)}.
\end{multline*}
This may be checked by putting the fractions over a common denominator.
\end{proof}

\begin{cor}
If $i \neq \tau i$ then
\begin{align}\label{stillstuck}
&
\ .
\end{align}
\end{cor}

\begin{proof}
By \cref{trick}, we have that
\begin{align*}
%
\right)\right]_{\substack{u:-1\\v:-1}}.
\end{align*}
The corollary follows from this, \cref{coincidence,differentworld}.
\end{proof}

\begin{lem}\label{covered}
Assuming that $i \neq \tau i$, the following hold in $\xUUi$:
\begin{align}
\left[\ 
%
\ 
\right]_{\!u:-1}\!\!\!\!.\label{rectified}
\end{align}
\end{lem}

\begin{proof}
We have that
\begin{align*}
%
\ .
\end{align*}
The identity
\cref{minerals} follows from this, and
then \cref{rectified} follows by rotating \cref{minerals}.
\end{proof}

\begin{lem}\label{hardplace1}
Still assuming that $i \neq \tau i$, for any $\lambda \in X^\imath$, 
we have that
\begin{multline*}
\left[\!%
\right]_{\substack{u:-1\\v:-1}}
\!\!\!\!.
\end{multline*}
\end{lem}

\begin{proof}
After some easy manipulations, the right hand side of the identity to be proved is
\begin{multline}\label{hardplace2}
\left[
\zeta_{\tau i}(\lambda)^{-1}\!
%
,
$$
where $n := \del_{\tau i} - \lambda_{\tau i}$.
We use this,
together with \cref{fancyicurl} for the first two,
to expand the four terms on the right hand side of \cref{hardplace2} (we will take the $u^{-1}$-coefficients shortly):
\begin{align*}
\zeta_{\tau i}(\lambda)^{-1}\!\!
%
\ .
\end{align*}
The $u^{-1}$-coefficient of the first term on the right hand side of the first equation here is the left hand side of \cref{hardplace1}.
It remains to show that the $u^{-1}$-coefficient of the sum of the four summations in the above  equations is 0. To check this, 
since all four summations have 
$\zeta_{\tau i}(\lambda)^{-1}$ in front, we can ignore this scale factor.
Also in the second summand the $u^{-1}$-coefficient is clearly 0 if
$t > n-s-1$, so this summation only needs to be taken over
$0 \leq t \leq n-s-1$ like the other three.
Then the result follows because for $r,s\geq 0$ and $0 \leq t \leq n-s-1$ we have that
\begin{multline*}
\left[\tfrac{u^{n-s-t-1}}{(u+x)^2}\left(u^{-r-1}
\begin{tikzpicture}[iQ,anchorbase,scale=1.2]
\draw (-1.06,0) arc(360:0:0.17);
\node at (-1.53,0) {\strandlabel{\tau i}};
\node at (-.69,0) {$(-u)$};
\region{0,0}{\lambda-\alpha_i};
\end{tikzpicture}\right)-
\tfrac{u^{n-s-t-1}}{(u+x)^2}\left[u^{-r-1}
\begin{tikzpicture}[iQ,anchorbase,scale=1.2]
\draw (-1.06,0) arc(360:0:0.17);
\node at (-1.53,0) {\strandlabel{\tau i}};
\node at (-.69,0) {$(-u)$};
\region{0,0}{\lambda-\alpha_i};
\end{tikzpicture}\right]_{u:<0}\right]_{u:-1}=\\
\left[\left(\tfrac{1}{u+x}
\left[t\tfrac{v^{n-s-t-1}}{(v+x)^2}\right]_{v:-1}
+\tfrac{1}{(u+x)^2}
\left[\tfrac{v^{n-s-t-1}}{v+x}\right]_{v:-1}\right)
\left(u^{-r-1}\begin{tikzpicture}[iQ,anchorbase,scale=1.2]
\draw (-1.06,0) arc(360:0:0.17);
\node at (-1.53,0) {\strandlabel{\tau i}};
\node at (-.69,0) {$(-u)$};
\region{0,0}{\lambda-\alpha_i};
\end{tikzpicture}\right)\right]_{u:-1}.
\end{multline*}
To prove this identity, note that the left hand side simplifies to
$$
\left[\tfrac{u^{n-s-t-1}}{(u+x)^2}\left[u^{-r-1}\!\!\!
\begin{tikzpicture}[iQ,anchorbase,scale=1.2]
\draw (-1.06,0) arc(360:0:0.17);
\node at (-1.53,0) {\strandlabel{\tau i}};
\node at (-.69,0) {$(-u)$};
\region{0,0}{\lambda-\alpha_i};
\end{tikzpicture}\!\!\right]_{u:\geq 0}\right]_{u:-1}.
$$
It is easy to see that this is equal to the right hand side when $t=n-s-1$.
The case $t=n-s-2$ is only slightly harder; it follows because
$$
\tfrac{u}{(u+x)^2}=
\tfrac{1}{u+x}
\left[\tfrac{v}{(v+x)^2}\right]_{v:-1}
+
\tfrac{1}{(u+x)^2}
\left[\tfrac{v}{v+x}\right]_{v:-1}.
$$
Finally suppose that $0 \leq t \leq n-s-3$ and set
$k := n-s-t-1$, so $k \geq 2$.
The identity to be proved is equivalent to
\begin{multline*}
\sum_{m \geq 0} 
\left[\tfrac{1}{(u+x)^2}\right]_{u:-m-k-1}\left[u^{-r-1}\!\!\!
\begin{tikzpicture}[iQ,anchorbase,scale=1.2]
\draw (-1.06,0) arc(360:0:0.17);
\node at (-1.53,0) {\strandlabel{\tau i}};
\node at (-.69,0) {$(-u)$};
\region{0,0}{\lambda-\alpha_i};
\end{tikzpicture}\!\!\right]_{u:m}
=\\
\sum_{m \geq 0}
\left[\tfrac{1}{u+x}
\tfrac{1}{(v+x)^2}+\tfrac{1}{(u+x)^2}
\tfrac{1}{v+x}\right]_{\substack{u:-m-1\\v:-k-1}}
\left[u^{-r-1}\!\!\!
\begin{tikzpicture}[iQ,anchorbase,scale=1.2]
\draw (-1.06,0) arc(360:0:0.17);
\node at (-1.53,0) {\strandlabel{\tau i}};
\node at (-.69,0) {$(-u)$};
\region{0,0}{\lambda-\alpha_i};
\end{tikzpicture}\!\!\right]_{u:m}.
\end{multline*}
This follows since by binomial expansions
the coefficient of 
$\left[u^{-r-1}\!\!\!
\begin{tikzpicture}[iQ,anchorbase,scale=1.2]
\draw (-1.06,0) arc(360:0:0.17);
\node at (-1.53,0) {\strandlabel{\tau i}};
\node at (-.69,0) {$(-u)$};
\region{0,0}{\lambda-\alpha_i};
\end{tikzpicture}\!\!\right]_{u:m}$
on the left hand side is
$(m+k)(-x)^{m+k-1}$,
and on the right hand side it is $(-x)^m\cdot k(-x)^{k-1}
+ m(-x)^{m-1}\cdot (-x)^k$.
\end{proof}

\begin{lem}
If $i \neq \tau i$ then
\begin{multline}
\left[
\begin{tikzpicture}[iQ,centerzero,scale=.98]
\draw[-] (0.4,0.6)  to [out=down,in=down,looseness=2.5] (0,.6)\toplabel{\tau i};
\draw[fakebubble,-] (.95,.1) arc(-360:0:0.18);
\bubblelabel{.77,.1}{u};
\stringlabel{.5,.1}{i};
\draw[-] (-.4,-.6) \botlabel{i}to [out=90,in=-90,looseness=1] (.2,-.1) to [out=90,in=0]
(0,0.22) to [out=180,in=up,looseness=1] (-.4,-.2) to [out=down,in=up] (0,-.6);
\draw[-] (-.4,.6)  to (-.4,.4) to [out=down,in=up] (.4,-.4) to (.4,-.6)\botlabel{i};
\circledbar{.33,.4}{u};
\circledbar{.13,.14}{u};
\end{tikzpicture}-
\begin{tikzpicture}[iQ,centerzero,scale=.98]
\draw[-] (-0.4,-0.6) \botlabel{i} to [out=up,in=up,looseness=2.5] (0,-.6);
\draw[fakebubble,-] (-0.8,-.1) arc(-180:180:0.18);
\bubblelabel{-.62,-.1}{u};
\stringlabel{-.99,-.1}{\tau i};
\draw[-] (.4,.6) to [out=-90,in=90,looseness=1] (-.2,.1) to [out=-90,in=180]
(0,-0.22) to [out=0,in=down,looseness=1] (.4,.2) to [out=up,in=down] (0,.6)\toplabel{\tau i};
\draw[-] (.4,-.6) \botlabel{i} to (.4,-.4) to [out=up,in=down] (-.4,.4) to (-.4,.6);
\circled{-.33,-.4}{u};
\circled{-.13,-.14}{u};
\end{tikzpicture}\right]_{u:-1}\!\!\!\!\!
=
\left[\begin{tikzpicture}[iQ,centerzero,scale=1.25]
\draw[-] (-0.2,-0.4) \botlabel{i} to [looseness=2.5,out=90,in=90] (0.2,-0.4);
\draw[-] (0.2,0.4) to [looseness=2.5,out=-90,in=-90] (-0.2,0.4)\toplabel{i};
\draw[-] (.84,-.4) \botlabel{i} to (.84,.4);
\circled{-.17,.2}{u};
\circled{-.17,-.2}{u};
\circled{.84,0}{u};
\draw[fakebubble] (0.33,0) arc(180:-180:0.13);
\node at (0.2,0) {\strandlabel{\tau i}};
\bubblelabel{0.46,0}{u};
\pinpin{.17,.23}{.85,.23}{1.6,.23}{Q_{i,\tau i}(y,x)};
\end{tikzpicture}-\begin{tikzpicture}[iQ,centerzero,scale=1.25]
\draw[-] (0.2,0.4)  to [looseness=2.5,out=-90,in=-90] (-0.2,0.4)\toplabel{\tau i};
\draw[-] (-0.2,-0.4) \botlabel{\tau i}to [looseness=2.5,out=90,in=90] (0.2,-0.4);
\draw[-] (-.84,.4)  to (-.84,-.4)\botlabel{i};
\circledbar{.17,-.2}{u};
\circledbar{.17,.2}{u};
\circledbar{-.84,0}{u};
\draw[fakebubble] (-0.3,0) arc(0:-360:0.13);
\node at (-0.62,0) {\strandlabel{i}};
\bubblelabel{-0.43,0}{u};
\pinpin{-.17,-.23}{-.85,-.23}{-1.6,-.23}{Q_{i,\tau i}(x,y)};
\end{tikzpicture}\right]_{u:-1}\\\label{thelaststraw}
+\begin{tikzpicture}[iQ,centerzero,scale=1.6]
\draw[-] (-0.3,0.3) -- (0.3,-0.3)\botlabel{i};
\draw[-] (0,0.3) \toplabel{\tau i}[out=-90,in=-90,looseness=2.5] to  (0.3,0.3);
\draw[-] (-0.3,-0.3) \botlabel{i}[out=90,in=90,looseness=2.5] to (0,-0.3);
\pinpinpinmid{.08,.1}{0,0}{-.08,-.1}{-1.25,0}{
\frac{Q_{\tau i, i}(x,y)-Q_{\tau i, i}(z,y)}{x-z}};
\pinpin{0.02,.2}{.28,.2}{.9,.2}{Q_{i,\tau i}(y,x)};
\end{tikzpicture}+
\begin{tikzpicture}[iQ,centerzero,scale=1.2]
\draw[-] (0,0.4) \toplabel{\tau i} to[out=-90,in=-90,looseness=2.5] (0.4,0.4);
\draw[-] (-.4,-0.4) \botlabel{i} to[out=90,in=90,looseness=2.5] (0,-0.4);
\draw[-] (0.4,-0.4) \botlabel{i} -- (-0.4,0.4);
\pinpin{-0.01,-.25}{.25,-.25}{1.1,-.25}{Q_{i,\tau i}(y,x)};
\pinpinpinmid{.1,.14}{0,0}{-.1,-.14}{-1.8,0}{
\frac{Q_{\tau i,i}(x,-x)-Q_{\tau i,i}(z,-z)}{x-z}};
\end{tikzpicture}\ .
\end{multline}
\end{lem}

\begin{proof}
Take any $\lambda \in X^\imath$.
We have that
\begin{align}
&\begin{tikzpicture}[iQ,centerzero,scale=1.6]
\draw[-] (-0.3,0.3) -- (0.3,-0.3)\botlabel{i};
\draw[-] (0,0.3) \toplabel{\tau i}[out=-90,in=-90,looseness=2.5] to  (0.3,0.3);
\draw[-] (-0.3,-0.3) \botlabel{i}[out=90,in=90,looseness=2.5] to (0,-0.3);
\pinpinpinmid{.08,.1}{0,0}{-.08,-.1}{-1.35,0}{
\frac{Q_{\tau i, i}(x,y)-Q_{\tau i, i}(z,y)}{x-z}};
\pinpin{0.02,.2}{.28,.2}{.9,.2}{Q_{i,\tau i}(y,x)};
\region{.8,-.1}{\lambda};
\end{tikzpicture}+
\begin{tikzpicture}[iQ,centerzero,scale=1.2]
\draw[-] (0,0.4) \toplabel{\tau i} to[out=-90,in=-90,looseness=2.5] (0.4,0.4);
\draw[-] (-.4,-0.4) \botlabel{i} to[out=90,in=90,looseness=2.5] (0,-0.4);
\draw[-] (0.4,-0.4) \botlabel{i} -- (-0.4,0.4);
\pinpinpinmid{.1,.14}{0,0}{-.1,-.14}{-1.9,0}{
\frac{Q_{\tau i,i}(x,-x)-Q_{\tau i,i}(z,-z)}{x-z}};
\pinpin{-0.01,-.25}{.25,-.25}{1.15,-.25}{Q_{i,\tau i}(y,x)};
\region{1,.2}{\lambda};
\end{tikzpicture}\notag\\
&\qquad=\begin{tikzpicture}[iQ,centerzero,scale=1.6]
\draw[-] (-0.3,0.3) -- (0.3,-0.3)\botlabel{i};
\draw[-] (0,0.3) \toplabel{\tau i}[out=-90,in=-90,looseness=2.5] to  (0.3,0.3);
\draw[-] (-0.3,-0.3) \botlabel{i}[out=90,in=90,looseness=2.5] to (0,-0.3);
\pinpinpinmid{.08,.1}{0,0}{-.08,-.1}{-3.2,0}{
\frac{(Q_{\tau i, i}(x,y)-Q_{\tau i, i}(z,y))Q_{\tau i, i}(z,-z)}{x-z}+\frac{Q_{\tau i, i}(x,y)(Q_{\tau i, i}(x,-x)-Q_{\tau i, i}(z,-z))}{x-z}};
\region{.35,0}{\lambda};
\end{tikzpicture}\notag\\
&\qquad=\begin{tikzpicture}[iQ,centerzero,scale=1.6]
\draw[-] (-0.3,0.3) -- (0.3,-0.3)\botlabel{i};
\draw[-] (0,0.3) \toplabel{\tau i}[out=-90,in=-90,looseness=2.5] to  (0.3,0.3);
\draw[-] (-0.3,-0.3) \botlabel{i}[out=90,in=90,looseness=2.5] to (0,-0.3);
\pinpinpinmid{.08,.1}{0,0}{-.08,-.1}{-2.2,0}{
\frac{
Q_{\tau i, i}(x,y)Q_{\tau i, i}(x,-x)-Q_{\tau i, i}(z,y)Q_{\tau i, i}(z,-z)}{x-z}};
\region{.35,0}{\lambda};
\end{tikzpicture}\stackrel{\cref{differentworld}}{=}
\left[
\begin{tikzpicture}[iQ,centerzero,scale=1.2]
\draw[-] (0,0.4) \toplabel{\tau i} to[out=-90,in=-90,looseness=3] (0.4,0.4);
\draw[-] (-.4,-0.4) \botlabel{i} to[out=90,in=90,looseness=3] (0,-0.4);
\draw[-] (0.4,-0.4) \botlabel{i} -- (-0.4,0.4);
\circled{-.07,-.2}{u};
\circled{.07,.2}{u};
\pin{0,0}{-1.8,.3}{Q_{\tau i, i}(u,x) Q_{\tau i,i}(u,-u)};
\region{.5,0}{\lambda};
\end{tikzpicture}
\right]_{\!u:-1}\!\!\!\!.\label{atthepool}
\end{align}
Now we observe that
\begin{align*}
Q_{\tau i, i}(u,-u)&\ \stackrel{\cref{being}}{=}\ (-1)^{\del_i+1}\zeta_i \zeta_{\tau i}u^{-a_{i,\tau i}} \ \substack{\cref{howthingswork}\\{\textstyle =}\\\cref{childhood}}\  -
\zeta_i(\lambda)\zeta_{\tau i}(\lambda)(-1)^{\del_i-\lambda_i}u^{-a_{i,\tau i}},\\
Q_{\tau i, i}(u,x) &\stackrel{\cref{boshgosh}}{=} 
t_{\tau i, i} R_{\tau i, i}(u,x)
\stackrel{\cref{luckya}}{=} c_{\tau i}(\alpha_i) R_{\tau i, i}(u,x)
\substack{\cref{purer}\\{\textstyle=}\\\cref{childhood}} 
(-1)^{a_{i,\tau i}}\frac{\zeta_i(\lambda-\alpha_i)}{\zeta_i(\lambda)} 
R_{\tau i, i}(u,x).
\end{align*}
So we can rewrite \cref{atthepool} to obtain
\begin{multline*}
- \zeta_i(\lambda-\alpha_i)
\zeta_{\tau i}(\lambda)
\left[ (-1)^{a_{i,\tau i}} \begin{tikzpicture}[iQ,centerzero,scale=1.2]
\draw[-] (-0.4,0.4)  -- (0.4,-0.4)\botlabel{i};
\draw[-] (0,-0.4) [out=90,in=90,looseness=2.5] to  (-0.4,-0.4)\botlabel{i};
\draw[-] (0.4,0.4) [out=-90,in=-90,looseness=2.5] to 
(0,.4)\toplabel{\tau i};
\circled{.05,.23}{u};
\circled{-.05,-.23}{u};
\draw (.6,.21) arc(180:-180:0.17);
\node at (.5,.21) {\strandlabel{i}};
\draw (.6,-.2) arc(180:-180:0.17);
\node at (.44,-.2) {\strandlabel{\tau i}};
\node at (1.18,.21) {$(u)$};
\node at (1.29,-.2) {$(-u)$};
\pin{-.2,.2}{-1.1,.2}{R_{\tau i, i}(u,x)};
\end{tikzpicture}\right]_{u:-1}\\
\stackrel{\cref{iinfgrass}}{=}
-\zeta_i(\lambda-\alpha_i)\zeta_{\tau i}(\lambda)
\left[
\begin{tikzpicture}[iQ,centerzero,scale=1.2]
\draw[-] (-0.4,0.4)  -- (0.4,-0.4)\botlabel{i};
\draw[-] (0,-0.4) [out=90,in=90,looseness=2.5] to  (-0.4,-0.4)\botlabel{i};
\draw[-] (0.4,0.4) [out=-90,in=-90,looseness=2.5] to 
(0,.4)\toplabel{\tau i};
\circled{.05,.23}{u};
\circled{-.05,-.23}{u};
\draw (.76,0) arc(360:0:0.17);
\node at (.35,0) {\strandlabel{i}};
\node at (1.01,0) {$(u)$};
\draw (-1.05,0) arc(360:0:0.17);
\node at (-1.53,0) {\strandlabel{\tau i}};
\node at (-.69,0) {$(-u)$};
\circledbar{.25,-.25}{u};
\circledbar{-.25,.25}{u};
\region{1.35,0}{\lambda};
\end{tikzpicture}\right]_{u:-1}\!\!\!\!.
\end{multline*}
Then we expand this
as a sum of four terms using \cref{hardplace1}:
it equals
\begin{align*}
&\left[
\begin{tikzpicture}[iQ,centerzero,scale=1.25]
\draw[-] (-0.4,0.4)  -- (0.4,-0.4)\botlabel{i};
\draw[-] (0,-0.4) [out=90,in=90,looseness=2.5] to  (-0.4,-0.4)\botlabel{i};
\draw[-] (0.4,0.4) [out=-150,in=180,looseness=1.4] to 
(.23,-0.05) to [out=0,in=-30,looseness=1.4] (0.06,.4)\toplabel{\tau i};
\circledbar{.32,.1}{u};
\circled{-.38,-.25}{u};
\draw[fakebubble] (-.53,0) arc(360:0:0.17);
\bubblelabel{-.7,0}{u};
\node at (-1.01,0) {\strandlabel{\tau i}};
\circled{-.2,.2}{u};
\region{.6,0}{\lambda};
\end{tikzpicture}-
\begin{tikzpicture}[iQ,centerzero,scale=1.25]
\draw[-] (0.4,-0.4) \botlabel{i} -- (-0.4,0.4);
\draw[-] (0,0.4) \toplabel{\tau i} [out=-90,in=-90,looseness=2.5] to  (0.4,0.4);
\draw[-] (-0.4,-0.4) \botlabel{i}[out=30,in=0,looseness=1.4] to 
(-.23,0.05) to [out=180,in=150,looseness=1.4] (-0.06,-.4);
\circled{-.32,-.1}{u};
\circledbar{.38,.25}{u};
\draw[fakebubble] (.78,0) arc(360:0:0.17);
\node at (.38,0) {\strandlabel{i}};
\bubblelabel{.61,0}{u};
\circledbar{.2,-.2}{u};
\circledbar{-.2,.2}{u};
\region{1,0}{\lambda};
\end{tikzpicture}\right]_{\!u:-1}+
\left[\begin{tikzpicture}[iQ,centerzero,scale=1.25]
\draw[-] (0.4,-0.4) \botlabel{i} -- (-0.4,0.4);
\draw[-] (0,0.4) \toplabel{\tau i} [out=-90,in=-90,looseness=2.5] to  (0.4,0.4);
\draw[-] (-0.4,-0.4) \botlabel{i}[out=90,in=90,looseness=2.5] to (0,-.4);
\circled{-.38,-.25}{v};
\circledbar{.38,.25}{u};
\draw[fakebubble] (.78,0) arc(360:0:0.17);
\node at (.38,0) {\strandlabel{i}};
\bubblelabel{.61,0}{u};
\draw[fakebubble] (-.53,0) arc(360:0:0.17);
\node at (-1.01,0) {\strandlabel{\tau i}};
\bubblelabel{-.69,0}{v};
\circledbar{-.2,.2}{u};
\circled{0,0}{v};
\circledbar{.2,-.2}{u};
\region{1,0}{\lambda};
\end{tikzpicture}-\begin{tikzpicture}[iQ,centerzero,scale=1.25]
\draw[-] (0.4,-0.4) \botlabel{i} -- (-0.4,0.4);
\draw[-] (0,0.4) \toplabel{\tau i} [out=-90,in=-90,looseness=2.5] to  (0.4,0.4);
\draw[-] (-0.4,-0.4) \botlabel{i}[out=90,in=90,looseness=2.5] to (0,-.4);
\circled{-.38,-.25}{u};
\circledbar{.38,.25}{v};
\draw[fakebubble] (.78,0) arc(360:0:0.17);
\bubblelabel{.61,0}{v};
\node at (.38,0) {\strandlabel{i}};
\draw[fakebubble] (-.53,0) arc(360:0:0.17);
\node at (-1.01,0) {\strandlabel{\tau i}};
\bubblelabel{-.69,0}{u};
\circledbar{0,0}{v};
\circled{.2,-.2}{u};
\circled{-.2,.2}{u};
\region{1,0}{\lambda};
\end{tikzpicture}\right]_{\substack{u:-1\\v:-1}}
\end{align*}
Now the lemma follows by taking 
the difference of the equations
\cref{minerals,rectified} and comparing with the identity just proved.
\end{proof}

\begin{proof}[Proof of \cref{necessary}]
We start from the ordinary braid relation \cref{ibraid} with strings labelled $(\tau i,i,i)$ (its right hand side is zero). Adding crossings to the bottom left pair of strings and to the top right pair of strings produces
the relation
\begin{align*}
\begin{tikzpicture}[iQ,centerzero,scale=1.1]
\draw[-] (-0.4,-0.8) \botlabel{i} to [out=up,in=down] (0,-.4)  to [out=up,in=down] (-.4,0) to [out=up,in=down] (.4,.8);
\draw[-] (0,-.8) \botlabel{\tau i} to [out=up,in=down] (-.4,-.4) to [out=up,in=down] (.4,.2) to [out=up,in=down] (0,.8);
\draw[-] (.4,-.8) \botlabel{i} to (.4,-.4) to [out=up,in=down] (-.4,.4) to (-.4,.8);
\end{tikzpicture}-
\begin{tikzpicture}[iQ,centerzero,scale=1.1]
\draw[-] (0.4,0.8) to [out=down,in=up] (0,.4)  to [out=down,in=up] (.4,0) to [out=down,in=up] (-.4,-.8) \botlabel{i};
\draw[-] (0,.8) to [out=down,in=up] (.4,.4) to [out=down,in=up] (-.4,-.2) to [out=down,in=up] (0,-.8) \botlabel{\tau i};
\draw[-] (-.4,.8) to (-.4,.4) to [out=down,in=up] (.4,-.4) to (.4,-.8)  \botlabel{i};
\end{tikzpicture}&=0.
\end{align*}
Then we use \cref{iquadratic} to expand the double crossings in this equation to deduce that
$$
\begin{tikzpicture}[iQ,centerzero,scale=1.5]
\draw[-] (-0.4,-0.4) \botlabel{i} -- (0.4,0.4);
\draw[-] (0,-0.4) \botlabel{\tau i} to[out=135,in=down] (-0.32,0) to[out=up,in=225] (0,0.4);
\draw[-] (0.4,-0.4) \botlabel{i} -- (-0.4,0.4);
\pinpin{-0.12,.3}{.3,.3}{.95,.3}{Q_{i,\tau i}(y,x)};
\end{tikzpicture}
-\begin{tikzpicture}[iQ,centerzero,scale=1.5]
\draw[-] (-0.4,-0.4) \botlabel{i} -- (0.4,0.4);
\draw[-] (0,-0.4) \botlabel{\tau i} to[out=45,in=down] (0.32,0) to[out=up,in=-45] (0,0.4);
\draw[-] (0.4,-0.4) \botlabel{i} -- (-0.4,0.4);
\pinpin{0.12,-.3}{-.3,-.3}{-1,-.3}{Q_{i,\tau i}(x,y)};
\end{tikzpicture}
=
\left[
\begin{tikzpicture}[iQ,centerzero,scale=1]
\draw[-] (0.4,0.6)  to [out=down,in=down,looseness=2.5] (0,.6)\toplabel{\tau i};
\draw[fakebubble,-] (.95,.1) arc(-360:0:0.18);
\bubblelabel{.77,.1}{u};
\stringlabel{.5,.1}{i};
\draw[-] (-.4,-.6) \botlabel{i}to [out=90,in=-90,looseness=1] (.2,-.1) to [out=90,in=0]
(0,0.22) to [out=180,in=up,looseness=1] (-.4,-.2) to [out=down,in=up] (0,-.6);
\draw[-] (-.4,.6)  to (-.4,.4) to [out=down,in=up] (.4,-.4) to (.4,-.6)\botlabel{i};
\circledbar{.33,.4}{u};
\circledbar{.13,.13}{u};
\end{tikzpicture}-\begin{tikzpicture}[iQ,centerzero,scale=1]
\draw[-] (-0.4,-0.6) \botlabel{i} to [out=up,in=up,looseness=2.5] (0,-.6);
\draw[fakebubble,-] (-0.8,-.1) arc(-180:180:0.18);
\bubblelabel{-.62,-.1}{u};
\stringlabel{-.99,-.1}{\tau i};
\draw[-] (.4,.6) to [out=-90,in=90,looseness=1] (-.2,.1) to [out=-90,in=180]
(0,-0.22) to [out=0,in=down,looseness=1] (.4,.2) to [out=up,in=down] (0,.6)\toplabel{\tau i};
\draw[-] (.4,-.6) \botlabel{i} to (.4,-.4) to [out=up,in=down] (-.4,.4) to (-.4,.6);
\circled{-.33,-.4}{u};
\circled{-.13,-.13}{u};
\end{tikzpicture}\right]_{u:-1}\!\!\!.
$$
The identity \cref{dandy} follows by adding \cref{stillstuck} to this equation
then expanding the right hand side of the result using \cref{thelaststraw}.
\end{proof}

%% file: appendix-b.tex
\appendix
\setcounter{section}{1}

\section{Properties of internal bubbles}\label{Appendix-B}

Notation in this appendix is as in \cref{newgrading}. 

\begin{lem}\label{wombles}
For $i \in I$ with $i \neq \tau i$ (equivalently, $\del_i \geq 0$), the following hold in $\UU(\del,\zeta)$:
\begin{align}
\begin{tikzpicture}[Q,centerzero,scale=1]
\draw[to-] (0.68,0) arc(0:360:0.4);
\node at (.28,-.5) {\strandlabel{i}};
\anticlockwiseinternalbubble{-.02,-.22};
\Circled{-.02,.25}{u};
\end{tikzpicture} 
&=
\zeta_i \left[
u^{\del_i}\ \begin{tikzpicture}[Q,centerzero,scale=1.4]
\draw[to-] (0.68,0) arc(0:360:0.25);
\stringlabel{.43,-.35}{i};
\stringlabel{-.3,-.3}{\tau i};
\Circled{.18,0}{u};
\filledanticlockwisebubble{-.3,0};	
\bubblelabel{-.3,0}{-u};
\end{tikzpicture} 
-u^{\del_i}\ \begin{tikzpicture}[Q,centerzero,scale=1.4]
\draw[to-] (0.68,0) arc(0:360:0.25);
\stringlabel{.43,-.35}{i};
\stringlabel{-.4,-.3}{\tau i};
\Circled{.23,.15}{u};
\anticlockwisebubble{-.4,0};	
\teleporter{.2,-.11}{-.24,-.11};
\end{tikzpicture} \ 
\right]_{u:<0},\label{puzzle1}\\\label{puzzle2}
\begin{tikzpicture}[Q,centerzero,scale=1]
\draw[to-] (-0.68,0) arc(180:-180:0.4);
\node at (-.28,-.5) {\strandlabel{\tau i}};
\clockwiseinternalbubble{0.02,.22};
\Circledbar{0.02,-.25}{u};
\end{tikzpicture} 
&=
-\zeta_i \left[
u^{\del_i}\ \begin{tikzpicture}[Q,centerzero,scale=1.4]
\draw[to-] (-0.68,0) arc(180:-180:0.25);
\stringlabel{-.43,-.35}{\tau i};
\stringlabel{.3,-.3}{i};
\Circledbar{-.18,0}{u};
\filledclockwisebubble{.3,0};	
\bubblelabel{.3,0}{u};
\end{tikzpicture} 
-u^{\del_i}\ \begin{tikzpicture}[Q,centerzero,scale=1.4]
\draw[to-] (-0.68,0) arc(180:-180:0.25);
\stringlabel{-.43,-.35}{\tau i};
\stringlabel{.4,-.3}{i};
\Circledbar{-.23,-.15}{u};
\clockwisebubble{.4,0};	
\teleporter{-.2,.11}{.24,.11};
\end{tikzpicture} \ 
\right]_{u:<0}.
\end{align}
\end{lem}

\begin{proof}
For the first identity, we expand the definition \cref{internal2} of the internal bubble to see that
\begin{equation*}
\begin{tikzpicture}[Q,centerzero,scale=1]
\draw[to-] (0.68,0) arc(0:360:0.4);
\node at (.28,-.5) {\strandlabel{i}};
\anticlockwiseinternalbubble{-0.02,-.22};
\Circled{-0.02,.25}{u};
\end{tikzpicture} 
=\zeta_i \left[
\begin{tikzpicture}[Q,centerzero,scale=1.4]
\draw[to-] (0.68,0) arc(0:360:0.25);
\stringlabel{.43,-.35}{i};
\stringlabel{-.4,-.3}{\tau i};
\Circled{.23,.18}{u};
\Circledbar{.23,-.18}{v};
\filledanticlockwisebubble{-.4,0};	
\bubblelabel{-.4,0}{v};
\multopendot{.18,0}{east}{\:\del_i\!};
\end{tikzpicture}\ \right]_{v:-1}
-\zeta_i\ \begin{tikzpicture}[Q,centerzero,scale=1.4]
\draw[to-] (0.68,0) arc(0:360:0.25);
\stringlabel{.43,-.35}{i};
\stringlabel{-.4,-.3}{\tau i};
\Circled{.23,.15}{u};
\anticlockwisebubble{-.4,0};
\multopendot{.23,-.16}{east}{\:\del_i\!};
\teleporter{.18,-.03}{-.19,-.03};
\end{tikzpicture}\ .
\end{equation*}
The second term is equal to the second term in \cref{puzzle1} thanks to \cref{trick}, and the first terms are equal by a similar trick.
Similarly, for the second identity, the internal bubble expands to give
$$
\begin{tikzpicture}[Q,centerzero,scale=1]
\draw[to-] (-0.68,0) arc(180:-180:0.4);
\node at (-.28,-.5) {\strandlabel{\tau i}};
\clockwiseinternalbubble{0.02,.22};
\Circledbar{0.02,-.25}{u};
\end{tikzpicture} 
=
(-1)^{\del_i+1}\zeta_i \left[
\ \begin{tikzpicture}[Q,centerzero,scale=1.4]
\draw[to-] (-0.68,0) arc(180:-180:0.25);
\stringlabel{-.43,-.35}{\tau i};
\stringlabel{.4,-.3}{i};
\multopendot{-.17,0}{west}{\!\del_i};
\Circledbar{-.26,.19}{v};
\Circledbar{-.26,-.19}{u};
\filledclockwisebubble{.4,0};	
\bubblelabel{.4,0}{v};
\end{tikzpicture}\  
\right]_{v:-1}
-(-1)^{\del_i+1}\zeta_i\ \begin{tikzpicture}[Q,centerzero,scale=1.4]
\draw[to-] (-0.68,0) arc(180:-180:0.25);
\stringlabel{-.43,-.35}{\tau i};
\multopendot{-.24,.17}{west}{\!\del_i};
\stringlabel{.4,-.3}{i};
\Circledbar{-.23,-.15}{u};
\clockwisebubble{.4,0};	
\teleporter{-.18,0.04}{.2,0.02};
\end{tikzpicture}\ ,
$$
which is equal to the right hand side of \cref{puzzle2}.
\end{proof}

\begin{rem}\label{awkwardness}
It is useful to be able to exploit the symmetries $\bar\Omega$ and 
$\Sigma$ from \cref{Omegainv,Sigmainv} extended to the localization to prove relations involving internal bubbles. We forget the grading for this since $\bar\Omega$ and $\Sigma$ do not respect the adapted grading.
More problematic is that $\bar\Omega$ and $\Sigma$ do not 
interchange clockwise and counterclockwise internal bubbles in the obvious way. However, the internal bubbles depend
implicitly on $\del$ and $\zeta$, and we can view
$\bar\Omega$ and $\Sigma$ instead as $\kk$-linear 2-functors
$\bar\Omega:\UU(\del,\zeta)^\op\rightarrow \UU({'}\del, {'}\zeta)$
and $\Sigma:\UU(\del,\zeta)^\rev\rightarrow \UU({'}\del, {'}\zeta)$ where ${'}\del, {'}\zeta$ are the primed parameters from \cref{sinead}.
Then both $\bar\Omega$ and $\Sigma$ take the clockwise and counterclockwise internal bubbles in $\UU(\del,\zeta)$ to the counterclockwise and clockwise ones 
in $\UU({'}\del, {'}\zeta)$ defined using the primed parameters.
To illustrate, we can take \cref{puzzle1}
with $i$ replaced by $\tau i$:
$$
 \ 
\right]_{\!u:<0}\!\!\!.
$$
Applying $\Sigma$ gives the identity
$$
%
 \ 
\right]_{u:<0},
$$
in $\UU({'}\del, {'}\zeta)$.
Replacing $u$ by $-u$, ${'}\del_i$ by $\del_i$ and ${'}\zeta_i$ by $\zeta_i$ in this identity recovers \cref{puzzle2}.
\end{rem}

\begin{lem}\label{slushy}
If $i \neq \tau i$ then
the following hold in $\UU(\del,\zeta)$:
\begin{align}\label{awkwardness1}
(-1)^{\del_{\tau i}+1} \zeta_{\tau i}\ 
%
.
\end{align}
\end{lem}

\begin{proof}
Fix a weight $\lambda \in X$.
We start from the right hand side and calculate:
\begin{align*}
\left[%
\ .
\end{align*}
\end{proof}

\begin{lem}\label{bubblyinverses}
Clockwise and counterclockwise internal bubbles are inverses of each other:
\begin{align}\label{bubinv}
%
\ .
\end{align}
\end{lem}

\begin{proof}
This is proved in the case that $i=\tau i$ in \cite[Lem.~3.5]{BWWbasis}. Now suppose that $i \neq \tau i$. We just prove the first equality. The others then follow since internal bubbles commute and we can rotate through 180$^\circ$ (or apply $\bar\Omega$).
Using the definition \cref{internal1} and \cref{slushy}, we have that
$$
%
.
$$
Now we use
\begin{align*}
%
\right]_{u:-1}\ .
\end{align*}
Putting these two things together, we deduce that
\begin{align*}
%
\ .
\end{align*}
\end{proof}

\begin{lem}\label{drainsfree}
Assuming $\tau i \neq j$,
the following hold in $\UU(\del,\zeta)$:
\begin{align}
%
.\label{dylan}
\end{align}
\end{lem}

\begin{proof}
Pick $\lambda \in X$.
The following proves \cref{marley}:
\begin{align*}
%
.
\end{align*}
The proof of \cref{dylan} is similar.
\end{proof}

\begin{lem}\label{internalbubbleslide}
Assuming $i \neq j$ and $i \neq \tau j$,
the following hold in $\UU(\del,\zeta)$:
\begin{align}\label{sauerkraut}
%
\right]_{u:-1}
\end{align}
for any $\lambda \in X$.
\end{lem}

\begin{proof}
Suppose first that $i \neq j$ and $i \neq \tau j$. 
We use \cref{dotslide,teleslide} to slide the dots and teleporter downward past the crossing in \cref{marley} to obtain
$$
%
.
$$
This proves the first equality from \cref{sauerkraut}. The proof of the second one is similar, or it can be deduced from the first equality by applying $\Sigma$ using the technique of \cref{awkwardness} plus \cref{bubblyinverses}.

To prove \cref{bob}, we start as in the first line of the proof of \cref{drainsfree}:
\begin{align*}
%
\ .
\end{align*}
The second term is now more complicated to simplify since there is an extra term:
\begin{align*}
\zeta_i
%
\!.
\end{align*}
Substituting this into our previous equation yields the first equality in \cref{bob}. The second follows on applying $\Sigma$ and using \cref{bubblyinverses}.

The following calculation proves \cref{blocked}:
\begin{align*}
&%
\right]_{u:-1}\!\!\!.
\end{align*}
The proof of \cref{drains} is similar, starting from \cref{marley} instead of \cref{dylan}.
\end{proof}

\begin{cor}\label{monday}
If $i \neq j$ and $i \neq \tau j$ then
\begin{equation}\label{littlek}
%
\right]_{u:-1}\!\!\!.
\end{align}
\end{cor}

\begin{lem}\label{belltime}
For all $i,j \in I$, we have that
\begin{align}\label{toomuchbarking}
%
\ .
\end{align}
\end{lem}

\begin{proof}
If $i = j = \tau i$ then this is proved in \cite[Lem.~3.11]{BWWbasis}.
For the other cases, suppose to start with just that $i \neq \tau j$.
Then
\begin{align*}
r_{\tau j, i}Q_{\tau j, i}(-x,y)
&= r_{\tau j, i} Q_{i,\tau j}(y,-x)
\stackrel{\cref{tauantisymmetric}}{=} r_{i,\tau j}^{-1} r_{i,j}^{-1} r_{\tau j,\tau i}^{-1}Q_{\tau i, j}(-y,x).
\end{align*}
When $i = j \neq \tau i$, this shows that
$Q_{\tau i, i}(-x,y) = Q_{\tau i, i}(-y,x)$, so that it is symmetric in $x$ and $y$.

Now suppose that $i \neq j$ and $i \neq \tau j$. Then, using
the identity proved in the previous paragraph for the unlabelled equality, we have that
\begin{align*}
.
\end{align*}
This proves \cref{toomuchbarking} when $i \neq j$ and $i \neq \tau j$.

If $i = \tau j \neq \tau i$, we instead use the two identities in \cref{bob}:
\begin{align*}
%
.
\end{align*}

Finally suppose that $i = j \neq \tau j$, and take any $\lambda \in X$.
By \cref{drains,blocked}, we have that
\begin{align*}
%
\right]_{u:-1}\!\!\!\!.
\end{align*}
The right hand sides of these equations are equal because $Q_{\tau i,i}(-x,y) = Q_{\tau i, i}(-y,x)$.
\end{proof}

\begin{cor}\label{stateschool}
We have that
\begin{equation}\label{quickie}
%
\end{equation}
for all $i,j \in I$.
\end{cor}

\begin{lem}\label{bedtime}
For $i \neq \tau i$ and $\lambda \in X$, we have that
\begin{align}\label{sleep}
%
\right]_{u:-1}\!\!.
\end{align}
\end{lem}

\begin{proof}
We have that
\begin{align*}
&%
\right]_{u:-1}\!\!.
\end{align*}
The identity \cref{sleep} follows from this and \cref{internal1}.
\end{proof}

\begin{lem}\label{newlem}
For $i \neq \tau i$ and $\lambda \in X$, we have that
\begin{align}\label{comfyhotel}
%
\right]_{u:-1}.
\end{align}
\end{lem}

\begin{proof}
We begin with
\begin{align*}
%
\right]_{u:-1}\!\!\!\!.
\end{align*}
Applying $\bar\Omega$ to this identity with $i$ replaced by $\tau i$
and $\lambda$ replaced with $-\lambda$
gives also that
\begin{align*}
-%
\right]_{u:-1}\!\!\!\!.
\end{align*}
Now we replace $u$ by $-u$ in the final term of this last expression then add it to the previous equation to obtain
\begin{align*}
%
\right]_{u:-1}\!\!\!.
\end{align*}
\end{proof}

\begin{lem}\label{goodnight}
For $i \in I$ with $i \neq \tau i$ and $\lambda \in X$, we have that
\begin{align}\label{prettyhard}
%
\ .
\end{align}
\end{lem}

\begin{proof}
We start by calculating
\begin{align*}
&%
\!
\right]_{\!u:-1}\!\!\!\!.
\end{align*}
Let us call the five terms in this final expression $A, B, C, D$ and $E$ in order from left to right, so it is $A+B+C+D+E$.
We have that
\begin{align*}
B &\ \substack{\cref{bubblegeneratingfunctions}\\{\textstyle=}}\;-\zeta_ic_{\tau i}(\lambda)^{-1}\!\left[ %
\right]_{\!u:-1}\!\!\!.
\end{align*}
    The first term in the expansion of $B$ cancels with the first term in the expansion of $C$,
    $D$ cancels with the third term in the expansion of $C$, the final term in the expansion of $C$ can be combined with $A$, and by \cref{bubblegeneratingfunctions} the first term of $E$ cancels with some of the second term of $C$ leaving just some fake bubble polynomial. We deduce that
    $A+B+C+D+E$ equals
   $$
-\frac{\zeta_ic_i(\lambda)}{c_{\tau i}(\lambda)}
\left[(-u)^{\del_i} %
\,\right]_{\!u:-1}\!\!\!\!\!\!.
$$
Now we replace $u$ by $-u$ (remembering the overall minus sign as we are taking the $u^{-1}$-coefficient)  to deduce that our expression equals
$$
\frac{\zeta_i c_i(\lambda)}{c_{\tau i}(\lambda)}
\left[u^{\del_i} %
\,\right]_{\!u:-1}\!\!\!\!\!\!.
$$
It just remains to simplify the last term:
$$
\zeta_i\!
\left[u^{\del_i} %
.
$$
\end{proof}

%% file: appendix-c.tex
\appendix
\setcounter{section}{2}

\section{Proof of
\texorpdfstring{\cref{boxingday}}{}}\label{Appendix-C}

There are ten relations to be checked to complete the proof.

\vspace{2mm}
\noindent
{\underline{Relation \cref{ibubblerel}}}:
In view of \cref{whybeta}, we need to check that
$$
\left[(-1)^{h_{\tau i}(\hat\lambda)} u^{\del_i}\ 
%
\right]_{u:\geq \del_i - \lambda_i}
= u^{\del_i - \lambda_i}
\id_{\one_{\hat\lambda}}.
$$
As $\lambda_i = h_i(\hat\lambda) -h_{\tau i}(\hat\lambda)$,
this follows 
by the following identities, which are consequences of \cref{infgrass0,infgrass1}:
\begin{align*}
\left[\ 
%
\right]_{u:\geq -h_{i}(\hat\lambda)}
\!\!\!\!\!\!\!\!\!\!\!\!&=
 u^{-h_i(\hat\lambda)}
\id_{\one_{\hat\lambda}}.
\end{align*}

\vspace{2mm}
\noindent
{\underline{Relation \cref{iinfgrass}}}:
By \cref{whybeta}, we need to show that
$$
(-1)^{h_i(\hat\lambda)+h_{\tau i}(\hat\lambda)} u^{\del_{\tau i}}(-u)^{\del_i}\ 
%
= 
(-1)^{\del_i-\lambda_i} u^{-a_{i,\tau i}}
\id_{\one_{\hat\lambda}}.
$$
This follows from \cref{infgrass,recap}.

\vspace{2mm}
\noindent
{\underline{Relation \cref{ibubslide}}}:
Using \cref{whybeta}, this reduces to checking the following two identities:
\begin{align*}
%
\ .
\end{align*}
The first of these follows from the image of
\cref{bubslide} under $\bar\Omega$.
For the second one, \cref{bubslide} gives instead that
$$
%
\ .
$$
Now we use the definitions \cref{boshgosh,tauantisymmetric}, which show that 
$R_{\tau i, \tau j}(-u,x) = (-1)^{a_{\tau i,\tau j}}
R_{i,j}(u,-x)$
and $R_{\tau i, j}(-u,-x) = (-1)^{a_{\tau i, j}}R_{i,\tau j}(u,x)$.

\vspace{2mm}
\noindent
{\underline{Relation \cref{izigzag}}}:
This follows from \cref{rightadj,leftpivots} using also \cref{bubblyinverses}.

\vspace{2mm}
\noindent
{\underline{Relations \cref{ipivotal}}}:
These follow from \cref{lotsmore}---this is the reason for the sign in \cref{psi1}.

\vspace{2mm}
\noindent
{\underline{Relations \cref{ipitchfork}}}:
It suffices to prove the first relation involving the cup. The other one then follows by rotation, i.e. we attach caps to the top left and top right strings and use the zig-zag identity. 

If $i = j = \tau i$ the relation is checked in the proof of \cite[Th.~4.2]{BWWbasis} (it is the relation referred to there as ``the second relation from (2.4)'').

Suppose next that $i \neq j$ and $i \neq \tau j$.
We have that
\begin{align*}
\Xi^\imath\left(\ \ 
\begin{tikzpicture}[iQ,baseline=-5pt,scale=1.2]
\draw[-] (-0.25,0.15) \toplabel{i} to[out=-90,in=-90,looseness=3] (0.25,0.15);
\draw[-] (-0.3,-0.4) to[out=up,in=down] (0,0.15)\toplabel{j};
\region{.4,0}{\hat\lambda};
\end{tikzpicture}\right)_{\hat\lambda}
&=
-r_{j,i}^{-1}
\begin{tikzpicture}[Q,baseline=-5pt,scale=1.7]
\draw[-to] (-0.25,0.15) \toplabel{i} to[out=-90,in=-90,looseness=3] (0.25,0.15);
\draw[to-] (-0.3,-0.4) to[out=up,in=down] (0,0.15)\toplabel{j};
\region{.4,0}{\hat\lambda};
\end{tikzpicture}
-r_{\tau i, \tau j}r_{\tau j,i}
\begin{tikzpicture}[Q,baseline=-5pt,scale=1.7]
\draw[to-] (-0.25,0.15) \toplabel{\tau i} to[out=-90,in=-90,looseness=3] (0.25,0.15);
\draw[-to] (-0.3,-0.4) to[out=up,in=down] (0,0.15)\toplabel{\tau j};
\anticlockwiseinternalbubble{.2,-.1};
\region{.4,0}{\hat\lambda};
\end{tikzpicture}
-r_{\tau j,i}
\begin{tikzpicture}[Q,baseline=-5pt,scale=1.7]
\draw[-to] (-0.25,0.15) \toplabel{i} to[out=-90,in=-90,looseness=3] (0.25,0.15);
\draw[-to] (-0.3,-0.4) to[out=up,in=down] (0,0.15)\toplabel{\tau j};
\region{.4,0}{\hat\lambda};
\end{tikzpicture}
-
r_{i,j}^{-1}r_{\tau i, j}^{-1}\begin{tikzpicture}[Q,baseline=-5pt,scale=1.7]
\draw[to-] (-0.25,0.15) \toplabel{\tau i} to[out=-90,in=-90,looseness=3] (0.25,0.15);
\draw[to-] (-0.3,-0.4) to[out=up,in=down] (0,0.15)\toplabel{j};
\anticlockwiseinternalbubble{-.24,-0.03};
\clockwiseinternalbubble{-0.1,-.26};
\anticlockwiseinternalbubble{.2,-.1};
\region{.45,0}{\hat\lambda};
\end{tikzpicture}\ .
\end{align*}
We need to show that this equals
\begin{align*}
\Xi^\imath\left(\ r_{i,j}^{-1}\ 
\begin{tikzpicture}[iQ,baseline=-5pt,scale=1.2]
\draw[-] (-0.25,0.15)\toplabel{i}  to[out=-90,in=-90,looseness=3] (0.25,0.15);
\draw[-] (0.3,-0.4)to[out=up,in=down] (0,0.15)\toplabel{j};
\region{.4,0}{\hat\lambda};
\end{tikzpicture}\right)_{\hat\lambda}
&=
-\ 
\begin{tikzpicture}[Q,baseline=-5pt,scale=1.7]
\draw[-to] (-0.25,0.15)\toplabel{i}  to[out=-90,in=-90,looseness=3] (0.25,0.15);
\draw[to-] (0.3,-0.4)to[out=up,in=down] (0,0.15)\toplabel{j};
\region{.4,0}{\hat\lambda};
\end{tikzpicture}
-\ 
r_{i,j}^{-1} r_{j,\tau i}^{-1} r_{\tau j, \tau i}^{-1}
\begin{tikzpicture}[Q,baseline=-5pt,scale=1.7]
\draw[to-] (-0.25,0.15)\toplabel{\tau i}  to[out=-90,in=-90,looseness=3] (0.25,0.15);
\draw[-to] (0.3,-0.4)to[out=up,in=down] (0,0.15)\toplabel{\tau j};
\anticlockwiseinternalbubble{-.2,-.1};
\anticlockwiseinternalbubble{.09,-.07};
\clockwiseinternalbubble{.26,-.27};
\region{.4,0}{\hat\lambda};
\end{tikzpicture}
-\ 
r_{\tau j, i}\begin{tikzpicture}[Q,baseline=-5pt,scale=1.7]
\draw[-to] (-0.25,0.15)\toplabel{i}  to[out=-90,in=-90,looseness=3] (0.25,0.15);
\draw[-to] (0.3,-0.4)to[out=up,in=down] (0,0.15)\toplabel{\tau j};
\region{.4,0}{\hat\lambda};
\end{tikzpicture}
-\ 
r_{i,j}^{-1} r_{\tau i, j}^{-1}\begin{tikzpicture}[Q,baseline=-5pt,scale=1.7]
\draw[to-] (-0.25,0.15)\toplabel{\tau i}  to[out=-90,in=-90,looseness=3] (0.25,0.15);
\draw[to-] (0.3,-0.4)to[out=up,in=down] (0,0.15)\toplabel{j};
\region{.4,0}{\hat\lambda};
\anticlockwiseinternalbubble{-.2,-.1};
\end{tikzpicture}\ .
\end{align*}
That the first and third terms are equal follows directly from \cref{ruby}, and the fourth terms are equal by \cref{ruby} plus an application of \cref{bubblyinverses}.
To see that the second terms are equal, we use \cref{belltime}:
$$
r_{i,j}^{-1} r_{j,\tau i}^{-1} r_{\tau j, \tau i}^{-1}
\begin{tikzpicture}[Q,baseline=-5pt,scale=2]
\draw[to-] (-0.25,0.15)\toplabel{\tau i}  to[out=-90,in=-90,looseness=3] (0.25,0.15);
\draw[-to] (0.3,-0.4)to[out=up,in=down] (0,0.15)\toplabel{\tau j};
\anticlockwiseinternalbubble{-.2,-.1};
\anticlockwiseinternalbubble{.09,-.07};
\clockwiseinternalbubble{.26,-.27};
\region{.4,0}{\hat\lambda};
\end{tikzpicture}
\stackrel{\cref{ruby}}{=}
r_{i,j}^{-1} r_{j,\tau i}^{-1} r_{\tau j, \tau i}^{-1}
r_{\tau i, \tau j}\begin{tikzpicture}[Q,baseline=-5pt,scale=2]
\draw[to-] (-0.25,0.15)\toplabel{\tau i}  to[out=-90,in=-90,looseness=3] (0.25,0.15);
\draw[-to] (-0.3,-0.4)to[out=up,in=down] (0,0.15)\toplabel{\tau j};
\anticlockwiseinternalbubble{-.23,0};
\anticlockwiseinternalbubble{-0.02,0};
\clockwiseinternalbubble{-.27,-.27};
\region{.4,0}{\hat\lambda};
\end{tikzpicture}
\ \substack{\cref{toomuchbarking}\\{\textstyle=}\\\cref{bubinv}}\ 
r_{\tau i, \tau j}
r_{\tau j,i}
\begin{tikzpicture}[Q,baseline=-5pt,scale=2]
\draw[to-] (-0.25,0.15)\toplabel{\tau i}  to[out=-90,in=-90,looseness=3] (0.25,0.15);
\draw[-to] (-0.3,-0.4)to[out=up,in=down] (0,0.15)\toplabel{\tau j};
\anticlockwiseinternalbubble{-.1,-.25};
\region{.4,0}{\hat\lambda};
\end{tikzpicture}\ .
$$

Next suppose that $i = j \neq \tau i$.
We must show that
\begin{align*}
\Xi^\imath\left(\ \ 
\begin{tikzpicture}[iQ,baseline=-5pt,scale=1.15]
\draw[-] (-0.25,0.15) \toplabel{i} to[out=-90,in=-90,looseness=3] (0.25,0.15);
\draw[-] (-0.3,-0.4) to[out=up,in=down] (0,0.15)\toplabel{i};
\region{.4,0}{\hat\lambda};
\end{tikzpicture}\right)_{\hat\lambda}\!
&=
-
\begin{tikzpicture}[Q,baseline=-5pt,scale=1.6]
\draw[-to] (-0.25,0.15) \toplabel{i} to[out=-90,in=-90,looseness=3] (0.25,0.15);
\draw[to-] (-0.3,-0.4) to[out=up,in=down] (0,0.15)\toplabel{i};
\region{.4,0}{\hat\lambda};
\end{tikzpicture}
-
\begin{tikzpicture}[Q,baseline=-5pt,scale=1.6]
\draw[to-] (-0.25,0.15) \toplabel{\tau i} to[out=-90,in=-90,looseness=3] (0.25,0.15);
\draw[-to] (-0.3,-0.4) to[out=up,in=down] (0,0.15)\toplabel{\tau i};
\anticlockwiseinternalbubble{.2,-.1};
\region{.4,0}{\hat\lambda};
\end{tikzpicture}
-\sgn(i)
\begin{tikzpicture}[Q,baseline=-5pt,scale=1.6]
\draw[-to] (-0.25,0.15) \toplabel{i} to[out=-90,in=-90,looseness=3] (0.25,0.15);
\draw[-to] (-0.3,-0.4) to[out=up,in=down] (0,0.15)\toplabel{\tau i};
\region{.4,0}{\hat\lambda};
\end{tikzpicture}
\!+\sgn(i)
\begin{tikzpicture}[Q,baseline=-5pt,scale=1.6]
\draw[to-] (-0.25,0.15) \toplabel{\tau i} to[out=-90,in=-90,looseness=3] (0.25,0.15);
\draw[to-] (-0.3,-0.4) to[out=up,in=down] (0,0.15)\toplabel{i};
\anticlockwiseinternalbubble{-.24,-0.03};
\clockwiseinternalbubble{-0.1,-.26};
\anticlockwiseinternalbubble{.2,-.1};
\region{.45,0}{\hat\lambda};
\end{tikzpicture}
\!-
\begin{tikzpicture}[Q,baseline=-5pt,scale=1.6]
\draw[to-] (-0.25,0.15) \toplabel{\tau i} to (-.25,-.4);
\draw[-to] (0,.15) \toplabel{i}to (0,-.1)  to [out=-90,in=-90,looseness=2] (.25,-.1) to (.25,.15);
\teleporter{-.25,0}{0,0};
\region{.4,0}{\hat\lambda};
\end{tikzpicture}+
\begin{tikzpicture}[Q,baseline=-5pt,scale=1.6]
\draw[-to] (-0.25,0.15) \toplabel{i} to (-.25,-.4);
\draw[to-] (0,.15) \toplabel{\tau i}to (0,-.1)  to [out=-90,in=-90,looseness=2] (.25,-.1) to (.25,.15);
\anticlockwiseinternalbubble{0.02,-.15};
\teleporter{-.25,0}{0,0};
\region{.4,0}{\hat\lambda};
\end{tikzpicture}
\end{align*}
equals
\begin{align*}
\Xi^\imath\left(\ 
\begin{tikzpicture}[iQ,baseline=-5pt,scale=1.15]
\draw[-] (-0.25,0.15)\toplabel{i}  to[out=-90,in=-90,looseness=3] (0.25,0.15);
\draw[-] (0.3,-0.4)to[out=up,in=down] (0,0.15)\toplabel{i};
\region{.4,0}{\hat\lambda};
\end{tikzpicture}\right)_{\hat\lambda}\!
&=
-
\begin{tikzpicture}[Q,baseline=-5pt,scale=1.6]
\draw[-to] (-0.25,0.15)\toplabel{i}  to[out=-90,in=-90,looseness=3] (0.25,0.15);
\draw[to-] (0.3,-0.4)to[out=up,in=down] (0,0.15)\toplabel{i};
\region{.4,0}{\hat\lambda};
\end{tikzpicture}\!
-
\begin{tikzpicture}[Q,baseline=-5pt,scale=1.6]
\draw[to-] (-0.25,0.15)\toplabel{\tau i}  to[out=-90,in=-90,looseness=3] (0.25,0.15);
\draw[-to] (0.3,-0.4)to[out=up,in=down] (0,0.15)\toplabel{\tau i};
\anticlockwiseinternalbubble{-.2,-.1};
\anticlockwiseinternalbubble{.09,-.07};
\clockwiseinternalbubble{.26,-.27};
\region{.4,0}{\hat\lambda};
\end{tikzpicture}\!
+\sgn(\tau i)
\begin{tikzpicture}[Q,baseline=-5pt,scale=1.6]
\draw[-to] (-0.25,0.15)\toplabel{i}  to[out=-90,in=-90,looseness=3] (0.25,0.15);
\draw[-to] (0.3,-0.4)to[out=up,in=down] (0,0.15)\toplabel{\tau i};
\region{.4,0}{\hat\lambda};
\end{tikzpicture}\!
-\sgn(\tau i)
\begin{tikzpicture}[Q,baseline=-5pt,scale=1.6]
\draw[to-] (-0.25,0.15)\toplabel{\tau i}  to[out=-90,in=-90,looseness=3] (0.25,0.15);
\draw[to-] (0.3,-0.4)to[out=up,in=down] (0,0.15)\toplabel{i};
\region{.4,0}{\hat\lambda};
\anticlockwiseinternalbubble{-.2,-.1};
\end{tikzpicture}
\!-\begin{tikzpicture}[Q,anchorbase,scale=.75]
 	\draw[to-,thin] (-1,.6) \toplabel{\tau i}to (-1,-.25) to [out=-90,in=-90,looseness=1.5] (-0.4,-.25) to[out=90,in=180,looseness=1] (-0.1,0) to[out=0,in=90,looseness=1] (0.25,-.3) to (0.25,-.6);
 	\draw[-to,thin] (-0.45,.6) \toplabel{i}to[out=-90,in=180] (-.1,.1) to[out=0,in=-90] (0.25,.6);
\teleporter{-.4,.3}{-.4,-.2};
     \clockwiseinternalbubble{.15,-.22};
          \anticlockwiseinternalbubble{-.92,-.35};
  \region{.5,0.1}{\hat\lambda};
\end{tikzpicture}
\!+\begin{tikzpicture}[Q,anchorbase,scale=.75]
\draw[-to,thin] (-1,.6) \toplabel{i}to (-1,-.25) to [out=-90,in=-90,looseness=1.5] (-0.4,-.25) to[out=90,in=180,looseness=1] (-0.1,0) to[out=0,in=90,looseness=1] (0.25,-.3) to (0.25,-.6);
 	\draw[to-,thin] (-0.45,.6)\toplabel{\tau i}to[out=-90,in=180] (-.1,.1) to[out=0,in=-90] (0.25,.6);
\teleporter{.18,.3}{.23,-.2};
\anticlockwiseinternalbubble{-.35,.25};
  \region{.6,0}{\hat\lambda};
\end{tikzpicture}.
\end{align*}
The last two terms are easily seen to be equal using the zig-zag identities and \cref{bubblyinverses}. The equality of the first four terms
follows like in the previous paragraph, using also that
$\sgn(\tau i) = -\sgn(i)$.

The final case $i = \tau j \neq \tau i$ follows in a similar way to the case $i=j\neq \tau i$ just treated.

\vspace{2mm}
\noindent
{\underline{Relation \cref{idotslide}}}:
In the case $i = j=\tau i$, this is checked in the proof of \cite[Th.~4.2]{BWWbasis} (it is the relation referred to there as ``the first relation from (2.5)''). 

Suppose that $i = \tau j \neq \tau i$. 
We prove the first equality in \cref{idotslide}. 
We have that
\begin{align*}
\Xi^\imath\left(\ \begin{tikzpicture}[iQ,centerzero]
\closeddot{-.15,-.15};
\draw[-] (-0.3,-0.3) \botlabel{i} -- (0.3,0.3);
\draw[-] (0.3,-0.3) \botlabel{\tau i} -- (-0.3,0.3);
\region{0.35,0}{\hat\lambda};
\end{tikzpicture}-\begin{tikzpicture}[iQ,centerzero]
\closeddot{.15,.15};
\draw[-] (-0.3,-0.3) \botlabel{i} -- (0.3,0.3);
\draw[-] (0.3,-0.3) \botlabel{\tau i} -- (-0.3,0.3);
\region{0.35,0}{\hat\lambda};
\end{tikzpicture}\right)_{\hat\lambda}\ &=
-\sgn(i)\begin{tikzpicture}[Q,anchorbase,scale=.7]
\draw[to-] (0.6,-.6)\botlabel{\tau i} to (-0.6,.6);
	\draw[to-,thin] (-0.6,-.6)\botlabel{i} to (0.6,.6);
  \region{.5,0}{\hat\lambda};
\opendot{-.3,-.3};
\end{tikzpicture}
-\sgn(i)\begin{tikzpicture}[Q,anchorbase,scale=.7]
	\draw[-to,thin] (0.6,-.6)\botlabel{ i} to (-0.6,.6);
	\draw[-to,thin] (-0.6,-.6)\botlabel{\tau i} to (0.6,.6);
  \region{.5,0}{\hat\lambda};
\opendot{-.3,-.3};
\end{tikzpicture}
+\begin{tikzpicture}[Q,anchorbase,scale=.7]
	\draw[to-,thin] (0.6,-.6)\botlabel{\tau i} to (-0.6,.6);
	\draw[-to,thin] (-0.6,-.6)\botlabel{\tau i} to (0.6,.6);
  \region{.5,0}{\hat\lambda};
\opendot{-.3,-.3};
\end{tikzpicture}
-\begin{tikzpicture}[Q,anchorbase,scale=.7]
	\draw[-to,thin] (0.6,-.6)\botlabel{i} to (-0.6,.6);
	\draw[to-,thin] (-0.6,-.6)\botlabel{i} to (0.6,.6);
\anticlockwiseinternalbubble{-.3,.3};
\clockwiseinternalbubble{.3,-.3};
\opendot{-.3,-.3};
\region{.7,0.1}{\hat\lambda};
\end{tikzpicture}
-\begin{tikzpicture}[Q,anchorbase,scale=.8]
 	\draw[to-,thin] (-0.45,-.6)\botlabel{i} to[out=90,in=180] (-.1,-.1) to[out=0,in=90] (0.25,-.6);
 	\draw[-to,thin] (-0.45,.6) \toplabel{\tau i}to[out=-90,in=180] (-.1,.1) to[out=0,in=-90] (0.25,.6);
\teleporter{-.4,.3}{-.4,-.3};
     \clockwiseinternalbubble{.15,-.25};
  \region{.5,0.1}{\hat\lambda};
\opendot{-.25,-.15};
\end{tikzpicture}-\begin{tikzpicture}[Q,anchorbase,scale=.8]
 	\draw[-to,thin] (-0.45,-.6)\botlabel{\tau i} to[out=90,in=180] (-.1,-.1) to[out=0,in=90] (0.25,-.6);
 	\draw[to-,thin] (-0.45,.6)\toplabel{i}to[out=-90,in=180] (-.1,.1) to[out=0,in=-90] (0.25,.6);
\teleporter{.19,.3}{.19,-.3};
\anticlockwiseinternalbubble{-.35,.25};
  \region{.55,0}{\hat\lambda};
\opendot{-.25,-.15};
\end{tikzpicture}\\
&\:\:\:\:\,+\sgn(i)\begin{tikzpicture}[Q,anchorbase,scale=.7]
\draw[to-] (0.6,-.6)\botlabel{\tau i} to (-0.6,.6);
	\draw[to-] (-0.6,-.6)\botlabel{i} to (0.6,.6);
  \region{.5,0}{\hat\lambda};
\opendot{.3,.3};
\end{tikzpicture}
+\sgn(i)\begin{tikzpicture}[Q,anchorbase,scale=.7]
	\draw[-to] (0.6,-.6)\botlabel{i} to (-0.6,.6);
	\draw[-to] (-0.6,-.6)\botlabel{\tau i} to (0.6,.6);
  \region{.5,0}{\hat\lambda};
\opendot{.3,.3};
\end{tikzpicture}
-\begin{tikzpicture}[Q,anchorbase,scale=.7]
	\draw[to-] (0.6,-.6)\botlabel{\tau i} to (-0.6,.6);
	\draw[-to] (-0.6,-.6)\botlabel{\tau i} to (0.6,.6);
  \region{.5,0}{\hat\lambda};
\opendot{.3,.3};
\end{tikzpicture}
+\begin{tikzpicture}[Q,anchorbase,scale=.7]
	\draw[-to] (0.6,-.6)\botlabel{i} to (-0.6,.6);
	\draw[to-] (-0.6,-.6)\botlabel{i} to (0.6,.6);
\anticlockwiseinternalbubble{-.3,.3};
\clockwiseinternalbubble{.3,-.3};
\opendot{.3,.3};
\region{.7,0.1}{\hat\lambda};
\end{tikzpicture}
-\begin{tikzpicture}[Q,anchorbase,scale=.8]
 	\draw[to-] (-0.45,-.6)\botlabel{i} to[out=90,in=180] (-.1,-.1) to[out=0,in=90] (0.25,-.6);
 	\draw[-to] (-0.45,.6) \toplabel{\tau i}to[out=-90,in=180] (-.1,.1) to[out=0,in=-90] (0.25,.6);
\teleporter{-.4,.3}{-.4,-.3};
     \clockwiseinternalbubble{.15,-.25};
  \region{.5,0.1}{\hat\lambda};
\opendot{.05,.15};
\end{tikzpicture}-\begin{tikzpicture}[Q,anchorbase,scale=.8]
 	\draw[-to] (-0.45,-.6)\botlabel{\tau i} to[out=90,in=180] (-.1,-.1) to[out=0,in=90] (0.25,-.6);
 	\draw[to-] (-0.45,.6)\toplabel{i}to[out=-90,in=180] (-.1,.1) to[out=0,in=-90] (0.25,.6);
\teleporter{.19,.3}{.19,-.3};
\anticlockwiseinternalbubble{-.35,.25};
  \region{.55,0}{\hat\lambda};
\opendot{.05,.15};
\end{tikzpicture}\\
&\ \substack{\cref{dotslide}\\{\textstyle=}\\\cref{partialfractions}}\ 
-\begin{tikzpicture}[Q,anchorbase,scale=.7]
	\draw[-to] (-0.4,-.6)\botlabel{\tau i} to[looseness=2.3,out=90,in=90] (0.4,-.6);
	\draw[-to] (-0.4,.6)\toplabel{\tau i} to[out=-90,in=-90,looseness=2.3] (0.4,.6);
  \region{.7,0}{\hat\lambda};
\end{tikzpicture}-
\begin{tikzpicture}[Q,anchorbase,scale=.7]
	\draw[to-] (-0.4,-.6)\botlabel{i} to[looseness=2.3,out=90,in=90] (0.4,-.6);
	\draw[to-] (-0.4,.6)\toplabel{i} to[out=-90,in=-90,looseness=2.3] (0.4,.6);
\anticlockwiseinternalbubble{-.3,.25};
\clockwiseinternalbubble{.3,-.25};
\region{.8,0}{\hat\lambda};
\end{tikzpicture}
-\begin{tikzpicture}[Q,anchorbase,scale=.7]
	\draw[to-] (-0.4,-.6)\botlabel{i} to[looseness=2.3,out=90,in=90] (0.4,-.6);
	\draw[-to] (-0.4,.6)\toplabel{\tau i} to[out=-90,in=-90,looseness=2.3] (0.4,.6);
\clockwiseinternalbubble{.3,-.25};
\region{.8,0}{\hat\lambda};
\end{tikzpicture}-
\begin{tikzpicture}[Q,anchorbase,scale=.7]
	\draw[-to] (-0.4,-.6)\botlabel{\tau i} to[looseness=2.3,out=90,in=90] (0.4,-.6);
	\draw[to-] (-0.4,.6)\toplabel{i} to[out=-90,in=-90,looseness=2.3] (0.4,.6);
\anticlockwiseinternalbubble{-.3,.25};
\region{.7,0}{\hat\lambda};
\end{tikzpicture}
=\Xi^\imath\left(-
\begin{tikzpicture}[iQ,centerzero,scale=1.2]
\draw[-] (-0.2,-0.3) \botlabel{i} to [looseness=2,out=90,in=90] (0.2,-0.3);
\draw[-] (0.2,0.3)  to [looseness=2,out=-90,in=-90] (-0.2,0.3)\toplabel{\tau i};
\region{.35,0}{\hat\lambda};
\end{tikzpicture}\right)_{\hat\lambda}.
\end{align*}
The second equality in \cref{idotslide} 
follows by similar considerations.

The case $i = j \neq \tau i$ follows by similar considerations.

The case that $i \neq j$ and $i \neq \tau j$ is even a bit easier.

\vspace{2mm}
\noindent
{\underline{Relation \cref{icurl}}}:
If $i = \tau i$, this is dealt with in the proof of \cite[Th.~4.2]{BWWbasis} (where it is ``the first relation from (2.4)'').
Now assume that $i \neq \tau i$. 
The fake bubble polynomial on the right hand side of \cref{icurl}
can be replaced by the bubble generating function multiplied by
$\zeta_{\tau i}(\lambda)$ according to \cref{ibubblegeneratingfunction},
and this scalar equals $(-1)^{h_i(\hat\lambda)}\frac{\zeta_{\tau i} c_{\tau i}(\hat\lambda)}{c_i(\hat\lambda)}$
by the definitions \cref{purer,childhood}.
Applying $\Xi^\imath$ to \cref{icurl} using \cref{whybeta}
reduces the proof 
to checking the following two identities:
\begin{align*}
-\begin{tikzpicture}[anchorbase,scale=1.2,Q]
\draw[-to] (0,-0.5)\botlabel{\tau i} to (0,0.5);
\stringlabel{.5,-.3}{i};
\anticlockwisebubble{.5,0};
\teleporter{0,0}{.3,0};
\clockwiseinternalbubble{.67,0};
\region{1,0}{\hat\lambda};
\end{tikzpicture}
-
\begin{tikzpicture}[anchorbase,scale=1.2,Q]
\draw[-to] (0,-0.5)\botlabel{\tau i} to[out=up,in=180] (0.3,0.2) to[out=0,in=up] (0.45,0) to[out=down,in=0] (0.3,-0.2) to[out=180,in=down] (0,0.5);
\anticlockwiseinternalbubble{.44,0};
\region{.7,0}{\hat\lambda};
\end{tikzpicture} &=
\frac{\zeta_{\tau i} c_{\tau i}(\hat\lambda)}{c_{i}(\hat\lambda)}
\left[\ 
u^{\del_{\tau i}}\ 
\begin{tikzpicture}[anchorbase,scale=1.2,Q]
\draw[-to] (-.1,-0.5)\botlabel{\tau i} to (-.1,0.5);
\Circled{-.1,0}{u};
\stringlabel{.37,-.3}{i};
\anticlockwisebubble{.37,0};
\node at (.92,0) {$(-u)$};
\clockwisebubble{1.5,0};
\stringlabel{1.5,-.3}{\tau i};
\node at (1.93,0) {$(u)$};
\region{2.3,0}{\hat\lambda};
\end{tikzpicture}
\right]_{u:-1},\\
\begin{tikzpicture}[anchorbase,scale=1.2,Q]
\draw[to-] (0,-0.5)\botlabel{i} to (0,0.5);
\stringlabel{.5,-.3}{\tau i};
\clockwisebubble{.5,0};
\teleporter{0,0}{.3,0};
\anticlockwiseinternalbubble{.7,0};
\region{1,0}{\hat\lambda};
\end{tikzpicture}
-
\begin{tikzpicture}[anchorbase,scale=1.2,Q]
\draw[to-] (0,-0.5)\botlabel{i} to[out=up,in=180] (0.3,0.2) to[out=0,in=up] (0.45,0) to[out=down,in=0] (0.3,-0.2) to[out=180,in=down] (0,0.5);
\clockwiseinternalbubble{.44,0};
\region{.7,0}{\hat\lambda};
\end{tikzpicture} &=\ 
\frac{\zeta_{\tau i} c_{\tau i}(\hat\lambda)}{c_i(\hat\lambda)}
\left[\ 
u^{\del_{\tau i}}\ 
\begin{tikzpicture}[anchorbase,scale=1.2,Q]
\draw[to-] (-.1,-0.5)\botlabel{i} to (-0.1,0.5);
\Circledbar{-0.1,0}{u};
\stringlabel{.37,-.3}{i};
\anticlockwisebubble{.37,0};
\node at (.92,0) {$(-u)$};
\clockwisebubble{1.5,0};
\stringlabel{1.5,-.3}{\tau i};
\node at (1.93,0) {$(u)$};
\region{2.3,0}{\hat\lambda};
\end{tikzpicture}
\right]_{u:-1}.
\end{align*}
The first identity is \cref{sleep} with $i$ replaced by $\tau i$.
The second identity follows from the first one on applying $\bar\Omega$ as explained in \cref{awkwardness}.

\vspace{2mm}
\noindent
{\underline{Relation \cref{iquadratic}}}:
If $i = j = \tau i$, this is checked in the proof of \cite[
Th.~4.2]{BWWbasis} (it is ``the first relation from (2.2)'').

If $i \neq j$ and $i \neq \tau j$,
the relation to be checked expands to the following four relations (we have also used \cref{bubblyinverses} to redistribute some internal bubbles):
\begin{align*}r_{i,j}^{-1}r_{j,i}^{-1}\begin{tikzpicture}[centerzero,scale=.94,Q]
\draw[-,thin] (0.28,0) to[out=90,in=-90] (-0.28,.75);
\draw[-,thin] (-0.28,0) to[out=90,in=-90] (0.28,.75);
\draw[to-,thin] (0.28,-.75)\botlabel{j} to[out=90,in=-90] (-0.28,0);
\draw[to-,thin] (-0.28,-.75)\botlabel{i} to[out=90,in=-90] (0.28,0);
\region{0.5,0}{\hat\lambda};
\end{tikzpicture}
&=
\begin{tikzpicture}[centerzero,scale=.94,Q]
\draw[to-,thin] (0.2,-.75)\botlabel{j} to (.2,.75);
\draw[to-,thin] (-0.3,-.75)\botlabel{i} to (-.3,.75); 
\Pinpin{-.3,0}{.2,0}{1.3,0}{Q_{i,j}(x,y)};
\region{2.3,0}{\hat\lambda};
\end{tikzpicture}\ ,&
r_{\tau i, \tau j} r_{\tau j, \tau i} r_{\tau i, j} r_{\tau j, i}\begin{tikzpicture}[centerzero,scale=.94,Q]
\draw[-to,thin] (0.28,0) to[out=90,in=-90] (-0.28,.75);
\draw[-to,thin] (-0.28,0) to[out=90,in=-90] (0.28,.75);
\draw[-,thin] (0.28,-.75)\botlabel{\tau j} to[out=90,in=-90] (-0.28,0);
\draw[-,thin] (-0.28,-.75)\botlabel{\tau i} to[out=90,in=-90] (0.28,0);
\region{0.5,0}{\hat\lambda};
\end{tikzpicture}
&=
\begin{tikzpicture}[centerzero,scale=.94,Q]
\draw[-to,thin] (0.2,-.75)\botlabel{\tau j} to (.2,.75);
\draw[-to,thin] (-0.3,-.75)\botlabel{\tau i} to (-.3,.75); 
\Pinpin{-.3,0}{.2,0}{1.5,0}{Q_{i,j}(-x,-y)};
\region{2.7,0}{\hat\lambda};
\end{tikzpicture}\ ,\\
r_{j,i}^{-1}\ 
\begin{tikzpicture}[centerzero,scale=.94,Q]
\draw[-,thin] (0.28,0) to[out=90,in=-90] (-0.28,.75);
\draw[-to,thin] (-0.28,0) to[out=90,in=-90] (0.28,.75);
\draw[-,thin] (0.28,-.75)\botlabel{\tau j} to[out=90,in=-90] (-0.28,0);
\draw[to-,thin] (-0.28,-.75)\botlabel{i} to[out=90,in=-90] (0.28,0);
\anticlockwiseinternalbubble{-.24,0};
\region{0.5,0}{\hat\lambda};
\end{tikzpicture}
&=
\begin{tikzpicture}[centerzero,scale=.94,Q]
	\draw[-to,thin] (0.2,-.75)\botlabel{\tau j} to (.2,.75);
	\draw[to-,thin] (-0.3,-.75)\botlabel{i} to (-.3,.75); 
    \Pinpin{-.3,0}{.2,0}{1.4,0}{Q_{i,j}(x,-y)};
\anticlockwiseinternalbubble{.2,-.4};
\region{2.5,0}{\hat\lambda};
\end{tikzpicture}\ ,&
r_{i,j}^{-1}\ 
\begin{tikzpicture}[centerzero,scale=.94,Q]
	\draw[-to,thin] (0.28,0) to[out=90,in=-90] (-0.28,.6);
	\draw[-,thin] (-0.28,0) to[out=90,in=-90] (0.28,.6);
	\draw[to-,thin] (0.28,-.75)\botlabel{j} to[out=90,in=-90] (-0.28,0);
	\draw[-,thin] (-0.28,-.75)\botlabel{\tau i} to[out=90,in=-90] (0.28,0);
    \clockwiseinternalbubble{.24,0};
    \region{0.55,0}{\hat\lambda};
\end{tikzpicture}
&=
\begin{tikzpicture}[centerzero,scale=.94,Q]
	\draw[to-,thin] (0.2,-.75)\botlabel{j} to (.2,.75);
	\draw[-to,thin] (-0.3,-.75)\botlabel{\tau i} to (-.3,.75); 
    \Pinpin{-.3,0}{.2,0}{1.4,0}{Q_{i,j}(-x,y)};
\clockwiseinternalbubble{-.3,.4};
\region{2.5,0}{\hat\lambda};
\end{tikzpicture}\ .
\end{align*}
The first of these follows by \cref{downwardscrossing}, and
the second by \cref{quadratic} and \cref{tauantisymmetric}.
The third and the fourth relations follow using
\cref{sauerkraut} and \cref{Altquadratic}.

Next suppose that $i = j \neq \tau i$. Using that $Q_{i,j}(x,y) = 0$, the relation to be checked expands to six relations:
\begin{align*}
\begin{tikzpicture}[centerzero,scale=.94,Q]
	\draw[-,thin] (0.28,0) to[out=90,in=-90] (-0.28,.75);
	\draw[-,thin] (-0.28,0) to[out=90,in=-90] (0.28,.75);
	\draw[to-,thin] (0.28,-.75)\botlabel{i} to[out=90,in=-90] (-0.28,0);
	\draw[to-,thin] (-0.28,-.75)\botlabel{i} to[out=90,in=-90] (0.28,0);
\region{0.5,0}{\hat\lambda};
\end{tikzpicture}
&=0,&
\begin{tikzpicture}[centerzero,scale=.94,Q]
	\draw[-to,thin] (0.28,0) to[out=90,in=-90] (-0.28,.75);
	\draw[-to,thin] (-0.28,0) to[out=90,in=-90] (0.28,.75);
	\draw[-,thin] (0.28,-.75)\botlabel{\tau i} to[out=90,in=-90] (-0.28,0);
	\draw[-,thin] (-0.28,-.75)\botlabel{\tau i} to[out=90,in=-90] (0.28,0);
\region{0.5,0}{\hat\lambda};
\end{tikzpicture}
&=
0,\\
-\sgn(i)\begin{tikzpicture}[centerzero,scale=.94,Q]
	\draw[to-,thin] (0.4,-.75)\botlabel{i} to (-0.4,.75);
	\draw[-to,thin] (-0.4,-.75)\botlabel{\tau i} to (0.4,.75);
	\teleporter{-.25,.48}{.25,.48};
    \region{0.5,0}{\hat\lambda};
\end{tikzpicture}
+\sgn(i)
\begin{tikzpicture}[centerzero,scale=.94,Q]
	\draw[to-,thin] (0.4,-.75)\botlabel{i} to (-0.4,.75);
	\draw[-to,thin] (-0.4,-.75)\botlabel{\tau i} to (0.4,.75);
	\teleporter{-.25,-.48}{.25,-.48};
        \region{0.5,0}{\hat\lambda};
\end{tikzpicture}&=0,&
-\sgn(i)\begin{tikzpicture}[centerzero,scale=.94,Q]
	\draw[-to,thin] (0.4,-.75)\botlabel{\tau i} to (-0.4,.75);
	\draw[to-,thin] (-0.4,-.75)\botlabel{i} to (0.4,.75);
	\teleporter{-.25,.48}{.25,.48};
        \region{0.5,0}{\hat\lambda};
\end{tikzpicture}
+\sgn(i)
\begin{tikzpicture}[centerzero,scale=.94,Q]
	\draw[-to,thin] (0.4,-.75)\botlabel{\tau i} to (-0.4,.75);
	\draw[to-,thin] (-0.4,-.75)\botlabel{i} to (0.4,.75);
	\teleporter{-.25,-.48}{.25,-.48};
        \region{0.5,0}{\hat\lambda};
\end{tikzpicture}
&=0,
\\
-\begin{tikzpicture}[centerzero,scale=.94,Q]
	\draw[-,thin] (0.28,0) to[out=90,in=-90] (-0.28,.75);
	\draw[-to,thin] (-0.28,0) to[out=90,in=-90] (0.28,.75);
	\draw[-,thin] (0.28,-.75)\botlabel{\tau i} to[out=90,in=-90] (-0.28,0);
	\draw[to-,thin] (-0.28,-.75)\botlabel{i} to[out=90,in=-90] (0.28,0);
    \anticlockwiseinternalbubble{-.24,0};
        \region{0.5,0}{\hat\lambda};
\end{tikzpicture}
+\begin{tikzpicture}[centerzero,scale=.94,Q]
\draw[to-,thin] (-.2,-.75)\botlabel{i} to (-.2,.75);
\draw[-to,thin] (.3,-.75)\botlabel{\tau i} to (.3,.75);
\anticlockwiseinternalbubble{.28,0};
\teleporter{-.2,.35}{.3,.35};
\teleporter{-.2,-.35}{.3,-.35};
    \region{0.6,0}{\hat\lambda};
\end{tikzpicture}
&=0,&
-\begin{tikzpicture}[centerzero,scale=.94,Q]
	\draw[-to,thin] (0.28,0) to[out=90,in=-90] (-0.28,.75);
	\draw[-,thin] (-0.28,0) to[out=90,in=-90] (0.28,.75);
	\draw[to-,thin] (0.28,-.75)\botlabel{i} to[out=90,in=-90] (-0.28,0);
	\draw[-,thin] (-0.28,-.75)\botlabel{\tau i} to[out=90,in=-90] (0.28,0);
    \clockwiseinternalbubble{.24,0};
        \region{0.6,0}{\hat\lambda};
\end{tikzpicture}
+\begin{tikzpicture}[centerzero,scale=.94,Q]
\draw[-to,thin] (-.2,-.75)\botlabel{\tau i} to (-.2,.75);
\draw[to-,thin] (.3,-.75)\botlabel{i} to (.3,.75);
\clockwiseinternalbubble{-.2,0};
\teleporter{-.2,.35}{.3,.35};
\teleporter{-.2,-.35}{.3,-.35};
    \region{0.55,0}{\hat\lambda};
\end{tikzpicture}
&=0.
\end{align*}
The first four of these follow easily from \cref{downwardscrossing,quadratic,dotslide}.
The fifth one follows because
\begin{align}\label{gordonbrown}
\begin{tikzpicture}[centerzero,scale=1.3,Q]
	\draw[-,thin] (0.28,0) to[out=90,in=-90] (-0.28,.6);
	\draw[-to,thin] (-0.28,0) to[out=90,in=-90] (0.28,.6);
	\draw[-,thin] (0.28,-.6)\botlabel{i} to[out=90,in=-90] (-0.28,0);
	\draw[to-,thin] (-0.28,-.6)\botlabel{\tau i} to[out=90,in=-90] (0.28,0);
    \anticlockwiseinternalbubble{-.24,0};
        \region{0.5,0}{\hat\lambda};
\end{tikzpicture}&\stackrel{\cref{bob}}{=}
\begin{tikzpicture}[centerzero,scale=1.3,Q]
	\draw[-,thin] (0.28,0) to[out=90,in=-90] (-0.28,.6);
	\draw[-to,thin] (-0.28,0) to[out=90,in=-90] (0.28,.6);
	\draw[-,thin] (0.28,-.6)\botlabel{i} to[out=90,in=-90] (-0.28,0);
	\draw[to-,thin] (-0.28,-.6)\botlabel{\tau i} to[out=90,in=-90] (0.28,0);
    \anticlockwiseinternalbubble{.16,-.42};
    \teleporter{-.23,-.43}{-.23,-.14};
        \teleporter{-.13,-.37}{-.13,-.22};
            \region{0.5,0}{\hat\lambda};
\end{tikzpicture}
\ \substack{\cref{altquadratic}\\{\textstyle=}\\\cref{teleslide}}\ 
\begin{tikzpicture}[centerzero,Q,scale=1.3]
\draw[to-,thin] (-.2,-.6)\botlabel{\tau i} to (-.2,.6);
\draw[-to,thin] (.3,-.6)\botlabel{i} to (.3,.6);
\anticlockwiseinternalbubble{.28,0};
\teleporter{-.2,.35}{.3,.35};
\teleporter{-.2,-.35}{.3,-.35};
    \region{0.55,0}{\hat\lambda};
\end{tikzpicture}\ .
\end{align}
Then the final one follows from the fifth one by applying $\Sigma$.

Finally suppose that $i = \tau j \neq \tau i$.
This is the first time that the additional sign from \cref{modifiedsigns} comes into play.
Again, there are six relations to check (after redistributing some internal bubbles):
\begin{align*}
-\begin{tikzpicture}[centerzero,scale=.94,Q]
	\draw[-,thin] (0.28,0) to[out=90,in=-90] (-0.28,.75);
	\draw[-,thin] (-0.28,0) to[out=90,in=-90] (0.28,.75);
	\draw[to-,thin] (0.28,-.75)\botlabel{\tau i} to[out=90,in=-90] (-0.28,0);
	\draw[to-,thin] (-0.28,-.75)\botlabel{i} to[out=90,in=-90] (0.28,0);
    \region{0.5,0}{\hat\lambda};
    \end{tikzpicture}
&=-
\begin{tikzpicture}[centerzero,scale=.94,Q]
	\draw[to-,thin] (0.2,-.75)\botlabel{\tau i} to (.2,.75);
	\draw[to-,thin] (-0.3,-.75)\botlabel{i} to (-.3,.75); 
    \Pinpin{-.3,0}{.2,0}{1.3,0}{Q_{i,\tau i}(x,y)};
        \region{2.4,0}{\hat\lambda};
\end{tikzpicture}\ ,\\
-\begin{tikzpicture}[centerzero,scale=.94,Q]
	\draw[-to,thin] (0.28,0) to[out=90,in=-90] (-0.28,.75);
	\draw[-to,thin] (-0.28,0) to[out=90,in=-90] (0.28,.75);
	\draw[-,thin] (0.28,-.75)\botlabel{i} to[out=90,in=-90] (-0.28,0);
	\draw[-,thin] (-0.28,-.75)\botlabel{\tau i} to[out=90,in=-90] (0.28,0);
    \region{0.5,0}{\hat\lambda};
    \end{tikzpicture}\,
&=-
\begin{tikzpicture}[centerzero,scale=.94,Q]
	\draw[-to,thin] (0.2,-.75)\botlabel{i} to (.2,.75);
	\draw[-to,thin] (-0.3,-.75)\botlabel{\tau i} to (-.3,.75); 
    \Pinpin{-.3,0}{.2,0}{1.5,0}{Q_{i,\tau i}(-x,-y)};
        \region{2.8,0}{\hat\lambda};
\end{tikzpicture}\ ,\\
-\begin{tikzpicture}[anchorbase,scale=.94,Q]
\draw[-to,thin] (-.35,-1.8)\botlabel{\tau i} to [out=90,in=180] (-.1,-1.4) to [out=0,in=90] (.15,-1.8);
	\draw[-,thin] (-0.35,-.3)\toplabel{i}  to [out=-90,in=90] (.15,-.9) to [out=-90,in=0] (-0.1,-1.2);
 \draw[to-,thin] (.15,-.3)  to [out=-90,in=90] (-.35,-.9) to [out=-90,in=180] (-.1,-1.2);
\teleporter{.13,-1.56}{.13,-1.04};
\anticlockwiseinternalbubble{-.32,-.9};
    \region{0.5,-1.1}{\hat\lambda};
\end{tikzpicture}
+
\begin{tikzpicture}[anchorbase,scale=.94,Q]
\draw[-to,thin] (-.35,1.8)\toplabel{i} to [out=-90,in=180] (-.1,1.4) to [out=0,in=-90] (.15,1.8);
	\draw[-,thin] (-0.35,.3) \botlabel{\tau i} to [out=90,in=-90] (.15,.9) to [out=90,in=0] (-0.1,1.2);
 \draw[to-,thin] (.15,.3)  to [out=90,in=-90] (-.35,.9) to [out=90,in=180] (-.1,1.2);
\teleporter{-.32,1.56}{-.32,1.04};
\clockwiseinternalbubble{.12,1};
    \region{0.5,1.1}{\hat\lambda};
\end{tikzpicture}\ &=\frac{\zeta_i c_i(\hat\lambda)}{c_{\tau i}(\hat\lambda)}\left[u^{\del_i}
\ \begin{tikzpicture}[anchorbase,scale=.94,Q]
\draw[-to,thin] (-.4,.8)\toplabel{i} to [out=-90,in=180] (-.1,.25) to [out=0,in=-90] (.2,.8);
\draw[-to,thin] (-.4,-.8)\botlabel{\tau i} to [out=90,in=180] (-.1,-.25) to [out=0,in=90] (.2,-.8);
\stringlabel{-1.3,-.3}{\tau i};
\anticlockwisebubble{-1.3,0};
\Circled{-.33,.5}{u};
\Circledbar{-.33,-.5}{u};
\node at (-.64,0) {$(-u)$};
\clockwisebubble{.5,0};
\stringlabel{.5,-.3}{i};
\node at (1,0) {$(u)$};
    \region{1.5,0}{\hat\lambda};
\end{tikzpicture}
\right]_{u:-1},
\\
-\begin{tikzpicture}[anchorbase,scale=.94,Q]
\draw[to-,thin] (-.35,1.8)\toplabel{\tau i} to [out=-90,in=180] (-.1,1.4) to [out=0,in=-90] (.15,1.8);
	\draw[to-,thin] (-0.35,.3) \botlabel{i} to [out=90,in=-90] (.15,.9) to [out=90,in=0] (-0.1,1.2);
 \draw[-,thin] (.15,.3)  to [out=90,in=-90] (-.35,.9) to [out=90,in=180] (-.1,1.2);
\teleporter{.12,1.56}{.12,1.04};
\anticlockwiseinternalbubble{-.32,1};
    \region{.5,1.1}{\hat\lambda};
\end{tikzpicture}+
\begin{tikzpicture}[anchorbase,scale=.94,Q]
\draw[to-,thin] (-.35,-1.8)\botlabel{i} to [out=90,in=180] (-.1,-1.4) to [out=0,in=90] (.15,-1.8);
	\draw[to-,thin] (-0.35,-.3)\toplabel{\tau i}  to [out=-90,in=90] (.15,-.9) to [out=-90,in=0] (-0.1,-1.2);
 \draw[-,thin] (.15,-.3)  to [out=-90,in=90] (-.35,-.9) to [out=-90,in=180] (-.1,-1.2);
\teleporter{-.33,-1.56}{-.33,-1.04};
\clockwiseinternalbubble{.12,-.9};
    \region{.5,-1}{\hat\lambda};
    \end{tikzpicture}
&=\frac{\zeta_i c_i(\hat\lambda)}{c_{\tau i}(\hat\lambda)}\left[u^{\del_i}
\ \begin{tikzpicture}[centerzero,scale=.94,Q]
\draw[to-,thin] (-.4,.8)\toplabel{\tau i} to [out=-90,in=180] (-.1,.25) to [out=0,in=-90] (.2,.8);
\draw[to-,thin] (-.4,-.8)\botlabel{i} to [out=90,in=180] (-.1,-.25) to [out=0,in=90] (.2,-.8);
\stringlabel{-1.3,-.3}{\tau i};
\anticlockwisebubble{-1.3,0};
\node at (-.64,0) {$(-u)$};
\clockwisebubble{.5,0};
\Circledbar{-.33,.5}{u};
\Circled{-.33,-.5}{u};
\stringlabel{.5,-.3}{i};
\node at (1,0) {$(u)$};
    \region{1.5,0}{\hat\lambda};
\end{tikzpicture}
\right]_{u:-1},
\\
\begin{tikzpicture}[centerzero,scale=.94,Q]
	\draw[-,thin] (0.28,0) to[out=90,in=-90] (-0.28,.75);
	\draw[-to,thin] (-0.28,0) to[out=90,in=-90] (0.28,.75);
	\draw[-,thin] (0.28,-.75)\botlabel{i} to[out=90,in=-90] (-0.28,0);
	\draw[to-,thin] (-0.28,-.75)\botlabel{i} to[out=90,in=-90] (0.28,0);
    \anticlockwiseinternalbubble{-.24,0};
        \region{.5,0}{\hat\lambda};
\end{tikzpicture}\ +\ 
\begin{tikzpicture}[centerzero,scale=.94,Q]
	\draw[-,thin] (0.3,-0.75) to[out=90, in=0] (0,-0.1);
	\draw[-to,thin] (0,-0.1) to[out = 180, in = 90] (-0.3,-.75)\botlabel{i};
    \stringlabel{.5,0}{\tau i};
  \draw[to-,thin] (1,0)++(-.28,0) arc(180:-180:0.28);
	\draw[to-,thin] (0.3,.75) to[out=-90, in=0] (0,0.1);
	\draw[-,thin] (0,0.1) to[out = -180, in = -90] (-0.3,.75)\toplabel{i};
    \clockwiseinternalbubble{1.24,0};
\teleporter{.28,.4}{.86,.25};
\teleporter{.28,-.4}{.86,-.25};
    \region{1.6,0}{\hat\lambda};
\end{tikzpicture}\ 
&=
-\begin{tikzpicture}[centerzero,scale=.94,Q]
	\draw[-to,thin] (0.2,-.75)\botlabel{i} to (.2,.75);
	\draw[to-,thin] (-0.3,-.75)\botlabel{i} to (-.3,.75); 
    \Pinpin{-.3,0}{.2,0}{1.4,0}{Q_{i,\tau i}(x,-y)};
\anticlockwiseinternalbubble{.2,-.4};
    \region{2.5,0}{\hat\lambda};
\end{tikzpicture}+\frac{\zeta_i c_i(\hat\lambda)}{c_{\tau i}(\hat\lambda)}
\left[u^{\del_i}
\ \begin{tikzpicture}[centerzero,scale=.94,Q]
\draw[-to,thin] (-.4,.8)\toplabel{i} to [out=-90,in=180] (-.1,.25) to [out=0,in=-90] (.2,.8);
\draw[to-,thin] (-.4,-.8)\botlabel{i} to [out=90,in=180] (-.1,-.25) to [out=0,in=90] (.2,-.8);
\stringlabel{-1.3,-.3}{\tau i};
\anticlockwisebubble{-1.3,0};
\Circled{-.33,.5}{u};
\Circled{-.33,-.5}{u};
\node at (-.64,0) {$(-u)$};
\clockwisebubble{.5,0};
\stringlabel{.5,-.3}{i};
\node at (1,0) {$(u)$};
    \region{1.5,0}{\hat\lambda};
    \end{tikzpicture}
\right]_{u:-1},\\
\begin{tikzpicture}[centerzero,scale=.94,Q]
	\draw[-to,thin] (0.28,0) to[out=90,in=-90] (-0.28,.75);
	\draw[-,thin] (-0.28,0) to[out=90,in=-90] (0.28,.75);
	\draw[to-,thin] (0.28,-.75)\botlabel{\tau i} to[out=90,in=-90] (-0.28,0);
	\draw[-,thin] (-0.28,-.75)\botlabel{\tau i} to[out=90,in=-90] (0.28,0);
    \clockwiseinternalbubble{.28,0};
        \region{.6,0}{\hat\lambda};
\end{tikzpicture}\ +\ 
\begin{tikzpicture}[centerzero,scale=.94,Q]
	\draw[to-,thin] (0.3,-0.75) to[out=90, in=0] (0,-0.1);
	\draw[-,thin] (0,-0.1) to[out = 180, in = 90] (-0.3,-.75)\botlabel{\tau i};
    \stringlabel{-.55,0}{i};
  \draw[to-,thin] (-1,0)++(.28,0) arc(0:360:0.28);
	\draw[to-,thin] (-0.3,.75) to[out=-90, in=180] (0,0.1);
	\draw[-,thin] (0,0.1) to[out = 0, in = -90] (0.3,.75)\toplabel{\tau i};
    \anticlockwiseinternalbubble{-1.24,0};
\teleporter{-.28,.4}{-.86,.25};
\teleporter{-.28,-.4}{-.86,-.25};
    \region{.5,0}{\hat\lambda};
\end{tikzpicture}
&=
-\begin{tikzpicture}[centerzero,scale=.94,Q]
	\draw[to-,thin] (0.2,-.75)\botlabel{\tau i} to (.2,.75);
	\draw[-to,thin] (-0.3,-.75)\botlabel{\tau i} to (-.3,.75); 
    \Pinpin{-.3,0}{.2,0}{1.4,0}{Q_{i,\tau i}(-x,y)};
\clockwiseinternalbubble{-.3,.4};
    \region{2.5,0}{\hat\lambda};
\end{tikzpicture}+\frac{\zeta_i c_i(\hat\lambda)}{c_{\tau i}(\hat\lambda)}
\left[u^{\del_i}
\ \begin{tikzpicture}[centerzero,scale=.94,Q]
\draw[to-,thin] (-.4,.8)\toplabel{\tau i} to [out=-90,in=180] (-.1,.25) to [out=0,in=-90] (.2,.8);
\draw[-to,thin] (-.4,-.8)\botlabel{\tau i} to [out=90,in=180] (-.1,-.25) to [out=0,in=90] (.2,-.8);
\stringlabel{-1.3,-.3}{\tau i};
\anticlockwisebubble{-1.3,0};
\Circledbar{-.33,.5}{u};
\Circledbar{-.33,-.5}{u};
\node at (-.64,0) {$(-u)$};
\clockwisebubble{.5,0};
\stringlabel{.5,-.3}{i};
\node at (1,0) {$(u)$};
    \region{1.5,0}{\hat\lambda};
\end{tikzpicture}
\right]_{u:-1}.
\end{align*}
The first follows from \cref{downwardscrossing}, the second from \cref{quadratic,tauantisymmetric},
the third from \cref{newlem},
the fourth follows from the third (with $i$ replaced by $\tau i$) on applying $\Sigma$ as in \cref{awkwardness},
and the fifth relation holds by \cref{goodnight}.
The final relation follows from the fifth one on composing 
on top with 
$\begin{tikzpicture}[Q,anchorbase,scale=1.2]
\draw[to-] (0,0) to (0,.5);
\draw[-to] (.4,0) to (.4,.5);
\anticlockwiseinternalbubble{0,.25};
\clockwiseinternalbubble{.4,.25}
\end{tikzpicture}\ $, 
cancelling internal bubbles with \cref{bubblyinverses,belltime}, 
then replacing $i$ by $\tau i$ and applying $\bar\Omega$.

\vspace{2mm}
\noindent
{\underline{Relation \cref{ibraid}}}:
Let $(i,j,k)$ be the labels of the strings at the bottom of the braids.
We noted in \cref{teapot} that the relation \cref{ibraid} is invariant under partial rotations. In view of the relations already proved, we can exploit this to reduce to proving the braid relation for just one triple $(i,j,k)$ from each of the rotation classes
$$
\left\{(i,j,k), (j,k,\tau i), (k, \tau i, \tau j),
(\tau i, \tau j, \tau k), 
(\tau j, \tau k, i), (\tau k, i, j)\right\}.
$$
Representatives for these classes are as follows:
\begin{itemize}
\item $(i,j,k)$ for $i,j,k \in I$ with $i \neq j, i \neq \tau j, j \neq k, j \neq \tau k, i \neq k$ and $i \neq \tau k$.
\item $(i,j,\tau i)$ for $i,j \in I$ with $i \neq \tau i$, $i \neq j$ and $i \neq \tau j$.
\item $(i,j,i)$ for $i,j \in I$ with $i \neq \tau i$, $i \neq j$ and $i \neq \tau j$.
\item $(i,j,i)$ for $i,j \in I$ with $i = \tau i$ and $i \neq j$.
\item $(i,i,i)$ for $i \in I$ with $i = \tau i$.
\item $(i,i,i)$ for $i \in I$ with $i \neq \tau i$.
\item $(i, \tau i, i)$ for $i \in I$ with $i \neq \tau i$.
\end{itemize} 
We will check these one by one.

\vspace{1mm}
\noindent
{\em Case $(i,j,k)$}:
This is the easiest case. Applying $\Xi^\imath$ 
to the relation \cref{ibraid} for $i,j,k$ 
with $i \neq j, i \neq \tau j, j \neq k, j \neq \tau k, i \neq k, i \neq \tau k$ (so that the right hand side is 0), 
also cancelling the constants which are the same on both sides and some pairs of oppositely oriented internal bubbles, this reduces to checking the following eight equalities:
\begin{align*}
.
\end{align*}
The bottom four equalities follow from the top four on applying $\bar\Omega$; for the last one, this also uses \cref{stateschool}.
Thus, we are reduced to checking the top four equalities. The first three follow directly from \cref{braid} and its rotations.
For the last one, by \cref{littlek,dotslide},
the identity to be checked is equivalent to
$$
%
\ ,
$$
which follows by \cref{altbraid,dotslide}.

\vspace{1mm}
\noindent
{\em Case $(i,j,\tau i)$ for $i \neq \tau i$}:
This is quite similar to Case $(i,j,k)$
taking $k = \tau i$. We obtain
the same eight equalities appearing as before, plus there are four additional equalities 
arising from the additional terms in \cref{psi5c} compared to \cref{psi5a}.
The first eight equalities follow by the same argument as in the previous case. 
After cancelling some inverse 
pairs of internal bubbles, the remaining four are as follows:
\begin{align*}
r_{j,\tau i}^{-1}r_{i,j}^{-1}
%
.
\end{align*}
The first two follow using \cref{teleslide,lotsmore,ruby,wax}.
For the last two, one has to use \cref{stateschool} first to replace the terms with three internal bubbles with terms with just one, then they follow using 
\cref{teleslide,lotsmore,ruby,wax} once again. 

\vspace{1mm}
\noindent
{\em Case $(i,j,i)$ for $i \neq \tau i$}:
This expands to 12 relations. 
The six in which the $j$ string is downward are
\begin{align*}
-r_{i,j}^{-1}r_{j,i}^{-1}%
\ .
\end{align*}
The first of these is \cref{downwardscrossing}. 
To prove the second, we use \cref{quickie} to see that the
left hand side equals
\begin{multline*}
%
.
\end{multline*}
The third and
the fourth equalities follow from the rotation of \cref{quadratic}.
By \cref{quickie,downwardscrossing}, the fifth identity to be proved is equivalent to
$$-\ 
%
,
$$
which follows using \cref{altquadratic,downwardscrossing} to simplify
the left hand side.
The proof of the sixth equation is similar to the proof of the fifth one.

The other six cases with upward middle string are
\begin{align*}
-r_{j,i}^{-1}%
\ .
\end{align*}
These may be verified in a similar way to the first six equalities, using also \cref{stateschool} to cancel internal bubbles in the fourth one.
Alternatively, they can be deduced from identities obtained from the first six by applying $\bar\Omega$.

\vspace{1mm}
\noindent
{\em Case $(i,j,i)$ for $i = \tau i$}:
This expands to 12 equations, eight of which (the first four of the downward six and the first four of the upward six) are the same as in the previous case,
with the same proofs. However, the remaining four identities are more complicated compared to the previous case
since a couple more terms arise from the extra two terms in \cref{psi5d} compared to \cref{psi5b}.
So we now need to 
show that
\begin{align*}
-r_{i,j}^{-1}
.
\end{align*}
In fact, the extra terms with cups and caps on the left hand sides of these equations cancel, leaving the same four identities that were already checked in the previous case.
For example, in the third equality here, we have that
$$
r_{j,i}^{-1}\ %
\ .
$$
The other three cases are similar.

\vspace{1mm}
\noindent
{\em Case $(i,i,i)$ for $i = \tau i$}:
This difficult case is checked in the proof of \cite[Th.~4.2]{BWWbasis} (where it is ``the second relation from (2.2)'').

\vspace{1mm}
\noindent
{\em Case $(i,i,i)$ for $i \neq \tau i$}:
There are 20 relations to be checked. We just list the first 10, with the rest essentially being obtained by applying $\bar\Omega$.
\begin{align*}
-%
\ .
\end{align*}
The first three of these follow from \cref{downwardsbraid} plus some rotation.
To establish the fourth equation,
we apply \cref{nomoretiktok} to remove the internal bubbles, also cancelling a couple of teleporters on both sides, to reduce to proving the identity
$$
%
\ .
$$
This is true because
\begin{align*}
%
\ .
\end{align*}
For the fifth relation, we expand the left hand side:
\begin{align*}
%
\ .
\end{align*}
This is what we want because
\begin{align*}
%
.
\end{align*}
The proof of the sixth relation is similar.
To prove the seventh equation, a couple of applications of \cref{teleslide} shows that both sides equal
$$
%
\ .
$$
The eighth equation follows in a similar way.
Finally, for the ninth equation, we first apply \cref{nomoretiktok} to see that it is equivalent to
$$
-%
,
$$
and this last equation 
follows by applying \cref{teleslide} to the top teleporter in the first term on the right hand side.
The tenth equation is proved similarly.

\vspace{1mm}
\noindent
{\em Case $(i,\tau i, i)$ for $i \neq \tau i$}:
This is the ibraid relation.
A more complicated expansion leads to 20 relations to be checked, which we split into two sets of 10. The second set is essentially the image of the first set under $\Omega$ up to redistributing some internal bubbles.
The first 10 relations are
\begin{align*}
%
\!\!
\right]_{\!u:-1}\!\!\!\!.
\end{align*}
We proceed to check these 10 relations, starting with the fourth one which we found was beyond our computational ability to check directly. 
Instead, we will prove it indirectly using \cref{necessary} and the invertibility of the morphisms \cref{bensidea} (this is why we included those).

At this point we have checked all of the defining relations except for 
ibraid, 
so we can appeal to \cref{necessary} to deduce that 
the image of the ibraid relation composed with
$%
$ holds in $\widehat{\UUloc}(\del,\zeta)$. 
In particular, the fourth equation follows because
$%
$ is invertible in 
$\widehat{\UUloc}(\del,\zeta)$ as we have inverted all of the 2-morphisms \cref{bensidea}.
It is not hard to verify the second, fifth, seventh, eighth and tenth relations directly, but there is no need to do that since they follow using the same trick as for the fourth one.

It remains to prove the first, third, sixth and ninth relations.
The first one is \cref{downwardsbraid}.
The third follows from \cref{downwardsbraid} transformed by a partial rotation.
The sixth relation follows from the fifth one by 
adding a clockwise internal bubble to the top right string and a counterclockwise internal bubble to the bottom right string, 
cancelling internal bubbles with \cref{bubblyinverses}, then applying $\Sigma$.
The ninth relation follows from the fifth one by attaching a rightward cap to the top right string and a leftward cup to the bottom left string.